\documentclass[a4paper,12pt]{article}
\usepackage{amssymb}
\usepackage{latexsym}
\newtheorem{Th}{Theorem}
\newtheorem{Prop}{Proposition}
\newtheorem{Co}{Corollary}
\newtheorem{Lm}{Lemma}
\newtheorem{Dfi}{Definition}
\newtheorem{Rm}{Remark}
\newtheorem{Con}{Conjecture}
\newcommand{\be}{\begin{equation}}
\newcommand{\ee}{\end{equation}}
\newcommand{\R}{\mathbb{R}}
\newcommand{\N}{\mathbb{N}}
\newcommand{\C}{\mathbb{C}}
\newcommand{\Z}{\mathbb{Z}}

\newcommand\res{\mathop{\hbox{\vrule height 7pt width .5pt depth 0pt
\vrule height .5pt width 6pt depth 0pt}}\nolimits}

\newcommand{\reset}{\setcounter{equation}{0}\setcounter{Th}{0}\setcounter{Prop}{0}\setcounter{Co}{0}
\setcounter{Lm}{0}\setcounter{Rm}{0}\setcounter{Con}{0}}

\def\ti{\tilde}
\def\lf{\left}
\def\rg{\right}

\def\al{\alpha}
\def\la{\lambda}
\def\e{\varepsilon}
\def\ep{\varepsilon}
\def\ds{\displaystyle}
\def\ov{\overline}
\def\Om{\Omega}
\def\om{\omega}
\def\p{\partial}
\def\bn{\vec{n}}
\def\bh{\mathbf{h}}
\def\bbe{\vec{e}}
\def\bH{\vec{H}}

\def\bp{\vec{\phi}}
\def\bP{\vec{\Phi}}

\def\res{\mathop{\hbox{\vrule height 7pt width .5pt 
depth 0pt\vrule height .5pt width 6pt depth 0pt}}\nolimits}
\begin{document}
\large
\title{\bf Conformally Invariant Variational Problems.}
\date{ }
\maketitle

\vskip 2 cm

\centerline{\bf \Large Tristan Rivi\`ere\footnote{ \large Department of Mathematics, ETH Zentrum
CH-8093 Z\"urich, Switzerland.}}

\vskip 3cm


\vspace{0.7cm}


\vspace{0.7cm}


\newpage

\section{Introduction}

These lecture notes form the cornerstone between two areas of Mathematics: calculus of variations and conformal invariance theory. 

\medskip

Conformal invariance plays a significant role in many areas of Physics, such as conformal field theory, renormalization theory, turbulence, general relativity. Naturally, it also plays an important role in geometry: theory of Riemannian surfaces, Weyl tensors, $Q$-curvature, Yang-Mills fields, etc... We shall be concerned with the study of conformal invariance in \underbar{analysis}. More precisely, we will focus on the study of nonlinear PDEs arising from conformally invariant variational problems (e.g. harmonic maps, prescribed mean curvature surfaces, Yang-Mills equations, amongst others).

\medskip

A transformation is called conformal when it preserves angles infinitesimally, that is, when its differential is a similarity at every point. Unlike in higher dimensions, the group of conformal transformations in two dimensions is very large ; it has infinite dimension. In fact, it contains as many elements as there are holomorphic and antiholomorphic maps. This particularly rich feature motivates us to restrict our attention on the two-dimensional case. Although we shall not be concerned with higher dimension, the reader should know that many of the results presented in these notes can be generalized to any dimension.

\medskip





\section{The parametric Plateau problem.}
\reset

\subsection{Introduction to the parametric Plateau problem.}

The first historical instance in which calculus of variations encountered conformal invariance took place early in the twentieth century with the resolution of the Plateau problem. Originally posed by J.-L. Lagrange in 1760, it was solved independently over 150 years later by J. Douglas and T. Rad\'o. In recognition of his work, the former was bestowed the first Fields Medal in 1936 (jointly with L. Alhfors).

\medskip

{\bf Plateau Problem.} {\it Given a jordan curve $\Gamma$ in $\R^m$, that is an injective continuous image of $S^1$, does there exist an immersion $u$ of the unit-disk $D^2$ such that $\partial D^2$ is homeomorphically sent onto $\Gamma$ by $u$ and for which $u(D^2)$ has a minimal area in the class of such immersions ?}

\medskip

The most natural approach for solving the Plateau problem would be to consider the direct minimization of the area functional 
of $C^1$ immersions $u$ from the disc $D^2$ into ${\R}^m$, sending it's boundary homeomorphically into $\Gamma$, and given explicitly by
\[
A(u)=\int_{D^2}|\p_xu\wedge\p_y u|\ dx\wedge dy\quad.
\]
where, for any pair of vector $a$ and $b$ in ${\R}^m$, $a\wedge b$ denotes the $2-$ vector\footnote{ If $\ep_i$ denotes the canonical basis of ${\R}^m$ and if $a=\sum_{i=1}^ma_i\,\ep_i$
(resp. $b=\sum_{i=1}^n\,b_i$)
\[
a\wedge b:= \sum_{i<j} [a_i\,b_j-a_j\,b_i]\ ep_i\wedge\ep_j
\]
and in particular, for the scalar product on $\wedge^2{\R}^m$ induced by the canonical scalar product on ${\R}^m$
\[
|a\wedge b|^2= \sum_{i<j} |a_i\,b_j-a_j\,b_i|^2
\]}
of $\wedge^2{\R}^m$ obtained by wedging $a$ and $b$.

\medskip

This natural approach however has little chance of success since the area $A(u)$ does not control enough of the map : in other words it is not enough {\it coercive}.
This observation is first the consequence of the huge invariance group of the Lagrangian $A$ : the space of positive diffeomorphisms of the disc $D^2$ : $Diff^+(D^2)$.

Indeed, given two distinct positive parametrizations $(x_1,x_2)$ and $(x'_1,x'_2)$ of the unit-disk $D^2$, there holds, for each pair of functions $f$ and $g$ on $D^2$, the identity
\[
\begin{array}{rl}
\ds df\wedge dg&\ds=\lf(\p_{x_1}f\p_{x_2}g-\p_{x_2}f\p_{x_1}g\rg)\ d{x_1}\wedge dx_2\\[5mm]
\ds\quad\quad&\ds=\lf(\p_{x'_1}f\p_{x'_2}g-\p_{x'_2}f\p_{x'_1}g\rg)\ d{x'_1}\wedge d{x'_2}\quad,
\end{array}
\]
so that, owing to $dx_1\wedge d{x_2}$ and $d{x_1}\wedge d{x_2}$ having the same sign, we find
\[
|\p_{x_1}f\p_{x_2}g-\p_{x_2}f\p_{x_1}g|\ dx\wedge dy=|\p_{x_1'}f\p_{x_2'}g-\p_{x_2'}f\p_{x_1'}g|\ dx'_1\wedge dx'_2\quad.
\]
This implies that $A$ is invariant through composition with positive diffeomorphisms. Thus if we would take a minimizing sequence $u_n$ of the area $A$ 
in the suitable class of immersions we can always compose this minimizing sequence with a degenerating sequence $\psi_n(x,y)\in Diff^+(D^2)$ in such a way that  
$u_n\circ\psi_n$ would for instance converge weakly to a point !

\medskip

Despite this possible degeneracy of the parametrisation, another difficulty while minimizing directly $A$ comes from the little control it gives on the
image itself $u_n(D^2)$. Starting from some given minimizing sequence $u_n$, one could indeed always modify $u_n$ by adding for instance tentacles of any sort filling more and more the space but without any significant cost
for $A(u_n)$ and keeping the sequence minimizing. The ''limiting'' object $\lim_{n\rightarrow +\infty} u_n(D^2)$ would then be some ''monster'' dense in ${\R}^m$.

\medskip

In order to overcome this lack of coercivity of the area lagrangian $A$, Douglas and Rad\'o proposed instead to minimize the energy of the map $u$

\[
E(u)=\frac{1}{2}\int_{D^2} |\p_x u|^2+|\p_y u|^2\ dx\wedge dy\quad .
\]
In contrast with $A$, $E$ has good coercivity properties and lower semicontinuity in the weak topology of the Sobolev space $W^{1,2}(D^2,{\R}^m)$, unlike the area One crucial observation is the following pointwise inequality valid for all $u$ dans $W^{1,2}(D^2,{\R}^3)$,
\[
|\p_xu\times\p_y u|\le \frac{1}{2}\lf[ |\p_x u|^2+|\p_y u|^2\rg]\quad\quad\mbox{ a.e. in }D^2
\]
and integrating this pointwise inequality over the disc $D^2$ gives
\[
A(u)\le E(u)\quad,
\]
with equality if and only if $u$ is weakly conformal, namely:
\[
|\p_xu|=|\p_yu|\quad\mbox{ et }\quad\p_xu\cdot\p_yu=0\quad\mbox{a.e.}\quad.
\]

\medskip

The initial idea of Douglas and Rad\'o bears resemblance to the corresponding strategy in 1 dimension while trying to minimize the length among all immersions $\gamma$ of the segment $[0,1]$ joining  two arbitrary points $a=\gamma(0)$ and $b=\gamma(1)$ in a riemannian manifold $(M^m,g)$,
\[
L(\gamma):=\int_{[0,1]}|\dot{\gamma}|_g\ dt\quad.
\]
The direct minimization of $L$ is here also made difficult by the existence of a huge non compact group of invariance for $L$ : the group of positive diffeomorphisms of the segment $[0,1]$ :
for any $C^1$ function $t(s)$ satisfying $t'>0$, $t(0)=0$ and $t(1)=1$ one has
\[
L(\gamma\circ s)=L(\gamma)\quad.
\]
The classical strategy to remedy to this difficulty is to minimize instead the energy of the immersion $\gamma$
\[
E(\gamma):=\int_0^1|\dot{\gamma}|_g^2\ dt\quad,
\]
in the space $W^{1,2}([0,1], M^m)$. The two lagrangians being related by the following inequality for any immersion $\gamma$
\[
L(\gamma)=\int_{[0,1]}|\dot{\gamma}|\ dt\le\lf[ \int_0^1|\dot{\gamma}|_g^2\ dt\rg]^{1/2}
\]
with equality if and only if the curve is in normal parametrisation :
\[
|\dot{\gamma}|_g=cte\quad\quad\mbox{ a.e }
\]
A classical results asserts that the minimization of $E$ provides a minimizer of $L$ in normal parametrization.

\medskip

Back to the two dimensional situation, in an ideal scenario, one could then hope a symmetry between the one dimensional and the two dimensional
cases in which, normal parametrization is replaced by conformal one, and to obtain, by minimizing $E$, a minimizer of $A$ in conformal parametrization.
This is eventually what will happen at the end but the 2-dimensional situation is analytically more complex and requires more work to be described as we will see below.

\medskip

In one dimension the Dirichlet energy $E$ is invariant under a finite dimensional group : the constant speed reparametrization. This is the main advantage
of working with $E$ instead of $L$. Similarly, in two dimension, the energy functional $E$ shares the same advantage over the area functional $A$ : while $A$ is invariant under the action of the {\it infinite} group of diffeomorphisms of $D^2$ into itself, the Dirichlet energy is invariant under the group of positive conformal diffeomorphism group of the disc, the M\"obius group
${\mathcal M}^+(D^2)$ which is 3 dimensional as we saw in the previous subsection and given by the holomorphic maps of the form
\[
f(w):=e^{i\theta}\,\frac{w-a}{1-\ov{a}w}
\]
for some $\theta\in{\R}$ and $a\in D^2$.

\medskip

This invariance of $E$ under conformal transformations may easily be seen by working with the complex variable $z=x_1+ix_2$. Indeed, we note
 $$\p_z:=\frac{1}{2}\lf(\p_{x_1}-i\p_{x_2}\rg)$$ et $$\p_{\ov z}:=\frac{1}{2}\lf(\p_{x_1}+i\p_{x_2}\rg)$$
so that $du=\p_zu\,dz+\p_{\ov{z}}u\,d\ov{z}$, and thus
\[
E(u)=\frac{i}{2}\int_{D^2}|\p_{z}u|^2+|\p_{\ov{z}}u|^2\ dz\wedge d\ov{z}\quad.
\]
Accordingly, if we compose $u$ with a conformal transformation, i.e. holomorphic, $z=f(w)$, there holds for $\ti{u}(w)=u(z)$ the identities
\[
|\p_{w}\ti{u}|^2=|f'(w)|^2\ |\p_zu|^2\circ f \ \quad \mbox{ and }\quad |\p_{\ov{w}}\ti{u}|^2=|f'(w)|^2\ |\p_{\ov{z}}u|^2\circ f \ \quad.
\]
Moreover, $dz\wedge d\ov{z}=|f'(w)|^2\ dw\wedge d\ov{w}$. Bringing altogether these results yields the desired conformal invariance $E(u)=E(\ti{u})$.

\medskip

An heuristic argument shed some light on the reason to believe that the strategy of minimizing the energy $E$ should provide a minimizer of $A$. Assuming one moment that we would have a 
smooth minimizer $u_\Gamma$ of the Dirichlet energy $E$ among our space of $C^1$ immersions sending homeomorphically $\p D^2$ onto the given Jordan curve $\Gamma$, then we claim that $u_\Gamma$ is conformal and minimizes also $A$ in the class. Indeed first, if $u_\Gamma$ would not minimize $A$ in this class, there would then be another immersion $v$ such that
$$A(v)<E(u_\Gamma)\quad.$$
Let $g$ be pull-back metric on the disc $D^2$, $g:=v^\ast g_{{\R}^m}$, where $g_{{\R}^m}$ denotes the standard scalar product in ${\R}^m$. The uniformization theorem gives the existence of a diffeomorphism $\Psi$ in $Diff^+(D^2)$
such that $\Psi^\ast g$ is conformal :
\[
e^{2\la}\, [dx_1^2+dx_2^2]= \Psi^\ast g=\Psi^\ast v^\ast g_{{\R}^m}=(v\circ\Psi)^\ast g_{{\R}^m}\quad.
\]
In other words we have that $v\circ\Psi$ is conformal therefore the following strict inequality holds
\[
E(v\circ\Psi)=A(v\circ\Psi)=A(v)<E(u_\Gamma)\quad,
\]
which is a contradiction. 

For similar reasons $u_\Gamma$ has to be conformal. Indeed if this would not be the case we would again find a diffeomorphism $\Psi$ such that $u_\Gamma\circ\Psi$
is conformal and, under this assumption that $u_\Gamma$ is not conformal we would have
\[
E(u_\Gamma\circ\Psi)=A(u_\Gamma\circ\Psi)=A(u_\Gamma)<E(u_{\Gamma})\quad,
\]
which would be again a contradiction.

\medskip

Of course this heuristic argument is based on the hypothesis that we have found a minimizer and that this minimizer is a smooth immersion, since we applied the uniformization theorem
to the induced metric $u_\Gamma^\ast g_{{\R}^m}$. In order to use the ''nice'' functional properties of the lagrangian $E$, as we mentioned above we have to enlarge the
class of candidates for the minimization to the space of Sobolev maps in $W^{1,2}(D^2,{\R}^m)$, continuous at the boundary, and sending $\p D^2$ monotonically onto $\Gamma$.
This weakening of the regularity for the space of maps prevents a-priori to make use of the uniformization theorem to a weak object such as $u^\ast g_{{\R}^m}$. However,
the following theorem of Morrey gives an ''almost uniformization'' result that permits to overcome this difficulty of the lack of regularity (see theorem 1.2 of \cite{Mor1} about
{\it $\epsilon-$conformal} parametrization)
\begin{Th}
\label{th-n.I.1}\cite{Mor1}
Let $u$ be a map in the Sobolev space $C^0\cap W^{1,2}(\ov{D^2},{\R}^m)$ and let $\ep>0$, then there exists an homeomorphism $\Psi$ of
the disc  such that $\Psi\in W^{1,2}(D^2,D^2)$,
\[
u\circ\Psi\in C^0\cap W^{1,2}(\ov{D^2},{\R}^m)\quad,
\]
and
\[
E(u\circ\Psi)\le A(u\circ\Psi)+\ep= A(u)+\ep\quad.
\]
\hfill $\Box$
\end{Th}

\medskip

The main difficulty remains to find a $C^1$ immersion minimizing the Dirichlet energy $E$. Postponing to later the requirement for the map $u$ to realize
an immersion of the disc, one could first try to find a general minimizer of $E$ within the class of Sobolev maps in $W^{1,2}(D^2,{\R}^m)$, continuous at the boundary, and sending $\p D^2$ monotonically onto $\Gamma$. By monotonically we mean the following relaxation of the homeomorphism condition - which is too restrictive in the first approach.
\begin{Dfi}
\label{df-n.I.1}
Let $\Gamma$ be a Jordan closed curve in ${\R}^m$, i.e. a subset of ${\R}^m$ homeomorphic to $S^1$, and let $\gamma$ be an homeomorphism from $S^1$ into $\Gamma$.
We say that a continuous map $\psi$ from $S^1$ into $\Gamma$ is weakly monotonic, if there exists a non decreasing  continuous function $\tau$ : $[0,2\pi]\rightarrow {\R}$
with $\tau(0)=0$ and $\tau(2\pi)=2\pi$ such that
\[
\forall \theta\in[0,2\pi]\quad\quad \psi(e^{i\theta})=\gamma(e^{i\tau(\theta)})\quad.
\]
\hfill $\Box$
\end{Dfi}

We now give the class in which we will proceed to the minimization argument. This is the following subset of $W^{1,2}(D^2,{\R}^m)$
\[
{\mathcal C}(\Gamma):=\lf\{
\begin{array}{c}
\ds u\in W^{1,2}\cap C^0(\ov{D^2},{\R}^m)\quad; \quad u_{|_{\partial D^2}}\in C^0(\p D^2,\Gamma)\\[5mm]
\ds\quad\mbox{ and $u$ is weakly monot. on }\p D^2
\end{array}
\rg\}
\]
One should stress the fact at this stage that the boundary data is a ''free Dirichlet boundary'' data in the sense that the value of $u$ on $\p D^2$
is not prescribed in a pointwise way but only globally by requiring that $u$ covers $\Gamma$ monotonically on $\p D^2$. This large degree of freedom
 is required in order to hope the parametrization to adopt the conformal configuration - see proposition~\ref{pr-n.I.2} below - but, in the mean time this is also all the source of our
 difficulties. Indeed, in contrast  with a ''classical'' Dirichlet boundary condition that would easily pass to the limit in the minimization process due
 to the continuity of the trace operator from $W^{1,2}(D^2,{\R}^m)$ into $H^{1/2}(D^2,{\R}^m)$, the requirement  in the''free Dirichlet boundary''  for $u$ to be continous
 on $\p D^2$ for instance does not a-priori pass to the limit\footnote{ Indeed for any $0<\al<1/2$ the map 
 \[
 z\rightarrow|\log|z-1||^\al
 \]
 is in $W^{1,2}(D^2)$ but is trace $e^{i\theta}\rightarrow |\log|e^{i\theta}-1||^\al$ is not bounded in $L^\infty(\p D^2)$.} since
 \be
 \label{n.I.0}
 W^{1,2}(D^2,{\R}^m)\nsubseteq C^0(\p D^2,{\R}^m)\quad.
 \ee

\subsection{A proof of the Plateau problem for rectifiable curves.}

We shall now give a proof of the existence of a minimizer under the assumption that $\Gamma$ has a finite length : in other words if $\gamma$ is an homeomorphism sending 
$S^1$ onto $\Gamma$ we assume that
\be
\label{n.I.2}
L(\Gamma):=\sup_{t_0=0<t_1<\cdots<t_N=2\pi}\sum_{k=1}^N|\gamma(e^{2\pi\, i\, t_k})-\gamma(e^{2\pi\, i\, t_{k+1}})|<+\infty
\ee
Such a Jordan curve is also simply called a {\it rectifiable closed curve}\footnote{This condition is implies that the curve is rectifiable in classical {\it Geometric measure theory} sense. Indeed one can pass to the limit in sequences of discretization of $S^1$ and use Federer-Fleming compactness theorem see\cite{Fe}. The finite length condition of our Jordan curve  also characterized by ${\mathcal H}^1(\Gamma)<+\infty$  where ${\mathcal H}^1$
denotes the 1-dimensional Hausdorff measure in ${\R}^m$ - see \cite{Fal} chapter 3.}

\begin{Th}
\label{th-I.0}
[Douglas-Rad\'o-Courant-Tonelli] Given a rectifiable closed curve $\Gamma$ in $\R^m$, there exists a  minimum $u$ for the energy $E$ within the space ${\mathcal C}(\Gamma)$ of $W^{1,2}(D^2,{\R}^m)$ functions mapping the boundary of the unit-disk $\p D^2$ onto $\Gamma$ continuously and monotonically. Any minimum $u$ of 
\be
\label{n.I.3}
\min_{u\in{\mathcal C}(\Gamma)}E(u)
\ee
satisfy
\[
u\quad\mbox{ minimizes $A$ in }{\mathcal C}(\Gamma)\quad,
\]
the following identities hold
\be
\label{I.1}
\left\{
\begin{array}{l}
\Delta u=0\quad\quad\mbox{ in }\quad D^2\\[5mm]
|\p_{x_1}u|^2-|\p_{x_2} u|^2-2i\,<\p_{x_1}u,\p_{x_2}u>=0\quad\quad\mbox{ in }\quad D^2\quad .
\end{array}\rg.
\ee
and moreover $u$ is in $C^\infty(D^2,{R}^m)\cap C^0(\ov{D^2},{\R}^m)$ and it realizes an homeomorphism from $\p D^2$ into $\Gamma$.
\hfill$\Box$
\end{Th} 
\begin{Rm}
\label{rm-n.I.1}
Assume $u$ is an immersion, the harmonicity and conformality condition exhibited in (\ref{I.1}) compared with the formula (???) imply that the mean curvature of $u(D^2)$ realizes a minimal surface. This
is not a surprise since $u$ has to minimize $A$ too in the class ${\mathcal C}(\Gamma)$ and  the expression of the  first area variation  (\ref{VI.63}) implies that the mean curvature of the immersion has to vanish for such an area minimizing immersion.
\end{Rm}

In order to fully solve the Plateau problem it remains to study whether the solutions given by the minimization problem (\ref{n.I.3}) is an immersion or not. 
A first observation in this direction can be made. The harmonicity of $u$ (\ref{I.1}) says that $div(\nabla u^k)=0$ on $D^2$ for each component $u^k$ of $u$. 
Applying Poincar\'e lemma, we can introduce the harmonic conjugates $v^k$ of $u^k$ satisfying
\[
(-\p_{x_2}v^k,\p_{x_1}v^k)=\nabla^\perp v^k:=\nabla u^k=(\p_{x_1}u^k,\p_{x_2}u^k)
\]
This implies that $f^k(z):=u^k-iv^k$ is holomorphic and
 $$
 |f'(z)|^2=4^{-1}\,[|\p_{x_1}u^k-\p_{x_2}v^k|^2+|\p_{x_2}u^k+\p_{x_1}v^k|^2]=|\nabla u^k|^2\quad.
 $$
 Since $u$ is conformal, it is an immersion at a given point (i.e $|du\wedge du|\ne0$) if and only if $|\nabla u|\ne 0$
 which is then equivalent to $|f'(z)|\ne 0$. Since $f'(z)$ is also holomorphic, we obtain that $u$ is a-priori  is an immersion
 away from isolated so called {\it branched points}. In other words $u$ is is what is called a {\it branched immersion}.
 
\medskip 

Finally one has to study the possibility for the branch points to exist or not. It is clear that in codimension larger than 1, such a branch point can exist as
one can see by taking a sub part of a complex algebraic curve in ${\R}^4\simeq {\C}^2$ such as 
\[
z\rightarrow (z^2,z^3)
\]
This disc is calibrated by the standard K\"ahler form of ${\C}^2$ and is then area minimizing for its boundary data, the curve $\Gamma$ given by $e^{i\theta}\rightarrow (e^{2i\theta},e^{3i\theta})$.

\medskip

In codimension 1, $m=3$, however the situation is much more constrained and the delicate analysis to study the possibility of the existence of branched points
is  beyond the scope of this chapter which is just intended to motivated the subsequent ones on general conformally invariant variational problems.
It has been proved by  R. Osserman that $u$ has no interior {\it true branched point} (see \cite{Oss}) and by R.Gulliver, R.Osserman and H.Royden that $u$ has no interior {\it false branched point} neither (see \cite{GOR}). Therefore, for $m=3$, any minimizer $u$ of (\ref{n.I.3}) is  an immersion in the interior of $D^2$.
A thorough presentation of the Plateau problem can be found for instance in \cite{DHKW1} and \cite{DHKW2}.

\medskip



We first establish some elementary properties of minimizers of (\ref{n.I.3}) assuming it exists. We are going to first establish (\ref{I.1}) postponing to later the existence question.

\medskip

We have the following elementary proposition
\begin{Prop}
\label{pr-n.I.0}
Let $u$ be the weak limit in $W^{1,2}$ of a minimizing sequence of $E$ in ${\mathcal C}(\Gamma)$. Then $u$ satisfies the Laplace equation 
\[
\forall\, i=1\cdots m\quad\quad\Delta u^i=0\quad\quad\mbox{ in }{\mathcal D}'(D^2)\quad.
\]
where $u^i$ are the components of $u$ in ${\R}^m$.\hfill $\Box$
\end{Prop}

\medskip

\noindent{\bf Proof of proposition~\ref{pr-n.I.0}.}

Let $u_n$ be a minimizing sequence. Modulo extraction of a subsequence we can always assume that $u_n$ weakly converges
to a limit $u\in W^{1,2}(D^2,{\R}^m)$. Let  $\phi\in C^\infty_0(D^2,{\R}^m)$ be a smooth function, compactly supported in $D^2$.

It is clear that for every $t\in {\R}$ and for every $n\in {\N}$ $$u_n+t\phi\in{\mathcal C}(\Gamma)\quad.$$ Hence we have for any $t\in{\R}$
\be
\label{n.s.I.0}
 \inf_{u\in{\mathcal C}(\Gamma)}E(u)=\lim_{n\rightarrow +\infty}E(u_n)\le\liminf_{n\rightarrow +\infty}E(u_n+t\phi)\quad.
\ee
Observe that we have
\[
\begin{array}{l}
E(u_n+t\phi)=E(u_n)\\[5mm]
\ds\quad\quad\quad+t\ \int_{D^2}<\nabla u_n\cdot\nabla\phi>\ dx+\frac{t^2}{2}\int_{D^2}|\nabla\phi|^2\ dx
\end{array}
\]
Since $u_n$ weakly converges to $u$ in $W^{1,2}(D^2)$ we have
\[
\int_{D^2}<\nabla u_n\cdot\nabla\phi>\ dx\quad\longrightarrow\quad \int_{D^2}<\nabla u\cdot\nabla\phi>\ dx\quad.
\]
This fact combined with the previous identity gives
\be
\label{n.s.I.1}
\begin{array}{l}
\ds\lim_{n\rightarrow+\infty}E(u_n+t\phi)=\inf_{u\in{\mathcal C}(\Gamma)}E(u)\\[5mm]
\ds\quad\quad\quad+t\ \int_{D^2}<\nabla u\cdot\nabla\phi>\ dx+\frac{t^2}{2}\int_{D^2}|\nabla\phi|^2\ dx\quad.
\end{array}
\ee
the inequality (\ref{n.s.I.0}) together with the identity (\ref{n.s.I.1}) imply
\[
\forall\, t\in {\R}\quad \quad t \int_{D^2}<\nabla u\cdot\nabla\phi>\ dx+\frac{t^2}{2}\int_{D^2}|\nabla\phi|^2\ dx\ge 0
\]
As a consequence of this inequality by dividing by $|t|$ and making $t$ tend to zero respectively from the right and from the left we have obtained
\be
\label{n.s.I.2}
\forall\, \phi\in C^\infty_0(D^2,{\R}^m)\quad\quad\int_{D^2}<\nabla u\cdot\nabla\phi>\ dx=0\quad.
\ee
which is the desired result and proposition~\ref{pr-n.I.0} is proved.\hfill $\Box$

\medskip

Observe that (\ref{n.s.I.2}) is equivalent to
\be
\label{n.I.4}
\forall\, \phi\in C^\infty_0(D^2,{\R}^m)\quad\quad\frac{d}{dt}E(u+t\phi)_{|_{t=0}}=0\quad.
\ee
This is the {\bf Euler-Lagrange equation} associated to our variational problem, the Laplace equation in the present case, and the proposition we just proved
can be summarized by saying that every weak limit
of our minimizing sequence  is  a critical point of the lagrangian
for {\bf variations in the target} : variations of the form $u+t\phi$. In order to prove the second line of (\ref{I.1}), the conformality of $u$, we will exploit instead that the minimizer $u$ - once we will prove it's existence -
has to be a critical point
 for {\bf variations in the domain} : variations of the form $u(id+tX)$ where $X\in C^\infty_0(D^2,{\R}^2)$. This condition is called the {\bf stationarity condition} :
$u$ satisfies
\be
\label{n.I.5}
\forall\, X\in C^\infty_0(D^2,{\R}^2)\quad\quad\frac{d}{dt}E(u(id+t X))_{|_{t=0}}=0\quad.
\ee
We shall now prove the following proposition.
\begin{Prop}
\label{pr-n.I.1}
A map $u$ in $W^{1,2}(D^2,{\R}^m)$ satisfies the stationarity condition (\ref{n.I.5}) if and only if its {\it Hopf differential}
\[
\begin{array}{l}
\ds h(u)=H(u)\ dz\otimes dz:=<\p_z u ,\p_zu>\ dz\otimes dz\\[5mm]
\ds\quad=4^{-1}\,\lf[|\p_{x_1}u|^2-|\p_{x_2}u|^2-2\, i\ <\p_{x_1}u,\p_{x_2}u>\rg] dz\otimes dz
\end{array}
\]
is holomorphic.\hfill$\Box$
\end{Prop}
\begin{Rm}
\label{rm-n.I.2}
The stationarity condition for general lagrangians in mathematical physics, equivalent to the conservation law 
\be
\label{n.I.5a}
\p_{\ov{z}}H(u)\equiv 0
\ee
for the Dirichlet energy- the so called $\sigma-$model - corresponds to the {\bf conservation of the stress-energy tensor} (see for instance (2.2) in II.2 of \cite{JaTa} for the Yang-Mills-Higgs lagrangian).\hfill $\Box$
\end{Rm}
\noindent{\bf Proof of proposition~\ref{pr-n.I.1}.} We denote $x_t$ the flow associated to $X$ such that $x(0)=x$ and $X=\sum_{i=1}^2X^i\p_{x_i}$. With this notation we apply the pointwise chain rule and obtain
\[
\p_{x_k}(u(x_t))=\sum_{i=1}^2\p_{x_i} u(x_t)\ \p_{x_k}x^i_t\quad,
\]
which gives in particular
\be
\label{n.I.5.b}
\begin{array}{l}
\ds\int_{D^2}|\nabla(u(x_t))|^2\ dx=\int_{D^2}|\nabla u|^2(x_t)\ dx\\[5mm]
\ds\quad\quad+2\ t\ \int_{D^2}(\p_{x_i}u)(x_t)\ (\p_{x_j}u)(x_t)\ \p_{x_i}X^j\ dx+o(t)
\end{array}
\ee
where $o(t)$ means Observe that for an $L^1$ function $f$ and any $\phi\in C^\infty(D^2)$ one has for $t$ small enough, since $X$ is compactly supported in $D^2$,
\[
\begin{array}{l}
\ds\int_{D^2}f(x_t)\ \phi(x)\ dx=\int_{x_t(D^2)} f(y)\ \phi(x_t^{-1})\ d(x_t^{-1}(y))\\[5mm]
\ds=\int_{D^2} f(y)\ (\phi(x)-t\nabla_X\phi+o(t))\ (1-t\, div\,X+o(t))\ dy
\end{array}
\]
Thus we obtain
\be
\label{n.I.5.c}
\frac{d}{dt}\lf(\int_{D^2}f(x_t)\ \phi(x)\ dx\rg)=-\int_{D^2} f(x)\ div(\phi\ X)\ dx
\ee
Hence we deduce from this identity and (\ref{n.I.5.b})
\be
\label{n.I.6}
\begin{array}{l}
\ds\frac{d}{dt}\lf(\int_{D^2}|\nabla(u(x_t))|^2\ dx\rg)_{|_{t=0}}\\[5mm]
\ds=-\int_{D^2} |\nabla u|^2\ div\,X\ dx+2\int_{D^2}\sum_{i,j=1}^2\p_{x_i}u\,\p_{x_j}u\ \p_{x_i}X^j\ dx
\end{array}
\ee
The assumption (\ref{n.I.5}) is then equivalent to
\be
\label{n.I.7}
\forall\, l=1,2\quad\quad\frac{\p}{\p x_l} |\nabla u|^2-2\sum_{k=1}^2\frac{\p}{\p x_k}\lf[\frac{\p u}{\p x_k}\,\frac{\p u}{\p x_l}\rg]=0\quad.
\ee
This reads
\be
\label{n.I.8}
\lf\{
\begin{array}{l}
\ds \frac{\p}{\p x_1}\lf(|\p_{x_1} u|^2-|\p_{x_2} u|^2\rg)+2\,\frac{\p}{\p x_2}\lf(\frac{\p u}{\p x_1}\ \frac{\p u}{\p x_2}\rg)=0\\[5mm]
\ds \frac{\p}{\p x_2}\lf(|\p_{x_1} u|^2-|\p_{x_2} u|^2\rg)-2\,\frac{\p}{\p x_2}\lf(\frac{\p u}{\p x_1}\ \frac{\p u}{\p x_2}\rg)=0\quad,
\end{array}
\rg.
\ee
which is also equivalent to (\ref{n.I.5a}) and proposition~\ref{pr-n.I.1} is proved \hfill $\Box$.

\begin{Rm}
\label{rm-n.I.3}
There are situations when, for a Lagrangian $L$, being a critical point {\bf for variations in the target} - i.e. solution to the {\bf Euler Lagrange Equation} -
\[
\forall\, \phi\in C^\infty_0(D^2,{\R}^m)\quad\quad\frac{d}{dt}L(u+t\phi)_{|_{t=0}}=0\quad,
\]
 implies that you are automatically a critical point for {\bf variations in the domain} - i.e. solution to the {\bf stationary equation}
 \[
 \forall\, X\in C^\infty_0(D^2,{\R}^2)\quad\quad\frac{d}{dt}L(u(id+t X))_{|_{t=0}}=0\quad.
 \]
 This happens in particular
 when the solution $u$ of the Euler Lagrange equation is {\bf smooth}. In that case indeed a variation of the form $u(id+tX)\simeq u(x)+t \nabla_Xu+o(t^2)$ 
 can be interpreted as being a variation in the target with $\phi$ of the form $\phi:=\nabla_Xu$.
 However this is {\bf not true in general} and there are situations where {\it weak solutions} to the Euler Lagrange Equations do not satisfy the stationarity condition
 (see \cite{Riv}).
 
 In our present case with the Dirichlet energy - i.e. $L=E$ - solutions to the Laplace equation are smooth and one obtains (\ref{n.I.8}) by multiplying
  $\Delta u =0$ respectively by $\p_{x_1} u$ and by $\p_{x_2} u$.\hfill $\Box$
\end{Rm}

\begin{Rm}
\label{rm-n.I.4}
This relation between the Euler Lagrange equation and the stationarity equation mentioned in the previous remark~\ref{rm-n.I.3} shed also some lights on the
reason why the stationarity equation is related to the conservation of the stress energy tensor. Take for instance the simplest system in classical mechanics of a single point
particle of mass $m$ moving in a potential $V$. The Lagrangian attached to this system is given by
\[
L(x(s)):=\int {m}\frac{\dot{x}^2(s)}{2}-V(x(s))\ ds
\]
and the {\bf law of motion} for this particle is given by the {\bf Euler-Lagrange} equation
\[
m\,\ddot{x}(s)+V'(x(s))=0
\]
{\bf Varitaions in the domain} correspond to perturbation of the form $x(s+t X)\simeq x(s)+t\,\dot{x}(s)+o(t^2)$. Multiplying the Euler Lagrange equation by the infinitesimal perturbation $\dot{x}(s)$ corresponding then to the variation in the domain gives the {\bf stationary equation}
\[
0=\dot{x}(s)\ [m\,\ddot{x}(s)+V'(x(s))]=\frac{d}{ds}\lf[\frac{m}{2}\ \dot{x}^2(s)+V(x(s))\rg]
\]
which is nothing but the {\bf conservation of energy}.
\hfill $\Box$
\end{Rm}

\medskip

We shall now exploit the specificity of the boundary condition imposed by the membership of the minimizer $u$ of $E$ in ${\mathcal C}(\Gamma)$, whose existence
is still assumed at this stage of the proof in order to get more information on the Hopf differential and the fact that $H(u)$ is identically zero.
This is the result of the {\it free Dirichlet} condition we are imposing. By imposing a {\it fixed Dirichlet} condition there would have been no reason for the holomorphic
Hopf differential to be identically zero and hence the minimizer $u$ to be conformal. Precisely we are now going to prove the following proposition.

\begin{Prop}
\label{pr-n.I.2}
Let $u$ be a map in $W^{1,2}(D^2,{\R}^m)$ satisfying
\be
\label{n.I.9}
\begin{array}{l}
\ds\forall\, X\in C^\infty(\ov{D^2},{\R}^2)\quad\mbox{ s. t. }\quad X\cdot x\equiv 0\quad\mbox{ on }\quad\p D^2\\[5mm]
\ds\quad \quad\frac{d}{dt}E(u(id+t X))_{|_{t=0}}=0\quad.
\end{array}
\ee
then
\[
|\p_{x_1}u|^2-|\p_{x_2}u|^2-2\, i\ <\p_{x_1}u,\p_{x_2}u>\equiv 0\quad\mbox{ on }D^2\quad.
\]
\hfill $\Box$
\end{Prop}
\noindent{\bf Proof of proposition~\ref{pr-n.I.2}.} Let $X$ be an arbitrary smooth vector-field on the disc $D^2$ satisfying
\be
\label{n.I.10}
X\cdot x\equiv 0\quad\quad\mbox{ on }\p D^2\quad.
\ee
Thus the flow $x_t$ of the vector-field $X$ preserves $D^2$ and (\ref{n.I.5.c}) still holds. Thus we have also (\ref{n.I.6}). The stationarity assumption (\ref{n.I.9})
implies then
\be
\label{n.I.11}
\lim_{r\rightarrow 1-}\int_{B^2_r(0)}\lf[-|\nabla u|^2\ div\,X+2\sum_{i,j=1}^2\p_{x_i}u\,\p_{x_j}u\ \p_{x_i}X^j\rg]\ dx=0\quad.
\ee
We adopt for $X$ the complex notation $X:=X^1+iX^2$ and we observe that (\ref{n.I.11}) becomes
\be
\label{n.I.12}
\lim_{r\rightarrow 1-}\int_{B^2_r(0)}\Re\lf[H(u)\ \frac{\p X}{\p\ov{z}}\rg]\ dx_1\wedge dx_2=0\quad,
\ee
which also reads
\be
\label{n.I.13}
\lim_{r\rightarrow 1-}\Re\lf(\frac{i}{2}\int_{B^2_r(0)}H(u)\ \frac{\p X}{\p\ov{z}}\ dz\wedge d\ov{z}\rg)=0
\ee
or equivalently
\be
\label{n.I.14}
-\lim_{r\rightarrow 1-}\Re\lf(\frac{i}{2}\int_{B^2_r(0)}H(u)\ dX\wedge dz\rg)=0
\ee
From proposition~\ref{pr-n.I.1} we know that $H(u)$ is holomorphic and then smooth in the interior of $D^2$. We can then integrate
by part and we obtain
\be
\label{n.I.15}
\lim_{r\rightarrow 1-}\Im\lf(\int_{B^2_r(0)}X\ dH(u)\wedge dz-\int_{\p B^2_r(0)} X\ H(u)\ dz\rg)=0\quad.
\ee
But $dH(u)\wedge dz=\p_{\ov{z}}H(u)\ d\ov{z}\wedge dz\equiv 0$, thus we finally obtain
\be
\label{n.I.16}
\lim_{r\rightarrow 1-}\Im\lf(\int_{\p B^2_r(0)} X\ H(u)(z)\  dz\rg)=0\quad.
\ee
We choose the vector field $X$ to be of the form $X= i\al\, z$ where $\al(e^{i\theta})$ is an arbitrary real function independent of $|z|$ in the neighborhood of $\p D^2$
(the boundary condition (\ref{n.I.10}) is clearly satisfied with this choice).
Observe that the restriction of $dz=d(r\,e^{i\theta})$ to $\p B^2_r(0)$ is equal to $i\,e^{i\theta}\,r\,d\theta=i\, z\, d\theta$. For this choice of $X$, (\ref{n.I.16}) becomes
\be
\label{n.I.17}
\lim_{r\rightarrow 1-}\Im\lf(\int_{\p B^2_r(0)} \al(\theta)\ H(u)(z)\,z^2\  d\theta\rg)=0\quad.
\ee
Let $z_0\in D^2$ and choose for $\al:=G(\theta,z_0)$ the Poisson Kernel such that, for any harmonic function $f$
\[
f(z_0)=\int_0^{2\pi} G(\theta,z_0)\ f(e^{i\theta})\ d\theta\quad.
\]
Since $H(u)$ is holomorphic, $H(u)(z)\, z^2$ is holomorphic and hence harmonic on $D^2$. This is also
of course the case for $H(u)(r\,z)\, r^2\,z^2$ for any $0<r<1$. We then obtain from (\ref{n.I.17}) 
\[
\begin{array}{l}
\ds 0=\lim_{r\rightarrow 1-}\Im\lf(\int_{0}^{2\pi} G(\theta,z_0)\ H(u)(r\,e^{i\theta})\,r^2 e^{2\,i\,\theta}\  d\theta\rg)\\[7mm]
\ds \quad=\lim_{r\rightarrow 1-}\Im(H(u)(r\,z_0)\, r^2\,z_0^2)=\Im(H(u)(z_0)\, z_0^2)\quad.
\end{array}
\]
This holds for any $z_0$ and therefore the holomorphic function $H(u)(z)\,z^2$, whose imaginary part vanishes, has to be 
identically equal to a real constant :
\[
H(u)(z)\, z^2\equiv c\quad.
\]
Since $H(u)$ is holomorphic without pole at the origin, $c$ has to be equal to zero. We have then established that $H(u)\equiv 0$ on $D^2$
which concludes the proof of proposition~\ref{pr-n.I.2}. \hfill $\Box$

\medskip

We will now concentrate our efforts for proving the existence of a minimizer of $E$ in ${\mathcal C}(\Gamma)$. While taking a minimizing
sequence $u_n$ we have already explained the risk for the boundary requirement $u_n\in C^0(\p D^2,\Gamma)$ not to be 
preserved at the limit do to the lack of Sobolev embedding (\ref{n.I.0}).

\medskip

Another difficulty lies in the remaining degree of freedom given by the invariance group of $E$ on ${\mathcal C}_\Gamma$, the so called gauge group
of our problem, which is here the M\"obius group ${\mathcal M}^+(D^2)$
which is three dimensional as we saw in section 1. The problem with this gauge group is non compact : by taking for instance a sequence $a_n\in D^2$ and $a_n\rightarrow (1,0)$
the sequence of maps
\[
\psi_n\ :\ z\ \longrightarrow\ \frac{z-a_n}{1-\ov{a_n} z}
\]
converges weakly to a constant map which is not in  ${\mathcal M}^+(D^2)$ anymore. Assuming we would have a sequence of minimizer $u_n$
converging to a solution of (\ref{n.I.3}) and satisfying the conclusions of theorem~\ref{th-I.0}, by composing $u_n$ with $\psi_n$ we still have a minimizing sequence since $E(u_n)=E(u_n\circ\psi_n)$. 
this new minimizing sequence of $E$ in ${\mathcal C}(\Gamma)$, $u_n\circ\psi_n$, converges then to a constant which cannot be a solution to the Plateau problem.

Thus all minimizing sequences cannot lead to a solution du in particular to the existence of a non compact gauge group ${\mathcal M}^+(D^2)$. This group however
is very small (in comparison with $Diff^+(D^2)$ in particular). 

 In order to break this gauge invariance it suffices to fix the images of 3 distinct points on the boundary. This is the {\bf three point normalization} method. Let $P_1$, $P_2$ and $P_3$ in $\p D^2$ and
three points in $\Gamma$ : $Q_1$, $Q_2$ and $Q_3$ in
the same order (with respect to the monotony given by definition~\ref{df-n.I.1}) and we introduce the following subspace of ${\mathcal C}_\Gamma$
\be
\label{n.z.I.1}
{\mathcal C}^\ast(\Gamma):=\lf\{u\in {\mathcal C}(\Gamma)\quad\mbox{s.t .} \quad\forall\, k=1,2,3\quad u(P_k)=Q_k\rg\}
\ee
The following elementary lemma, whose proof is left to the reader, garanties that
\be
\label{n.I.1}
\inf_{u\in{\mathcal C}(\Gamma)}E(u)=\inf_{u\in{\mathcal C}^\ast(\Gamma)}E(u)\quad.
\ee
\begin{Lm}
\label{lm-n.I.1}
Let $P_1$, $P_2$ and $P_3$ be 3 distinct points on $\p D^2$ indexed in a trigonometric order then there is a unique element $f\in{\mathcal M}^+(D^2)$, $$f(z)=e^{i\theta}\,(z-a)/(1-\ov{a}z)$$ such that
$f(e^{2\,i\,k\pi/3})=P_k$ for $k=1,2,3$.\hfill $\Box$
\end{Lm}

\medskip

We are now going to prove the {\bf closure of ${\mathcal C}^\ast(\Gamma)$ for the sequential weak $W^{1,2}$ topology}. This closure implies the existence of a minimizer
of $E$ in ${\mathcal C}^\ast(\Gamma)$ and hence in ${\mathcal C}(\Gamma)$ too. Precisely we have the following theorem.

\begin{Lm}
\label{lm-n.I.3}
For any positive constant $C\ge\inf_{{\mathcal C}^\ast(\Gamma)}E(u)$, the trace on $\p D^2$ of the subset of elements $u$ in ${\mathcal C}^\ast(\Gamma)$ satisfying $E(u)\le C$ is equicontinuous.  \hfill $\Box$
\end{Lm}
From an equicontinuous sequence one can always extract a subsequence that uniformly converges (Arzel\`a-Ascoli's theorem). Hence one deduces from this lemma the following corollary.
\begin{Co}
\label{co-n.I.1}
Let $u_n$ be a sequence in ${\mathcal C}^\ast(\Gamma)$ weakly converging to a map $u_\infty$ in $W^{1,2}$ then $u_\infty$ is continuous and monotone
on $\p D^2$.
\hfill $\Box$
\end{Co}
The main tool we shall use in order to prove the lemma~\ref{lm-n.I.3} is an argument  introduced in the framework
of the Plateau problem by R. Courant. The argument is based on Fubini theorem in order to  extract a ''good slice'', an arc of circle on which the energy is controlled and on which we can
apply Sobolev embedding which is better in this one dimensional curve than on the whole 2-dimensional disc.

\begin{Lm}
\label{lm-n.I.4} {\bf [Courant Lemma]}
Let $u\in W^{1,2}(D^2,{\R}^m)$ and let $p\in \p D^2$. For any $0<\delta<1$. Then there exists $\rho\in[\delta,\sqrt{\delta}]$ such that $$\nabla u\in L^2(\p B_\rho(p)\cap D^2)\quad,$$ and
\be
\label{n.I.18}
\begin{array}{l}
\ds\|u(x)-u(y)\|_{L^\infty((\p B_\rho(p)\cap D^2)^2)}^2\le\lf[\int_{\p B_\rho(p)\cap D^2}|\nabla u|\ d\theta\rg]^2\\[7mm]
\ds\quad\quad\quad\quad\quad\le \frac{4\pi}{\log\frac{1}{\delta}}\ \int_{D^2}|\nabla u |^2\ dx
\end{array}
\ee
where $d\theta$ is the length form on $\p B_\rho(p)\cap D^2$.\hfill $\Box$
\end{Lm}

\noindent{\bf Proof of lemma~\ref{lm-n.I.4}.}
Using Fubini theorem we have
\[
\begin{array}{l}
\ds\int_{D^2}|\nabla u|^2\ dx\ge\int_{D^2\cap(B_{\sqrt{\delta}}(p)\setminus B_\delta(p))}|\nabla u|^2\ dx\\[7mm]
\ds\quad\quad\quad\ge\int_{\delta}^{\sqrt{\delta}} d\rho\int_{D^2\cap\,\p B_{\rho}(p)}|\nabla u|^2\ d\theta\\[7mm]
\ds\quad\quad\quad\ge ess\inf\lf\{\rho\int_{D^2\cap\,\p B_{\rho}(p)}|\nabla u|^2\ d\theta\rg\}\ \int_{\delta}^{\sqrt{\delta}}\frac{d\rho}{\rho}
\end{array}
\]
Thus there exists a radius $\rho\in[\delta,\sqrt{\delta}]$ such that
\be
\label{n.I.19}
\rho\ \int_{D^2\cap\,\p B_{\rho}(p)}|\nabla u|^2\ d\theta\le\frac{2}{\log\frac{1}{\delta}}\int_{D^2}|\nabla u|^2\ dx
\ee
Cauchy-Schwartz inequality gives
\be
\label{n.I.20}
\lf[\int_{D^2\cap\,\p B_{\rho}(p)}|\nabla u|^2\ d\theta\rg]^2\le |D^2\cap\,\p B_{\rho}(p)|\ \int_{D^2\cap\,\p B_{\rho}(p)}|\nabla u|^2\ d\theta\quad.
\ee
Since $|D^2\cap\,\p B_{\rho}(p)|\le 2\pi\rho$, combining (\ref{n.I.19}) and (\ref{n.I.20}) gives (\ref{n.I.18}) and lemma~\ref{lm-n.I.4} is proved.\hfill $\Box$

\medskip

\noindent{\bf Proof of lemma~\ref{lm-n.I.3}.} Let $C>0$ satisfying
\[
C\ge\inf_{{\mathcal C}^\ast(\Gamma)}E(u)\quad.
\]
and denote
\[
{\mathcal C}^\ast_C(\Gamma):=\lf\{u\in {\mathcal C}^\ast(\Gamma)\quad\mbox{ s. t. }\quad E(u)\le C\rg\}\quad.
\]
The lemma is equivalent to the following claim.
\be
\label{n.I.21}
\begin{array}{l}
\ds\forall\,\ep>0\quad\exists\,\delta>0\quad\mbox{ s.t. }\quad\forall\, u\in {\mathcal C}^\ast_C(\Gamma)\quad\forall\,p,q\,\in \p D^2\\[5mm]
\ds \quad\quad |p-q|<\delta\quad\Longrightarrow\quad |u(p)-u(q)|<\ep
\end{array}
\ee
Since $\Gamma$ is the image of a \underbar{continuous} and  \underbar{injective} map $\gamma$ from $S^1$ into ${\R}^m$ the following
claim holds
\be
\label{n.I.22}
\begin{array}{l}
\ds\forall\,\ep>0\quad\exists\,\eta>0\quad\mbox{ s.t. }\quad\forall\ 0<\theta_1<\theta_2\le 2\pi\\[5mm]
\ds |\gamma(e^{i\theta_2})-\gamma(e^{i\theta_1})|<\eta\ \Longrightarrow\ \|\gamma(e^{i\theta})-\gamma(e^{i\theta_1})\|_{L^\infty([\theta_1,\theta_2])}<\ep
\end{array}
\ee
This claim can be proved by contradiction and we leave the details of the argument to the reader.

\medskip

We are heading now to the proof of (\ref{n.I.21}). Let $\ep>0$  such that
\be
\label{n.I.23}
2\ep<\inf_{i\ne j}|Q_i-Q_j|
\ee
where the $Q_i$ are the 3 fixed points on $\Gamma$ appearing in the definition (\ref{n.z.I.1}) of ${\mathcal C}^\ast(\Gamma)$. We are then considering $\ep$ small enough
in such a way that each ball of radius $\epsilon$ contains at most one $Q_i$.

\medskip

$\ep>0$ being fixed and satisfying (\ref{n.I.23}), we consider $\eta>0$ given by (\ref{n.I.22}).

\medskip

Consider $1>\delta>0$ to be fixed later but satisfying at least 
\be
\label{n.I.24}
2\sqrt{\delta}<\inf_{i\ne j}|P_i-P_j|
\ee

Take an arbitrary pair of points $p$ and $q$ in $\p D^2$ such that $|p-q|<\delta$. Let $p_0$ be the point on the small arc ${p\,q}\subset\p D^2$ right in the middle of this arc :
$|p-p_0|=|q-p_0|<\delta/2$. Consider $\rho\in[\delta,\sqrt{\delta}]$ given by the {\bf Courant lemma}~\ref{lm-n.I.4} and satisfying
\[
\|u(x)-u(y)\|_{L^\infty((\p B_\rho(p_0)\cap D^2)^2)}^2\le \frac{4\pi}{\log\frac{1}{\delta}}\ \int_{D^2}|\nabla u |^2\ dx
\]
Let $p'$ and $q'$ be the two points given by the intersection between $\p D^2$ and $\p B_\rho(p_0)$. Since $u$ is continuous \underbar{up to the boundary} we deduce that
\be
\label{n.I.25}
|u(p')-u(q')|\le \frac{4\pi}{\log\frac{1}{\delta}}\ \int_{D^2}|\nabla u |^2\ dx\le \frac{4\pi\ C}{\log\frac{1}{\delta}}
\ee
We fix now $\delta$ in such a way that
\[
\frac{4\pi\ C}{\log\frac{1}{\delta}}<\eta
\]
Because of (\ref{n.I.22}) one of the two arcs in $\Gamma$ connecting $u(p')$ and $u(q')$ has to be contained in a $B^m_\ep-$ball. Since $\ep$ has been chosen small
enough satisfying (\ref{n.I.23}) this arc connecting $u(p')$ and $u(q')$ and contained in a $B^m_\ep-$ball can contain at most one of the $Q_i$. In the mean time
$\delta$ has been chosen small enough in such a way that $B_{\sqrt{\delta}}(p_0)\cap \p D^2$ contains also at most one of the $P_i$, thus,
since $u$ in monotonic on $\p D^2$, the arc connecting $u(p')$ and $u(q')$ and contained in a $B^m_\ep-$ball has to be $u(\p B_\rho(p_0)\cap D^2)$ and we have than proved that
\[
|u(p)-u(q)|\le |u(p')-u(q')|<\ep\quad.
\]
This concludes the proof of claim (\ref{n.I.21}) and lemma~\ref{lm-n.I.3} is proved.\hfill $\Box$

\medskip

In order to establish that $E$ posses a minimizer in ${\mathcal C}^\ast(\Gamma)$ we are going to establish the following result.

\medskip

\begin{Lm}
\label{lm-n.I.5}
Let $u$ be the weak limit of a minimizing sequence of $E$ in ${\mathcal C}^\ast(\Gamma)$, then $u\in C^0(\ov{D^2},{\R}^m)$.\hfill $\Box$
\end{Lm}
Combining this lemma with the equicontinuity proved in lemma~\ref{lm-n.I.3} we obtain that the weak limit of a minimizing sequence satisfies
all the conditions for the membership in ${\mathcal C}^\ast(\Gamma)$ and then realizes a minimizer of (\ref{n.I.3}).

\medskip

\noindent{\bf Proof of lemma~\ref{lm-n.I.5}.}

From proposition~\ref{pr-n.I.0} we know that the weak limit $u$ of our minimizing sequence satisfies the Laplace equation and is therefore smooth and hence continuous
in the interior of $D^2$. From the equicontinuity proved in lemma~\ref{lm-n.I.3}, using Arzel\`a-Ascoli's theorem, we obtain moreover that the restriction of $u$ to the boundary of $D^2$
is continuous. It remains then to study the continuity of $u$ at the boundary while ''approaching'' this boundary from the inside of the disc.

\medskip

Let $u_n$ be a weakly converging minimizing sequence. Introduce $v_n\in W^{1,2}(D^2,{\R}^m)$ to be the unique solution of 
\[
\lf\{
\begin{array}{l}
\ds\Delta v_n=0\quad\quad\mbox{ in }\quad D^2\\[5mm]
\ds  v_n=u_n\quad\quad\quad\mbox{ on }\quad\p D^2
\end{array}
\rg.
\]
Observe that $v_n$ fulfills the conditions for the membership in ${\mathcal C}^\ast(\Gamma)$ and moreover the Dirichlet principle
implies 
\[
\int_{D^2}|\nabla v_n|^2\ dx\le \int_{D^2}|\nabla u_n|^2\ dx\quad.
\]
Thus $v_n$ is also a minimizing sequence. Since $u_n$ uniformly converges on the boundary $\p D^2$ to $u$, $v_n$ weakly converges
in $W^{1,2}$ to $v$ the unique solution of
\[
\lf\{
\begin{array}{l}
\ds\Delta v=0\quad\quad\mbox{ in }\quad D^2\\[5mm]
\ds  v=u\quad\quad\quad\mbox{ on }\quad\p D^2
\end{array}
\rg.
\]
This system is satisfied by $u$, hence $u=v$ and we have that the new minimizing sequence $v_n$ also converge to $u$. We replace then $u_n$ by $v_n$.
Since the maps $v_n$ satisfy the Laplace equation, the Maximum Principle implies that
\[
\begin{array}{l}
\forall k\,l\, \in {\N}\quad\forall i=1\cdots m\quad\quad\\[5mm]
\quad\quad\quad\|v_k^i-v^i_l\|_{L^\infty(D^2)}\le\|v_k^i-v^i_l\|_{L^\infty(\p D^2)}=\|u_k^i-u^i_l\|_{L^\infty(\p D^2)}\quad.
\end{array}
\]
Thus, the uniform convergence of $u_n$ to $u$ on $\p D^2$ implies the uniform convergence of $v_n$ to $u$ on $\p D^2$ and this implies that $u$
is in $C^0(\ov{D^2},{\R}^m)$ and  lemma~\ref{lm-n.I.5} is proved.\hfill $\Box$

\medskip

One point has not been addressed at this stage in order to prove theorem~\ref{th-I.0}, this is the question to know whether ${\mathcal C}(\Gamma)$ is empty or not
while assuming $\Gamma$ to be a rectifiable closed curve only. We are going to answer positively to this question.

\medskip

\subsection{Existence of Parametric disc extensions to Jordan rectifiable curves.}

One point has not been addressed at this stage in order to prove theorem~\ref{th-I.0}, this is the question to know whether ${\mathcal C}(\Gamma)$ is empty or not
while assuming $\Gamma$ to be a rectifiable closed curve only. We are going to answer positively to this question.

\begin{Th}
\label{lm-n.I.6}
Let $\Gamma$ be a closed Jordan rectifiable curve then ${\mathcal C}(\Gamma)$ is not empty.\hfill $\Box$
\end{Th}
Our approach consists in approximating $\Gamma$ by smooth curves and pass to the limit but before to do so we shall first establish an isoperimetric
inequality for smooth conformal parametrization.

\newpage

\section{Conformally invariant coercice Lagrangians with quadratic growth, in dimension 2.}

\reset

The resolution of the Plateau problem proposed by Douglas and Rad\'o is an example of the use of a conformal invariant Lagrangian $E$ to approach an ``extrinsic" problem:  minimizing the area of a disk with fixed boundary. The analysis of this problem was eased by the high simplicity of the equation (\ref{I.1}) satisfied by the critical points of $E$. It is the Laplace equation. Hence, questions related to unicity, regularity, compactness, etc... can be handled with a direct application of the maximum principle. In the coming three chapters, we will be concerned with analogous problems (in particular regularity issues) related to the critical points to general conformally invariant, coercive Lagrangians with quadratic growth. As we will discover, the maximum principle no longer holds, and one must seek an alternate way to compensate this lack. The conformal invariance of the Lagrangian will generate a very peculiar type of nonlinearities in the corresponding Euler-Lagrange equations. We will see how the specific structure of these nonlinearities enable one to recast the equations in divergence form. This new formulation, combined to the results of {\it integration by compensation}, will provide the substrate to understanding a variety of problems, such as Willmore surfaces, poly-harmonic and $\al$-harmonic maps, Yang-Mills fields, Hermitte-Einstein equations, wave maps, etc...

\medskip

We consider a Lagrangian of the form
\be
\label{II.1}
L(u)=\int_{D^2}l(u,\nabla u)\ dx\ dy\quad,
\ee
where the integrand $l$ is a function of the variables $z\in {\R}^m$ and $p\in {\R}^2\otimes{\R}^m$, which satisfy the following coercivity and ``almost quadratic" conditions in $p$\,:
\be
\label{II.2}
C^{-1}\ |p|^2\le l(z,p)\le C\ |p|^2\quad ,
\ee
We further assume that $L$ is conformally invariant: for each positive conformal transformation $f$ of degree 1, and for each map $u\in W^{1,2}(D^2,{\R}^m)$, there holds
\be
\label{II.3}
\begin{array}{rl}
\ds L(u\circ f)&\ds=\int_{f^{-1}(D^2)}l(u\circ f,\nabla(u\circ f))\ dx'\ dy'\\[5mm]
 &\ds=\int_{D^2}l(u,\nabla u)\ dx\ dy=L(u)\quad .
 \end{array}
\ee
{\bf Example 1.} The Dirichlet energy described in the Introduction,
\[
E(u)=\int_{D^2}|\nabla u|^2\ dx\, dy\quad,
\]
whose critical points satisfy the Laplace equation (\ref{I.1}), which, owing to the conformal hypothesis, geometrically describes {\it minimal surfaces}. Regularity and compactness matters relative to this equation are handled with the help of the maximum principle.
 
 \medskip
 
 \noindent {\bf Example 2.} Let an arbitrary in $\R^m$ be given, namely $(g_{ij})_{i,j\in {\N}_m}\in C^1({\R}^m, {\mathcal S}^+_m)$, where ${\mathcal S}^+_m$ denotes the subset of $M_m({\R})$, comprising the symmetric positive definite $m\times m$ matrices. We make the following uniform coercivity and boundedness hypothesis:
 \[
\exists\:\: C>0\quad\quad \mbox{ s. t. }\quad C^{-1}\delta_{ij}\le g_{ij}\le C\delta_{ij}\quad\mbox{ on }{\R}^m.
\]
Finally, we suppose that
\[
\|\nabla g\|_{L^{\infty}({\R}^m)} <+\infty\quad.
\]
With these conditions, the second example of quadratic, coercive, conformally invariant Lagrangian is
\[
\begin{array}{rl}
E_g(u)&\ds=\frac{1}{2}\int_{D^2}\lf<\nabla u,\nabla u\rg>_g\ dx\,dy\\[5mm]
 &\ds=\frac{1}{2}\int_{D^2}\sum_{i,j=1}^mg_{ij}(u)\nabla u^i\cdot\nabla u^j\ dx\,dy\quad .
 \end{array}
\]
Note that Example 1 is contained as a particular case.\\
Verifying that $E_g$ is indeed conformally invariant may be done analogously to the case of the Dirichlet energy, via introducing the complex variable $z=x+iy$. No new difficulty arises, and the details are left to the reader as an exercise.\\
The weak critical points of $E_g$ are the functions $u\in W^{1,2}(D^2,{\R}^m)$ which satisfy
\[
\forall\xi\in C^\infty_0(D^2,{\R}^m)\quad\quad\frac{d\ }{dt}E_g(u+t\xi)_{|_{t=0}}=0\quad.
\]
An elementary computation reveals that $u$ is a weak critical point of $E_g$ if and only if the following Euler-Lagrange equation holds in the sense of distributions:
\be
\label{II.4}
\forall i=1\cdots m\quad\quad \Delta u^i+\sum_{k,l=1}^m\Gamma^i_{kl}(u)\,\nabla u^k\cdot\nabla u^l=0\quad.
\ee
Here, $\Gamma_{kl}^i$ are the Christoffel symbols corresponding to the metric $g$, explicitly given by
\be
\label{II.4za}
\Gamma^i_{kl}(z)=\frac{1}{2}\sum_{s=1}^mg^{is}\lf(\p_{z_l}g_{km}+\p_{z_k}g_{lm}-\p_{z_m}g_{kl}\rg)\quad,
\ee
where $(g^{ij})$ is the inverse matrix of $(g_{ij})$.

Equation (\ref{II.4}) bears the name {\it harmonic map equation}\footnote{One way to interpret (\ref{II.4}) as the two-dimensional equivalent of the geodesic equation in normal parametrization,
\[
\frac{d^2x^i}{dt^2}+\sum_{k,l=1}^m\Gamma_{kl}^i\frac{dx^k}{dt}\,\frac{dx^l}{dt}=0\quad .
\]} with values in $({\R}^m,g)$.

Just as in the flat setting, if we further suppose that $u$ is conformal, then (\ref{II.4}) is in fact equivalent to $u(D^2)$ being a minimal surface in $({\R}^m,g)$.

We note that $\Gamma^i(\nabla u,\nabla u):=\sum_{k,l=1}^m\Gamma^i_{kl}\nabla u^k\cdot\nabla u^l$, so that the harmonic map equation can be recast as
\be
\label{II.4b}
\Delta u+\Gamma(\nabla u,\nabla u)=0\quad .
\ee
This equation raises several analytical questions:
\begin{itemize}
\item[(i)] {\bf Weak limits} : Let $u_n$ be a sequence of solutions of (\ref{II.4b}) with uniformly bounded energy $E_g$. Can one extract a subsequence converging weakly in $W^{1,2}$ to a harmonic map ?

\item[(ii)] {\bf Palais-Smale sequences} :
 Let $u_n$ be a sequence of solutions of (\ref{II.4b}) in $W^{1,2}(D^2,{\R}^m)$ with uniformly bounded energy $E_g$, and such that
 \[
\Delta u_n+\Gamma(\nabla u_n,\nabla u_n)=\delta_n\rightarrow 0\quad\quad\mbox{ strongly in }
H^{-1}\quad.
\]
 Can one extract a subsequence converging weakly in $W^{1,2}$ to a harmonic map ?

\item[(iii)] {\bf Regularity of weak solutions} : Let $u$ be a map in $W^{1,2}(D^2,{\R}^m)$ which satisfies (\ref{II.4}) distributionally. How regular is $u$ ? Continuous, smooth, analytic, etc...
\end{itemize}

The answer to (iii) is strongly tied to that of (i) and (ii). We shall thus restrict our attention in these notes on regularity matters.

\medskip

Prior to bringing into light further examples of conformally invariant Lagrangians, we feel worthwhile to investigate deeper the difficulties associated with the study of the regularity of harmonic maps in two dimensions.

\medskip
The harmonic map equation (\ref{II.4b}) belongs to the class of elliptic systems with {\it quadratic growth}, also known as {\it natural growth}, of the form
\be
\label{II.6}
\Delta u=f(u,\nabla u)\quad ,
\ee
where $f(z,p)$ is an arbitrary continuous function for which there exists constants $C_0>0$ and $C_1>0$ satisfying
\be
\label{II.7}
\forall z\in {\R}^m\quad\forall p\in {\R}^2\otimes{\R}^m\quad\quad f(z,p)\le C_1|p|^2+C_0\quad.
\ee
In dimension two, these equations are critical for the Sobolev space $W^{1,2}$. Indeed,
\[
u\in W^{1,2}\Rightarrow\Gamma(\nabla u, \nabla u)\in L^1\Rightarrow \nabla u\in L^p_{loc}(D^2)\quad\forall p<2\quad.
\]
In other words, from the regularity standpoint, the demand that $\nabla u$ be square-integrable provides the information that\footnote{Actually, one can show that $\nabla u$ belongs to the weak-$L^2$ Marcinkiewicz space $L_{loc}^{2,\infty}$ comprising those measurable functions $f$ for which
\be
\label{II.5}
\sup_{\la>0}\la^2\ \lf|\{p\in\omega\ ;\ |f(p)|>\la\}\rg|<+\infty\quad,
\ee
where $|\cdot|$ is the standard Lebesgue measure. Note that $L^{2,\infty}$ is a slightly larger space than $L^2$. However, it possesses the same scaling properties.} $\nabla u$ belongs to $L^p_{loc}$ for all $p<2$. We have thus lost a little bit of information! Had this not been the case, the problem would be ``boostrapable", thereby enabling a successful study of the regularity of $u$. Therefore, in this class of problems, the main difficulty lies in the aforementioned slight loss of information, which we symbolically represent by $L^2\rightarrow L^{2,\infty}$.\\

There are simple examples of equations with quadratic growth in two dimensions for which the answers to the questions (i)-(iii) are all negative. Consider\footnote{This equation is conformally invariant. However, as shown by J. Frehse \cite{Fre}, it is also the Euler-Lagrange equation derived from a Lagrangian which is {\it not} conformally invariant:
\[
L(u)=\int_{D^2}\lf(1+\frac{1}{1+e^{12\, u}\ (\log1/|(x,y)|)^{-12}}\rg)\ |\nabla u|^2(x,y)\ dx\  dy\quad.
\]
}
\be
\label{II.8}
\Delta u+|\nabla u|^2=0\quad.
\ee
This equation has quadratic growth, and it admits a solution in $W^{1,2}(D^2)$ which is unbounded in $L^\infty$, and thus discontinuous. It is explicitly given by
\[
u(x,y):=\log\log\frac{2}{\sqrt{x^2+y^2}}\quad.
\]
The regularity issue can thus be answered negatively. Similarly, for the equation (\ref{II.8}), it takes little effort to devise counter-examples to the weak limit question (i), and thus to the question (ii). To this end, it is helpful to observe that $C^2$ maps obey the general identity
\be
\label{II.8b}
\Delta e^u=e^u\ \lf[\Delta u+|\nabla u|^2\rg]\quad .
\ee
One easily verifies that if $v$ is a positive solution of
$$\Delta v=-2\pi\sum_i\la_i\,\delta_{a_i}\quad,$$  where
$\la_i>0$ and $\delta_{a_i}$ are isolated Dirac masses, then $u:=\log{v}$ provides a solution\footnote{Indeed, per (\ref{II.8b}), we find $\Delta u+|\nabla u|^2=0$ away from the points $a_i$. Near these points, $\nabla u$ asymptotically behaves as follows:
\[
|\nabla u|=|v|^{-1}\,|\nabla v|\simeq \big(|(x,y)-a_i|\ \log|(x,y)-a_i|\big)^{-1}\ \in L^2\quad.
\]
Hence, $|\nabla u|^2\in L^1$, so that $\Delta u+|\nabla u|^2$ is a distribution in $H^{-1}+L^{1}$ supported on the isolated points $a_i$. From this, it follows easily that 
\[
\Delta u+|\nabla u|^2=\sum_i\mu_i\ \delta_{a_i}\quad .
\]
Thus, $\Delta u$ is the sum of an $L^1$ function and of Dirac masses. But because $\Delta u$ lies in $H^{-1}$, the coefficients $\mu_i$ must be zero. Accordingly, $u$ does belong to
$W^{1,2}$.} in $W^{1,2}$ of (\ref{II.8}). 
We then select a strictly positive regular function $f$ with integral equal to 1, and supported on the ball of radius $1/4$ centered on the origin. There exists a sequence of atomic measures with positive weights $\la_i^n$ such that
\be
 f_n=\sum_{i=1}^{n}\la_i^n\ \delta_{a_i^n}\qquad\textnormal{and}\qquad \sum_{i=1}^n\la_i^n=1\quad,
\ee
which converges as Radon measures to $f$. We next introduce
\[
u_n(x,y):=\log\lf[\sum_{i=1}^n\la_i^n\ \log\frac{2}{|(x,y)-a_i^n|}\rg]\quad.
\]
On $D^2$, we find that
\be
\label{II.8c}
v_n=\sum_{i=1}^n\la_i^n\ \log\frac{2}{|(x,y)-a_i^n|}>\sum_{i=1}^n\la_i^n\ \log\frac{8}{5}=\log\frac{8}{5}\quad.
\ee
On the other hand, there holds
\[
\begin{array}{l}
\ds\int_{D^2}|\nabla u_n|^2=-\int_{D^2}\Delta u_n=-\int_{\p D^2}\frac{\p u_n}{\p r}\\[5mm]
\ds\quad\le
\int_{\p D^2} \frac{|\nabla v_n|}{|v_n|}\le \frac{1}{\log\frac{8}{5}}\int_{\p D^2}|\nabla v_n|\le C
\end{array}
\]
for some constant $C$ independent of $n$ . Hence, $(u_n)_n$ is a sequence of solutions to  (\ref{II.8}) uniformly bounded in $W^{1,2}$. Since the sequence $(f_n)$ converges as Radon measures to $f$, it follows that for any $p<2$, the sequence $(v_n)$ converges strongly in $W^{1,p}$ to
\[
v:=\log\frac{2}{r}\ast f\quad.
\]
The uniform upper bounded (\ref{II.8c}) paired to the aforementioned strong convergence shows that for each $p<2$, the sequence $u_n=\log v_n$ converges strongly in $W^{1,p}$ to
\[
u:=\log\lf[\log\frac{2}{r}\ast f\rg]
\]
From the hypotheses satisfied by $f$, we see that $\Delta(e^u)=-2\pi\  f\ne 0$. As $f$ is regular, so is thus $e^u$, and therefore, owing to (\ref{II.8b}), $u$ cannot fullfill (\ref{II.8}).

Accordingly, we have constructed a sequence of solutions to (\ref{II.8}) which converges weakly in $W^{1,2}$ to a map that is {\it not} a solution to (\ref{II.8}).

\medskip

\noindent{\bf Example 3.}
We consider a map $(\om_{ij})_{i,j\in{\N}_m}$ in 
$C^1({\R}^m,so(m))$, where $so(m)$ is the space antisymmetric square $m\times m$ matrices. We impose the following uniform bound
$$
\|\nabla \om\|_{L^\infty(D^2)}<+\infty\quad .
$$
For maps $u\in W^{1,2}(D^2,\R^m)$, we introduce the Lagrangian
\be
\label{II.9}
E^\om(u)=\frac{1}{2}\int_{D^2}|\nabla u|^2+\sum_{i,j=1}^m\om_{ij}(u)\p_xu^i\p_yu^j-\p_yu^i\p_xu^j\ dx\,dy
\ee
The conformal invariance of this Lagrangian arises from the fact that $E^\om$ is made of the conformally invariant Lagrangian $E$ to which is added the integral over $D^2$ of the 2-form $\omega=\omega_ij dz^i\wedge dz^j$ {\it pulled back} by $u$. Composing $u$ by an arbitrary positive diffeomorphism of $D^2$ will not affect this integral, thereby making $E^\om$ into a conformally invariant Lagrangian.\\
The Euler-Lagrange equation deriving from (\ref{II.9}) for variations of the form $u+t\xi$, where $\xi$ is an arbitrary smooth function with compact support in $D^2$, is found to be 
\be
\label{II.10}
\Delta u^i-2\sum_{k,l=1}^m H^i_{kl}(u)\ \nabla^\perp u^k\cdot\nabla u^l=0\qquad\forall\:\:i=1,\dots,m.
\ee
Here, $\nabla^\perp u^l=(-\p_yu^k,\p_xu^k)$ \footnote{in our notation, $\nabla^\perp u^k\cdot\nabla u^l$ is the Jacobian
\[
\nabla^\perp u^k\cdot\nabla u^l=\p_xu^k\p_yu^l-\p_yu^k\p_xu^l\quad.
\]
} 
while $H^i_{kl}$ is antisymmetric in the indices $k$ et $l$. It is the coefficient of the $\R^m$-valued two-form $H$ on $\R^m$ 
\[
H^i(z):=\sum_{k,l=1}^mH^i_{kl}(z)\ dz^k\wedge dz^l\quad .
\]
The form $H$ appearing in the Euler-Lagrange equation (\ref{II.10}) is the unique solution of
\[
\begin{array}{rl}
\forall z\in{\R}^m\quad\forall U,V,W\in {\R}^m &\ \\[5mm]
 \ds d\om_z(U,V,W)&\ds=4\,U\cdot H(V,W)\\[5mm]
 &\ds=4\sum_{i=1}^mU^i\,
H^i(V,W)\quad.
\end{array}
\]
For instance, in dimension three, $d\om$ is a 3-form which can be identified with a function on ${\R}^m$. More precisely, there exists
$H$ such that $d\om=4H\ dz^1\wedge dz^2\wedge dz^3$. In this notation (\ref{II.10}) may be recast, for each $i\in\{1,\ldots,m\}$, as
\be
\label{II.11}
\Delta u^i=2H(u)\ \p_xu^{i+1}\p_yu^{i-1}-\p_xu^{i-1}\p_yu^{i+1}\quad ,
\ee
where the indexing is understood in ${\Z}_3$. The equation (\ref{II.11}) may also be written
$$
\Delta u=2H(u)\ \p_xu\times \p_yu\quad ,
$$
which we recognize as (\ref{I.2}), the {\it prescribed mean curvature equation.} 

In a general fashion, the equation (\ref{II.10}) admits the following geometric interpretation. Let $u$ be a conformal solution of (\ref{II.10}), so that $u(D^2)$ is a surface whose mean curvature vector at the point $(x,y)$ is given by
\be
\label{II.12}
e^{-2\la}\ u^\ast H=\lf(e^{-2\la}\ \sum_{k,l=1}^m H^i_{kl}(u)\ \nabla^\perp u^k\cdot\nabla u^l\rg)_{i=1\cdots m}\quad,
\ee
where $e^\la$ is the conformal factor $e^\la=|\p_x u|=|\p_y u|$.
As in Example 2, the equation (\ref{II.10}) forms an elliptic system with quadratic growth, thus critical in dimension two for the $W^{1,2}$ norm. The analytical difficulties relative to this nonlinear system are thus, {\it a priori}, of the same nature as those arising from the {\it harmonic map equation}.

\medskip

\noindent{\bf Example 4.}
In this last example, we combine the settings of Examples 2 and 3 to produce a mixed problem. Given on $\R^m$ a metric $g$ and a two-form $\om$, both $C^1$ with uniformly bounded Lipschitz norm, consider the Lagrangian
\[
E^\om_g(u)=\frac{1}{2}\int_{D^2}\lf<\nabla u,\nabla u\rg>_g\ dx\,dy+u^\ast\om\quad.
\]
As before, it is a coercive conformally invariant Lagrangian with quadratic growth. Its critical points satisfy the Euler-Lagrangian equation
\be
\label{II.13}
\Delta u^i+\sum_{k,l=1}^m\Gamma^i_{kl}(u)\nabla u^k\cdot\nabla u^l
-2\sum_{k,l=1}^mH^i_{kl}(u)\nabla^\perp u^k\cdot\nabla u^l=0\quad,
\ee
for $i=1\cdots m$.\\
Once again, this elliptic system admits a geometric interpretation which generalizes the ones from Examples 2 and 3. Whenever a conformal map $u$ satisfies (\ref{II.13}), then $u(D^2)$ is a surface in $(\R^m,g)$ whose mean curvature vector is given by (\ref{II.12}). 
The equation (\ref{II.13}) also forms an elliptic system with quadratic growth, and critical in dimension two for the $W^{1,2}$ norm.

\medskip
Interestingly enough, M. Gr\"uter showed that {\it any} coercive conformally invariant Lagrangian with quadratic growth is of the form $E^\om_g$ for some appropriately chosen $g$ and $\om$. 

\begin{Th} 
\label{th-I.1}
\cite{Gr}
Let $l(z,p)$ be a real-valued function on ${\R}^m\times {\R}^2\otimes{\R}^m$, which is $C^1$ in its first variable and $C^2$ in its second variable. Suppose that $l$ obeys the coercivity and quadratic growth conditions
\be
\label{Ia.13}
\begin{array}{l}
\ds\exists C>0\quad\mbox{t.q.}\quad\forall z\in{\R}^m\quad\forall p\in {\R}^{2}\otimes{\R}^m\\[5mm]
\quad C^{-1}|p|^2\le l(X,p)\le C|p|^2\quad.
\end{array}
\ee
Let $L$ be the Lagrangian
\be
\label{Ia.14}
{ L}(u)=\int_{D^2} l(u,\nabla u)(x,y)\ dx\ dy
\ee
acting on $W^{1,2}(D^2,\R^m)$-maps $u$. We suppose that $L$ is conformally invariant: for every conformal application $\phi$ positive and of degree 1, there holds
\be
\label{Ia.15}
{L}(u\circ\phi)=\int_{\phi^{-1}(D^2)}l(u\circ\phi,\nabla (u\circ\phi))(x,y)\ dx\ dy
={ L}(u)\quad .
\ee
Then there exist on $\R^m$ a $C^1$ metric $g$ and a $C^1$ two-form $\omega$ such that
\be
\label{Ia.16}
{L}=E^\omega_g\quad.
\ee
\end{Th}

\medskip

\noindent{\bf Maps taking values in a submanifold of $\R^m$.}

\medskip

Up to now, we have restricted our attention to maps from $D^2$ into a manifold with only one chart $({\R}^n,g)$. More generally, it is possible to introduce the Sobolev space $W^{1,2}(D^2,N^n)$, where $(N^n,g)$ is an oriented $n$-dimensional $C^2$-manifold. When this manifold is compact without boundary (which we shall henceforth assume, for the sake of simplicity), a theorem by Nash guarantees that it can be isometrically immersed into Euclidean space $\R^m$, for $m$ large enough. We then define
\[
W^{1,2}(D^2, N^n):=\lf\{u\in W^{1,2}(D^2,{\R}^m)\ ;\ u(p)\in N^n\ \mbox{ a.e. }p\in D^2\rg\}
\]
Given on $N^n$ a $C^1$ two-form $\om$, we may consider the Lagrangian
\be
\label{Ia.17}
E^\om(u)=\frac{1}{2}\int_{D^2}|\nabla u|^2\ dx\ dy+u^\ast\om
\ee
acting on maps $u\in W^{1,2}(D^2,N^n)$. The critical points of $E^\om$ are defined as follows. Let $\pi_N$ be the orthogonal projection on $N^n$ which to each point in a neighborhood of $N$ associates its nearest orthogonal projection on $N^n$. For points sufficiently close to $N$, the map $\pi_N$ is regular. We decree that $u\in W^{1,2}(D^2,N^n)$ is a critical point of $E^\om$ whenever there holds
\be
\label{Ia.18}
\frac{d}{dt}E^\om(\pi_N(u+t\xi))_{t=0}=0\quad,
\ee
for all $\xi\in C^\infty_0(D^2,{\R}^m)$.\\
It can be shown\footnote{in codimension 1, this is done below.} that (\ref{Ia.18}) is satisfied by $u\in C^\infty_0(D^2,{\R}^m)$ if and only if $u$ obeys the Euler-Lagrange equation
\be
\label{Ia.19}
\Delta u+A(u)(\nabla u,\nabla u)= H(u)(\nabla ^\perp u,\nabla u)\quad,
\ee
where $A\,(\equiv A_z)$ is the second fundamental form at the point $z\in N^n$ corresponding to the immersion of $N^n$ into $\R^m$. To a pair of vectors in $T_zN^n$, the map $A_z$ associates a vector orthogonal to $T_zN^n$. In particular, at a point $(x,y)\in D^2$, the quantity $A_{(x,y)}(u)(\nabla u, \nabla u)$ is the vector of ${\R}^m$ given by
\[
A_{(x,y)}(u)(\nabla u,\nabla u):=A_{(x,y)}(u)(\p_xu,\p_xu)+A_{(x,y)}(u)(\p_yu,\p_yu)\quad.
\]
For notational convenience, we henceforth omit the subscript $(x,y)$. \\ 
Similarly, $H(u)(\nabla^\perp u,\nabla u)$ at the point $(x,y)\in D^2$ is the vector in $\R^m$ given by
\[
\begin{array}{rl}
\ds H(u)(\nabla^\perp u,\nabla u)&\ds:=H(u)(\p_x u,\p_y u)-H(u)(\p_y u,\p_x u)\\[5mm]
   &\ds=2H(u)(\p_xu,\p_yu)\quad,
\end{array}
\]
where $H\,(\equiv H_z)$ is the $T_zN^n$-valued alternating two-form on $T_zN^n$\,: 
\[
\forall\  U,V,W\ \in T_zN^n\quad\quad d\om(U,V,W):=U\cdot H_z(V,W)\quad.
\]
Note that in the special case when $\om=0$, the equation (\ref{Ia.19}) reduces to\be
\label{Ia.20}
\Delta u+A(u)(\nabla u,\nabla u)=0\quad ,
\ee
which is known as the {\it $N^n$-valued harmonic map equation}.

\medskip

We now establish (\ref{Ia.19}) in the codimension 1 case. Let $\nu$ be the normal unit vector to $N$. The form $\om$ may be naturally extended on a small neighborhood of $N^n$ via the pull-back $\pi_N^\ast\om$ of the projection $\pi_N$. Infinitesimally, to first order, considering variations for $E^\om$ of the form $\pi_N(u+t\xi)$ is tantamount to considering variations of the kind $u+t\,d\pi_N(u)\xi$, which further amounts to focusing on variations of the form $u+tv$, where $v\in W^{1,2}(D^2,{\R}^m)\cap L^\infty$ satisfies $v\cdot \nu(u)=0$ almost everywhere.\\
Following the argument from Example 3, we obtain that $u$ is a critical point of $E^\om$ whenever for all $v$ with $v\cdot \nu(u)=0$ a.e., there holds 
\[
\int_{D^2}\sum_{i=1}^m\lf[\Delta u^i-2\sum_{k,l=1}^m H^i_{kl}(u)\ \nabla^\perp u^k\cdot\nabla u^l\rg]\ v^i\ dx\ dy=0\quad,
\]
where $H$ is the vector-valued two-form on $\R^m$ given for $z$ on $N^n$ by 
\[
\forall\  U,V,W\ \in {\R}^m\quad\quad d\pi_N^\ast\om(U,V,W):=U\cdot H_z(V,W)\quad.
\]
In the sense of distributions, we thus find that
\be
\label{Ia.21}
\lf[\Delta u-H(u)(\nabla^\perp u,\nabla u)\rg]\wedge\nu(u)=0\quad.
\ee
Recall, $\nu\circ u\in L^\infty\cap W^{1,2}(D^2,{\R}^m)$. Accordingly (\ref{Ia.21}) does indeed make sense in ${\mathcal D}'(D^2)$.

Note that if any of the vectors $U$, $V$, and $W$ is normal to $N^n$, i.e. parallel to $\nu$, then $d\pi_N^\ast\om(U,V,W)=0$, so that
\[
\nu_z\cdot H_z(V,W)=0\qquad\forall\:\:V\,,W\,\in\,\R^m\:.
\]
Whence,
\be
\label{Ia.22}
\begin{array}{l}
\ds\lf[\Delta u-H(u)(\nabla^\perp u,\nabla u)\rg]\cdot\nu(u)=\Delta u\cdot\nu(u)\\[5mm]
 \ds=div(\nabla u\cdot \nu(u))-\nabla u\cdot\nabla(\nu(u))=-\nabla u\cdot\nabla(\nu(u))
 \end{array}
 \ee
 where we have used the fact that $\nabla u\cdot \nu(u)=0$ holds almost everywhere, since $\nabla u$ is tangent to $N^n$. \\

Altogether, (\ref{Ia.21}) and (\ref{Ia.22}) show that $u$ satisfies in the sense of distributions the equation
\be
\label{Ia.23}
\Delta u-H(u)(\nabla^\perp u,\nabla u)=-\nu(u)\ \nabla(\nu(u))\cdot\nabla u\quad .
\ee
In codimension 1, the second fundamental form acts on a pair of vectors $(U,V)$ in $T_zN^n$ via
\be
\label{Ia.24}
A_z(U,V)=\nu(z)\ <d\nu_zU,V>\quad ,
\ee
so that, as announced, (\ref{Ia.23}) and (\ref{Ia.19}) are identical.

\newpage

We close this section by stating a conjecture formulated by Stefan d'Hildebrandt in the late 1970s. 

\begin{Con} \cite{Hil} \cite{Hil2}
The critical points with finite energy of a coercive conformally invariant Lagrangian with quadractic growth are H\"older continuous.
\end{Con}

The remainder of these lecture notes shall be devoted to establishing this conjecture. Although its resolution is closely related to the compactness questions (i) and (ii) previously formulated on page 9, for lack of time, we shall not dive into the study of this point.\\
Our proof will begin by recalling the first partial answers to Hildebrandt's conjecture provided by H. Wente and F. H\'elein, and the importance in their approach of the r\^ole played by {\it conservations laws} and {\it integration by compensation}.\\
Then, in the last section, we will investigate the theory of linear elliptic systems with antisymmetric potentials, and show how to apply it to the resolution of Hildebrandt's conjecture.

\newpage

\section{Integrability by compensation theory applied to some conformally invariant Lagrangians}
\reset

\subsection{Constant mean curvature equation (CMC)}

Let $H\in{\R}$ be constant. We study the analytical properties of solutions in $W^{1,2}(D^2,\R^3)$ of the equation
\be
\label{III.1}
\Delta u-2H\ \p_xu\times\p_y u=0\quad.
\ee
The Jacobian structure of the right-hand side enable without much trouble, inter alia, to show that {\bf Palais-Smale sequences} converge weakly:

Let $F_n$ be a sequence of distributions converging to zero in $H^{-1}(D^2,{\R}^3)$, and let $u_n$ be a sequence of functions uniformly bounded in $W^{1,2}$ and satisfying the equation
\[
\Delta u_n-2H\ \p_xu_n\times \p_yu_n=F_n\rightarrow 0\mbox{ strongly in }H^{-1}(D^2)\quad.
\]
We use the notation 
\be
\label{III.2}
\begin{array}{l}
(\p_xu_n\times\p_yu_n)^i=\p_x u^{i+1}_n\p_y u^{i-1}_n-\p_xu^{i-1}_n\p_yu^{i+1}_n\\[5mm]
\quad=\p_x(u^{i+1}_n\ \p_yu^{i-1}_n)-\p_y(u^{i+1}_n\p_xu^{i-1}_n)\quad.
\end{array}
\ee
The uniform bounded on the $W^{1,2}$-norm of $u_n$ enables the extraction of a subsequence $u_{n'}$ weakly converging in $W^{1,2}$ to some limit $u_\infty$. With the help of the Rellich-Kondrachov theorem, we see that the sequence $u_n$ is strongly compact in $L^2$. In particular, we can pass to the limit in the following quadratic terms
$$
u^{i+1}_n\ \p_yu^{i-1}_n\rightarrow u^{i+1}_\infty\ \p_yu^{i-1}_\infty\quad \quad\mbox{ in }{\mathcal D}'(D^2)\quad
$$
and
$$
u^{i+1}_n\ \p_xu^{i-1}_n\rightarrow u^{i+1}_\infty\ \p_xu^{i-1}_\infty\quad \quad\mbox{ in }{\mathcal D}'(D^2)\quad.
$$
Combining this to (\ref{III.2}) reveals that $u_\infty$ is a solution of the CMC equation (\ref{III.1}).

\medskip

Obtaining information on the regularity of weak $W^{1,2}$ solutions of the CMC equation (\ref{III.2}) requires some more elaborate work. More precisely, a result from the theory of integration by compensation due to H. Wente is needed.

\begin{Th}
\label{th-III.1}
\cite{We}
Let $a$ and $b$ be two functions in $W^{1,2}(D^2)$, and let $\phi$ be the unique solution in $W^{1,p}_0(D^2)$ - for $1\le p<2$ - of the equation
\be
\label{III.3}
\lf\{
\begin{array}{l}
\ds-\Delta\phi=\p_xa\,\p_yb-\p_xb\,\p_ya\quad\quad\mbox{ in }D^2\\[5mm]
\ds\varphi=0\quad\quad\quad\mbox{ on }\p D^2\quad.
\end{array}
\rg.
\ee
Then $\phi$ belongs to $C^0\cap W^{1,2}(D^2)$ and
\be
\label{III.4}
\|\phi\|_{L^\infty(D^2)}+\|\nabla\phi\|_{L^2(D^2)}\le C_0\ \|\nabla a\|_{L^2(D^2)}\ \|\nabla b\|_{L^2(D^2)}\quad .
\ee
where $C_0$ is a constant independent of $a$ and $b$.\footnote{Actually, one shows that theorem~\ref{th-III.1} may be generalized to arbitrary oriented Riemannian surfaces, with a constant $C_0$ \underbar{independent of the surface}, which is quite a remarkable and useful fact. For more details, see \cite{Ge} and \cite{To}.}
\hfill$\Box$
\end{Th}
{\bf Proof of theorem~\ref{th-III.1}.}
We shall first assume that $a$ and $b$ are smooth, so as to legitimize the various manipulations which we will need to perform. The conclusion of the theorem for general $a$ and $b$ in $W^{1,2}$ may then be reached through a simple density argument. In this fashion, we will obtain the continuity of $\phi$ from its being the uniform limit of smooth functions.

\medskip

Observe first that integration by parts and a simple application of the Cauchy-Schwarz inequality yields the estimate
\[
\begin{array}{rl}
\ds\int_{D^2}|\nabla\phi|^2=-\int_{D^2}\phi\,\Delta\phi&\ds\le \|\phi\|_\infty \ \|\p_xa\,\p_yb-\p_xb\,\p_ya\|_1\\[5mm]
 &\le2\,\|\phi\|_\infty \|\nabla a\|_2\ \|\nabla b\|_2\quad.
 \end{array}
\]
Accordingly, if $\phi$ lies in $L^\infty$, then it automatically lies in $W^{1,2}$.

{\bf {Step 1}.} given two functions $\ti{a}$ and $\ti{b}$ in $C^\infty_0({\C})$, which is dense in $W^{1,2}({\C})$, we first establish the estimate (\ref{III.4}) for
\be
\label{III.4a}
\ti{\phi}:=\frac{1}{2\pi}\log\, \frac{1}{r}\ast \lf[\p_x\ti{a}\,\p_y\ti{b}-\p_x\ti{b}\,\p_y\ti{a}\rg]\quad .
\ee
Owing to the translation-invariance, it suffices to show that
\be
\label{III.5}
|\ti{\phi}(0)|\le C_0\ \|\nabla\ti{a}\|_{L^2({\C})}\ \|\nabla\ti{b}\|_{L^2({\C})}\quad .
\ee
We have
\[
\begin{array}{l}
\ds\ti{\phi}(0)=-\frac{1}{2\pi}\int_{{\R}^2}\log r\ \p_x\ti{a}\,\p_y\ti{b}-\p_x\ti{b}\,\p_y\ti{a}\\[5mm]
\ds\quad=-\frac{1}{2\pi}\int_0^{2\pi}\int_0^{+\infty}\log\,r\ \frac{\p}{\p r}\lf(\ti{a}\ \frac{\p\ti{b}}{\p \theta}\rg)-\frac{\p}{\p \theta}\lf(\ti{a}\ \frac{\p\ti{b}}{\p r}\rg)\ dr\ d\theta\\[5mm]
\ds \quad=\frac{1}{2\pi}\int_0^{2\pi}\int_0^{+\infty}\ti{a}\ \frac{\p\ti{b}}{\p \theta}\ \frac{dr}{r}\ d\theta
\end{array}
\]
Because $\int_0^{2\pi}\frac{\p\ti{b}}{\p \theta}\ d\theta=0$, we may deduct from each circle $\p B_r(0)$ a constant \`a $\ti{a}$ chosen to have average $\ov{\ti{a}}_r$ on $\p B_r(0)$. Hence, there holds
\[
\ds\ti{\phi}(0)=\frac{1}{2\pi}\int_0^{2\pi}\int_0^{+\infty}[\ti{a}-\ov{\ti{a}}_r]\ \frac{\p\ti{b}}{\p \theta}\ \frac{dr}{r}\ d\theta\quad .
\]
Applying successively the Cauchy-Schwarz and Poincar\'e inequalities on the circle $S^1$, we obtain
\[
\begin{array}{l}
\ds|\ti{\phi}(0)|\le\frac{1}{2\pi}\int_0^{+\infty}\frac{dr}{r}\lf(\int_0^{2\pi}|\ti{a}-\ov{\ti{a}}_r|^2\rg)^\frac{1}{2}\ \lf(\int_0^{2\pi}\lf|\frac{\p\ti{b}}{\p\theta}\rg|^2\rg)^\frac{1}{2}\\[5mm]
\ds\quad\le\frac{1}{2\pi}\int_0^{+\infty}\frac{dr}{r}\lf(\int_0^{2\pi}\lf|\frac{\p\ti{a}}{\p\theta}\rg|^2\rg)^\frac{1}{2}\ \lf(\int_0^{2\pi}\lf|\frac{\p\ti{b}}{\p\theta}\rg|^2\rg)^\frac{1}{2}
\end{array}
\]
The sought after inequality (\ref{III.5}) may then be inferred from the latter via applying once more the Cauchy-Schwarz inequality.

\medskip

Returning to the disk $D^2$, the Whitney extension theorem yields the existence of
$\ti{a}$ and $\ti{b}$ such that
\be
\label{III.6}
\int_{{\C}}|\nabla\ti{a}|^2\le C_1\ \int_{D^2}|\nabla{a}|^2\quad,
\ee
and
\be
\label{III.6a}
\int_{{\C}}|\nabla\ti{b}|^2\le C_1\ \int_{D^2}|\nabla{b}|^2\quad.
\ee
Let $\ti{\phi}$ be the function in (\ref{III.4a}). The difference $\phi-\ti{\phi}$ satisfies the equation
\[
\lf\{
\begin{array}{l}
\ds\Delta(\phi-\ti{\phi})=0\quad\quad\mbox{ in }D^2\\[5mm]
\ds\phi-\ti{\phi}=-\ti{\phi}\quad\quad\quad\mbox{ on }\p D^2
\end{array}
\rg.
\]
The {\it maximum principle} applied to the inequalities (\ref{III.5}), (\ref{III.6}) and (\ref{III.6a}) produces
\[
\|\phi-\ti{\phi}\|_{L^\infty(D^2)}\le\|\ti{\phi}\|_{L^\infty(\p D^2)}\le C\|\nabla a\|_2\ \|\nabla b\|_2\quad .
\]
With the triangle inequality $|\|\phi\|_\infty-\|\ti{\phi}\|_\infty|\le \|\phi-\ti{\phi}\|_\infty$ and the inequality (\ref{III.5}), we reach the desired $L^\infty$-estimate of $\phi$, and therefore, per the above discussion, the theorem is proved. \hfill $\Box$ 

\newpage

{\bf Proof of the regularity of the solutions of the CMC equation.} 

\medskip

Our first aim will be to establish the existence of a positive constant $\al$ such that
\be
\label{III.7}
\sup_{\rho<1/4,\  p\in B_{1/2}(0)}\rho^{-\al}\ \int_{B_\rho(p)}|\nabla u|^2<+\infty\quad .
\ee
Owing to a classical result from Functional Analysis\footnote{See for instance \cite{Gi}.}, the latter implies that $u\in C^{0,\al/2}(B_{1/2}(0))$ . From this, we deduce that $u$ is locally H\"older continuous in the interior of the disk $D^2$. We will then explain how to obtain the smoothness of $u$ from its H\"older continuity.

Let $\ep_0>0$. There exists some radius $\rho_0>0$ such that for every $r<\rho_0$ and every point $p$ in $B_{1/2}(0)$
\[
\int_{B_r(p)}|\nabla u|^2<\ep_0\quad.
\]
We shall in due time adjust the value $\ep_0$ to fit our purposes. In the sequel, $r<\rho_0$. On $B_r(p)$, we decompose $u=\phi+v$ in such a way that 
\[
\lf\{
\begin{array}{l}
\ds\Delta\phi=H\ \p_xu\times \p_yu\quad\quad\quad\mbox{ in }\quad B_r(p)\\[5mm]
\ds \phi=0\quad\quad\quad\mbox{ on }\quad\p B_r(p)
\end{array}
\rg.
\]
Applying theorem~\ref{th-III.1} to $\phi$ yields
\be
\label{III.7aa}
\begin{array}{rl}
\ds\int_{B_r(p)}|\nabla\phi|^2&\ds\le C_0|H|\ \int_{B_r(p)}|\nabla u|^2\ \int_{B_r(p)}|\nabla u|^2\\[5mm]
 &\ds\le C_0|H|\ \ep_0\ \int_{B_r(p)}|\nabla u|^2\quad.
 \end{array}
\ee
The function $v=u-\phi$ is harmonic. To obtain useful estimates on $v$, we need the following result.

\begin{Lm}
\label{lm-III.1}
Let $v$ be a harmonic function on $D^2$. For every point $p$ in $D^2$, the function
\[
\rho\longmapsto \frac{1}{\rho^2}\int_{B_\rho(p)}|\nabla v|^2
\]
is increasing.
\hfill$\Box$
\end{Lm}
{\bf Proof.}
Note first that
\be
\label{III.7a}
\frac{d}{d\rho}\lf[\frac{1}{\rho^2}\int_{B_\rho(p)}|\nabla v|^2\rg]=-\frac{2}{\rho^3}\int_{B_\rho(p)}|\nabla v|^2+\frac{1}{\rho^2}\int_{\p B_\rho(p)}|\nabla v|^2\quad.
\ee
Denote by $\ov{v}$ the average of $v$ on $\p B_\rho(p)$ : $\ov{v}:=|\p B_\rho(p)|^{-1}\int_{\p B_\rho(p)}\,v$. Then, there holds
\[
0=\int_{B_\rho(p)}(v-\ov{v})\ \Delta v=-\int_{B_\rho(p)}|\nabla v|^2+\int_{\p B_\rho(p)}(v-\ov{v})\ \frac{\p v}{\p \rho}\quad.
\]
This implies that
\be
\label{III.8}
\frac{1}{\rho}\int_{B_\rho(p)}|\nabla v|^2\le\lf(\frac{1}{\rho^2}\int_{\p B_\rho(p)}|v-\ov{v}|^2\rg)^\frac{1}{2}\ \lf(\int_{\p B_\rho(p)}\lf|\frac{\p v}{\p \rho}\rg|^2\rg)^\frac{1}{2}\quad .
\ee
In Fourier space, $v$ satisfies $v=\sum_{n\in{\Z}}a_n\,e^{in\theta}$ and $v-\ov{v}=\sum_{n\in{\Z}^\ast}a_n\,e^{in\theta}$. Accordingly,
\[
\frac{1}{2\pi\rho}\int_{\p B_\rho(p)}|v-\ov{v}|^2=\sum_{n\in{\Z}^\ast}|a_n|^2\le\sum_{n\in{\Z}^\ast}|n|^2|a_n|^2\le\frac{1}{2\pi}\int_0^{2\pi}\lf|\frac{\p v}{\p\theta}\rg|^2\ d\theta\quad.
\]
Combining the latter with (\ref{III.8}) then gives
\be
\label{III.9}
\frac{1}{\rho}\int_{B_\rho(p)}|\nabla v|^2\le\lf(\int_{\p B_\rho(p)}\lf|\frac{1}{\rho}\frac{\p v}{\p\theta}\rg|^2\rg)^\frac{1}{2}\ \lf(\int_{\p B_\rho(p)}\lf|\frac{\p v}{\p \rho}\rg|^2\rg)^\frac{1}{2}\quad.
\ee
If we multiply the Laplace equation throughout by $(x-x_p)\,\p_xv+(y-y_p)\,\p_yv$, and then integrate by parts over $B_\rho(p)$, we reach the {\bf Pohozaev identity} :
\be
\label{III.10}
2\int_{\p B_\rho(p)}\lf|\frac{\p v}{\p\rho}\rg|^2=\int_{\p B_\rho(p)}|\nabla v|^2\quad .
\ee
Altogether with (\ref{III.9}), this identity implies that the right-hand side of (\ref{III.7a}) is positive, thereby concluding the proof \footnote{ Another proof of lemma~\ref{lm-III.1} goes as follows :
if $v$ is harmonic then $f:=|\nabla v|^2$ is sub-harmonic - $\Delta|\nabla v|^2\ge 0$ - and an elementary calculation shows that for any non negative  subharmonic
function $f$  in ${\R}^n$ one has $d/dr(r^{-n}\int_{B_r}f)\ge 0$.}. \hfill $\Box$

\medskip

We now return to the proof of the regularity of the solutions of the CMC equation. Per the above lemma, there holds
\be
\label{III.11}
\int_{B_{\rho/2}(p)}|\nabla v|^2\le\frac{1}{4}\int_{B_\rho(p)}|\nabla v|^2\quad.
\ee
Since $\Delta v=0$ on $B_\rho(p)$, while $\phi=0$ on $\p B_\rho(p)$, we have
\[
\int_{B_\rho(p)}\nabla v\cdot\nabla\phi=0\quad .
\]
Combining this identity to the inequality in (\ref{III.11}), we obtain
\be
\label{III.12}
\begin{array}{rl}
\ds \int_{B_{\rho/2}(p)}|\nabla(v+\phi)|^2&\ds\le\frac{1}{2}\int_{B_\rho(p)}|\nabla(v+\phi)|^2\\[5mm]
 &\ds\ +3\int_{B_\rho(p)}|\nabla\phi|^2\quad .
\end{array}
\ee
which, accounting for (\ref{III.7aa}), yields
\be
\label{III.13}
\int_{B_{\rho/2}(p)}|\nabla u|^2\le \lf(\frac{1}{2}+3\ C_0\ |H|\ \ep_0\rg)\int_{B_\rho(p)}|\nabla u|^2\quad .
\ee
If we adjust $\ep_0$ sufficiently small as to have $3\ C_0\ |H|\ \ep_0<1/4$, it follows that
\be
\label{III.14}
\int_{B_{\rho/2}(p)}|\nabla u|^2\le\frac{3}{4}\int_{B_\rho(p)}|\nabla u|^2\quad .
\ee
Iterating this inequality gives the existence of a constant $\al>0$ such that for all $p\in B_{1/2}(0)$ and all $r<\rho$, there holds
\[
\int_{B_r(p)}|\nabla u|^2\le \lf(\frac{r}{\rho_0}\rg)^\al\ \int_{D^2}|\nabla u|^2\quad,
\]
which implies (\ref{III.7}). Accordingly, the solution $u$ of the CMC equation is H\"older continuous.\\
Next, we infer from (\ref{III.7}) and (\ref{III.1}) the bound
\be
\label{III.15}
\sup_{\rho<1/2,\  p\in B_{1/2}(0)}\rho^{-\al}\ \int_{B_\rho(p)}|\Delta u|<+\infty\quad .
\ee
A classical estimate on Riesz potentials gives
$$
|\nabla u|(p)\le C\frac{1}{|x|}\ast \chi_{B_{1/2}}\ |\Delta u|+C\qquad\forall\:\:p\in B_{1/4}(0)\quad,
$$
where $\chi_{B_{1/2}}$ is the characteristic function of the ball $B_{1/2}(0)$. Together with injections proved by Adams in \cite{Ad}, the latter shows that $u\in W^{1,q}(B_{1/4}(0))$ for any $q>(2-\al)/(1-\al)$. Substituted back into (\ref{III.1}), this fact implies that $\Delta u\in L^r$ for some $r>1$. The equation the becomes subcritical, and a standard bootstrapping argument eventually yields that $u\in C^\infty$. This concludes the proof of the regularity of solutions of the CMC equation.
 
\medskip

\subsection{Harmonic maps with values in the sphere $S^n$}

When the target manifold $N^n$ has codimension 1, the {\it harmonic map equation} (\ref{Ia.20}) becomes (cf. (\ref{Ia.24}))
\be
\label{III.16}
-\Delta u=\nu(u)\ \nabla(\nu(u))\cdot\nabla u\quad,
\ee
where $u$ still denotes the normal unit-vector to the submanifold $N^n\subset\R^{n+1}$. In particular, if $N^n$ is the sphere $S^n$, there holds $\nu(u)=u$, and the equation reads
\be
\label{III.17}
-\Delta u=u\ |\nabla u|^2\quad .
\ee
Another characterization of (\ref{III.17}) states that the function $u\in W^{1,2}(D^2,S^n)$ satisfies (\ref{III.17}) if and only if 
\be
\label{III.18}
u\wedge \Delta u=0\quad\quad\mbox{ in }{\mathcal D}'(D^2)\quad.
\ee
Indeed, any $S^n$-valued map $u$ obeys
\[
0=\Delta \frac{|u|^2}{2}=div(u\,\nabla u)=|\nabla u|^2+u\,\Delta u
\]
so that $\Delta u$ is parallel to $u$ as in (\ref{III.18}) if and only if the proportionality is 
$-|\nabla u|^2$. This is equivalent to (\ref{III.17}). Interestingly enough, J. Shatah \cite{Sha} observed that (\ref{III.18}) is tantamount to
\be
\label{III.19}
\forall i,j=1\cdots n+1\quad\quad div(u^i\,\nabla u_j-u_j\,\nabla u^i)=0\quad .
\ee
This formulation of the equation for $S^n$-valued harmonic maps enables one to pass to the weak limit, just as we previously did in the CMC equation.

\medskip

The {\bf regularity of $S^n$-valued harmonic maps} was obtained by F.H\'elein, \cite{He}. It is established as follows. 

For each pair of indices $(i,j)$ in $\{1\cdots n+1\}^2$, the equation (\ref{III.19}) reveals that the vector field $u^i\,\nabla u^j-u^j\,\nabla u^i$ forms a curl term, and hence there exists 
$B^i_j\in W^{1,2}$ with
\[
\nabla^\perp B^i_j=u^i\,\nabla u_j-u_j\,\nabla u^i\quad .
\]
In local coordinates, (\ref{III.17}) may be written
\be
\label{III.20}
-\Delta u^i=\sum_{j=1}^{n+1}u^i\,\nabla u_j\cdot\nabla u^j\quad.
\ee
We then make the field $\nabla^\perp B^i_j$ appear on the right-hand side by observing that
\[
\sum_{j=1}^{n+1}u_j\,\nabla u^i\cdot\nabla u^j=\nabla u^i\cdot\nabla\lf(\sum_{j=1}^{n+1}|u^j|^2/2\rg)=\nabla u^i\cdot\nabla|u|^2/2=0\quad .
\]
Deducting this null term from the right-hand side of (\ref{III.20}) yields that for all $i=1\cdots n+1$, there holds
\be
\label{III.21}
\begin{array}{rl}
\ds-\Delta u^i&\ds=\sum_{j=1}^{n+1}\nabla^\perp B^i_j\cdot\nabla u^j\\[5mm]
 &\ds=\sum_{j=1}^{n+1}\p_xB^i_j\,\p_y u^j-\p_yB^i_j\,\p_xu^i\quad.
\end{array} 
\ee
We recognize the same Jacobian structure which we previously employed to establish the regularity of solutions of the CMC equation. It is thus possible to adapt mutatis mutandis our argument to (\ref{III.21}) so as to infer that $S^n$-valued harmonic maps are regular.

\subsection{H\'elein's moving frames method and the regularity of harmonic maps mapping into a manifold.}

When the target manifold is no longer a sphere (or, more generally, when it is no longer homogeneous), the aforementioned Jacobian structure disappears, and the techniques we employed no longer seem to be directly applicable. 

\medskip

To palliate this lack of structure, and thus extend the regularity result to harmonic maps mapping into an arbitrary manifold, F. H\'elein devised the {\it moving frames method}. The divergence-form structure being the result of the global symmetry of the target manifold, H\'elein's idea consists in expressing the harmonic map equation in preferred moving frames, called {\it Coulomb frames}, thereby compensating for the lack of global symmetry with ``infinitesimal symmetries".

\medskip

This method, although seemingly unnatural and rather mysterious, has subsequently proved very efficient to answer regularity and compactness questions, such as in the study of nonlinear wave maps (see \cite{FMS}, \cite{ShS}, \cite{Tao1}, \cite{Tao2}). For this reason, it is worthwhile to dwell a bit more on H\'elein's method.
  
\medskip
  
\noindent We first recall the main result of F. H\'elein.
  
\begin{Th}
\label{th-III.2} \cite{He}
Let $N^n$ be a closed $C^2$-submanifold of ${\R}^m$. Suppose that $u$ is a harmonic map in $W^{1,2}(D^2,N^n)$ that weakly satisfies the harmonic map equation (\ref{Ia.20}). Then $u$ lies in $C^{1,\alpha}$ for all $\alpha<1$. 
\end{Th}

\noindent{\bf Proof of theorem~\ref{th-III.2} when $N^n$ is a two-torus.}
 
 The notion of {\it harmonic coordinates} has been introduced first in general relativity by Yvonne Choquet-Bruaht in the early fifties. She discovered that the formulation of Einstein equation in these coordinates simplifies in a spectacular way. This idea of searching optimal charts among all possible ''gauges'' has also been very efficient for harmonic maps into manifolds.
 Since the different works of Hildebrandt, Karcher, Kaul, J\"ager, Jost, Widman...etc in the seventies it was known that the intrinsic harmonic map system (\ref{II.13}) becomes for instance almost
 ''triangular'' in  harmonic coordinates $(x^\al)_\al$ in the target which are minimizing the Dirichlet energy $\int_U |dx^\al|^2_g\ dvol_g$. The drawback of this approach is that
 working with harmonic coordinates requires to localize in the target and to restrict only to maps taking values into  a \underbar{single chart} in which such coordinates exist !
 While looking at regularity question this assumption is very restrictive as long as we don't know that the harmonic map $u$ is continuous for instance. It is not excluded
 {\it a priori} that the weak harmonic map $u$ we are considering ''covers the whole target'' even locally in the domain.
 
 \medskip
 
 The main idea of Frederic H\'elein was to extend the notion of {\it harmonic coordinates} of Choquet-Bruhat by replacing it with the more flexible harmonic or {\it Coulomb orthonormal frame},
 notion for which no localization in the target is needed anymore.
 Precisely the Coulomb orthonormal frames are mappings $e=(e_1,\cdots,e_n)$ from the domain $D^2$ into the orthonormal basis of $T_uN^n$ minimizing 
 the Dirichlet energy but for the covariant derivatives $D_g$ in the target :
 \[
 \int_{D^2}\sum_{i=1}^n|D_ge_i|^2\ dx\,dy=\int_{D^2}\sum_{i,k=1}^n|(e_k, \nabla e_i)|^2\ dx\,dy\quad.
\] 
 where $(\cdot,\cdot)$ denotes the canonical scalar product in ${\R}^m$.
 
 \medskip

We will consider the case when $N^n$ is a two-dimensional parallelizable manifold (i.e. admitting a global basis of tangent vectors for the tangent space), namely a torus $T^2$ arbitrarily immersed into Euclidean space $\R^m$, for $m$ large enough. The case of the two-torus is distinguished. Indeed, in general, if a harmonic map $u$ takes its values in an immersed manifold $N^n$, then it is possible to lift $u$ to a harmonic map $\tilde{u}$ taking values in a parallelizable torus $(S^1)^q$ of higher dimension. Accordingly, the argument which we present below can be analogously extended to a more general setting\footnote{although the lifting procedure is rather technical. The details are presented in Lemma 4.1.2 from \cite{He}.}.

 \medskip
 
Let $u\in W^{1,2}(D^2, T^2)$ satisfy weakly (\ref{Ia.20}). We equip $T^2$ with a global, regular, positive orthonormal tangent frame field $(\ep_1,\ep_2)$. Let $\ti{e}:=(\ti{e}_1,\ti{e}_2)\in W^{1,2}(D^2,{\R}^m\times{\R}^m)$ be defined by the composition
\[
\ti{e}_i(x,y):=\ep_i(u(x,y))\quad.
\]
The map $(\ti{e})$ is defined on $D^2$ and it takes its values in the tangent frame field to $T^2$. Define the energy
\be
\label{III.22}
\min_{\psi\in W^{1,2}(D^2,{\R})}\int_{D^2}|(e_1,\nabla e_2)|^2\ dx\, dy\quad,
\ee
where $(\cdot,\cdot)$ is the standard scalar product on ${\R}^m$, and
\[
e_1(x,y)+ie_2(x,y):= e^{i\psi(x,y)}\,(\ti{e}_1(x,y)+i\ti{e}_2(x,y))\quad .
\]
We seek to optimize the map $(\ti{e})$ by minimizing this energy over the $W^{1,2}(D^2)$-maps taking values in the space of rotations of the plane ${\R}^2\simeq T_{u(x,y)}T^2$. Our goal is to seek a frame field as regular as possible in which the harmonic map equation will be recast. The variational problem (\ref{III.22}) is well-posed, and it further admits a solution in $W^{1,2}$. Indeed, there holds
\[
|(e_1,\nabla e_2)|^2=|\nabla\psi+(\ti{e}_1,\nabla\ti{e}_2)|^2\quad.
\]
Hence, there exists a unique minimizer in $W^{1,2}$ which satisfies
\be
\label{III.23}
0=div\lf(\nabla\psi+(\ti{e}_1,\nabla\ti{e}_2)\rg)=div((e_1,\nabla e_2))\quad.
\ee
A priori, $(e_1,\nabla e_2)$ belongs to $L^2$. But actually, thanks to the careful selection brought in by the variational problem (\ref{III.22}), we shall discover that the frame field $(e_1,\nabla e_2)$ over $D^2$ lies in $W^{1,1}$, thereby improving the original $L^2$ belongingness\footnote{Further yet, owing to a result of Luc Tartar \cite{Tar2}, we know that $W^{1,1}(D^2)$ is continuously embedded in the Lorentz space $L^{2,1}(D^2)$, whose dual is the Marcinkiewicz weak-$L^2$ space $L^{2,\infty}(D^2)$, whose definition was recalled in (\ref{II.5}). A measurable function $f$ is an element of $L^{2,1}(D^2)$ whenever
\be
\label{III.24}
\int_0^{+\infty}\lf|\lf\{p\in D^2\ ;\ |f(p)|>\la\rg\}\rg|^\frac{1}{2}\ d\la\quad.
\ee}.
Because the vector field $(e_1,\nabla e_2)$ is divergence-free, there exists some function  $\phi\in W^{1,2}$ such that
\be
\label{III.24a}
(e_1,\nabla e_2)=\nabla^\perp \phi\quad.
\ee
On the other hand, $\phi$ satisfies by definition
\be
\label{III.25}
-\Delta \phi=(\nabla e_1,\nabla^\perp e_2)=\sum_{j=1}^{m}\p_ye_1^j\p_xe^j_2-\p_xe_1^j\p_ye_2^j\quad .
\ee
The right-hand side of this elliptic equation comprises only Jacobians of elements of $W^{1,2}$. This configuration is identical to those previously encountered in our study of the constant mean curvature equation and of the equation of $S^n$-valued harmonic maps. In order to capitalize on this particular structure, we call upon an extension of Wente's theorem~\ref{th-III.1} due to Coifman, Lions, Meyer, and Semmes.
\begin{Th}
\label{th-III.3}
\cite{CLMS}
Let $a$ and $b$ be two functions in $W^{1,2}(D^2)$, and let $\phi$ be the unique solution in $W^{1,p}_0(D^2)$, for $1\le p<2$ , of the equation
\be
\label{III.26}
\lf\{
\begin{array}{l}
\ds-\Delta\phi=\p_xa\,\p_yb-\p_xb\,\p_ya\quad\quad\mbox{ in }D^2\\[5mm]
\ds\quad\phi=0\quad\quad\quad\quad\mbox{ on }\p D^2\quad.
\end{array}
\rg.
\ee
Then $\phi$ lies in $W^{2,1}$ and
\be
\label{III.27}
\|\nabla^2\phi\|_{L^1(D^2)}\le C_1\ \|\nabla a\|_{L^2(D^2)}\ \|\nabla b\|_{L^2(D^2)}\quad .
\ee
where $C_1$ is a constant independent of $a$ and $b$.\footnote{Theorem~\ref{th-III.1} is a corollary of theorem~\ref{th-III.3} owing to the Sobolev embedding
$W^{2,1}(D^2)\subset W^{1,2}\cap C^0$ . In the same vein, theorem~\ref{th-III.3} was preceded by two intermediary results. The first one, by Luc Tartar \cite{Tar1}, states that the Fourier transform of $\nabla\phi$ lies in the Lorentz space $L^{2,1}$, which also implies theorem~\ref{th-III.1}. The second one, due to Stefan M\"uller, obtains the statement of theorem~\ref{th-III.3} under the additional hypothesis that the Jacobian  $\p_xa\,\p_yb-\p_xb\,\p_ya$ be positive.}
\hfill$\Box$
\end{Th}
Applying this result to the solution $\phi$ of (\ref{III.25}) then reveals that $(e_1,\nabla e_2)$ is indeed an element of $W^{1,1}$. 

\medskip

We will express the harmonic map equation (\ref{Ia.20}) in this particular Coulomb frame field, distinguished by its increased regularity. Note that (\ref{Ia.20}) is equivalent to
\be
\label{III.28}
\lf\{
\begin{array}{l}
\ds(\Delta u,e_1)=0\\[5mm]
\ds(\Delta u,e_2)=0
\end{array}
\rg.
\ee
Using the fact that
\[
\begin{array}{l}
\ds\p_x u, \p_yu\in T_uN^n=\mbox{vec}\{e_1,e_2\}\\[5mm]
\ds(\nabla e_1, e_1)=(\nabla e_2, e_2)=0\\[5mm]
\ds(\nabla e_1,e_2)+(e_1,\nabla e_2)=0
\end{array}
\]
we obtain that (\ref{III.28}) may be recast in the form
\be
\label{III.29}
\lf\{
\begin{array}{l}
\ds div((e_1,\nabla u))= -(\nabla e_2,e_1)\cdot (e_2,\nabla u)\\[5mm]
\ds div((e_2,\nabla u))=(\nabla e_2,e_1)\cdot (e_1,\nabla u)
\end{array}
\rg.
\ee
On the other hand, there holds
\be
\label{III.30}
\lf\{
\begin{array}{l}
\ds rot((e_1,\nabla u))= -(\nabla^\perp e_2,e_1)\cdot (e_2,\nabla u)\\[5mm]
\ds rot((e_2,\nabla u))=(\nabla^\perp e_2,e_1)\cdot (e_1,\nabla u)
\end{array}
\rg.
\ee
We next proceed by introducing the Hodge decompositions in $L^2$ of the frames $(e_i,\nabla u)$, for $i\in\{1,2\}$. In particular, there exist four functions
$C_i$ and $D_i$ in $W^{1,2}$ such that
\[
(e_i,\nabla u)=\nabla C_i+\nabla^\perp D_i\quad .
\]
Setting $W:=(C_1,C_2,D_1,D_2)$, the identities (\ref{III.29}) et (\ref{III.30}) become
\be
\label{III.31}
-\Delta W=\Om\cdot\nabla W\quad ,
\ee
where $\Om$ is the vector field valued in the space of $4\times 4$ matrices defined by
\be
\label{III.32}
\Om=
\lf(
\begin{array}{cccc}
0 & -\nabla^\perp\phi & 0 & -\nabla\phi\\[5mm]
\nabla^\perp\phi & 0 & \nabla\phi & 0\\[5mm]
0 & \nabla\phi & 0 & -\nabla^\perp\phi\\[5mm]
-\nabla\phi & 0 &\nabla^\perp\phi & 0
\end{array}
\rg)
\ee

Since $\phi\in W^{2,1}$, the following theorem \ref{th-III.4} implies that $\nabla W$, and hence $\nabla u$, belong to $L^p$ for some $p>2$, thereby enabling the initialization of a bootstrapping argument analogous to that previously encountered in our study of the CMC equation. This procedure yields that $u$ lies in $W^{2,q}$ for all $q<+\infty$. Owing to the standard Sobolev embedding theorem, it follows that $u\in C^{1,\al}$, which concludes the proof of the desired theorem~\ref{th-III.2} in the case when the target manifold of the harmonic map $u$ is the two-torus. \hfill $\Box$

\bigskip

\begin{Th}
\label{th-III.4}
Let $W$ be a solution in $W^{1,2}(D^2,{\R}^n)$ of the linear system
\be
\label{III.33}
-\Delta W=\Om\cdot\nabla W\quad ,
\ee
where $\Om$ is a $W^{1,1}$ vector field on $D^2$ taking values in the space of $n\times n$ matrices. Then $W$ belongs to $W^{1,p}(B_{1/2}(0))$, for some $p>2$. In particular, $W$ is H\"older continuous\footnote{The statement of theorem~\ref{III.4} is optimal. To see this, consider $\,u=\log\log1/r=W$. One verifies easily that $u\in W^{1,2}(D^2, T^2)$ satisfies weakly (\ref{Ia.20}). Yet, $\Om\equiv\nabla u$ fails to be $W^{1,1}$, owing to
\[
\int_0^1\frac{dr}{r\log\frac{1}{r}}=+\infty\quad .
\]}
\footnote{The hypothesis $\Om\in W^{1,1}$ may be replaced by the condition that $\Om\in L^{2,1}$.}.
\hfill $\Box$
\end{Th}

\noindent{\bf Proof of theorem~\ref{th-III.4}.}

Just as in the proof of the regularity of solutions of the CMC equation, we seek to obtain a Morrey-type estimate via the existence of some constant $\al>0$ such that
\be
\label{III.34}
\sup_{p\in B_{1/2}(0)\;,\; 0<\rho<1/4}\rho^{-\al}\int_{B_\rho(p)}|\Delta W|<+\infty\quad.
\ee
The statement of the theorem is then a corollary of an inequality involving Riesz potentials (cf. \cite{Ad} and the CMC equation case on page 28 above).

Let $\ep_0>0$ be some constant whose size shall be in due time adjusted to fit our needs. There exists some radius $\rho_0$ such that for every $r<\rho_0$ and every point $p\in B_{1/2}(0)$, there holds
\[
\|\Om\|_{L^{2,1}(B_r(p))}<\ep_0\quad.
\]
Note that we have used the aforementioned continuous injection $W^{1,1}\subset L^{2,1}$.
\medskip

Henceforth, we consider $r<\rho_0$. On $B_r(p)$, we introduce the decomposition $W=\Phi+V$, with
\[
\lf\{
\begin{array}{l}
\ds\Delta\Phi=\Om\cdot \nabla W\quad\quad\mbox{ in }B_r(p)\\[5mm]
\ds \quad\Phi=0\quad\quad\quad\mbox{ on }\p B_r(p)\quad .
\end{array}
\rg.
\]
A classical result on Riesz potentials (cf. \cite{Ad}) grants the existence of a constant $C_0$ independent of $r$ and such that
\be
\label{III.35}
\begin{array}{rl}
\ds\|\nabla\Phi\|_{L^{2,\infty}(B_r(p))}&\ds\le C_0\int_{B_r(p)}|\Om\cdot\nabla W|\\[5mm]
 &\ds\le C_0\|\Om\|_{L^{2,1}(B_r(p))}\ \|\nabla W\|_{L^{2,\infty}(B_r(p))}\\[5mm]
  &\ds\le C_0\,\ep_0\ \|\nabla W\|_{L^{2,\infty}(B_r(p))}
  \end{array}
\ee  
As for the function $V$, since it is harmonic, we can call upon lemma~\ref{lm-III.1} to deduce that for every $0<\delta <1$ there holds
\be
\label{III.36}
\begin{array}{rl}
\ds\|\nabla V\|^2_{L^{2,\infty}(B_{\delta r}(p))}&\ds\le \|\nabla V\|^2_{L^{2}(B_{\delta r}(p))}\\[5mm]
 &\ds\le \lf(\frac{4\delta}{3}\rg)^2\ \|\nabla V\|^2_{L^{2}(B_{3 r/4}(p))}\\[5mm]
 &\ds\le C_1\ \lf(\frac{4\delta}{3}\rg)^2\ \|\nabla V\|^2_{L^{2,\infty}(B_{r}(p))}\quad ,
 \end{array}
 \ee
where $C_1$ is a constant independent of $r$. Indeed, the $L^{2,\infty}$-norm of a harmonic function on the unit ball controls all its other norms on balls of radii inferior to   $3/4$. 

\medskip

We next choose $\delta$ independent of $r$ and so small as to have $C_1\lf(\frac{4\delta}{3}\rg)^2<1/16$. We also adjust $\ep_0$ to satisfy $C_0\ep_0\ <1/8$. Then, combining 
(\ref{III.35}) and (\ref{III.36}) yields the following inequality
\be
\label{III.37}
\|\nabla W\|_{L^{2,\infty}(B_{\delta r}(p))}\le \frac{1}{2}\|\nabla W\|_{L^{2,\infty}(B_r(p))}\quad,
\ee
valid for all $r<\rho_0$ and all $p\in B_{1/2}(0)$.\\
Just as in the regularity proof for the CMC equation, the latter is iterated to eventually produce the estimate
\be
\label{III.38}
\sup_{p\in B_{1/2}(0)\;,\; 0<\rho<1/4}\rho^{-\al}\|\nabla W\|_{L^{2,\infty}(B_\rho(p))}<+\infty\quad.
\ee
Calling once again upon the duality $L^{2,1}-L^{2,\infty}$, and upon the upper bound on $\|\Om\|_{L^{2,1}(D^2)}$ provided in (\ref{III.38}), we infer that
\be
\label{III.39}
\sup_{p\in B_{1/2}(0)\;,\; 0<\rho<1/4}\rho^{-\al}\|\Om\cdot\nabla W\|_{L^1(B_\rho(p))}<+\infty\quad, 
\ee
thereby giving (\ref{III.34}). This concludes the proof of the desired statement.\hfill $\Box$

\newpage

\section{A proof of Heinz-Hildebrandt's regularity conjecture.}
\reset

The methods which we have used up to now to approach Hildebrandt's conjecture and obtain the regularity of $W^{1,2}$ solutions of the generic system
\be
\label{IV.1}
\Delta u+A(u)(\nabla u,\nabla u)= H(u)(\nabla ^\perp u,\nabla u)
\ee
rely on two main ideas:
\begin{itemize}
\item[i)] recast, as much as possible, quadratic nonlinear terms as linear combinations of Jacobians or as {\it null forms} ;
\item[ii)] project equation (\ref{IV.1}) on a {\it moving frame} $(e_1\cdots e_n)$ satisfying the {\it Coulomb gauge condition}
\[
\forall i,j=1\cdots m\quad\quad div((e_j,\nabla e_i))=0 \quad.
\]
\end{itemize}
Both approaches can be combined to establish the H\"older continuity of $W^{1,2}$ solutions of (\ref{IV.1}) when the target manifold $N^n$ is $C^2$, and when the prescribed mean curvature $H$ is Lipschitz continuous (see \cite{Bet1}, \cite{Cho}, and \cite{He}). Seemingly, these are the weakest possible hypotheses required to carry out the above strategy.

\medskip
However, to fully solve Heinz-Hildebrandt's conjecture, one must replace the Lipschitzean condition on $H$ by its being an element of $L^\infty$. This makes quite a difference\,!

\medskip
Despite its evident elegance and verified usefulness, H\'elein's moving frames method suffers from a relative opacity:\footnote{Yet another drawback of the moving frames method is that it lifts an $N^n$-valued harmonic map, with $n>2$, to another harmonic map, valued in a parallelizable manifold $(S^1)^q$ of higher dimension. This procedure requires that $N^n$ have a higher regularity than the ``natural" one (namely, $C^5$ in place of $C^2$). It is only under this more stringent assumption that the regularity of $N^n$-valued harmonic maps was obtained in \cite{Bet2} and \cite{He}. The introduction of Schr\"odinger systems with antisymmetric potentials in \cite{RiSt} enabled to improve these results.}
\noindent
what makes nonlinearities of the form
\[
A(u)(\nabla u,\nabla u)- H(u)(\nabla ^\perp u,\nabla u) \quad,
\]
so special and more favorable to treating regularity/compactness matters than seemingly simpler nonlinearities, such as
\[
|\nabla u|^2\quad \quad,
\]
which we encountered in Section 1\,?

\medskip
 
\noindent The moving frames method does not address this question. 

\medskip

We consider a weakly harmonic map $u$ with finite energy, on $D^2$ and taking values in a regular oriented closed submanifold $N^n\subset\R^{n+1}$ of codimension 1. We saw at the end of Section 2 that $u$ satisfies the equation
\be
\label{IV.2}
-\Delta u=\nu(u)\ \nabla(\nu(u))\cdot\nabla u\quad ,
\ee
where $\nu$ is the normal unit-vector to $N^n$ relative to the orientation of $N^n$. 

In local coordinates, (\ref{IV.2}) may be recast as
\be
\label{IV.3}
-\Delta u^i=\nu(u)^i\ \sum_{j=1}^{n+1}\nabla(\nu(u))_j\cdot\nabla u^j\qquad\forall\:\:i=1\cdots n+1\quad.
\ee
In this more general framework, we may attempt to adapt H\'elein's operation which changes (\ref{III.20}) into (\ref{III.21}). The first step of this process is easily accomplished. Indeed, since $\nabla u$ is orthogonal to $\nu(u)$, there holds
\[
\sum_{j=1}^{n+1}\nu_j(u)\,\nabla u^j=0\quad.
\]
Substituting this identity into (\ref{IV.4}) yields another equivalent formulation of the equation satisfies by $N^n$-valued harmonic maps, namely 
\be
\label{IV.4}
-\Delta u^i=\sum_{j=1}^{n+1}\lf(\nu(u)^i\ \nabla(\nu(u))_j-\nu(u)_j\ \nabla(\nu(u))^i\rg)\cdot\nabla u^j\quad .
\ee
On the contrary, the second step of the process can not {\it a priori} be extended. Indeed, one cannot justify that the vector field
\[
\nu(u)^i\ \nabla(\nu(u))_j-\nu(u)_j\ \nabla(\nu(u))^i
\]
is divergence-free. This was true so long as $N^n$ was the sphere $S^n$, but it fails so soon as the metric is ever so slightly perturbed. What remains however robust is the \underbar{antisymmetry} of the matrix
\be
\label{IV.5}
\Om:=\lf(\nu(u)^i\ \nabla(\nu(u))_j-\nu(u)_j\ \nabla(\nu(u))^i\rg)_{i,j=1\cdots n+1}\quad.
\ee
It turns out that the antisymmetry of $\Om$ lies in the heart of the problem we have been tackling in these lecture notes. The following result sheds some light onto this claim.
\begin{Th}
\label{th-IV.1}\cite{Riv1}
Let $\Om$ be a vector field in $ L^2(\wedge^1 D^2\otimes so(m))$, thus takings values in the space antisymmetric $m\times m$ matrices $so(m)$. Suppose that $u$ is a map in
$W^{1,2}(D^2,{\R}^m)$ satisfying the equation\footnote{In local coordinates, (\ref{IV.6}) reads
\[
-\Delta u^i=\sum_{j=1}^{m}\Om^i_j\cdot\nabla u^j\qquad\forall\:\:i=1\cdots m\quad.
\]}
\be
\label{IV.6}
-\Delta u=\Om\cdot\nabla u\quad\quad{ in }\quad{\mathcal D}'(D^2)\quad .
\ee
Then there exists some $p>2$ such that $u\in W^{1,p}_{loc}(D^2,{\R}^m)$. In particular, $u$ is H\"older continuous.\hfill$\Box$
\end{Th}

Prior to delving into the proof of this theorem, let us first examine some of its implications towards answering the questions we aim to solve.

\medskip

First of all, it is clear that theorem~\ref{th-IV.1} is applicable to the equation (\ref{IV.4}) so as to yield the regularity of harmonic maps taking values in a manifold of codimension 1.

\medskip

Another rather direct application of theorem~\ref{th-IV.1} deals with the solutions of the prescribed mean curvature equation in $\R^3$, 
\[
\Delta u=2H(u)\ \p_x u\times \p_y u\quad\quad\mbox{ dans }\quad{\mathcal D}'(D^2)\quad.
\]
This equation can be recast in the form
\[
\Delta u=H(u)\nabla^\perp u\times\nabla u\quad,
\]
Via introducing
\[
\Om:=H(u)\lf(
\begin{array}{ccc}
0 & -\nabla^\perp u_3 &\nabla^\perp u_2\\[5mm]
\nabla^\perp u_3 &0 &-\nabla^\perp u_1\\[5mm]
-\nabla^\perp u_2 &\nabla^\perp u_1& 0
\end{array}
\rg)
\]
we observe successively that $\Om$ is antisymmetric, that it belongs to $L^2$ whenever $H$ belongs to $L^\infty$, and that $u$ satisfies (\ref{IV.6}). The hypotheses of theorem~\ref{th-IV.1} are thus all satisfied, and so we conclude that that $u$ is H\"older continuous.

\medskip

This last example outlines clearly the usefulness of theorem~\ref{th-IV.1} towards solving
Heinz-Hildebrandt's conjecture. Namely, it enables us to weaken the Lipschitzean assumption on $H$ found in previous works (\cite{Hei1}, \cite{Hei2}, \cite{Gr2}, \cite{Bet1}, ...), by only requiring that $H$ be an element of $L^\infty$. This is precisely the condition stated in Hildebrandt's conjecture. By all means, we are in good shape.

\medskip

In fact, Hildebrandt's conjecture will be completely resolved with the help of the following result.

\begin{Th}
\label{th-IV.2}\cite{Riv1}
Let $N^n$ be an arbitrary closed oriented $C^2$-submanifold of ${\R}^m$, with $1\le n<m$, and let $\om$ be a $C^1$ two-form on $N^n$. Suppose that $u$ is a critical point in  $W^{1,2}(D^2,{N^n})$ of the energy
\[
E^\om(u)=\frac{1}{2}\int_{D^2}|\nabla u|^2(x,y)\ dx\ dy+ u^\ast\om\quad .
\]
Then $u$ fulfills all of the hypotheses of theoreme~\ref{th-IV.1}, and therefore is H\"older continuous. \hfill $\Box$
\end{Th}

\noindent{\bf Proof of theorem~\ref{th-IV.2}.}

\medskip

The critical points of $E^\om$ satisfy the equation (\ref{Ia.19}), which, in local coordinates, takes the form
\be
\label{IV.7}
\Delta u^i=-\sum_{j,k=1}^m H^i_{jk}(u)\ \nabla^\perp u^k\cdot\nabla u^j-\sum_{j,k=1}^mA^i_{jk}(u)\ \nabla u^k\cdot \nabla u^j\quad ,
\ee
for $i=1\cdots m$. Denoting by $(\ep_i)_{i=1\cdots m}$ the canonical basis of $\R^m$, we first observe that since
$$H^i_{jk}(z)=d\om_z(\ep_i,\ep_j\ep_k)$$ 
the antisymmetry of the 3-forme $d\om$ yields for every $z\in{\R}^m$ the identity $H^i_{jk}(z)=-H^j_{ik}(z)$. Then, (\ref{IV.7}) becomes
\be
\label{IV.8}
\Delta u^i=-\sum_{j,k=1}^m (H^i_{jk}(u)-H^j_{ik}(u))\ \nabla^\perp u^k\cdot\nabla u^j-\sum_{j,k=1}^mA^i_{jk}(u)\ \nabla u^k\cdot \nabla u^j\quad .
\ee
On the other hand, $A(u)(U,V)$ is orthogonal to the tangent plane for every choice of vectors $U$ et $V$\footnote{Rigorously speaking, $A$ is only defined for pairs of vectors which are tangent to the surface. Nevertheless, $A$ can be extended to all pairs of vectors in $\R^m$ in a neighborhood of $N^n$ by applying the pull-back of the projection on $N^n$. This extension procedure is analogous to that outlined on page 18.}. In particular, there holds
\be
\label{IV.9}
\sum_{j=1}^m A^j_{ik}\ \nabla u^j=0\qquad\forall\:\:i,k=1\cdots m\quad.
\ee
Inserting this identity into (\ref{IV.8}) produces
\be
\label{IV.10}
\begin{array}{rl}
\ds\Delta u^i&\ds=-\sum_{j,k=1}^m (H^i_{jk}(u)-H^j_{ik}(u))\ \nabla^\perp u^k\cdot\nabla u^j\\[5mm]
 &\ds\quad-\sum_{j,k=1}^m(A^i_{jk}(u)-A^j_{ik}(u))\ \nabla u^k\cdot \nabla u^j\quad.
 \end{array}
\ee 
The $m\times m$ matrix $\Om:=(\Om^i_j)_{i,j=1\cdots m} $ defined via
\[
\Om^i_j:=\sum_{k=1}^m(H^i_{jk}(u)-H^j_{ik}(u))\ \nabla^\perp u^k+\sum_{k=1}^m(A^i_{jk}(u)-A^j_{ik}(u))\ \nabla u^k\quad,
\]
is evidently antisymmetric, and it belongs to $L^2$. With this notation, (\ref{IV.10}) is recast in the form (\ref{IV.6}), and thus all of the hypotheses of theorem~\ref{th-IV.1} are fulfilled, thereby concluding the proof of theorem~\ref{th-IV.2}.\hfill$\Box$

\medskip

\noindent{\bf On the conservation laws for Schr\"odinger systems with antisymmetric potentials.}

\medskip

Per the above discussion, there only remains to establish theorem~\ref{th-IV.1} in order to reach our goal. To this end, we will express the Schr\"odinger systems with antisymmetric potentials in the form of {\it conservation laws}. More precisely, we have

\begin{Th}
\label{th-IV.3}\cite{Riv1}
Let $\Om$ be a matrix-valued vector field on $D^2$ in $L^2(\wedge^1D^2,so(m))$. Suppose that $A$ and $B$ are two $W^{1,2}$ functions on $D^2$ taking their values in the same of square $m\times m$ matrices which satisfy the equation
\be
\label{IV.11}
\nabla A-A\Om=-\nabla^\perp B\quad.
\ee
If $A$ is almost everywhere invertible, and if it has the bound
\be
\label{IV.11a}
\|A\|_{L^\infty(D^2)}+\|A^{-1}\|_{L^\infty(D^2)}<+\infty\quad,
\ee
then $u$ is a solution of the Schr\"odinger system (\ref{IV.6}) if and only if it satisfies the conservation law

\be
\label{IV.12}
div(A\nabla u-B\nabla^\perp u)=0\quad.
\ee
If (\ref{IV.12}) holds, then $u\in W^{1,p}_{loc}(D^2,{\R}^m)$ for any $1\le p<+\infty$, and therefore $u$ is H\"older continuous in the interior of $D^2$, $C^{0,\al}_{loc}(D^2)$ 
for any $\al<+\infty$
.\hfill $\Box$
\end{Th}

We note that the conservation law (\ref{IV.12}), when it exists, generalizes the conservation laws previously encountered in the study of problems with symmetry, namely:

\medskip

1) In the case of the constant mean curvature equation, the conservation law  (\ref{III.1}) is (\ref{IV.12}) with the choice
\[
A_{ij}=\delta_{ij}\quad,
\]
and
\[
B=\lf(
\begin{array}{ccc}
0&-H\,u_3 &H\, u_2\\[5mm]
H\,u_3&0 &-H\,u_1\\[5mm]
-H\, u_2& H\,u_1&0
\end{array}
\rg)
\]

\medskip

2) In the case of $S^n$-valued harmonic maps, the conservation law (\ref{III.21}) is (\ref{IV.12}) for
\[
A_{ij}=\delta_{ij}\quad,
\]
and $B=(B^i_j)$ with
\[
\nabla^\perp B^i_j=u^i\,\nabla u_j-u_j\,\nabla u^i\quad .
\]

\medskip

The ultimate part of this section will be devoted to constructing $A$ and $B$, for any given antisymmetric $\Om$, with sufficiently small $L^2$-norms (cf. theorem~\ref{th-IV.4} below). As a result, all coercive conformally invariant Lagrangians with quadratic growth will yield conservation laws written in divergence form. This is quite an amazing fact. Indeed, while in cases of the CMC and $S^n$-valued harmonic map equations the existence of conservation laws can be explained by Noether's theorem\footnote{roughly speaking, symmetries give rise to conservation laws. In both the CMC and $S^n$-harmonic map equations, the said symmetries are tantamount to the corresponding Lagrangians being invariant under the group of isometries of the target space $\R^m$.}, one may wonder {\bf which hidden symmetries yield the existence of the general divergence form (\ref{IV.12})\,?} This profound question shall unfortunately not be addressed here.

 \medskip
 
Prior to constructing $A$ and $B$ in the general case, we first establish theorem~\ref{th-IV.3}.
 
 \medskip
 
 \noindent{\bf Proof of theorem~\ref{th-IV.3}.}

The first part of the theorem is the result of the elementary calculation,
\[
\begin{array}{rl}
\ds div(A\,\nabla u-B\,\nabla^\perp u)&\ds=A\, \Delta u+\nabla A\cdot\nabla u-\nabla B\cdot\nabla^\perp u\\[5mm]
 &\ds=A\ \Delta u+(\nabla A+\nabla^\perp B)\cdot\nabla u\\[5mm]
 &\ds=A(\Delta u +\Om\cdot\nabla u)=0
 \end{array}
\]
Regularity matters are settled as follows. Just as in the previously encountered problems, we seek to employ a Morrey-type argument via the existence of some constant $\al>0$ such that
\be
\label{IV.13}
\sup_{p\in B_{1/2}(0)\;,\; 0<\rho<1/4}\rho^{-\al}\int_{B_\rho(p)}|\Delta u|<+\infty\quad.
\ee
The fact that $\nabla u$ belongs to $L^p_{loc}(D^2)$ for some $p>2$ is then deduced through calling upon the inequalities in \cite{Ad}, exactly in the same manner as we previously outlined.
Finally once we know that $\nabla u$ belongs to $L^p_{loc}(D^2)$ for some $p>2$ we deduce the whole regularity result stated in the theorem by using the following lemma.

\begin{Lm}
\label{lm-IV.ax}
Let $m\in {\N}\setminus\{0\}$ and $u\in W^{1,p}_{loc}(D^2,{\R}^m)$ for some $p>2$ satisfying
\[
-\Delta u=\Om\cdot\nabla u
\]
where\footnote{Observe that in this lemma no antisymmetry assumption is made for $\Omega$ which is an arbitrary $m\times m-$matrix valued $L^2-$vectorfield on $D^2$.} $\Om\in L^2(D^2,M_m({\R})\otimes{\R}^2)$ then $u\in W^{1,q}_{loc}(D^2,{\R}^m)$ for any $q<+\infty$.\hfill$\Box$
\end{Lm}

\medskip

Let $\ep_0>0$ be some constant whose value will be adjusted in due time to fit our needs. There exists a radius $\rho_0$ such that for every $r<\rho_0$ and every point $p$ dans $B_{1/2}(0)$, there holds
\be
\label{IV.14a}
\int_{B_r(p)}|\nabla A|^2+|\nabla B|^2<\ep_0\quad.
\ee

\medskip
\noindent
Henceforth, we consider only radii $r<\rho_0$.

\medskip

Note that $A\nabla u$ satisfies the elliptic system
\[
\lf\{
\begin{array}{rl}
\ds div(A\nabla u)&\ds=\nabla B\cdot\nabla^\perp u=\p_yB\,\p_xu-\p_xB\,\p_yu\\[5mm]
\ds rot(A\nabla u)&\ds=-\nabla A\cdot\nabla ^\perp u=\p_x A\,\p_yu-\p_yA\,\p_x u
\end{array}
\rg.
\]
We proceed by introducing on $B_r(p)$ the linear Hodge decomposition in $L^2$ of $A\nabla u$. Namely, there exist two functions $C$ and $D$, unique up to additive constants, elements of $W^{1,2}_0(B_r(p))$ and $W^{1,2}(B_r(p))$ respectively, and such that
\be
\label{IV.14}
A\nabla u=\nabla C+\nabla^\perp D\quad .
\ee
To see why such $C$ and $D$ do indeed exist, consider first the equation
\be
\label{IV.15}
\lf\{
\begin{array}{l}
\ds\Delta C=div(A\nabla u) =\p_yB\,\p_xu-\p_xB\,\p_yu\\[5mm]
\ds\quad C=0\quad .
\end{array}
\rg.
\ee
Wente's theorem (\ref{th-III.1}) guarantees that $C$ lies in $W^{1,2}$, and moreover
\be
\label{IV.16}
\int_{D^2}|\nabla C|^2\le C_0\ \int_{D^2}|\nabla B|^2\ \int_{D^2}|\nabla u|^2\quad.
\ee 
By construction, $div(A\nabla u-\nabla C)=0$. Poincar\'e's lemma thus yields the existence of $D$ in $W^{1,2}$ with $\nabla^\perp D:=A\nabla u-\nabla C$, and
\be
\label{IV.17}
\begin{array}{l}
\ds\int_{D^2}|\nabla D|^2\ds\le 2\int_{D^2}|A\nabla u|^2+|\nabla C|^2\\[5mm]
 \ds\quad\le 2\|A\|_\infty\int_{D^2}|\nabla u|^2+2C_0\ \int_{D^2}|\nabla B|^2\ \int_{D^2}|\nabla u|^2\quad .
\end{array}
\ee 
The function $D$ satisfies the identity
\[
\Delta D=-\nabla A\cdot\nabla ^\perp u=\p_x A\,\p_yu-\p_yA\,\p_x u\quad.
\]
Just as we did in the case of the CMC equation, we introduce the decomposition
$D=\phi+v$, with $\phi$ fulfilling
\be
\label{IV.18}
\lf\{
\begin{array}{l}
\ds\Delta\phi=\p_x A\,\p_yu-\p_yA\,\p_x u\quad\quad\mbox{ in }\quad B_r(p)\\[5mm]
\ds \quad\phi=0\quad\quad\quad\mbox{ on }\quad\p B_r(p)\quad ,
\end{array}
\rg.
\ee
and with $v$ being harmonic. Once again, Wente's theorem~\ref{th-III.1} gives us the estimate
\be
\label{IV.19}
\int_{B_r(p)}|\nabla\phi|^2\le C_0\int_{B_r(p)}|\nabla A|^2\ \int_{B_r(p)}|\nabla u|^2\quad.
\ee
The arguments which we used in the course of the regularity proof for the CMC equation may be recycled here so as to obtain the analogous version of (\ref{III.12}), only this time on the ball $B_{\delta r}(p)$, where $0<\delta<1$ will be adjusted in due time. More precisely, we find
\be
\label{IV.20}
\begin{array}{rl}
\ds \int_{B_{\delta r}(p)}|\nabla D|^2&\ds\le 2\delta^2\int_{B_r(p)}|\nabla D|^2\\[5mm]
 &\ds\ +3\int_{B_r(p)}|\nabla\phi|^2\quad .
\end{array}
\ee
Bringing altogether (\ref{IV.14a}), (\ref{IV.16}), (\ref{IV.17}), (\ref{IV.19}) et (\ref{IV.20}) produces
\be
\label{IV.21}
\begin{array}{rl}
\ds\int_{B_{\delta r}(p)}|A\,\nabla u|^2&\ds\le3\delta^2\int_{B_r(p)}|A\,\nabla u|^2\\[5mm]
 &\ds \quad+C_1\,\ep_0\ \int_{B_r(p)}|\nabla u|^2
 \end{array}
 \ee
Using the hypotheses that $A$ and $A^{-1}$ are bounded in $L^\infty$, it follows from (\ref{IV.21}) that for all  $1>\delta>0$, there holds the estimate
\be
\label{IV.22}
\begin{array}{rl}
\ds\int_{B_{\delta\, r}(p)}|\nabla u|^2&\ds\le3\|A^{-1}\|_\infty \,\|A\|_\infty\delta^2\int_{B_r(p)}|\nabla u|^2\\[5mm]
 &\ds \quad+C_1\,\|A^{-1}\|_\infty\ep_0\ \int_{B_r(p)}|\nabla u|^2\quad.
 \end{array}
 \ee
Next, we choose $\ep_0$ and $\delta$ strictly positive, independent of $r$ et $p$, and such that
\[
3\|A^{-1}\|_\infty \,\|A\|_\infty\delta^2+C_1\,\|A^{-1}\|_\infty\ep_0=\frac{1}{2}\quad.
\]
For this particular choice of $\delta$, we have thus obtained the inequality
\[
\int_{B_{\delta\, r}(p)}|\nabla u|^2\le \frac{1}{2}\int_{B_r(p)}|\nabla u|^2\quad.
\]
Iterating this inequality as in the previous regularity proofs yields the existence of some constant $\al>0$ for which
\[
\sup_{p\in B_{1/2}(0)\;,\; 0<\rho<1/4}\rho^{-2\al}\int_{B_\rho(p)}|\nabla u|^2<+\infty\quad.
\]
Since $|\Delta u|\le|\Om|\,|\nabla u|$, the latter gives us (\ref{IV.13}), thereby concluding the proof of theorem~\ref{th-IV.3}.\hfill $\Box$

\bigskip

There only now remains to establish the existence of the functions $A$ and $B$ in $W^{1,2}$ satisfying the equation (\ref{IV.11}) and the hypothesis (\ref{IV.11a}).

\medskip

{\bf The construction of conservation laws for systems with antisymmetric potentials, and the proof of theorem~\ref{th-IV.1}.}

\medskip

The following result, combined to theorem~\ref{th-IV.3}, implies theorem~\ref{IV.1}, itself yielding theorem~\ref{th-IV.2}, and thereby providing a proof of Hildebrandt's conjecture, as we previously explained.

\medskip

\begin{Th}
\label{th-IV.4}\cite{Riv1}
There exists a constant $\ep_0(m)>0$ depending only on the integer $m$, such that for every vector field $\Om\in L^2(D^2,so(m))$ with
\be
\label{IV.23}
\int_{D^2}|\Om|^2<\ep_0(m)\quad,
\ee
it is possible to construct $A\in L^\infty(D^2,Gl_m({\R}))\cap
W^{1,2}$ and $B\in W^{1,2}(D^2,{ M}_m({\R}))$ with the properties
\begin{itemize}
\item[i)] 
\be
\label{IV.24}
\int_{D^2}|\nabla A|^2 +\|dist(A,SO(m))\|_{L^\infty(D^2)}\le C(m)\int_{D^2}|\Om|^2\quad,
\ee
\item[ii)]
\be
\label{IV.26}
div\lf(\nabla_\Om A\rg):=div\lf(\nabla A- A \Om\rg)=0\quad,
\ee
\end{itemize}
where $C(m)$ is a constant depending only on the dimension $m$.\hfill $\Box$
\end{Th}

\medskip

Prior to delving into the proof of theorem~\ref{th-IV.4}, a few comments and observations are in order. 

\medskip

\noindent Glancing at the statement of the theorem, one question naturally arises: {\bf why is the antisymmetry of $\Om$ so important\,?}

\medskip 

\noindent It can be understood as follow. 

\medskip

\noindent 
In the simpler case when $\Om$ is {\bf divergence-free}, we can write $\Om$ in the form
\[
\Om=\nabla^\perp\xi\quad,
\]
for some $\xi\in W^{1,2}(D^2,so(m))$. In particular, the statement of theorem~\ref{th-IV.4} is settled by choosing
\be\label{buzuc}
A_{ij}\;=\;\delta_{ij}\qquad{and}\qquad B_{ij}\;=\;\xi_{ij}\quad.
\ee

\medskip

Accordingly, it seems reasonable in the general case to seek a solution pair $(A,B)$ which comes as ``close" as can be to (\ref{buzuc}). A first approach consists in performing a {\bf linear Hodge decomposition} in $L^2$ of $\Om$. Hence, for some $\xi$ and $P$ in $W^{1,2}$, we write
\be
\label{IV.26a}
\Om= \nabla^\perp\xi-\nabla P\quad.
\ee
In this case, we see that if $A$ exists, then it must satisfy the equation
\be
\label{IV.27}
\Delta A=\nabla A\cdot\nabla^\perp\xi-div(A\nabla P)\quad .
\ee
This equation is critical in $W^{1,2}$. The first summand $\nabla A\cdot\nabla^\perp\xi$ on the right-hand side of (\ref{IV.27}) is a Jacobian. This is a desirable feature with many very good analytical properties, as we have previously seen. In particular, using integration by compensation (Wente's theorem~\ref{III.1}), we can devise a bootstrapping argument beginning in $W^{1,2}$. On the other hand, the second summand $div(A\nabla P)$ on the right-hand side of (\ref{IV.27}) displays no particular structure. All which we know about it, is that $A$ should {\it a-priori} belong to $W^{1,2}$. But this space is not embedded in $L^\infty$, and so we cannot a priori conclude that  $A\nabla P$ lies in $L^2$, thereby obstructing a successful analysis...

\medskip
However, not all hope is lost for the {\bf antisymmetric structure} of $\Om$ still remains to be used. The idea is to perform a {\bf nonlinear} Hodge decomposition\footnote{which is tantamount to a {\bf change of gauge}.} in $L^2$ of $\Om$. Thus, let $\xi\in W^{1,2}(D^2,so(m))$ and $P$ be a $W^{1,2}$ map taking values in the group $SO(m)$ of proper rotations of $\R^m$, such that
\be
\label{IV.28}
\Om=P\,\nabla^\perp\xi \,P^{-1}-\nabla P\, P^{-1}\quad.
\ee
At first glance, the advantage of (\ref{IV.28}) over (\ref{IV.27}) is not obvious. If anything, it seems as though we have complicated the problem by having to introduce left and right multiplications by $P$ and $P^{-1}$. On second thought, however, since rotations are always bounded, the map $P$ in (\ref{IV.28}) is an element of $W^{1,2}\cap L^\infty$, whereas in (\ref{IV.27}), the map $P$ belonged only to $W^{1,2}$. This slight improvement will actually be sufficient to successfully carry out our proof. Furthermore, (\ref{IV.28}) has  yet another advantage over (\ref{IV.27}). Indeed, whenever $A$ and $B$ are solutions of (\ref{IV.26}), there holds
\[
\begin{array}{l}
\ds\nabla_{\nabla^\perp\xi}(AP)=\nabla(AP)-(AP)\,\nabla^\perp\xi\\[5mm]
\ds\quad\quad=  \nabla A\, P+A\,\nabla P-A\,P\ (P^{-1}\Om P+P^{-1}\nabla P)\\[5mm]
\ds\quad\quad=(\nabla_\Om A)P=-\nabla^\perp B\ P\quad.
\end{array}
\]
Hence, via setting $\ti{A}:=AP$, $\ti{A}$, we find
\be
\label{IV.29}
\Delta \ti{A}=\nabla\ti{A}\cdot\nabla^\perp\xi+\nabla^\perp B\cdot\nabla P\quad.
\ee
Unlike (\ref{IV.27}), the second summand on the right-hand side of (\ref{IV.29}) is a linear combination of Jacobians of terms which lie in $W^{1,2}$. Accordingly, calling upon theorem~\ref{th-III.1}, we can control $\ti{A}$ in $L^\infty\cap W^{1,2}$. This will make a  bootstrapping argument possible. 

\medskip

One point still remains to be verified. Namely, that the nonlinear Hodge decomposition (\ref{IV.28}) does exist. This can be accomplished with the help of a result of Karen Uhlenbeck\footnote{In reality, this result, as it is stated here, does not appear in the original work of Uhlenbeck. In \cite{Riv1}, it is shown how to deduce theorem~\ref{th-IV.5} from Uhlenbeck's approach.}. 
\begin{Th}
\label{th-IV.5}\cite{Uhl}, \cite{Riv1}
Let $m\in{\N}$. There are two constants $\ep(m)>0$ and $C(m)>0$, depending only on $m$, such that for each vector field $\Om\in L^2(D^2,so(m))$ with
\[
\int_{D^2}|\Om|^2<\ep(m)\quad ,
\]
there exist $\xi\in W^{1,2}(D^2,so(m))$ and $P\in W^{1,2}(D^2,SO(m))$ satisfying
\be
\label{IV.30}
\Om=P\,\nabla^\perp\xi\, P^{-1}-\nabla P\, P^{-1}\quad,
\ee
\be
\label{IV.30a}
\xi=0\quad\quad\quad\mbox{ on }\quad \p D^2\quad,
\ee
and
\be
\label{IV.31}
\int_{D^2}|\nabla\xi|^2+\int_{D^2}|\nabla P|^2\le C(m)\ \int_{D^2}|\Om|^2\quad .
\ee
\hfill$\Box$
\end{Th}

\noindent{\bf Proof of theorem~\ref{th-IV.4}.}

Let $P$ and $\xi$ be as in  theorem~\ref{th-IV.5}. To each $A\in L^\infty\cap W^{1,2}(D^2,M_m({\R}))$ we associate $\ti{A}=AP$. Suppose that $A$ and $B$ are solutions of (\ref{IV.26}). Then $\ti{A}$ and $B$ satisfy the elliptic system
\be
\label{IV.32}
\lf\{
\begin{array}{l}
\ds\Delta \ti{A}=\nabla\ti{A}\cdot\nabla^\perp\xi+\nabla^\perp B\cdot\nabla P\\[5mm]
\ds \Delta B=-div(\ti{A}\,\nabla\xi\ P^{-1})+\nabla^\perp\ti{A}\cdot\nabla P^{-1}\quad.
\end{array}
\rg.
\ee
We first consider the invertible elliptic system
\be
\label{IV.32z}
\lf\{
\begin{array}{l}
\ds\Delta \ti{A}=\nabla\hat{A}\cdot\nabla^\perp\xi+\nabla^\perp \hat{B}\cdot\nabla P\\[5mm]
\ds \Delta B=-div(\hat{A}\,\nabla\xi\ P^{-1})+\nabla^\perp\hat{A}\cdot\nabla P^{-1}\\[5mm]
\ds \frac{\p \ti{A}}{\p \nu}=0\quad\mbox{ and }\quad B=0\quad\mbox{ on }\quad\p D^2\\[5mm]
\ds \int_{D^2}\ti{A}=\pi^2\ Id_m
\end{array}
\rg.
\ee
where $\hat{A}$ and $\hat{B}$ are arbitrary functions in $L^\infty\cap W^{1,2}$ and in $W^{1,2}$ respectively. An analogous version\footnote{whose proof is left as an exercise.} of theorem~\ref{th-III.1} with Neuman boundary conditions in place of Dirichlet conditions, we deduce that the unique solution $(\ti{A},B)$ of (\ref{IV.32}) satisfies the estimates
\be
\label{IV.33}
\begin{array}{rl}
\ds\int_{D^2}|\nabla \ti{A}|^2+\|\ti{A}-Id_m\|^2_\infty&\ds\le C\ \int_{D^2}|\nabla\hat{A}|^2\ \int_{D^2}|\nabla\xi|^2\\[5mm]
 &\ds\ +C\ \int_{D^2}|\nabla\hat{B}|^2\ \int_{D^2}|\nabla P|^2\quad ,
\end{array}
\ee
and
\be
\label{IV.34}
\begin{array}{rl}
\ds\int_{D^2}|\nabla (\ti{B}-B_0)|^2&\ds\le C\ \|\hat{A}-Id_m\|^2_\infty\ \int_{D^2}|\nabla\xi|^2\ \\[5mm]
 &\ds\ + C\ \int_{D^2}|\nabla\hat{A}|^2\ \int_{D^2}|\nabla P|^2\quad,
\end{array}
\ee
where $B_0$ is the solution in $W^{1,2}$ of
\be
\label{IV.35}
\lf\{
\begin{array}{l}
\Delta B_0=-div(\nabla \xi\ P^{-1})\quad\quad\mbox{ in }\quad D^2\\[5mm]
\quad B_0=0\quad\quad\quad\mbox{on }\quad\p D^2
\end{array}
\rg.
\ee  
Hence, if
\[
\int_{D^2}|\nabla P|^2+|\nabla\xi|^2
\]
is sufficiently small (this can always be arranged owing to (\ref{IV.31}) and the hypothesis (\ref{IV.23})), then a standard fixed point argument in the space $\big(L^\infty\cap W^{1,2}(D^2,M_m({\R}))\big)\times W^{1,2}(D^2,M_m({\R}))$ yields the existence of the solution $(\ti{A},B)$ of the system
\be
\label{IV.36}
\lf\{
\begin{array}{l}
\ds\Delta \ti{A}=\nabla\ti{A}\cdot\nabla^\perp\xi+\nabla^\perp {B}\cdot\nabla P\\[5mm]
\ds \Delta B=-div(\ti{A}\,\nabla\xi\ P^{-1})+\nabla^\perp\ti{A}\cdot\nabla P^{-1}\\[5mm]
\ds \frac{\p \ti{A}}{\p \nu}=0\quad\mbox{ and }\quad B=0\quad\mbox{ on }\quad\p D^2\\[5mm]
\ds \int_{D^2}\ti{A}=\pi^2\ Id_m
\end{array}
\rg.
\ee
By construction, this solution satisfies the estimate (\ref{IV.24}) with $A=\ti{A}\, P^{-1}$.

\medskip
The proof of theorem~\ref{th-IV.4} will then be finished once it is established that $(A,B)$ is a solution of (\ref{IV.26}). 

\medskip
To do so, we introduce the following linear Hodge decomposition in $L^2$\,:
\[
\nabla \ti{A}-\ti{A}\nabla^\perp\xi+\nabla^\perp B\ P=\nabla C +\nabla^\perp D
\]
where $C=0$ on $\p D^2$. The first equation in (\ref{IV.36}) states that $\Delta C=0$, so that $C\equiv 0$ sur $D^2$. The second equation in (\ref{IV.36}) along with the boundary conditions imply that $D$ satisfies
\be
\label{IV.37}
\lf\{
\begin{array}{l}
\ds div(\nabla D\ P^{-1})=0\quad\quad\mbox{ in }\quad\p D^2\\[5mm]
\ds\quad D=0\quad\quad\mbox{ on }\quad\p D^2\quad .
\end{array}
\rg.
\ee
Thus, there exists $E\in W^{1,2}(D^2,M_n({\R}))$ such that 
\be
\label{IV.38}
\lf\{
\begin{array}{l}
\ds -\Delta E=\nabla^\perp D\cdot\nabla P^{-1}\quad\quad\mbox{ in }\quad D^2\\[5mm]
\ds \quad \frac{\p E}{\p \nu}=0\quad\quad\mbox{ on }\quad\p D^2
\end{array}
\rg.
\ee
The analogous version of theorem~\ref{th-III.1} with Neuman boundary conditions yields the estimate
\be
\label{IV.39}
\int_{D^2}|\nabla E|\le C_0\ \int_{D^2}|\nabla D|^2\ \int_{D^2}|\nabla P^{-1}|^2\quad.
\ee
Moreover, because $\nabla D=\nabla^\perp E\ P$, there holds $|\nabla D|\le |\nabla E|$. Put into (\ref{IV.39}), this shows that if $\int_{D^2}|\nabla P|^2$ is chosen sufficiently small (i.e. for $\ep_0(m)$ in (\ref{IV.23}) small enough), then $D\equiv 0$. Whence, we find 
\[
\nabla \ti{A}-\ti{A}\nabla^\perp\xi+\nabla^\perp B\ P=0\quad\quad\mbox{ in }\quad D^2\quad,
\]
thereby ending the proof of theorem~\ref{th-IV.4}. \hfill $\Box$

\newpage

\section{A PDE version of the constant variation method for Schr\"odinger Systems with anti-symetric potentials.}

In this part we shall look at various, a-priori critical, elliptic systems with antisymmetric potentials and
extend the approach we developed in the previous section in order to establish their hidden sub-critical
nature. Before to do so we will look at what we have done so far from a different perspective
than the one suggested by the ''gauge theoretic'' type arguments we used.

In 2 dimension, in order to prove the sub-criticality of the system
\be
\label{V.1}
-\Delta u=\Om\cdot \nabla u\quad,
\ee
where $u\in W^{1,2}(D^2,{\R}^m)$ and $\Omega\in L^2(\wedge^1D^2,so(m))$ - sub-criticality meaning
that the fact that $u$ solves (\ref{V.1}) imply that it is in fact more regular than the initial assumption $u\in W^{1,2}$ - we proceeded as follows : we first constructed a solution of the equation
\be
\label{V.2}
div(\nabla P\,P^{-1})=div(P\,\Omega\, P^{-1})\quad,
\ee
and multiplying $\nabla u$ by the rotation $P$ in the equation (\ref{V.1}) we obtain that it is
equivalent to
\be
\label{V.3}
-div(P\,\nabla u)=\nabla^\perp\xi\cdot P\,\nabla u\quad,
\ee
where
$$\nabla^\perp\xi:=-\nabla P\, P^{-1}+P\,\Omega\, P^{-1}$$
Then in order to have a pure Jacobian in the r.h.s to (\ref{V.3}) we looked for a replacement of $P$ by a perturbation of the form
$A:=(id+\epsilon)\ P$ where $\epsilon$ is a $m\times m$ matrix valued map hopefully
small in $L^\infty\cap W^{1,2}$. This is obtained by solving the following well posed problem in $W^{1,2}\cap L^\infty$ - due to Wente's theorem - 
\be
\label{V.4}
\left\{
\begin{array}{l}
\ds div(\nabla\ep\, P)=div((id+\ep)\ \nabla^\perp\xi\, P)\quad\quad\mbox{ in }D^2\\[5mm]
\ds \ep=0\quad\quad\quad\mbox{ on }\p D^2
\end{array}
\right.
\ee
Posing $\nabla^\perp B:=\nabla\ep\, P-(id+\ep)\ \nabla^\perp\xi\, P$ equation (\ref{V.1})
becomes equivalent to
\be
\label{V.5}
-div(A\,\nabla u)=\nabla^\perp B\cdot\nabla u\quad.
\ee 

Trying now to reproduce this procedure in one space dimension leads to the following. We aim 
to solve

\be
\label{V.6}
-u''=\Omega\, u'
\ee
where $u\ :\ [0,1]\rightarrow\ {\R}^m$ and $\Om\ :\ [0,1]\rightarrow\ so(m)$. We then construct a solution
$P$ to
(\ref{V.3}) which in one dimension becomes
\[
(P'\, P^{-1})'=(P\,\Omega\, P^{-1})'\quad .
\]
A special solution is given by $P$ solving
\be
\label{V.7}
P^{-1}\ P'=\Omega
\ee
Computing now $(Pu')'$ gives the 1-D analogue of (\ref{V.3}) which is
\be
\label{V.8}
(P u')'=0\quad,
\ee
indeed the curl operator in one dimension is trivial, and jacobians too, since curl-like vector fields corresponds to functions $f$ satisfying
$div\, f=f'=0$ and then can be taken equal to zero. At this stage it is not necessary to go to the ultimate step and perturb $P$ into $(id+\ep)P$ since the r.h.s of (\ref{V.8}) is already a pure 1D jacobian (that is
zero !) and since we have succeeded in writing equation (\ref{V.6}) in conservative form. \\ 

The reader has then noticed that the 1-D analogue of our approach to write equation (\ref{V.1})
in conservative form is the well known {\bf  variation of the constant method} : a solution is constructed
to some {\bf auxiliary equation} - (\ref{V.2}) or (\ref{V.7}) - which has been carefully chosen
in order to absorb the ''worst'' part of the r.h.s. of the original equation while comparing our
given solution to the constructed solution of the auxiliary equation.\\

We shall then in the sequel forget the geometrical interpretation of the auxiliary equation (\ref{IV.30})
in terms of gauge theory that we see as being specific to the kind of equation (\ref{V.1}) we were
looking at and keep the general philosophy of the classical variation of the constant method for Ordinary Differential Equations that we are now extending to other classes of Partial Differential Equations different from (\ref{V.1}).\\

We establish the following result.

\begin{Th}
\label{th-V.1}\cite{Riv2}
Let $n>2$ and $m\ge 2$ there exists $\ep_0>0$ and $C>0$ such that for any $\Omega\in L^{n/2}(B^n,so(m))$
there exists $A\in L^\infty\cap W^{2,{n/2}}(B^n, Gl_m({\R})$ satisfying
\begin{itemize}
\item[i)]\be
\label{zV.9}
\|A\|_{W^{2,n/2}(B^n)}\le C\ \|\Om\|_{L^{n/2}(B^n)}\quad,
\ee
\item[ii)]
\be
\label{zV.10}
\Delta A+A\Om=0\quad.
\ee
\end{itemize}
Moreover for any map $v$ in $L^{n/(n-2)}(B^n,{\R}^m)$ 
\be
\label{V.9}
-\Delta v=\Omega\, v\quad\quad \Longleftrightarrow\quad\quad div\lf(A\,\nabla v-\nabla A\,v\rg)=0
\ee
and we deduce that $v\in L^\infty_{loc}(B^n)$.\hfill $\Box$
\end{Th}

\begin{Rm}
\label{rm-V.1}
Again the assumptions $v\in L^{n/(n-2)}(B^n,{\R}^m)$ and $\Omega\in L^{n/2}(B^n,so(m))$ make equation (\ref{V.9}) critical
in dimension $n$ : Inserting this information in the r.h.s. of (\ref{V.9}) gives $\Delta v\in L^1$ which
implies in return $v\in L^{n/(n-2),\infty}_{loc}$, which corresponds to our definition of being critical for
an elliptic system.
\end{Rm}

\begin{Rm}
\label{rm-V.2}
We have then been able to write critical systems of the kind $-\Delta v=\Om\,v$ in conservative form whenever $\Om$ is antisymmetric.
This ''factorization of the divergence'' operator is obtain through the construction of a solution $A$ to the auxiliary equation $\Delta A+A\,\Om=0$
exactly like in the {\bf constant variation method in 1-D}, a solution to the auxiliary equation $A''+A\,\Om=0$ permits to factorize the derivative
in the ODE given by $ -z''=\Om\, z$ which becomes, after multiplication by $A$ : $(A\, z'-A'\,z)'=0$.
\end{Rm}

\noindent{\bf Proof of theorem~\ref{th-V.1} for $n=4$.}

The goal again here, like in the previous sections, is to establish a Morrey type estimate for $v$ that could be re-injected in the equation and converted into an $L^q_{loc}$ due to Adams result in \cite{Ad}.\\

We shall look for some  map $P$ from $B^4$ into the space $SO(m)$ solving some ad-hoc auxiliary
equation. Formal computation - we still don't kow which regularity for $P$ we should assume at this
stage - gives
\be
\label{V.10}
\begin{array}{rl}
\ds-\Delta(P\, v)&\ds=-\Delta P\, v-P\Delta v -2\nabla P\cdot\nabla v\\[5mm]
\ds &\ds=\lf(\Delta P\, P^{-1}+P\,\Omega\, P^{-1}\rg)\ Pv \\[5mm]
\ds &\ds\quad-2\, div(\nabla P\, P^{-1}\ Pv)\quad.
\end{array}
\ee
In view of the first term in the  r.h.s. of equation (\ref{V.10}) it is natural to look for
$\Delta P$ having the same regularity as $\Omega$ that is $L^2$. Hence we are looking for $P\in W^{2,2}(B^4,SO(m))$. Under such an assumption the second term in the r.h.s. is not problematic while
working in the function space $L^2$ for $v$ indeed, standard elliptic estimates give 
\be
\label{V.11}
\|\Delta_0^{-1}\lf(div(\nabla P\, P^{-1}\ Pv)\rg)\|_{L^2(B^4)}\le \|\nabla P\|_{L^4(B^4)}\ \|v\|_{L^2(B^4)}
\ee 
where $\Delta_0^{-1}$ is the map which to $f$ in $W^{-1,4/3}$ assigns the function $u$ in $W^{1,4/3}(B^4)$ satisfying $\Delta u=f$ and equal to zero on the boundary of $B^4$. Hence, if we localize in space
that ensures that $\|\nabla P\|_{L^4(B^4)}$ is small enough, the contribution of the second term
of the r.h.s. of (\ref{V.10}) can be ''absorbed'' in the l.h.s. of (\ref{V.10}) while working with the $L^2$ norm
of $v$.\\

The first term in the r.h.s. of (\ref{V.10}) is however problematic while intending to work with the $L^2$
norm of $v$. It is indeed only in $L^1$ and such an estimate  does not give in return an $L^2$ 
control of $v$. It is then tempting to look for a map $P$ solving an auxiliary equation in such a way that this term vanishes. Unfortunately such a hope cannot be realized due to the fact that when $P$ is a map
into the rotations
$$
P\,\Omega\, P^{-1}\in so(m)\quad\mbox{ but a-priori }\quad \Delta P\, P^{-1}\notin so(m)\quad .
$$
The idea is then to cancel ''as much as we can'' in the first term of the r.h.s. of (\ref{V.10}) by looking at a solution to the following auxiliary equation\footnote{
which does not have, to our knowledge, a geometric relevant interpretation similar to the Coulomb Gauge extraction we used in the previous section.} :
\be
\label{V.12}
ASym(\Delta P\, P^{-1})+P\,\Omega\, P^{-1}=0
\ee
where $ASym(\Delta P\, P^{-1})$ is the antisymmetric part of $\Delta P\, P^{-1}$ given by
\[
ASym(\Delta P\, P^{-1}):=\frac{1}{2}\lf(\Delta P\, P^{-1}-P\,\Delta P^{-1}\rg)\quad.
\]
Precisely the following proposition holds
\begin{Prop}
\label{pr-V.1}\label{DR2}
There exists $\ep_0>0$ and $C>0$ such that for any $\Om\in L^2(B^4,so(m))$ satisfying $\|\Om\|_{L^2}<\ep_0$ there exists
$P\in W^{2,2}(B^4,SO(m))$ such that
\[
ASym(\Delta P\, P^{-1})+P\,\Omega\, P^{-1}=0\quad,
\]
and
\be
\label{V.12a}
\|\nabla P\|_{W^{1,2}(B^4)}\le C\ \|\Om\|_{L^2(B^4)}\quad.
\ee
\hfill $\Box$
\end{Prop}
Taking the gauge $P$ given by the previous proposition and denoting $w:=P\,v$ the system (\ref{V.10}) becomes
\be
\label{V.13}
{\mathcal L}_{P}w:=-\Delta w -(\nabla P\ P^{-1})^2\ w+2\ div(\nabla P\ P^{-1}\ w)=0\quad,
\ee
where  we have used that 
$$
\begin{array}{l}
Symm(\Delta P\, P^{-1})=2^{-1}\lf(\Delta P\, P^{-1}+P\,\Delta P^{-1}\rg)\\[5mm]
\quad\quad=2^{-1}div\lf(\nabla P\, P^{-1}+P\,\nabla P^{-1}\rg)+\nabla P\cdot\nabla P^{-1}\\[5mm]
\quad\quad=-(\nabla P\ P^{-1})^2\quad.
\end{array}
$$
Denote for  any $Q\in W^{2,2}(B^4,M_m({\R}))$ ${\mathcal L}^\ast_PQ$ the ''formal adjoint'' to ${\mathcal L}_{P}$ acting on $Q$ 
\[
{\mathcal L}_P^\ast Q:=-\Delta Q-2\ \nabla Q\cdot\nabla P\ P^{-1}-Q\ (\nabla P\ P^{-1})^2\quad.
\]
In order to factorize the divergence operator in (\ref{V.13}) it is natural to look for $Q$ satisfying ${\mathcal L}_P^\ast Q=0$.
One has indeed
\be
\label{V.14}
\begin{array}{l}
\ds 0=QL_Pw-wL^\ast_PQ\\[5mm]
\ds \quad={div( Q\nabla w-[\nabla Q+2\nabla P\ P^{-1}]\ w)}\quad.
\end{array}
\ee
It is not so difficult to construct $Q$ solving ${\mathcal L}^\ast_PQ=0$ for a $Q\in W^{2,p}(B^4,M_m({\R})$  (for $p<2$). However,
in order to give a meaning to (\ref{V.14}) we need at least $Q\in W^{2,2}$ and, moreover the invertibility of the matrix $Q$ almost everywhere is also
needed for the conservation law.
In the aim of producing a $Q\in W^{2,2}(B^4,Gl_m({\R}))$ solving ${\mathcal L}^\ast_PQ=0$  it  is important to observe first that  $-(\nabla P\ P^{-1})^2$ is a non-negative symmetric matrix \footnote{Indeed it is a sum of non-negative
symmetric matrices : each of the matrices $\p_{x_j} P\ P^{-1}$ is antisymetric
and hence its square $(\p_{x_j} P\ P^{-1})^2$ is symmetric non-positive.} but that it is in a smaller space than $L^2$ : the space $L^{2,1}$. Combining
the improved Sobolev embeddings \footnote{$L^{4,2}(B^4)$ is the Lorentz space of measurable functions $f$ such that the decreasing rearangement 
$f^\ast$ of $f$ satisfies $\int_0^\infty t^{-1/2} (f^\ast)^2(t)\ dt<+\infty$.}
space $W^{1,2}(B^4)\hookrightarrow L^{4,2}(B^4)$ and the fact that the product of two functions in $L^{4,2}$ is in $L^{2,1}$ (see \cite{Tar2}) we deduce from 
(\ref{V.12a}) that
\be
\label{V.14zz}
\|(\nabla P\ P^{-1})^2\|_{L^{2,1}(B^4)}\le C\ \|\Om\|^2_{L^2(B^4)}\quad.
\ee
Granting these three important properties for $-(\nabla P\ P^{-1})^2$ (symmetry, positiveness and improved integrability) one can prove the following 
result (see \cite{Riv3}).
\begin{Th}
\label{th-V.2}
There exists $\ep>0$ such that
\[
\forall P\in W^{2,2}(B^4,SO(m))\quad\mbox{ satisfying }\quad \|\nabla P\|_{W^{1,2}(B^4)}\le\ \ep
\]
there exists a unique $Q\in W^{2,2}\cap L^\infty(B^4,Gl_m({\R}))$ satisfying
\be
\label{V.15}
\lf\{
\begin{array}{l}
-\Delta Q-2\ \nabla Q\cdot\nabla P\ P^{-1}-Q\ (\nabla P\ P^{-1})^2=0\\[5mm]
Q=id_m
\end{array}
\rg.
\ee
and
\be
\label{V.16}
\|Q-id_m\|_{L^\infty\cap W^{2,2}}\le C\ \|\nabla P\|_{W^{1,2}(B^4)}^2\quad.
\ee
\hfill$\Box$
\end{Th}
Taking $A:=P\,Q$ we have constructed a solution to
\be
\label{zzV.10}
\Delta A+A\Om=0\quad.
\ee
such that
\be
\label{zzV.11}
\|dist(A,SO(m))\|_\infty+\|A-Id\|_{W^{2,2}}\le\,C\ \|\Om\|_{L^2}
\ee
and then for any map $v$ in $L^2(B^4,{\R}^m)$ the following equivalence holds 
\be
\label{zzV.12}
-\Delta v=\Omega\, v\quad\quad \Longleftrightarrow\quad\quad div\lf(A\nabla v- \nabla A\,v\rg)=0
\ee
Having now the equation $-\Delta v=\Omega\, v$ in the form $$div\lf(\nabla w-2\nabla A\, A^{-1}\ w\rg)=0$$
where $w:=A\,v$ permits
to obtain easily the following  Morrey estimate : $\forall\,\rho<1$
\be
\label{zzV.13}
\forall\,\rho<1\quad\quad\sup_{x_0\in B_{\rho}(0),\ r<1-\rho}r^{-\nu}\int_{B_r(x_0)}|w|^2<+\infty\quad,
\ee
for some $\nu>0$. As in the previous sections one deduces using Adams embeddings  that $w\in L^q_{loc}$ for some $q>2$. Bootstarping this information in the equation gives $v\in L^\infty_{loc}$
(see \cite{Riv3} for a complete description of these arguments).

\subsection{Concluding remarks.}

More jacobian structures or anti-symmetric structures have been discovered in other conformally invariant problems such as  Willmore surfaces \cite{Riv1},  bi-harmonic maps
into manifolds  \cite{LaRi},
 1/2-harmonic maps into manifolds \cite{DR1} and \cite{DR2}...etc. 
 Applying then integrability by compensation results in the spirit of what has been presented above, analysis questions such as the regularity of weak solutions,
 the behavior of sequences of solutions or the compactness of Palais-Smale sequences....have been solved in these works.
 
Moreover, beyond the conformal dimension, while considering the same problems, but in dimension larger than the conformal one, similar approaches can be very efficient  and the same strategy of proofs can sometimes be developed  successfully (see for instance \cite{RiSt}).

\newpage

\section{The Willmore Functional}

\reset

The first appearance of the so-called {\it Willmore functional} goes back to the work of Sophie Germain  on elastic surfaces. Following the main lines of
J. Bernoulli and Euler's studies of the mechanics of sticks in the first half of the XVIII-th century she formulated what she called the {\it fundamental hypothesis} : {\it at one point of the surface
the elastic force which counterbalances the external forces is proportional to the sum of the principal curvature at this point} i.e. what we call the {\it Mean curvature} today.

In  modern elasticity theory (see for instance \cite{LaLi}) the infinitesimal free energy of an elastic membrane at a point  is a product of the area element with a weighted sum of the square of the mean curvature (the Willmore integrand) and the total Gauss curvature (whose integral is a topological invariant for a closed surface).  This has been originally proposed on
physical grounds by G Kirchhoff in 1850 \cite{Ki}. The equation of equilibrium in the absence of external
forces are the critical points to the integral of the free energy that happens to be the {\it Willmore surfaces} since the Gauss curvature times the area element  is locally a jacobian 
and it's integral a null lagrangian.

\subsection{The Willmore Energy of a Surface in ${\R}^3$.} 

Let $S$ be a 2-dimensional submanifold of ${\R}^3$. We assume $S$ to be oriented and we denote by $\vec{n}$ the associated {\it Gauss map} : the unit normal giving this orientation. 

\medskip

The {\it first fundamental form} is the induced metric on $\Sigma$ that we denote by $g$.
\[
\forall p\in S\quad\forall \vec{X},\vec{Y}\in T_pS\quad\quad g(\vec{X},\vec{Y}):=\lf<\vec{X},\vec{Y}\rg>
\]
where $\lf<\cdot,\cdot\rg>$ is the canonical scalar product in ${\R}^3$. The volume form associated to $g$ on $S$ at the point $p$ is given 
by
\[
dvol_g:= \sqrt{ det (g(\p_{x_i},\p_{x_j}))}\ dx_1\wedge dx_2\quad,
\]
where $(x_1,x_2)$ are arbitrary local positive coordinates\footnote{Local coordinates, denoted $(x_1,x_2)$ is
a diffeomorphism $x$ from an open set in ${\R}^2$ into an open set in $\Sigma$. For any point $q$ in this open set of $S$ we shall denote $x_i(q)$ 
 the canonical coordinates in ${\R}^2$ of $x^{-1}(q)$. Finally $\p_{x_i}$ is the vector-field on $S$ given by $\p x/\p x_i$.}

\medskip

The second fundamental form at $p\in S$ is the bilinear map which assigns to a pair of vectors $\vec{X}$,  $\vec{Y}$ in $T_pS$ an orthogonal vector to $T_pS$ that we shall
denote $\vec{\mathbb I}(\vec{X},\vec{Y})$. This normal vector 
''expresses'' how much the Gauss map varies along these directions $\vec{X}$ and $\vec{Y}$. Precisely it is given by
\be
\label{VI.0}
\begin{array}{cccc}
\ds\vec{\mathbb I}_p\ :&\ds T_pS\times T_pS & \ds\longrightarrow &\ds N_pS\\[5mm]
 &\ds (\vec{X},\vec{Y}) &\ds\longrightarrow &\ds - \lf<d\vec{n}_p\cdot\vec{X},\vec{Y}\rg>\ \vec{n}(p)
 \end{array}
 \ee
Extending smoothly $\vec{X}$ and $\vec{Y}$ locally first on $S$ and then in a neighborhood of $p$ in ${\R}^3$, since $<\vec{n}, Y>=0$ on $S$ one has
\be
\label{VI.1}
\begin{array}{rl}
\ds\vec{\mathbb I}_p(\vec{X},\vec{Y}):= \lf<\vec{n}_p\cdot,d\vec{Y}_p\cdot\vec{X}\rg>\ \vec{n}&\ds= d\vec{Y}_p\cdot X-\nabla_{\vec{X}}\vec{Y}\\[5mm]
 &\ds =\ov{\nabla}_{\vec{X}}{\vec{Y}}-\nabla_{\vec{X}}{\vec{Y}}\quad.
\end{array}
\ee
where $\nabla$ is the Levi-Civita connection on $S$ generated by $g$, it is given by $\pi_T(d\vec{Y}_p\cdot X)$ where $\pi_T$ is the orthogonal projection 
onto $T_pS$, and $\ov{\nabla}$ is the the Levi-Civita connection associated to the flat metric and is simply given by $\ov{\nabla}_{\vec{X}}{\vec{Y}}=d\vec{Y}_p\cdot X$.

An elementary but fundamental property of the {\it second fundamental form} says that it is symmetric\footnote{This can be seen combining equation (\ref{VI.1}) with the fact that the two Levi Civita connections $\nabla$ and $\ov{\nabla}$ are symmetric
and hence we have respectively
\[
\nabla_{\vec{X}}\vec{Y}-\nabla_{\vec{Y}}\vec{X}=[X,Y]
\]
and
\[
\ov{\nabla}_{\vec{X}}\vec{Y}-\ov{\nabla}_{\vec{Y}}\vec{X}=[X,Y]\quad.
\]

}. It can then be diagonalized in an orthonormal basis
and the two eigenvalues $\kappa_1$ and $\kappa_2$ are called the {\it principal curvatures} of the surface at $p$. The {\it mean curvature} is then given by
\[
H:=\frac{\kappa_1+\kappa_2}{2} 
\]
and the {\it mean curvature vector} is given by
\be
\label{VI.2}
\vec{H}:=H\ \vec{n}=\frac{1}{2} tr(g^{-1}\ \vec{\mathbb I})=\frac{1}{2} \sum_{ij=1}^2g^{ij} \vec{\mathbb I}(\p_{x_i},\p_{x_j})\quad,
\ee
where $(x_1,x_2)$ are arbitrary local coordinates and $(g^{ij})_{ij}$ is the inverse matrix to $(g(\p_{x_i},\p_{x_j}))$. In particular if $(\vec{e}_1,\vec{e}_2)$ is an orthonormal basis of $T_pS$, (\ref{VI.2})
becomes 
\be
\label{VI.3}
\vec{H}=\frac{\vec{\mathbb I}(\vec{e}_1,\vec{e}_2)+\vec{\mathbb I}(\vec{e}_2,\vec{e}_2)}{2}\quad.
\ee
The Gauss curvature is given by
\be
\label{VI.4}
K:=\frac{\ds det\lf(\vec{n}\cdot\vec{\mathbb I}(\p_{x_i},\p_{x_j})\rg)}{det (g_{ij})}=\kappa_1\ \kappa_2\quad.
\ee
The {\it Willmore functional} of the surface $\Sigma$ is given by
\[
W(S)=\int_{S}|\vec{H}|^2\ dvol_g=\frac{1}{4}\int_{S}|\kappa_1+\kappa_2|^2\ dvol_g
\]
One can rewrite this energy in various ways and get various interpretations of this energy. Assume first that $\Sigma$ is closed (compact without boundary)
the {\it Gauss-Bonnet theorem}\footnote{See for instance \cite{doC1}} asserts that the integral of $K\ dvol_g$ is proportional to a topological invariant of $S$ : $\chi(S)$ the
Euler class of $S$. Precisely one has
\be
\label{VI.5}
\begin{array}{rl}
\ds\int_{S}K\ dvol_g=\int_{S}\kappa_1\ \kappa_2\ dvol_g&=2\pi\ \chi(S)\\[5mm]
 &=4\pi\ (1-g(S))\quad,
 \end{array}
\ee
where $g(S)$ denotes the genus of $S$. Combining the definition of $W$ and this last identity one obtains\footnote{The last identity comes from the fact that, at a point $p$, taking an orthonormal basis $(\vec{e}_1,\vec{e}_2)$ of $T_pS$ one has, since $<d\vec{n},\vec{n}>=0$ :
\[
|d\vec{n}|_g^2=\sum_{i,j=1}^2<d\vec{n}\cdot\vec{e_i},\vec{e}_j>^2=\sum_{i,j=1}^2|\vec{\mathbb I}(\vec{e}_i,\vec{e}_j)|^2=|\vec{\mathbb I}|^2\quad.
\]}
\be
\label{VI.6}
\begin{array}{rl}
\ds W(S)-\pi\ \chi(S)&\ds=\frac{1}{4}\int_{S}(\kappa_1^2+\kappa_2^2)\ dvol_g\\[5mm]
 &\ds=\frac{1}{4}\int_{S}|\vec{\mathbb I}|^2\ dvol_g
 \ds=\frac{1}{4}\int_{S}|d\vec{n}|^2_g\ dvol_g\quad.
 \end{array}
 \ee
Hence modulo the addition of a topological term, the Willmore energy corresponds to the {\it Sobolev homogeneous $\dot{H}^1-$energy of the Gauss map} for the induced
metric $g$.

\medskip

Consider now, again for a closed surface, the following identity based again on Gauss-Bonnet theorem (\ref{VI.5}) :
\be
\label{VI.7}
\begin{array}{rl}
W(S)&\ds=\frac{1}{4}\int_{S}(\kappa_1-\kappa_2)^2\ dvol_g+\int_{S}\kappa_1\ \kappa_2\ dvol_g\\[5mm]
 &\ds=\frac{1}{4}\int_{S}(\kappa_1-\kappa_2)^2\ dvol_g+2\pi\ \chi(S)\quad.
 \end{array}
\ee
Hence, modulo this time the addition of the topological invariant of the surface $2\pi\chi(S)$, the Willmore energy identifies with an energy that penalizes
the lack of {\it umbilicity} when the two principal curvatures differ. The Willmore energy in this form is commonly called {\it Umbilic Energy}.

\medskip

Finally there is an interesting last expression of the Willmore energy that we would like to give here. Consider a conformal parametrization $\vec{\Phi}$ of the surface
$S$ for the conformal structure induced by the metric $g$ (it means that we take  a Riemann surface $\Sigma^2$ and a conformal diffeomorphism $\vec{\Phi}$ from $\Sigma^2$ into $S$).
In local conformal coordinates $(x_1,x_2)$ one has $\p_{x_1}\vec{\Phi}\cdot\p_{x_2}\vec{\Phi}=0$, $|\p_{x_1}\vec{\Phi}|^2=|\p_{x_2}\vec{\Phi}|^2=e^{2\la}$ and the mean-curvature vector is given 
by
\[
\vec{H}=\frac{e^{-2\la}}{2}\Delta\vec{\Phi}=\frac{1}{2}\Delta_g\vec{\Phi}\quad,
\]
where $\Delta$ denotes the negative ''flat'' laplacian, $\Delta=\p^2_{x_1^2}+\p^2_{x_2^2}$, and $\Delta_g$ is the intrinsic negative Laplace-Beltrami operator. With these
notations the Willmore energy becomes
\be
\label{VI.8}
W(\Sigma^2)=\frac{1}{4}\int_{\Sigma^2}|\Delta_g\vec{\Phi}|^2\ dvol_g\quad.
\ee
Hence the Willmore energy identifies to $1/4-$th of the {\it Bi-harmonic Energy}\footnote{ This formulation has also the advantage to show clearly why  Willmore energy is a 4-th order
elliptic problem of the parametrization.} 
of any conformal parametrization $\vec{\Phi}$. 

\subsection{The role of Willmore energy in different areas of sciences and technology.}

As mentioned in the beginning of the present chapter, Willmore energy has been first considered in the framework of the modelization of elastic surfaces and was 
proposed in particular by Poisson \cite{Poi} in 1816 as a Lagrangian from which the equilibrium states of such elastic surfaces could be derived. 
Because of it's simplicity and the fundamental properties it satisfies the Willmore energy has appeared in many area of science for two centuries already.
One can quote the following fields in which the Willmore energy plays an important role.

\begin{itemize}
\item[-] {\it Conformal geometry :} Because of it's conformal invariance (see the next section), the Willmore energy has been introduced in the early XX-th century - see the book of Blaschke \cite{Bla3} - as a fundamental tool in Conformal or {\it M\"obius geometry} of submanifolds.
\item[-] {\it General relativity :} The {\it Willmore energy} arises as being the main term in the expression of the {\it Hawking Mass} of a closed surface in space which measures  the bending of ingoing and outgoing rays of light that are orthogonal to $S$ surrounding the region of space whose mass is to be defined. It is given by
\[
m_{H}(S):=\frac{1}{64\ \pi^{3/2}}\ |S|^{1/2}\ \lf(16\pi-\int_{S}|\vec{H}|^2\ dvol_g\rg)\quad,
\]
where $|S|$ denotes the area of the surface $S$.
\item[-] {\it Cell Biology :} The {\it Willmore Energy} is the main term in the so called {\it Helfrich Energy} in cell biology modelizing the free elastic energy of lipid bilayers membranes
(see \cite{Helf}). The
{\it Helfrich Energy} of a membrane $S$ is given by 
\[
F_H(S):=\int_{S}(2H+C_0)^2\ dvol_g+C_1\ \int_{S}\ dvol_g\quad.
\]
\item[-]{\it Non linear elasticity.} As we mentioned above this is the field where Willmore functional was first introduced. In {\it plate theory }, see for instance \cite{LaLi}, the infinitesimal free energy is a linear combination of the mean curvature to the square and  the Gauss curvature. In a recent mathematical work the Willmore functional is derived as
a Gamma limit of 3-dimensional bending energies for thin materials. There is then a mathematically rigorous consistency between the modelization
of 3 dimensional elasticity and the modelization of the free energy of thin materials given by the Willmore functional see \cite{FJM}.  

\item[-] {\it Optics and Lens design} : The {\it Willmore Energy} in the umbilic form (\ref{VI.7}) arises in methods for the design of multifocal optical elements (see for instance \cite{KaRu}).
\end{itemize}
\subsection{ The Willmore Energy of immersions into ${\R}^m$ and more generalizations.}

\subsubsection{ The Willmore Energy of immersions into ${\R}^m$.}

We extend the definition of the Willmore energy to the class of smooth (at least $C^2$) immersions into ${\R}^m$. Let $\Sigma^2$ be an abstract 2-dimensional oriented manifold
(oriented abstract surface) and let $\vec{\Phi}$ be a $C^2$ immersion of $\Sigma^2$. 
 
 \medskip
 
The {\it first fundamental form} associated to the immersion at the point $p$ is the scalar product $g$ on $T_p\Sigma$ given by $g:=\vec{\Phi}^\ast g_{{\R}^m}$ where $g_{{\R}^m}$
is the canonical metric on ${\R}^m$. Precisely one has
\[
\forall p\in \Sigma^2\quad\forall {X},{Y}\in T_p\Sigma^2\quad\quad g({X},{Y}):=\lf<d\vec{\Phi}\cdot{X},d\vec{\Phi}\cdot Y\rg>
\]
where $\lf<\cdot,\cdot\rg>$ is the canonical scalar product in ${\R}^m$. The volume form associated to $g$ on $\Sigma^2$ at the point $p$ is given 
by
\[
dvol_g:= \sqrt{ det (g(\p_{x_i},\p_{x_j}))}\ dx_1\wedge dx_2\quad,
\]
where $(x_1,x_2)$ are arbitrary local positive coordinates.

\medskip

We shall denote by $\vec{e}$ the map which to a point in $\Sigma^2$ 
assigns the oriented 2-plane\footnote{We denote $\ti{G}_p({\R}^m)$, the Grassman
space of oriented p-planes in ${\R}^m$ that we interpret as the space of unit  simple p-vectors in ${\R}^m$ which is included in the grassmann algebra $\wedge^p{\R}^m$.
} given by the push-forward by $\vec{\Phi}$ of the oriented tangent space $T_p\Sigma^2$. Using a positive orthonormal basis $(\vec{e}_1,\vec{e}_2)$ of $\vec{\Phi}_\ast T_p\Sigma^2$,
an explicit expression of $\vec{e}$ is given by
\[
\vec{e}=\vec{e}_1\wedge\vec{e}_2\quad.
\]
With these notations the Gauss Map which to every point $p$ assigns the oriented $m-2-$orthogonal plane to $\vec{\Phi}_\ast T_p\Sigma^2$ is given by
\[
\vec{n}=\star\vec{e}=\vec{n}_1\wedge\cdots\wedge\vec{n}_{m-2}\quad,
\]
where $\star$ is the Hodge operator\footnote{The Hodge operator on ${\R}^m$ is the linear map from $\wedge^p{\R}^m$ into $\wedge^{m-p}{\R}^m$ which to a $p-$vector $\alpha$
assigns the $m-p-$vector $\star\al$ on ${\R}^m$ such that for any $p-$vector $\beta$ in $\wedge^p{\R}^m$ the following identity holds
\[
\beta\wedge\star\al=<\beta,\al>\ \ep_1\wedge\cdots\wedge\ep_m\quad,
\]
where $(\ep_1,\cdots,\ep_m)$ is the canonical orthonormal basis of ${\R}^m$ and $<\cdot,\cdot>$ is the canonical scalar product on $\wedge^p{\R}^m$.} from $\wedge^2{\R}^m$  into
$\wedge^{m-2}{\R}^m$ and $\vec{n}_1\cdots\vec{n}_{m-2}$ is a positive orthonormal basis of the oriented normal plane to $\vec{\Phi}_\ast T_p\Sigma^2$ (this implies in particular that
$(\vec{e}_1,\vec{e}_2,\vec{n}_1,\cdots,\vec{n}_{m-2})$ is an orthonormal basis of ${\R}^m$).

\medskip

We shall denote by $\pi_{\vec{n}}$ the orthogonal projection onto the $m-2-$plane at $p$ given by $\vec{n}(p)$. Denote by  

\medskip

The {\it second fundamental form}
associated to the immersion $\vec{\Phi}$ is the following map
 \[
\begin{array}{cccl}
\ds\vec{\mathbb I}_p :&\ds T_p\Sigma^2\times T_p\Sigma^2 & \ds\longrightarrow &\ds (\vec{\Phi}_\ast T_p\Sigma^2)^\perp\\[5mm]
 &\ds ({X},{Y}) &\ds\longrightarrow &\ds \vec{\mathbb I}_p({X},{Y}):=\pi_{\vec{n}}(d^2\vec{\Phi}(X,Y))
 \end{array}
 \]
where $X$ and $Y$ are extended smoothly into local smooth vector-fields around $p$. One easily verifies that, though $d^2\vec{\Phi}(X,Y)$ might depend
on these extensions, $\pi_{\vec{n}}(d^2\vec{\Phi}(X,Y))$ does not depend on these extensions and we have then defined a tensor.

\medskip

Let $\vec{X}:=d\vec{\Phi}\cdot X$ and $\vec{Y}:=d\vec{\Phi}\cdot Y$. Denote also by $\pi_T$ the orthogonal projection onto $\vec{\Phi}_\ast T_q\Sigma^2$.
\be
\label{VI.9}
\begin{array}{rl}
\ds\pi_{\vec{n}}(d^2\vec{\Phi}(X,Y))&\ds=d(d\vec{\Phi}_p\cdot X)\cdot Y-\pi_T(d(d\vec{\Phi}_p\cdot X)\cdot Y)\\[5mm]
 &\ds=d\vec{X}\cdot \vec{Y}-\nabla_YX\\[5mm]
 &\ds=\ov{\nabla}_{\vec{Y}}\vec{X}-\nabla_YX\quad.
 \end{array}
\ee
where $\ov{\nabla}$ is the Levi-Civita connection in ${\R}^m$ for the canonical metric and $\nabla$ is the Levi-Civita connection on $T\Sigma^2$ induced by the metric $g$.
Here again, as in the 3-d case in the previous section, from the fact that Levi-Civita connections are symmetric we can deduce the symmetry of the second fundamental form.

\medskip

Similarly to the 3-d case, the {\it mean curvature vector}\footnote{observe that the notion of mean curvature $H$ does not make sense any more in codimension
larger than $1$ unless a normal direction is given.} is given by
\be
\label{VI.10}
\vec{H}:=\frac{1}{2} tr(g^{-1}\ \vec{\mathbb I})=\frac{1}{2} \sum_{ij=1}^2g^{ij} \vec{\mathbb I}(\p_{x_i},\p_{x_j})\quad,
\ee
where $(x_1,x_2)$ are arbitrary local coordinates in $\Sigma^2$ and $(g^{ij})_{ij}$ is the inverse matrix to $(g(\p_{x_i},\p_{x_j}))$.

\medskip

We can now give the general formulation of the {\it Willmore energy} of an immersion $\vec{\Phi}$ in ${\R}^m$ of an abstract surface $\Sigma^2$ :
\[
W(\vec{\Phi}):=\int_{\Sigma^2}|\vec{H}|^2\ dvol_g\quad.
\]

A fundamental theorem by Gauss gives an expression of the intrinsic Gauss curvature at a point $p\in\Sigma^2$ in terms of the second fundamental form
of any immersion of the surface in ${\R}^m$. Precisely this theorem says (see theorem 2.5 chapter 6 of \cite{doC2})
\be
\label{VI.11}
K(p)=\lf<\vec{\mathbb I}(e_1,e_1),\vec{\mathbb I}(e_2,e_2)\rg>-\lf<\vec{\mathbb I}(e_1,e_2),\vec{\mathbb I}(e_1,e_2)\rg>
\ee
where $(e_1,e_2)$ is an arbitrary orthonormal basis of $T_p\Sigma^2$. From this identity we deduce easily
\be
\label{VI.12}
|\vec{\mathbb I}|^2=4|\vec{H}|^2-2K\quad.
\ee
Hence, using Gauss bonnet theorem, we obtain the following expression of the Willmore energy of an immersion into ${\R}^m$
of an arbitrary closed surface 
\be
\label{VI.13}
W(\vec{\Phi})=\frac{1}{4}\int_{\Sigma^2}|\vec{\mathbb I}|^2\ dvol_g+\pi\ \chi(\Sigma^2)\quad.
\ee
Let us take locally about $p$ a {\it normal frame} : a smooth map $(\vec{n}_1,\cdots,\vec{n}_{m-2})$ from a neighborhood $U\subset\Sigma^2$ into $(S^{m-1})^{m-2}$ such that
for any point $q$ in $U$ $(\vec{n}_1(q),\cdots,\vec{n}_{m-2}(q))$ realizes a positive orthonormal basis of $(\vec{\Phi}_\ast T_q\Sigma^2)^\perp$. Then
\[
\pi_{\vec{n}}(d^2\vec{\Phi}(X,Y))=\sum_{\al=1}^{m-2} \lf<d^2\vec{\Phi}(X,Y),\vec{n}_\al\rg>\ \vec{n}_\al\quad,
\]
from which we deduce the following expression - which is the natural extension of (\ref{VI.0}) -
\be
\label{VI.9a}
\vec{\mathbb I}_p({X},{Y})\ds=-\sum_{\al=1}^{m-2} \lf< d\vec{n}_\al\cdot X,\vec{Y}\rg>\ \vec{n}_\al\quad,
\ee
where we denote $\vec{Y}:=d\vec{\Phi}\cdot Y$. Let $(e_1,e_2)$ be an orthonormal basis of $T_p\Sigma^2$, the previous expression of the second fundamental form implies
\be
\label{VI.13a}
|\vec{\mathbb I}_p|^2=\sum_{i,j=1}^2\sum_{\al=1}^{m-2}|<d\vec{n}_\al\cdot e_i,\vec{e}_j>|^2\sum_{i=1}^2 \sum_{\al=1}^{m-2} |<d\vec{n}_\al,\vec{e}_i>|^2
\ee
Observe that 
\be
\label{VI.14}
\begin{array}{rl}
d\vec{n}&\ds = \sum_{\al=1}^{m-2} (-1)^{\al-1} d\vec{n}_\al\wedge_{\beta\ne\al}\vec{n}_\beta=d\vec{n}\\[5mm]
 &\ds =\sum_{i=1}^2 \sum_{\al=1}^{m-2} (-1)^{\al-1} <d\vec{n}_\al,\vec{e}_i>\ \vec{e}_i\wedge_{\beta\ne\al}\vec{n}_\beta
 \end{array}
\ee
$(\vec{e}_i\wedge_{\beta\ne\al}\vec{n}_\beta)$ for $\al=1\cdots m-2$ and $i=1,2$ realizes a free family of $2\, (m-2)$ orthonormal vectors in $\wedge^{m-2}{\R}^m$. Hence
\be
\label{VI.15}
|d\vec{n}|^2_g=\sum_{i=1}^2 \sum_{\al=1}^{m-2} |<d\vec{n}_\al,\vec{e}_i>|^2=|\vec{\mathbb I}_p|^2\quad.
\ee
Combining (\ref{VI.13}), (\ref{VI.13a}) and (\ref{VI.15}) we obtain 
\be
\label{VI.16}
W(\vec{\Phi})=\frac{1}{4}\int_{\Sigma^2}|d\vec{n}|_g^2\ dvol_g+\pi\ \chi(\Sigma^2)\quad,
\ee
which generalizes to arbitrary immersions of closed 2-dimensional surfaces the identity (\ref{VI.6}).

\subsubsection{The Willmore Energy of immersions into a riemannian Manifold $(M^m,\ov{h})$.}

Let $\Sigma^n$ be an abstract $n-$dimensional oriented manifold. 
Let $(M^m,\ov{g})$ be an arbitrary riemannian manifold of dimension larger or equal to $n+1$. 
 The Willmore energy of an immersion $\vec{\Phi}$ into  $(M^m,\ov{g})$
can be defined in a similar way as in the previous subsection by formally replacing the exterior differential $d$ with the Levi-Civita connection $\ov{\nabla}$ on $M$ induced
by the ambient metric $\ov{g}$. Precisely we denote still by $g$ the pull-back of the ambient metric $\ov{g}$ by $\vec{\Phi}$.
\[
\forall p\in \Sigma^2\quad\forall {X},{Y}\in T_p\Sigma^2\quad\quad g({X},{Y}):=\ov{g}(d\vec{\Phi}\cdot{X},d\vec{\Phi}\cdot Y)
\]
The volume form associated to $g$ on $\Sigma^2$ at the point $p$ is still given 
by
\[
dvol_g:= \sqrt{ det (g(\p_{x_i},\p_{x_j}))}\ dx_1\wedge\cdots\wedge dx_n\quad,
\]
where $(x_1,\cdots,x_n)$ are arbitrary local positive coordinates. The second fundamental form associated to the immersion $\vec{\Phi}$ at a point $p$
of $\Sigma^n$ is the following map
 \[
\begin{array}{cccl}
\ds\vec{\mathbb I}_p :&\ds T_p\Sigma^n\times T_p\Sigma^n & \ds\longrightarrow &\ds (\vec{\Phi}_\ast T_p\Sigma^n)^\perp\\[5mm]
 &\ds ({X},{Y}) &\ds\longrightarrow &\ds \vec{\mathbb I}_p({X},{Y}):=\pi_{\vec{n}}(\ov{\nabla}_{\vec{Y}}(d\vec{\Phi}\cdot X))
 \end{array}
 \]
where $\pi_{\vec{n}}$ denotes the orthogonal projection from $T_{\vec{\Phi}(p)}M$ onto the space orthogonal to $\vec{\Phi}_\ast(T_p\Sigma^n)$ with respect to the $\ov{g}$
metric and  that we have denoted by $ (\vec{\Phi}_\ast T_p\Sigma^n)^\perp$. As before we have also used the following notation $\vec{Y}=d\vec{\Phi}\cdot Y$.

\medskip

With this definition of $\vec{\mathbb I}$, the {\it mean curvature vector} is given by 
\be
\label{VI.16a}
\vec{H}=\frac{1}{n}tr(g^{-1}\vec{\mathbb I})=\frac{1}{n}\sum_{i,j=1}^ng^{ij}\vec{\mathbb I}(\p_{x_i},\p_{x_j})\quad,
\ee
where we are using local coordinates $(x_1\cdots x_n)$ on $\Sigma^n$. Hence we have now every elements in order to define
the Willmore energy of $\vec{\Phi}$ which is given by
\[
W(\vec{\Phi}):=\int_{\Sigma^n}|\vec{H}|^2\ dvol_g\quad.
\]
\subsection{Fundamental properties of the Willmore Energy.}

The ''universality'' of Willmore energy which appears in various area of sciences is due to the simplicity of it's expression but also to
the numerous fundamental properties it satisfies. One of the most important properties of this Lagrangian is it's conformal invariance
that we will present in the next subsections. Other important properties relate the amount of Willmore energy of an immersion to the
''complexity'' of this immersion. For instance, in the second subsection we will present the {\it Li-Yau $8\pi$ threshold} below which every immersion of a closed
surface is in fact an embedding. We will end up this section by raising the fundamental question in calculus of variation regarding the existence
of minimizers of the Willmore energy under various constraints. At this occasion we will state the so-called {\it Willmore conjecture}.

\subsubsection{Conformal Invariance of the Willmore Energy.}

The conformal invariance of the Willmore energy was known since the work of Blaschke \cite{Bla3} in 3 dimension, in the general case it is a consequence of the following theorem due to Bang Yen Chen \cite{Che}.

\begin{Th}
\label{th-VI.1}
Let $\vec{\Phi}$ be the immersion of an $n-$dimensional manifold $\Sigma^n$ into a riemannian manifold $(M^m,g)$. Let $\mu$ be a smooth function
in $M^m$ and let $h$ be the conformally equivalent metric given by $h:=e^{2\mu}\ g$. We denote by $\vec{H}^{g}$ and
$\vec{H}^{{h}}$ the mean curvature vectors of the immersion $\vec{\Phi}$ respectively in $(M^m,{g})$ and $(M^m,{h})$. We also denote by $K^{{g}}$
and $K^{{h}}$ the extrinsic scalar curvatures respectively of $(\Sigma^n,\vec{\Phi}^\ast{g})$ and $(\Sigma^n,\vec{\Phi}^\ast{h})$. With the previous notations the following identity holds
\be
\label{VI.17}
e^{2\mu}\ \lf(|\vec{H}^h|^2_h-K^h+\ov{K}^h\rg)=|\vec{H}^g|^2_g-K^g+\ov{K}^g\quad.
\ee
where $\ov{K}^g$ (resp. $\ov{K}^h$) is the sectional curvature of the subspace $\vec{\Phi}_\ast T_p\Sigma^2$ in the manifold $(M^m,g)$ (resp. $(M^m,h)$). 
$K^g-\ov{K}^g$ ( resp. $K^h-\ov{K}^h$) is also\hfill $\Box$
\end{Th}
\begin{Rm}
\label{rm-VI.1b}
$K^g-\ov{K}^g$ ( resp. $K^h-\ov{K}^h$) is also called the extrinsic scalar curvature of the immersion $\vec{\Phi}$ os $\Sigma^n$ into $(M^m,g)$ (resp. into $(M^m,h)$).
\end{Rm}
{\bf Proof of theorem~\ref{th-VI.1}.}

Let $\ov{\nabla}^g$ (resp. $\ov{\nabla}^h$) be the Levi-Civita connection induced by the metric $g$ (resp. $h$) on $M^m$. Let $\nabla^g$ (resp. $\nabla^h$) be the
Levi-Civita connection induced by the restriction of the metric $g$ (resp. $h$) on $\Sigma^n$. By an abuse of notation we shall still write $g$ (resp. $h$) for
the pull back by $\vec{\Phi}$ on $\Sigma^n$ of the restriction of the metric $g$ (resp. $h$). As in the previous section for any vector $X\in T\Sigma^n$ we denote by $\vec{X}$ the push forward by $\vec{\Phi}$ of $X$ : $\vec{X}:=d\vec{\Phi}\cdot X$. With these notations, the second fundamental form $\vec{\mathbb I}^g$ of the immersion $\vec{\Phi}$ of $\Sigma^n$ into $(M^m,g)$
at a point $p\in\Sigma^n$
is defined as follows
\be
\label{VI.17a}
\forall X,Y\in T_p\Sigma^n\quad\quad\vec{\mathbb I}^g(X,Y)=\ov{\nabla}^g_{\vec{Y}}\vec{X}-{\nabla^g_YX}\quad,
\ee
and similarly the second fundamental form $\vec{\mathbb I}^g$ of the immersion $\vec{\Phi}$ of $\Sigma^n$ into $(M^m,g)$
at a point $p\in\Sigma^n$
is given by
\be
\label{VI.17b}
\forall X,Y\in T_p\Sigma^n\quad\quad\vec{\mathbb I}^h(X,Y)=\ov{\nabla}^h_{\vec{Y}}\vec{X}-\nabla^h_YX\quad.
\ee
Denote by $\Gamma^{g,k}_{ij}$ (resp. $\Gamma^{h,k}_{ij}$) the Kronecker symbols of the metric $g$ (resp. $h$) in some local
coordinates $(x_1,\cdots,x_n)$ in a neighborhood of a point $p$ in $\Sigma^n$.
Using the explicit form (\ref{II.4za}) of these Kronecker symbols that we already recalled in the first part of the course  we have, for all $i,j,k\in\{1,\cdots,n\}$,
\be
\label{VI.18}
\begin{array}{rl}
\Gamma_{ij}^{h,k}&\ds=\frac{1}{2}\sum_{s=1}^nh^{ks}\ \lf[\p_{x_j}h_{si}+\p_{x_i}h_{sj}-\p_{x_s}h_{ij}\rg]\\[5mm]
 &\ds=\sum_{s=1}^n\frac{1}{2}g^{ks}\ \lf[\p_{x_j}g_{si}+\p_{x_i}g_{sj}-\p_{x_s}g_{ij}\rg]\\[5mm]
 &\quad\ds+\sum_{s=1}^ng^{ks}\ \lf[\p_{x_j}\mu\ g_{si}+\p_{x_i}\mu\ g_{sj}-\p_{x_s}\mu\ g_{ij}\rg]\quad.
 \end{array}
\ee
Using the fact that $\sum_{s=1}^ng^{ks}g_{si}=\delta^k_i$ we deduce from the previous identity (\ref{VI.18}) that for any triple $i,j,k$ in $\{1\cdots n\}$
\be
\label{VI.19}
\Gamma_{ij}^{h,k}-\Gamma_{ij}^{g,k}=\delta^k_i\p_{x_j}\mu+\delta^k_j\p_{x_i}\mu-\sum_{s=1}^n\p_{x_s}\mu\ g^{ks}\,g_{ij}\quad.
\ee
Using the expression of the Levi-Civita in local coordinates by the mean of Christoffel symbols\footnote{In local coordinates
one has
\[
\nabla^g_YX=\sum_{k=1}^n\ \lf(\sum_{i=1}^nY^i\ \p_{x_i}X^k+\sum_{i,j=1}^n\Gamma^{g,k}_{ij}\ X^i\ X^j\rg)\ \p_{x_k}\quad.
\]
} we have
\be
\label{VI.20}
\nabla^{h}_YX-\nabla^g_YX=\sum_{k=1}^n\,\lf[\sum_{i,j=1}^n\lf(\Gamma_{ij}^{h,k}-\Gamma_{ij}^{g,k}\rg)\,X^i\,Y^j\rg]\,\p_{x_k}\quad.
\ee
Combining (\ref{VI.19}) and (\ref{VI.20}) we obtain
\be
\label{VI.21}
\nabla^{h}_YX-\nabla^g_YX=d\mu\cdot X\ Y+d\mu\cdot Y\ X-g(X,Y)\ U\quad,
\ee
where $U\in T_p\Sigma^n$ is the dual vector to the 1-form $d\mu$ for the metric $g$ which is given by
\be
\label{VI.22}
\forall Z\in T_p\Sigma^n\quad\quad g(U,Z)=d\mu\cdot Z\quad.
\ee
Similarly we have
\be
\label{VI.23}
\ov{\nabla}^{h}_{\vec{Y}}\vec{X}-\vec{\nabla}^g_{\vec{Y}}\vec{X}=d\mu\cdot \vec{X}\ \vec{Y}+d\mu\cdot \vec{Y}\ \vec{X}-g(\vec{X},\vec{Y})\ \vec{U}\quad,
\ee
where $\vec{V}\in T_{\vec{\Phi}(p)}M^m$ is the dual vector to the 1-form $d\mu$ for the metric $g$ which is given by
\be
\label{VI.24}
\forall \vec{Z}\in T_{\vec{\Phi}(p)}M^m\quad\quad g(\vec{V},\vec{Z})=d\mu\cdot \vec{Z}\quad.
\ee
Combining (\ref{VI.17a}), (\ref{VI.17b}), (\ref{VI.22}) and (\ref{VI.24}) we obtain
\be
\label{VI.25}
\vec{\mathbb I}^h(X,Y)-\vec{\mathbb I}^g(X,Y)=-g(X,Y)\,\vec{W}\quad .
\ee
where $\vec{W}:=\vec{V}-\vec{U}$. From (\ref{VI.22}) and (\ref{VI.24}) we obtain that
\[
\forall \vec{Z}\in \vec{\Phi}_\ast (T_p\Sigma^n)\quad\quad g(\vec{Z},\vec{W})=0\quad.
\]
hence $\vec{U}$ is the orthogonal projection of $\vec{V}$ onto $\vec{\Phi}_\ast (T_p\Sigma^n)$ and $\vec{W}$
is the orthogonal projection onto the normal space to $\vec{\Phi}(\Sigma^n)$. 

\medskip

Let $\vec{\xi}$ be a unit normal vector to $\Sigma^n$ in $(M^m,g)$. For any such a vector $\vec{\xi}$, the $g-$scalar product between
$\vec{\xi}$ and $\vec{\mathbb I}^g$ at a point $p\in\Sigma^n$ produces a symmetric bilinear two-form $g(\vec{\xi},\vec{\mathbb I}^g(\cdot,\cdot))$ on $(T_p\Sigma^n,g)$.
We denote by $\kappa_1^g(\xi)$ and $\kappa_2^g(\xi)$ the eigenvalues of this symmetric two form for the $g$ scalar product.
Let $({e}_1,\cdots,e_n)$ be an orthonormal basis to $(T_p\Sigma^n,g)$, we have
\be
\label{VI.26}
\begin{array}{rl}
\ds g\lf(\vec{\xi},\sum_{i=1}^n\vec{\mathbb I}^g(e_i,e_i)\rg)&\ds=tr_g\lf(g(\vec{\xi},\vec{\mathbb I}^g(\cdot,\cdot))\rg)\\[5mm]
 &\ds=\kappa_1(\vec{\xi})+\cdots +\kappa_n(\vec{\xi})\quad.
 \end{array}
\ee
Let $(\vec{n}_1,\cdots,\vec{n}_{m-2})$ be an orthonormal basis for the metric $g$ of the normal space to $\vec{\Phi}_\ast(T_p\Sigma^n)$. Combining the expression
of the mean curvature vector (\ref{VI.16a}) together with (\ref{VI.26}) we obtain the following expression of the mean curvature vector $\vec{H}_g$
of the immersion $\vec{\Phi}$ into $(M^m,g)$ :
\be
\label{VI.27}
\vec{H}^g=\frac{1}{n}\sum_{\al=1}^{m-2}\sum_{i=1}^n\kappa_i^g(\vec{n}_\al)\ \vec{n}_\al\quad.
\ee
Hence 
\be
\label{VI.27a}
\begin{array}{rl}
|\vec{H}^g|_g^2&\ds=\frac{1}{n^2}\sum_{\al=1}^{m-2}\lf|\sum_{i=1}^n\kappa_i^g(\vec{n}_\al)\rg|^2\\[5mm]
 &\ds=\frac{1}{n^2}\sum_{\al=1}^{m-2}\sum_{i<j}\frac{\ds\lf(\kappa_i^g(\vec{n}_\al)-\kappa_j^g(\vec{n}_\al)\rg)^2}{n-1}\\[5mm]
  &\ds\quad+\frac{2}{n\,(n-1)}\sum_{\al=1}^{m-2}\sum_{i<j}\kappa_i^g(\vec{n}_\al)\ \kappa_j^g(\vec{n}_\al)
\end{array}
\ee
The Gauss theorem (theorem 2.5 chap.6 in \cite{doC2}) gives an extrinsic expression in terms of the second fundamental form
 of the sectional curvature of the manifold $(\Sigma^n,g)$ at an oriented tangent two-plane given by a pair $(e_i,e_j)$ of normal vectors in $(T_p\Sigma^n,g)$
orthogonal to each other :
\be
\label{VI.27b}
S^g(e_i,e_j)-\ov{S}^g(e_i,e_j)=\vec{\mathbb I}(e_i,e_i)\cdot\vec{\mathbb I}(e_j,e_j)-|\vec{\mathbb I}(e_i,e_j)|^2\quad,
\ee
where $\ov{S}^g(e_i,e_j)$ is the sectional curvature of $(M^m,g)$ at the oriented two plane given by the orthonormal pair $(d\vec{\Phi}\cdot e_i, d\vec{\Phi}\cdot e_j)$. The intrinsic formulation of $S^g(e_i,e_j)$ using the riemann tensor  $R^g$ of $(\Sigma^n,\vec{\Phi}^\ast g)$ is given by
\be
\label{VI.27c}
S^g(e_i,e_j)=g\lf(R^g(e_i,e_j)\, e_j,e_i)\rg)
\ee
The scalar curvature of $(\Sigma^n,g)$ is the following average of the sectional curvatures of the different 2-planes defined
by pairs of elements of a given orthonormal basis $(e_1,\cdots, e_n)$ of $(T_p\Sigma^n,g)$ :
\be
\label{VI.27d}
K^g:=\frac{2}{n\,(n-1)}\sum_{i<j}S^g(e_i,e_j)\quad.
\ee
Combining
the identities (\ref{VI.27b}) and (\ref{VI.27d}) and taking the orthonormal basis $(\vec{n}_\al)_{\al=1\cdots m-2}$ that we already fixed above,
we obtain
\be
\label{VI.27e}
\begin{array}{l}
\ds K^g- \ov{K}^g\\[5mm]
\ds=\frac{2}{n\,(n-1)}\sum_{\al=1}^{m-2}\sum_{i<j}g(\vec{n}_\al,\vec{\mathbb I}(e_i,e_i))\ g(\vec{n}_\al,\vec{\mathbb I}(e_j,e_j))-|g(\vec{n}_\al,\vec{\mathbb I}(e_i,e_j))|^2\quad.
\end{array}
\ee
It is a classical fact that for any matrix $\sum_{ij}a_{ii}\,a_{jj}-a_{ij}^2$ is invariant under the change of  basis hence for a given $\al$ we can choose $(e_1\cdots e_n)$
to be the orthonormal basis in which the symmetric bilinear form $g(\vec{n}_\al,\vec{\mathbb I}(\cdot,\cdot))$ is diagonal. We then obtain that
\be
\label{VI.27f}
\begin{array}{l}
\ds\sum_{i<j} g(\vec{n}_\al,\vec{\mathbb I}(e_i,e_i))\ g(\vec{n}_\al,\vec{\mathbb I}(e_j,e_j))-|g(\vec{n}_\al,\vec{\mathbb I}(e_i,e_j))|^2\\[5mm]
\ds\quad\quad=\sum_{i<j}\kappa_i(\vec{n}_\al)\ \kappa_j(\vec{n}_\al)\quad.
\end{array}
\ee
Thus, finally, we obtain combining (\ref{VI.27e}) and (\ref{VI.27f})
\be
\label{VI.27g}
\ds\quad\quad K^g-\ov{K}^g=\frac{2}{n\,(n-1)}\sum_{i<j}\kappa_i(\vec{\xi})\,\kappa_j(\vec{\xi})\quad.
\ee
Hence (\ref{VI.27a}) becomes
\be
\label{VI.27h}
|\vec{H}^g|_g^2-K^g+\ov{K}^g=\frac{1}{n^2}\sum_{\al=1}^{m-2}\sum_{i<j}\frac{\ds\lf(\kappa_i^g(\vec{n}_\al)-\kappa_j^g(\vec{n}_\al)\rg)^2}{n-1}\quad.
\ee

\medskip

Since $h=e^{2\mu}\, g$, the normal space to $\vec{\Phi}_\ast(T_p\Sigma^n)$ for $g$ coincides with the one for $h$, moreover 
$(e^{-\mu}\,\vec{n}_1,\cdots,e^{-\mu}\,\vec{n}_{m-2})$ realizes an orthonormal basis for the metric $g$. Hence the corresponding
expression to (\ref{VI.27}) for the mean curvature vector $\vec{H}^h$ of the same immersion $\vec{\Phi}$ but into $(M^m,h)$ reads
\be
\label{VI.28}
\vec{H}^h=\frac{e^{-\mu}}{n}\sum_{\al=1}^{m-2}\sum_{i=1}^n\kappa_i^h(e^{-\mu}\,\vec{n}_\al)\ \vec{n}_\al\quad.
\ee
A vector $e\in T_p\Sigma^2\setminus\{0\}$ is an eigenvector for $g(\vec{n}_\al,\vec{\mathbb I}^g(\cdot,\cdot))$ w.r.t. the metric $g$ if and only if there exists
a real number $\kappa$ such that
\be
\label{VI.29}
g(\vec{n}_\al,\vec{\mathbb I}^g(\cdot, e))=\kappa\ g(\cdot, e)\quad.
\ee
This implies that
\[
h(\vec{n}_\al,\vec{\mathbb I}^h(\cdot, e)+g(\cdot,e)\ \vec{W})=\kappa\ h(\cdot,e)\quad.
\]
In other words we have obtained
\be
\label{VI.30}
h(e^{-\mu}\,\vec{n}_\al,\vec{\mathbb I}^h(\cdot,e))=e^{-\mu}\ \lf[\kappa-g(\vec{n}_\al,\vec{W})\rg]\ h(\cdot,e)\quad.
\ee
This later identity says then that $e$ is also an eigenvector of $h(e^{-\mu}\,\vec{n}_\al,\vec{\mathbb I}^h(\cdot,\cdot))$ with eigenvalue
$e^{-\mu}\ \lf[\kappa-g(\vec{n}_\al,\vec{W})\rg]$. We then have 
\be
\label{VI.31}
\kappa_i^h(e^{-\mu}\,\vec{n}_\al)= e^{-\mu}\ \lf[\kappa_i^g(\vec{n}_\al)-g(\vec{n}_\al,\vec{W})\rg]\quad.
\ee
This implies that $\forall \al\in\{1\cdots m-2\}$ and $\forall i,j\in\{1\cdots n\}$
\be
\label{VI.32}
\lf|\kappa_i^g(\vec{n}_\al)-\kappa_j^g(\vec{n}_\al)\rg|^2=e^{2\mu}\, \lf|\kappa_i^h(e^{-\mu}\,\vec{n}_\al)-\kappa_j^h(e^{-\mu}\,\vec{n}_\al)\rg|^2\quad.
\ee
Combining (\ref{VI.27h}) (which is also valid for $h$ of course) with (\ref{VI.32}) gives (\ref{VI.17}) and theorem~\ref{th-VI.1} is proved.
\hfill $\Box$

\medskip

We will make use of the following corollary of theorem~\ref{th-VI.1}.

\medskip

\begin{Co}
\label{co-VI.2}
Let $\Sigma^2$ be a closed smooth oriented 2-dimensional manifold and let $\vec{\Phi}$ be an immersion 
of $\Sigma^2$ into an oriented riemannian manifold $(M^m,g)$. Let $\Psi$ be a positive conformal diffeomorphism
from $(M^m,g)$ into another riemannian oriented manifold $(N^m,k)$ then we have the following pointwise identity
everywhere on $\Sigma^2$
\be
\label{VI.33}
\begin{array}{l}
\ds\lf[|\vec{H}^{\vec{\Phi}^\ast g}|^2-K^{\vec{\Phi}^\ast g}+\ov{K}^g\rg]\ dvol_{\vec{\Phi}^\ast g}\\[5mm]
\ds\quad\quad=\lf[|\vec{H}^{(\Psi\circ\vec{\Phi})^\ast k}|^2-K^{(\Psi\circ\vec{\Phi})^\ast k}+\ov{K}^k\rg]\ dvol_{(\Psi\circ\vec{\Phi})^\ast k}
\end{array}
\ee
where $\ov{K}^g$ (resp. $\ov{K}^k$) is the sectional curvature of the $2-$plane $\vec{\Phi}_\ast T\Sigma^2$ in  $(M^m,g)$  (resp. of the two
plane $\Psi_\ast\vec{\Phi}_\ast T\Sigma^2$ in $(N^m,k)$). In particular the following equality holds
\be
\label{VI.34}
\begin{array}{l}
\ds W(\vec{\Phi})+\int_{\Sigma^2}\ov{K}^g\ dvol_{\vec{\Phi}^\ast g}\\[5mm]
\ds\quad\quad=W(\Psi\circ\vec{\Phi})+\int_{\Sigma^2}\ov{K}^k\ dvol_{(\Psi\circ\vec{\Phi})^{\ast} g}\quad.
\end{array}
\ee
\hfill $\Box$
\end{Co}
{\bf Proof of corollary~\ref{co-VI.2}.}

By definition $\Psi$ realizes an isometry between $(M^m,\Psi^\ast k)$ and $(N^m,k)$. Let $\mu\in {\R}$ such that $e^\mu\ g=\Psi^\ast k$. We can apply the previous
theorem with $h=\Psi^\ast k$ and we obtain
\be
\label{VI.35}
\lf[|\vec{H}^{\vec{\Phi}^\ast g}|^2-K^{\vec{\Phi}^\ast g}\rg]= e^{2\mu}\ \lf[|\vec{H}^{(\Psi\circ\vec{\Phi})^\ast k}|^2-K^{(\Psi\circ\vec{\Phi})^\ast k}\rg]\quad.
\ee
It is also clear that the volume form $dvol_g$ given by the restriction to $\Sigma^2$ of the metric $g$ is equal to $e^{-2\mu}\, dvol_{\Psi^\ast k}$
where $dvol_{\Psi^\ast k}$ is equal to the volume form given by the restriction to $\Sigma^2$ of the metric $\Psi^\ast k$. Hence combining this last fact with (\ref{VI.33}).
(\ref{VI.34}) is obtained by integrating (\ref{VI.33}) over $\Sigma^2$, the scalar curvature terms canceling each other on both sides of the identity due to Gauss Bonnet theorem.
Corollary~\ref{co-VI.2} is then proved. 
\hfill $\Box$
\subsubsection{Li-Yau Energy lower bounds and the Willmore conjecture.}

It is a well known fact that the integral of the curvature $\kappa_{\vec{\Phi}}$ of the immersion $\vec{\Phi}$ of a closed curved $\Gamma$ ($\p\Gamma=\emptyset$)  in ${\R}^m$ 
is always 
larger than $2\pi$ :
\be
\label{VI.36}
F(\vec{\Phi})=\int_{\Gamma}\kappa_{\vec{\Phi}}\ge 2\pi\quad,
\ee
and is equal to $2\pi$ if and only if $\Gamma\simeq S^1$ and it's image by $\vec{\Phi}$ realizes a convex planar curve. Using the computations in the previous section
one easily verifies that $F$ is conformally invariant $F(\Psi\circ\vec{\Phi})=F(\vec{\Phi}) $ and plays a bit the role of a 1-dimensional version of the Willmore functional.
As the immersion gets more complicated one can expect the lower bound (\ref{VI.36}) to increase. For instance a result by Milnor says that if the immersion $\vec{\Phi}$ is 
a \underbar{knotted} embedding then the lower bound (\ref{VI.36}) is multiplied by 2 :
\be
\label{VI.37}
\int_{\Gamma}\kappa_{\vec{\Phi}}\ge 4\pi\quad.
\ee
It is natural to think that this general philosophy can be transfered to the Willmore functional : if the surface $\Sigma^2$ and the immersion $\vec{\Phi}$ become ''more complicated''
then there should exist increasing general lower bounds for $W(\vec{\Phi})$. First we give the result corresponding to the first lower bound (\ref{VI.36}) for curves. It was first proved by 
Thomas Willmore for $m=3$ and by Bang-Yen Chen \cite{BYC2} in the general codimension case.

\begin{Th}
\label{th-VI.3}
Let $\Sigma^2$ be a closed surface and let $\vec{\Phi}$ be a smooth immersion of $\Sigma^2$ into ${\R}^m$. Then the following inequality holds
\be
\label{VI.38}
W(\vec{\Phi})=\int_{\Sigma^2}|\vec{H}_{\vec{\Phi}}|^2\ dvol_{\vec{\Phi}^\ast g_{{\R}^m}}\ge 4\pi\quad.
\ee 
Moreover equality holds if and only if $\Sigma^2$ is a 2-sphere and $\vec{\Phi}$ realizes - modulo translations and dilations - an embedding onto $S^2$ the unit 
sphere of ${\R}^3\subset{\R}^m$.\hfill$\Box$
\end{Th}
{\bf Proof of theorem~\ref{th-VI.3}}
 Denote by $N\Sigma^2$ the pullback by $\vec{\Phi}$ of the normal sphere bundle to $\vec{\Phi}(\Sigma^2)$ in ${\R}^m$ made of the unit normal vectors
 to $\vec{\Phi}_\ast T\Sigma^2$. Denote by $\vec{G}=(G_1,\cdots,G_m)$ the map from $N\Sigma^2$ into $S^{m-1}$ the unit sphere which to an element in $N\Sigma^2$ assigns the corresponding
 unit vector in $S^{m-1}$. For any $k\in {\N}$, $k\ge 1$, we denote
 \[
 \om_{S^{k-1}}=\frac{1}{|S^{k-1}|}\sum_{j=1}^k(-1)^{j-1}x_j\, \wedge_{l\ne j}dx_l\quad.
 \]
 Observe that 
 \[
 \int_{S^{k-1}}\om_{S^{k-1}}=1\quad.
 \]
 Locally on $\Sigma^2$ we can choose an orthonormal positively oriented frame of the normal plane : $(\vec{n}_1\cdots\vec{n}_{m-2})$. We also chose locally
 an orthonormal tangent frame $(\vec{e}_1,\vec{e}_2)$. This normal frame permits locally, over an open disk $U\subset\Sigma^2$, to trivialize $N\Sigma^2$ and the map $\vec{G}$ can be seen as a map
 from $U\times S^{m-3}$ into $S^{m-1}$ : 
 \[
 \vec{G}(p,s)=\sum_{\al=1}^{m-2}s_\al\ \vec{n}_\al(p)\quad,
 \]
 where $\sum_\al s^2_\al=1$. Let $p\in U$, in order to simplify the notations we may 
 assume that $(\vec{e}_1,\vec{e}_2,\vec{n}_1\cdots \vec{n}_{m-2})$ at $p$ coincides with the canonical basis of ${\R}^m$. Hence we have $G_1(p)=G_2(p)=0$ and
 $G_{\al+2}(p)=s_\al$ for $\al=1\cdots m-2$.
 \[
 dG_1(p)=\sum_{\al=1}^{m-2} s_\al\ <d\vec{n}_\al,\vec{e}_1>\quad,
 \]
 \[
 dG_2(p)=\sum_{\al=1}^{m-2} s_\al\ <d\vec{n}_\al,\vec{e}_2>\quad,
 \]
and
\[
\forall\ \al=1\cdots m-2\quad\quad dG_{\al+2}=ds_\al+\sum_{\beta=1}^{m-2}s_\beta\  <d\vec{n}_\beta,\vec{n}_\al>\quad.
\]
Thus
\be
\label{VI.38a}
\begin{array}{l}
\ds\vec{G}^\ast\om_{S^{m-1}}(p)=\frac{1}{|S^{m-1}|} \sum_{\al,\beta=1}^{m-2}s_\al\,s_\beta\  <d\vec{n}_\al,\vec{e}_1>\wedge<d\vec{n}_\al,\vec{e}_2>\\[5mm]
\ds\quad\quad\quad\quad\wedge\sum_{j=1}^{m-2}(-1)^{j+1}\ s_j\ 
\wedge_{l\ne j}ds_l\\[5mm]
\ds\quad\quad\quad=\frac{|S^{m-3}|}{|S^{m-1}|}\ <d\vec{e}_1,\vec{G}>\wedge<d\vec{e}_2,\vec{G}>\wedge\om_{S^{m-3}}(s)
\end{array}
\ee
Observe that
\be
\label{VI.38b} 
<d\vec{e}_1,\vec{G}>\wedge<d\vec{e}_2,\vec{G}> =det\lf(<\vec{G},\vec{\mathbb I}>\rg)\ dvol_{\vec{\Phi}^\ast g_{{\R}^m}}
\ee
and denoting $\Om$ the closed $m-3-$form on $N\Sigma^2$ which is invariant under the action of $SO(m-2)$ over the fibers and 
whose integral over each fiber is equal to one\footnote{This means that $\Om$ is the {\it Thom class} of the normal bundle $N\Sigma^2$ it coincides in particular with $\om_{S^{m-3}}$ in the coordinates $s$.}, combining
(\ref{VI.38a}) and (\ref{VI.38b}),
we have proved finally that\footnote{where we are using the fact that $|S^{m-3}|/|S^{m-1}|=(m-2)/2\pi$.}
\be
\label{VI.38c}
\vec{G}^\ast\om_{S^{m-1}}=\frac{m-2}{2\pi}det\lf(<\vec{G},\vec{\mathbb I}>\rg)\ dvol_{\vec{\Phi}^\ast g_{{\R}^m}}\wedge\Om
\ee
The form $<\vec{G},\vec{\mathbb I}>$ is bilinear symmetric, if $\kappa_1$, $\kappa_2$ are it's eigenvalues, since $4\kappa_1\ \kappa_2\le(\kappa_1+\kappa_2)^2$, we deduce the pointwise inequality
\be
\label{VI.38d}
\begin{array}{l}
\ds\vec{G}^\ast\om_{S^{m-1}}\le \frac{m-2}{2\pi}\lf|\frac{tr_g(<\vec{G},\vec{\mathbb I}>)}{2}\rg|^2\ dvol_{\vec{\Phi}^\ast g_{{\R}^m}}\wedge\Om\\[5mm]
\ds\quad\quad =\frac{m-2}{2\pi}\lf|<\vec{G},\vec{H}_{\vec{\Phi}}>\rg|^2\ dvol_{\vec{\Phi}^\ast g_{{\R}^m}}\wedge\Om
\end{array}
\ee
Assume that at $p$ the vector $\vec{H}_{\vec{\Phi}}$ is parallel to $\vec{n}_1$, one has, using the coordinates $s$ we introduced,
\[
\begin{array}{l}
\ds\int_{Fiber}\lf|<\vec{G},\vec{H}_{\vec{\Phi}}>\rg|^2\ \Om=\int_{S^{m-3}} s_1^2\ |\vec{H}_{\vec{\Phi}}|^2\ \om_{S^{m-3}}\\[5mm]
\ds\quad\quad=\frac{|\vec{H}_{\vec{\Phi}}|^2}{m-2}\ \sum_{j=1}^{m-2}\int_{S^{m-3}}s_j^2\ \om_{S^{m-3}}=\frac{|\vec{H}_{\vec{\Phi}}|^2}{m-2}\quad.
\end{array}
\]
Denote $N^+\Sigma^2$ the subset of $N\Sigma^2$ on which $det\lf(<\vec{G},\vec{\mathbb I}>\rg)$ is non-negative. We have then proved
\be
\label{VI.38e}
\int_{N^+\Sigma^2}\vec{G}^\ast\om_{S^{m-1}}\le \frac{1}{2\pi}\int_{\Sigma^2}|\vec{H}_{\vec{\Phi}}|^2\ \ dvol_{\vec{\Phi}^\ast g_{{\R}^m}}\quad.
\ee
Now we claim that each points in $S^{m-1}$ admits at least two preimages by $\vec{G}$ at which $det\lf(<\vec{G},\vec{\mathbb I}>\rg)\ge 0$ unless the immersed surface is included
in an hyperplane of ${\R}^{m}$. 

Indeed, let $\vec{\xi}\in {S^{m-1}}$ and consider the affine hyper-plane $\Xi_{a}$ given by $<x,\vec{\xi}>=a$
Let $$a_+=\max\{a\in{\R}\ ;\ \Xi_{a}\cap\vec{\Phi}(\Sigma^2)\ne\emptyset\}$$ and $$a_-=\min\{a\in{\R}\ ;\ \Xi_{a}\cap\vec{\Phi}(\Sigma^2)\ne\emptyset\}$$ 
Assume $a_+=a_-$ then the surfaced is immersed in the hyper-plane $\Xi_{a_+}=\Xi_{a_-}$ and thus the claim is proved. If $a_+>a_-$, it is clear that $\vec{\Phi}(\Sigma^2)$
is tangent to $\Xi_{a_+}$ at a point $\vec{\Phi}(p_+)$ and hence $\vec{\xi}$ belongs to the normal space to $\vec{\Phi}_\ast T_p\Sigma^2$ which means in other words 
that $\vec{\xi}\in \vec{G}(N_{p_+}\Sigma^2)$. Similarly,  $\vec{\Phi}(\Sigma^2)$
is tangent to $\Xi_{a_-}$ at a point $\vec{\Phi}(p_-)$ and hence $\vec{\xi}\in \vec{G}(N_{p_-}\Sigma^2)$. Since $a_-<a_+$, $\vec{\Phi}(p_-)\ne\vec{\Phi}(p_+)$ and then we have proved that $\vec{\xi}$ admits at least two prei-mages by $\vec{G}$. We claim now that $det\lf(<\vec{G},\vec{\mathbb I}>\rg)\ge 0$ at these points. Since the whole surface is contained in one of the 
half space given by the affine plane $(\vec{\Phi}(p_\pm),\vec{\xi})$ we have that $p\rightarrow<\vec{\xi},\vec{\Phi}(p)-\vec{\Phi}(p_\pm)>$ has an absolute maximum (for $+$) or minimum (for $-$) at $p=\pm p$. In both cases the determinant of the 2 by 2 Hessian has to be non-negative :
\[
det\lf(<\vec{\xi},\p^2_{x_i, x_j}\vec{\Phi}>\rg)\ge 0
\]
from which we deduce, since $\vec{\xi}$ is normal to $\Phi_\ast T_{p_\pm}\Sigma^2$,
\[
det\lf(<\vec{\xi},\vec{\mathbb I}>\rg)\ge 0 
\]
\medskip

Assume ${\R}^m$ is the smallest affine subspace in which $\vec{\Phi}({\Sigma}^2)$ is included. Then , using Federer's co-area formula, the claim implies
\be
\label{VI.38f}
2\le\int_{S^{m-1}}  \#\{\vec{G}^{-1}(\vec{\xi}),\ \xi\in N^+\Sigma^2\}\ \om_{S^{m-1}}(\vec{\xi})=\int_{N^+\Sigma^2}\vec{G}^\ast\om_{S^{m-1}}\quad.
\ee
Combining (\ref{VI.38e}) and (\ref{VI.38f}) gives (\ref{VI.38}) :
\[
4\pi\le\int_{\Sigma^2}|\vec{H}|^2\ \ dvol_{\vec{\Phi}^\ast g_{{\R}^m}}
\]
Assume this inequality is in fact an equality. Then all the inequalities above are equalities and in particular we have that 
\[
\lf|\frac{tr_g(<\vec{G},\vec{\mathbb I}>)}{2}\rg|^2=det\lf(<\vec{G},\vec{\mathbb I}>\rg)
\]
which implies the fact that the immersed surface is totally umbilic and it has then to be a translation of a sphere homothetic to $S^2$ the unit sphere of ${\R}^3\subset {\R}^m$.
This concludes the proof of theorem~\ref{th-VI.3}.    \hfill $\Box$ 

\medskip

The first part of theorem~\ref{th-VI.3} is in fact a special case of this more general result due to P. Li and S.T. Yau (see \cite{LiYa}).
\begin{Th}
\label{th-VI.4}
Let $\Sigma^2$ be a closed surface and let $\vec{\Phi}$ be a smooth immersion of $\Sigma^2$ into ${\R}^m$. Assume there exists a point $p\in {\R}^m$ with at least $k$ pre-images 
by $\vec{\Phi}$, then the following inequality holds
\be
\label{VI.39}
W(\vec{\Phi})=\int_{\Sigma^2}|\vec{H}_{\vec{\Phi}}|^2\ dvol_{\vec{\Phi}^\ast g_{{\R}^m}}\ge 4\pi\,k\quad.
\ee 
\hfill$\Box$
\end{Th}
An important corollary of the previous theorem is the following {\it Li-Yau $8\pi-$threshold} result\footnote{A weaker version of this result has been first established for spheres in codimension 2 (i.e. $m=4$)
by Peter Wintgen \cite{Win}.} .
\begin{Co}
\label{co-VI.5}
Let $\Sigma^2$ be a closed surface and let $\vec{\Phi}$ be a smooth immersion of $\Sigma^2$ into ${\R}^m$. If
\[
W(\vec{\Phi})=\int_{\Sigma^2}|\vec{H}_{\vec{\Phi}}|^2\ dvol_{\vec{\Phi}^\ast g_{{\R}^m}}<8\pi\quad,
\]
then $\vec{\Phi}$ is an embedding\footnote{We recall that embeddings are injective immersions. Images of manifolds
by embeddings are then submanifolds.}.\hfill $\Box$
\end{Co}

\noindent{\bf Proof of theorem~\ref{VI.4}.} Assume the origin $O$ of ${\R}^m$ admits $k-$pre-images by $\vec{\Phi}$. Let $\pi$ be inverse of the stereographic projection of ${\R}^m$ into $S^m$ which sends $O$ to the north pole of $S^m$ and $\infty$ to the south pole. Let $D_\la$ be the homothecy of center $O$ and factor $\la$. On every ball $B_R(O)$, the restriction of $D_\la\circ\vec{\Phi}$ converge in any $C^l$ norm to a union of $k$ planes restricted to $B_R(O)$ as $\la$ goes to $+\infty$. Hence $\pi\circ D_\la\circ\vec{\Phi}(\Sigma^2)$ restricted to any compact subset of $S^m\setminus\{south pole\}$ converges strongly in any $C^l$ norm to a union
of totally geodesic 2-spheres $S_1,\cdots, S_k$. It is well known that $\Psi_\la:=\pi\circ D_\la$ are conformal transformations. Applying then corollary~\ref{co-VI.2} we have
in one hand
\be
\label{VI.40}
W(\vec{\Phi})=\int_{\Sigma^2}|\vec{H}_{\Psi_\la\circ\vec{\Phi}}|^2\ dvol_{(\Psi_\la\circ\vec{\Phi})^\ast g_{S^m}}+\int_{\Sigma^2}\ dvol_{(\Psi_\la\circ\vec{\Phi})^\ast g_{S^m}}\quad.
\ee
where we used that the unit sphere is a constant sectional curvature space : for any two plane $\sigma$ in $TS^m$ $\ov{K}^{S^m}(\sigma)=1$. In the other hand the previous discussion leads to the following inequality
\be
\label{VI.41}
\begin{array}{l}
\ds\liminf_{\la\rightarrow +\infty}\int_{\Sigma^2}|\vec{H}_{\Psi_\la\circ\vec{\Phi}}|^2\ dvol_{(\Psi_\la\circ\vec{\Phi})^\ast g_{S^m}}+\int_{\Sigma^2}\ dvol_{(\Psi_\la\circ\vec{\Phi})^\ast g_{S^m}}\\[5mm]
\ds\quad\quad\ge\sum_{j=1}^k\int_{S_j}|\vec{H}|^2_{g_{S_j}}\ dvol_{g_{S_j}}+Area(S_j)\quad.
\end{array}
\ee
For each $S_j$ there is an isometry of $S^m$  sending $S_j$ to the canonical 2-sphere $S^2\subset{\R}^3\subset{\R}^m$. $S^2$ is minimal in $S^m$ ($\vec{H}\equiv 0$ ) and 
$Area(S^2)=4\pi$. Thus we deduce that
\be
\label{VI.42}
\forall j=1\cdots k\quad\quad\int_{S_j}|\vec{H}|^2_{g_{S_j}}\ dvol_{g_{S_j}}+Area(S_j)=4\pi\quad.
\ee
Combining (\ref{VI.40}), (\ref{VI.41}) and (\ref{VI.42}) gives Li-Yau inequality (\ref{VI.39}).\hfill $\Box$

Li and Yau  established moreover  a connection between  the Willmore energy of an immersed surface and it's the conformal class  that
provided new lower bounds. 

Let $\vec{\Phi}$ be a conformal parametrization of the immersion of a riemann surface \footnote{The pair $(\Sigma^2,c)$ denotes a closed 2-dimensional manifold $\Sigma^2$ together with a  
fixed conformal class on this surface.} $(\Sigma^2,c)$ into $S^m$. The $m-$conformal volume of
$\vec{\Phi}$ is the following quantity
\[
V_c(m,\vec{\Phi})=\sup_{\Psi\in Conf(S^m)}\int_{\Sigma^2}\ dvol_{\vec{\Phi}^\ast\Psi^\ast g_{S^m}}\quad.
\]
where $Conf(S^m)$  denotes the space of conformal diffeomorphism of $S^m$. Then Li and Yau define the {\it $m-$conformal volume} of a Rieman surface $(\Sigma^2,c)$ to be the following quantity
\[
V_c(m,(\Sigma^2,c))=\inf_{\vec{\Phi}} V_c(m,\vec{\Phi})
\]
where $\vec{\Phi}$ runs over all conformal immersions of $(\Sigma^2,c)$. 

\medskip

Let $\vec{\Phi}$ be a conformal immersion of a surface $\Sigma^2$ into ${\R}^m$ and $\pi$ be the inverse of the stereographic projection from ${\R}^m$ into $S^m$ (which is a conformal map).
Corollary~\ref{co-VI.2}  gives
\[
W(\vec{\Phi})=\int_{\Sigma^2}|\vec{H}_{\pi\circ\vec{\Phi}}|^2\ dvol_{(\pi\circ\vec{\Phi})^\ast g_{S^m}}+\int_{\Sigma^2}\ dvol_{(\pi\circ\vec{\Phi})^\ast g_{S^m}}\quad.
\]
from which one deduces the following lemma.
\begin{Lm}
\label{lm-VI.5}
Let $(\Sigma^2,c)$ be a closed Riemann surface and $\vec{\Phi}$ be a conformal immersion of this surface. Then the following
inequality holds
\[
\int_{\Sigma^2}|\vec{H}_{\vec{\Phi}}|^2\ dvol_{\vec{\Phi}^\ast g_{{\R}^m}}\ge V_c(m,\Sigma^2)\quad.
\]
Moreover equality holds if and only if $\vec{\Phi}(\Sigma^2)$ is the stereographic projection of a minimal surface of $S^m$.\hfill $\Box$ 
\end{Lm}
The main achievement of \cite{LiYa} is to provide lower bounds of $V_c(m,\Sigma^2)$ in terms of the conformal class of $\Sigma^2$. In particular they
establish the following result
\begin{Th}
\label{th-VI.5a}
Let $(T^2,c)$ be a torus equipped with the conformal class given by the flat torus $R^2/a{\Z}+b{\Z}$ where $a=(1,0)$ and $b=(x,y)$
where $0\le x\le 1/2$ and $\sqrt{1-x^2}\le y\le 1$, then for any $m\ge 3$
 \[
 2\pi^2\le V_c(m,(T^2,c))\quad.
 \]
\hfill $\Box$
\end{Th}
Combining lemma~\ref{lm-VI.5} and theorem~\ref{th-VI.5a} gives $2\pi^2$ as a lower bound to the Willmore energy of conformal
immersions of riemann surfaces in some sub-domain of Moduli space of the torus.  In fact the following statement
has been conjectured by T.Willmore in 1965
\begin{Con}
\label{con-VI.a}
Let $\vec{\Phi}$ be an immersion in ${\R}^m$ of the two dimensional torus $T^2$ then 
\[
\int_{T^2}|\vec{H}_{\vec{\Phi}}|^2\ dvol_{\vec{\Phi}^\ast g_{{\R}^m}}\ge 2\pi^2
\]
Equality should hold only for $\Phi(T^2)$ being equal to a Moebius transform of the stereographic projection into ${\R}^3$
of the Clifford torus $T^2_{cliff}:=\{1/{\sqrt{2}}\ (e^{i\theta},e^{i\phi})\in{\C}^2 ; (\theta,\phi)\in {\R}^2\}\subset {\R}^3$.\hfill $\Box$
\end{Con}
This conjecture has stimulated a lot of works. Recently F.C. Marques and A.Neves have submitted a proof of it in codimension 1 : $m=3$ (see \cite{MaNe}).

\subsection{The Willmore Surface Equations.}

In the previous subsection we have presented some lower bounds for the Willmore energy under various constraints. It is natural to 
look at the existence of the optimal surfaces for which the Willmore energy achieves it's lower bound. For instance minimal surfaces $S$
are absolute minima of the Willmore energy since they satisfy $\vec{H}=0$ and then $W(S)=0$. More generally  these optimal immersed surfaces 
under various constraints (fixed genus, boundary values...etc) will be  {\it critical points}  to the Willmore energy that are called {\it Willmore surfaces}. 
As we saw the Willmore energy in conformal coordinates identify with the $L^2$ norm of the immersion. It is therefore a 4-th order 
PDE generalizing the minimal surface equation $\vec{H}_{\vec{\Phi}}=0$ which is of second order\footnote{
Recall that in conformal coordinates $\vec{H}_{\vec{\Phi}}=0$ in ${\R}^m$ is equivalent to $\Delta\vec{\Phi}=0$.} (This is reminiscent for instance to the bi-harmonic map equation which is the 4th order 
generalization of the harmonic map equation which is of order 2 or, in the linear world, the Laplace equation being the 2nd order version of the Cauchy-Riemann
which is a 1st order PDE ...etc). This idea of having a 4-th order generalization of the minimal surface equation
was present in the original nomenclature given by Wilhelm Blaschke and his school
where {\it Willmore surfaces } were called {\it conformal minimal surfaces} see \cite{Bla3} and \cite{Tho}. Because of this strong link with minimal surface theory combined
with the fundamental conformal invariance property one can then naturally expect the family of {\it Willmore surfaces} to be of special interest in geometry .

\subsubsection{The Euler-Lagrange equation of Shadow, Thomsen and Weiner.} 

We first introduce the notion of Willmore surfaces.
\begin{Dfi}
\label{df-VI.6}
Let $\vec{\Phi}$ be a smooth immersion of a surface $\Sigma^2$ such that $W(\vec{\Phi})<+\infty$. $\vec{\Phi}$ is a critical point for $W$ if
\be
\label{VI.43}
\forall\vec{\xi}\in C^\infty_0({\Sigma^2},{\R}^m)\quad\quad\frac{d}{dt}W(\vec{\Phi}+t\vec{\xi})_{t=0}=0
\ee
Such an immersion is called {\it Willmore}.\hfill $\Box$
\end{Dfi}
Willmore immersions are characterized by an Euler Lagrange equation which has been discovered in dimension 3 by Shadow\footnote{See the comment by Whilhelm Blaschke in Ex. 7 ¤83 chapter 8 of \cite{Bla3}.} and also appear
in the PhD work of Gerhard Thomsen \cite{Tho}, student of Wilhelm Blaschke. The equation in general codimension, $m\ge 3$ arbitrary has been derived by Joel Weiner
in \cite{Wei}.
\begin{Th}
\label{th-VI.7}
A smooth immersion $\vec{\Phi}$ into ${\R}^m$ is Willmore (i.e. satisfies (\ref{VI.43})) if and only if it solves the following equation
\be
\label{VI.44}
\Delta_\perp\vec{H}_{\vec{\Phi}}-2|\vec{H}_{\vec{\Phi}}|^2\,\vec{H}_{\vec{\Phi}}+\ti{A}(\vec{H}_{\vec{\Phi}})=0
\ee
where $\Delta_{\perp}$ is the negative covariant laplacian operator for the connection\footnote{The assocated covariant derivative for any normal
vector-field $\vec{X}$ is given by
\[
D_Z\vec{X}:=\pi_{\vec{n}}(d\vec{X}\cdot Z)\quad.
\]} induced by the ambiant metric on the normal bundle
to $\vec{\Phi}(\Sigma^2)$ : for all $\vec{X}$ normal vector field to $\vec{\Phi}(\Sigma^2)$ one has
\[
\begin{array}{l}
\ds\Delta_{\perp}\vec{X}:=-\pi_{\vec{n}}\lf[d^{\ast_g}\lf(\pi_{\vec{n}}\lf[d\vec{X}\rg]\rg)\rg]\\[5mm]
\ds\quad\quad=\pi_{\vec{n}}\lf[(det\,g)^{-1/2}\p_{x_i}\lf(g^{ij}\, \sqrt{det\,g}\ \pi_{\vec{n}}\lf[\p_{x_j}\vec{X}\rg]\rg)\rg]
\end{array}
\]
where $d^{\ast_g}$ is the adjoint of $d$ for the induced scalar product $g:=\vec{\Phi}^\ast g_{{\R}^m}$ on $\Sigma^2$ and where we are using local coordinates on $\Sigma^2$ in the last line. Finally $\ti{A}$ is the following linear map
\be
\label{VI.44a}
\forall \vec{X}\in {\R}^m\quad\quad \ti{A}(\vec{X})=\sum_{i,j=1}^2\vec{\mathbb I}(e_i,e_j)\ <\vec{\mathbb I}(e_i,e_j),\vec{X}>\quad,
\ee
where $(e_1,e_2)$ is an orthonormal basis of $T\Sigma^2$ for the induced metric $g:=\vec{\Phi}^\ast g_{{\R}^m}$.\hfill $\Box$
\end{Th}
In the sequel we shall often omit the suscribt $\vec{\Phi}$ when there is no ambiguity and for instance $\vec{H}_{\vec{\Phi}}$
will sometimes be simply denoted $\vec{H}$.

\medskip

In dimension 3 the equation (\ref{VI.44}) takes a simpler form. $\vec{n}$ in this case is an $S^2$ valued vector and 
\[
\pi_{\vec{n}}(d\vec{H})=\lf<\vec{n},d\vec{H}\rg>\ \vec{n}=d H\ \vec{n}\quad.
\]
where $\ast_g$ is the Hodge operator for the  induced metric $g$ on $\Sigma^2$ pulled-back of the ambiant
metric in ${\R}^m$ by $\vec{\Phi}$. Hence
\[
d^{\ast_g}\lf(\pi_{\vec{n}}(d\vec{H})\rg)=d^{\ast_g}\lf(dH\ \vec{n}\rg)=d^{\ast_g}dH\ \vec{n}+\ast_{g}\lf(dH\wedge\ast_{g}d\vec{n}\rg)\quad.
\]
Thus
\be
\label{VI.45}
\Delta_\perp\vec{H}=-\pi_{\vec{n}}\lf[d^{\ast_g}\lf(\pi_{\vec{n}}(d\vec{H})\rg)\rg]=-d^{\ast_g}dH\ \vec{n}=\Delta_gH\ \vec{n}\quad,
\ee
where $\Delta_g$ is the negative Laplace Beltrami operator for the induced metric $g$ on $\Sigma^2$. Since $\vec{\mathbb I}$ takes value in the normal plane to $\vec{\Phi}_\ast T\Sigma^2$, we have
\[
\begin{array}{l}
\ds\ti{A}(\vec{H})=\sum_{i,j=1}^2\vec{\mathbb I}(e_i,e_j)\ \lf<\vec{\mathbb I}(e_i,e_j),\vec{H}\rg>\\[5mm]
\ds\quad=H\ \sum_{i,j=1}^2|\vec{\mathbb I}|^2(e_i,e_j)\ \vec{n}\quad.
\end{array}
\]
We take for $(e_1,e_2)$ an orthonormal basis of principal directions for $<\vec{n},{\mathbb I}>(\cdot,\cdot)$. Thus we have $|\vec{\mathbb I}|^2(e_i,e_j)=\kappa_i^e\ \delta_i^j$
and
\be
\label{VI.46}
\begin{array}{l}
\ti{A}(\vec{H})=H\ (\kappa_1^2+\kappa_2^2)\ \vec{n}\\[5mm]
\quad = H\ (4H^2-2K)\ \vec{n}
\end{array}
\ee
Combining (\ref{VI.44}) with (\ref{VI.45}) and (\ref{VI.46}), we obtain the following result which is a particular case 
of theorem~\ref{th-VI.7}.
\begin{Th}
\label{th-VI.8a}
A smooth immersion $\vec{\Phi}$ into ${\R}^3$ is Willmore (i.e. satisfies (\ref{VI.43})) if and only if it solves the following equation
\be
\label{VI.47}
\Delta_gH+2H\ (H^2-K)=0\quad,
\ee
where $H$ is the mean curvature of the immersed surface $\vec{\Phi}(\Sigma^2)$, $K$ the Gauss curvature and $\Delta_g$ the negative
Laplace-Beltrami operator for the induced metric $g$ obtained by pulling-back by $\vec{\Phi}$ the standard metric in ${\R}^3$
\hfill$\Box$
\end{Th}
{\bf Proof of theorem~\ref{th-VI.7}.}

Let
\[
\vec{\Phi} :\ [0,1]\times\Sigma^2\longrightarrow {\R}^m
\] 
be a smooth map such that $\vec{\Phi}(t,\cdot)$ (that we also denote $\vec{\Phi}_t$) is an immersion for every $t$ and $\vec{\Phi}_0$ is a Willmore surface.
Denote $\vec{H}(t,p)$ the mean-curvature vector of the surface $\vec{\Phi}_t(\Sigma^2)$ at the point $\vec{\Phi}_t(p)$. 

Let $T{\R}^m\res\vec{\Phi}([0,1]\times\Sigma^2)$ be the restriction of the tangent bundle to ${\R}^m$ over $\vec{\Phi}([0,1]\times\Sigma^2)$.

 The pull-back
bundle $\vec{\Phi}^{-1}(T{\R}^m\res\vec{\Phi}([0,1]\times\Sigma^2))$ can be decomposed into a direct bundle sum
\[
\begin{array}{rcl}
\vec{\Phi}^{-1}(T{\R}^m\res\vec{\Phi}([0,1]\times\Sigma^2))=T&\oplus& N\\[5mm]
 &\downarrow &\\[5mm]
 I&\times &\Sigma
 \end{array}
\]
where the fiber $T_{(t,p)}$ over the point $(t,p)\in [0,1]\times\Sigma^2$ is made of the vectors in ${\R}^m$ tangent to $\vec{\Phi}_t(\Sigma^2)$
and $N_{(t,p)}$ is made of the vectors in ${\R}^m$ normal to $\vec{\Phi}_\ast(T_{(t,p)}\Sigma^2)$ at $\vec{\Phi}(t,p)$. On $T$ we define the connection $\nabla$ as follows :
let $\sigma$ be a section of $T$ then we set
\[
\forall X\in T_{(t,p)}([0,1]\times\Sigma^2)\quad\quad\quad\nabla_X\sigma:=\pi_T(d\sigma\cdot X)\quad.
\]
where $\pi_T$ is the orthogonal projection onto $\vec{\Phi}_\ast(T_{(t,p)}\Sigma^2)$. On $N$ we define the connection $D$ as follows :
let $\tau$ be a section of $N$ then we set
\[
\forall X\in T_{(t,p)}([0,1]\times\Sigma^2)\quad\quad\quad D_X\sigma:=\pi_{\vec{n}}(d\sigma\cdot X)\quad.
\]
where $\pi_{\vec{n}}$ is the orthogonal projection onto the normal space to $\vec{\Phi}_\ast(T_{(t,p)}\Sigma^2)$.

\medskip

The first part of the proof consists in computing
\[
D_{\frac{\p}{\p t}}\vec{H}(0,p)\quad\quad\quad\forall p\in \Sigma^2\quad.
\]
To this purpose we introduce in a neighborhood of $(0,p)$ some special trivialization of $T$ and $N$.

Let $(\vec{e}_1,\vec{e}_2)$ be a positive orthonormal basis of $\vec{\Phi}_\ast(T_{(0,p)}\Sigma^2)$. Both $\vec{e}_1$ and $\vec{e}_2$ are points in the bundle $T$.
We first transport parallely $\vec{e}_1$ and $\vec{e}_2$, with respect to the connection $\nabla$, along the path $\{(t,p)\ ;\ \forall t\in [0,1]\}$. The resulting
parallely transported vectors are denoted $\vec{e}_i(t,p)$. The fact that this transport is parallel w.r.t. $\nabla$ implies that
 $(\vec{e}_1(t,p),\vec{e}_2(t,p))$ realizes a positive orthonormal basis\footnote{Indeed 
 \[
 \begin{array}{l}
\ds \nabla_{\frac{\p}{\p t}}\vec{e}_i=0\quad\Rightarrow\quad\pi_T\lf(d\vec{e}_i\cdot\frac{\p}{\p t}\rg)=0\\[5mm]
 \ds\quad\Rightarrow
 \quad\forall\,i,j\quad\lf<\frac{\p \vec{e}_i}{\p t},\vec{e}_j\rg>=0\quad\Rightarrow \quad\forall\,i,j\quad\frac{\p}{\p t}<\vec{e}_i,\vec{e}_j>=0\quad.
 \end{array}
 \]
 } of $\vec{\Phi}_\ast({t}\times T_p\Sigma^2)$.

 
Next we extend $\vec{e}_i(t,p)$ parallely and locally in $\{t\}\times \Sigma^2$ with respect again to $\nabla$ along geodesics in $\{t\}\times\Sigma^2$ starting
from $(t,p)$
for the induced metric $\vec{\Phi}_t^\ast g_{{\R}^m}$. 

We will denote also ${e}_i:=(\vec{\Phi}_t^{-1})_\ast\vec{e}_i$. By definition one has
\[
\ds\vec{H}(t,p):=\frac{1}{2}\sum_{s=1}^2\pi_{\vec{n}}(d\vec{e}_s\cdot{e}_s)=\frac{1}{2}\sum_{s=1}^2D_{{e}_s}\vec{e}_s\quad,
\]
where we have extended the use of the notation $D$ for any section of $T\oplus N$ in an obvious way\footnote{$D_X\sigma:=\pi_{\vec{n}}(d\sigma\cdot X)$. $D$ does not realizes a connection on the whole $T\oplus N$.}. Hence we have
\be
\label{VI.48}
\ds D_{\frac{\p}{\p t}}\vec{H}= \frac{1}{2}\sum_{s=1}^2D_{\frac{\p}{\p t}}(D_{{e}_s}\vec{e}_s)
\ee
 Since on $\vec{\Phi}([0,1]\times\{p\})$ one has $\nabla\vec{e}_s\equiv 0$ and since moreover $\nabla +D$ coincides with the
differentiation\footnote{ $d$ is the flat standard connection on $T{\R}^m$} $d$ in ${\R}^m$ one has
\be
\label{VI.49}
\begin{array}{l}
\ds D_{\frac{\p}{\p t}}(D_{{e}_s}\vec{e}_s)(0,p)=D_{\frac{\p}{\p t}}(d\vec{e}_s\cdot{e}_s)(0,p)\\[5mm]
\ds\quad=\pi_{\vec{n}}\lf(d\lf(\p_t \vec{e}_s\rg)\cdot e_s\rg)+\pi_{\vec{n}}\lf(d\vec{e}_s\cdot\p_t e_s\rg)
\end{array}
\ee 
Observe that for a fixed $q\in \Sigma^2$ $e_s(t,q)$ stays in $T_q\Sigma^2$ as $t$ varies and then there is no need
of a connection to define $\p_te_s$. Observe moreover that $\p_t e_s=[\p_t,e_s]$. Thus (\ref{VI.48}) together with (\ref{VI.49}) imply
\be
\label{VI.50}
 D_{\frac{\p}{\p t}}\vec{H}= \frac{1}{2}\sum_{s=1}^2D_{e_s}\lf(\p_t \vec{e}_s\rg)+ \frac{1}{2}\sum_{s=1}^2D_{[\p_t,e_s]}\vec{e}_s\quad.
\ee 
We denote $V:=\p_t\vec{\Phi}$. Observe that
\be
\label{VI.51}
\begin{array}{l}
\p_t\vec{e}_s=\p_t(d\vec{\Phi}\cdot e_s)=d\vec{V}\cdot e_s+d\vec{\Phi}_t\cdot\p_te_s\\[5mm]
\ds\quad\quad=d\vec{V}\cdot e_s+d\vec{\Phi}\cdot[\p_t,e_s]=d\vec{V}\cdot e_s+[\vec{V},\vec{e}_s]\quad.
\end{array}
\ee
Observe moreover that for two tangent vector-fields $X$ and $Y$ on $\Sigma^2$ one has\footnote{Indeed
in local coordinates $(x_1,x_2)$ in $\Sigma_2$ one has
\[
\begin{array}{l}
\ds D_X(d\vec{\Phi}\cdot Y)=\sum_{k,j=1}^2\pi_{\vec{n}}(\p_{x_k}(\p_{x_j}\vec{\Phi}\,Y_j)X_k)\\[5mm]
= \sum_{j,k=1}^2\pi_{\vec{n}}(\p_{x_j}(\p_{x_k}\vec{\Phi}\,X_k)Y_j)-\sum_{j,k=1}^2\pi_{\vec{n}}(\p_{x_k}\vec{\Phi} Y_j\, \p_{x_j}X_k)=D_Y(d\vec{\Phi}\cdot X)\quad,
\end{array}
\]
where we have used the fact for all $j,k=1,2$
\[
\vec{n}(\p_{x_k}\vec{\Phi}\, Y_j\, \p_{x_j}X_k)=0\quad\mbox{ and }\quad \vec{n}(\p_{x_j}\vec{\Phi}\, X_k\,\p_{x_k}Y_j)\quad.
\]}
\be
\label{VI.52}
D_X(d\vec{\Phi}\cdot Y)=D_Y(d\vec{\Phi}\cdot X)\quad.
\ee
This last identity implies in particular that 
\be
\label{VI.53}
D_{e_s}([\vec{V},\vec{e}_s])=D_{[\p_t,e_s]}\vec{e}_s\quad.
\ee
Combining (\ref{VI.50}), (\ref{VI.51}) and (\ref{VI.53}) gives
\be
\label{VI.54}
D_{\frac{\p}{\p t}}\vec{H}=\frac{1}{2}\sum_{s=1}^2D_{e_s}(d\vec{V}\cdot e_s)+\sum_{s=1}^2D_{e_s}[\vec{V},\vec{e}_s]\quad.
\ee
We have
\[
[\vec{V},\vec{e}_s]=d\vec{\Phi}\cdot [\p_t,e_s]=d\vec{e}_s\cdot\p_t-d\vec{V}\cdot e_s
\]
Since $[\vec{V},\vec{e}_s]$ is tangent to $\vec{\Phi}_t(\Sigma^2)$, by composing the previous identity
with the orthogonal projection $\pi_T$ on the tangent space to $\vec{\Phi}_t(\Sigma^2)$ we obtain
\be
\label{VI.55}
[\vec{V},\vec{e}_s]=\pi_T(d\vec{e}_s\cdot\p_t)-\pi_T(d\vec{V}\cdot e_s)=\nabla_{{\p_t}}\,\vec{e}_s-\nabla_{e_s}\,\vec{V}\quad.
\ee
Using (\ref{VI.52}) we have
\[
D_{e_s}(\nabla_{\p_t}\, \vec{e}_s)=D_{\nabla_{\p_t}\,\vec{e}_s}\vec{e}_s
\]
Since $\nabla_{\p_t}\,\vec{e}_s=0$ along the curve $[0,1]\times\{p\}$ one deduces from the previous identity
that
\be
\label{VI.56}
D_{e_s}(\nabla_{\p_t}\,\vec{e}_s)(t,p)\equiv 0\quad.
\ee
Combining (\ref{VI.54}) with (\ref{VI.55}) and (\ref{VI.56}) we obtain
\be
\label{VI.57}
D_{\frac{\p}{\p t}}\vec{H}=\frac{1}{2}\sum_{s=1}^2\lf[D_{e_s}(d\vec{V}\cdot e_s)-2D_{e_s}(\nabla_{e_s}\vec{V})\rg]\quad.
\ee
We decompose $\vec{V}=\vec{V}^N+\vec{V}^T$ where $\vec{V}^N=\pi_{\vec{n}}(\vec{V})$ and $\vec{V}^T=\pi_T(\vec{V})$ and  after writting
$d\vec{V}\cdot e_s=\nabla_{e_s}\vec{V}+D_{e_s}\vec{V},$ (\ref{VI.57}) becomes
\be
\label{VI.58}
\begin{array}{rl}
\ds D_{\frac{\p}{\p t}}\vec{H}&\ds=\frac{1}{2}\sum_{s=1}^2\lf[D_{e_s}(D_{e_s}\vec{V}^N)-D_{e_s}(\nabla_{e_s}\vec{V}^N)\rg]\\[5mm]
 &\ds\ + \frac{1}{2}\sum_{s=1}^2\lf[D_{e_s}(D_{e_s}\vec{V}^T)-D_{e_s}(\nabla_{e_s}\vec{V}^T)\rg]
\end{array}
\ee
Observe that
\be
\label{VI.59}
\ds\sum_{s=1}^2D_{e_s}(D_{e_s}\vec{V}^N)=\Delta_{\perp}\vec{V}^N\quad,
\ee
and 
\be
\label{VI.60}
\ds D_{e_s}(\nabla_{e_s}\vec{V}^N)=\vec{\mathbb I}(e_s,\nabla_{e_s}\vec{V}^N)\quad.
\ee
We have moreover
\[
\begin{array}{l}
\ds\nabla_{e_s}\vec{V}^N=\sum_{k=1}^2<\nabla_{e_s}\vec{V}^N,\vec{e}_k>\ \vec{e}_k=\sum_{k=1}^2<d\vec{V}^N\cdot e_s,\vec{e}_k>\ \vec{e}_k\\[5mm]
\ds\quad\quad=-\sum_{k=1}^2<\vec{V}^N,d\vec{e}_k\cdot e_s>\ \vec{e}_k=-\sum_{k=1}^2<\vec{V}^N,D_{e_s}\vec{e}_k>\ \vec{e}_k
\end{array}
\]
Combining this last identity with (\ref{VI.60}) gives
\be
\label{VI.61}
\ds\sum_{s=1}^2 D_{e_s}(\nabla_{e_s}\vec{V}^N)=-\sum_{s,k=1}^2 D_{e_s}\vec{e}_k\ \vec{\mathbb I}(\vec{e}_s,\vec{e}_k)=-\ti{A}(\vec{V}^N)\ .
\ee
Now, since $\vec{\Phi}_0+t\vec{V}^T$ preserves infinitesimally the surface $\vec{\Phi}_0(\Sigma^2)$ an since the Willmore energy
is independent of the parametrization,  one has
\[
\frac{d}{dt}\int_{\Sigma^2}|\vec{H}_{\vec{\Phi}_0+t\vec{V}^T}|^2\ dvol_{(\vec{\Phi}_0+t\vec{V}^T)^\ast g_{{\R}^m}}=0\quad.
\]
In other words, tangent variations do not change the Willmore energy at the first order. Therefore we can assume $\vec{V}=\vec{V}^N$
and combining (\ref{VI.58}) with (\ref{VI.59}) and (\ref{VI.61}) gives finally if $\vec{V}^T=0$
\be
\label{VI.62}
D_{\frac{\p}{\p t}}\vec{H}=\frac{1}{2}\lf[\Delta_\perp\vec{V}^N+\ti{A}(\vec{V}^N)\rg]\quad.
\ee
Let $(g_{ij}^t)_{ij}:=(\vec{\Phi}+t\vec{V})^\ast g_{{\R}^m}$, a straightforward computation gives in local coordinates
\[
g_{ij}^t=g_{ij}+t\ \lf[<\p_{x_i}\vec{V},\p_{x_j}\vec{\Phi}>+<\p_{x_j}\vec{V},\p_{x_i}\vec{\Phi}>\rg]+o(t)\quad,
\]
from which we deduce
\[
\ds \frac{det(g_{ij}^t)}{det(g_{ij})}-1=2\ \sum_{k,j=1}^2g^{jk}\ <\p_{x_k}\vec{V},\p_{x_j}\vec{\Phi}>+o(t)\ .
\]
where $(g^{ik})$ is the inverse matrix to $(g_{ij})$. Thus
\[
\frac{\frac{d}{dt}\sqrt{det(g_{ij}^t)}(0)}{\sqrt{det(g_{ij})}}=\sum_{k=1}^2g^{jk}\ <\p_{x_k}\vec{V},\p_{x_j}\vec{\Phi}>\ .
\]
Since we are led to consider only $\vec{V}$ orthogonal to $(\vec{\Phi}_0)_\ast(T\Sigma^2)$ we have
\[
\begin{array}{l}
\ds\sum_{k=1}^2g^{jk}\ <\p_{x_k}\vec{V},\p_{x_j}\vec{\Phi}>=-\sum_{k=1}^2g^{jk}\ <\vec{V},\p^2_{x_j\,x_k}\vec{\Phi}>\\[5mm]
\ds\quad\quad=-\sum_{k=1}^2g^{jk}\ <\vec{V},\pi_{\vec{n}}(\p^2_{x_j\,x_k}\vec{\Phi})>=-2<\vec{H},\vec{V}>\quad.
\end{array}
\]
Hence we have proved that
\be
\label{VI.63}
\frac{d}{dt}(dvol_{\vec{\Phi}_t^\ast g_{{\R}^m}})(0)=-2<\vec{H},\vec{V}>\ dvol_{\vec{\Phi}^\ast g_{{\R}^m}}\quad.
\ee
Combining (\ref{VI.62}) and (\ref{VI.63}) we obtain\footnote{These first variations of $\vec{H}$ and $dvol_g$ were known in codimension 1 probably even before W.Blaschke see ¤ 117 in \cite{Bla1}.}
\[
\begin{array}{l}
\ds\frac{d}{dt}\int_{\Sigma^2}|\vec{H}_t|^2\ dvol_{\vec{\Phi}_t^\ast g_{{\R}^m}}=\ds 2\int_{\Sigma^2} <D_{\p_t}\vec{H},\vec{H}>\ dvol_{\vec{\Phi}^\ast g_{{\R}^m}}\\[5mm]
\ds\quad\quad-2\int_{\Sigma^2}|\vec{H}|^2\ <\vec{H},\vec{V}>\ dvol_{\vec{\Phi}^\ast g_{{\R}^m}}\\[5mm]
 \ds\quad=\int_{\Sigma^2}\lf<\Delta_\perp\vec{V}+\ti{A}(\vec{V})-2|\vec{H}|^2\,\vec{V},\vec{H}\rg>\ dvol_{\vec{\Phi}^\ast g_{{\R}^m}}\\[5mm]
 \ds\quad=\int_{\Sigma^2}\lf<\Delta_\perp\vec{H}+\ti{A}(\vec{H})-2|\vec{H}|^2\,\vec{H},\vec{V}\rg>\ dvol_{\vec{\Phi}^\ast g_{{\R}^m}}
 \end{array}
 \]
The immersion $\vec{\Phi}$ is Willmore if and only if $$ \frac{d}{dt}\int_{\Sigma^2}|\vec{H}_t|^2\ dvol_{\vec{\Phi}_t^\ast g_{{\R}^m}}=0$$ for any perturbation $\vec{V}$
which is equivalent to (\ref{VI.44}) and theorem~\ref{th-VI.7} is proved.\hfill $\Box$

\medskip

As mentioned in the introduction of this course questions we are interested in are analysis questions for conformally invariant lagrangians.

 In this part of the course devoted to Willmore
Lagrangian we shall look at the following problems  :

\begin{itemize}
\item[i)] Does there exists a minimizer of Willmore functional among all smooth immersions for a fixed
2-dimensional surface $\Sigma^2$ ? and, if yes, can one estimate the energy and special
properties of such a minimizer ?
\item[ii)] Does there exists a minimizers of Willmore functional among a more restricted class of immersions such
as conformal immersions for a fixed chosen conformal class $c$ on $\Sigma$ ? or does there exist
a minimizing immersion of Willmore functional among all immersions into ${\R}^3$ enclosing a domain of given volume and
realizing  a fixed area... 
\item[iii)] What happens to a sequence of weak Willmore immersions of a surface $\Sigma^2$ having a uniformly bounded energy
at the limit ? does it convergence in some sense to a surface which is still Willmore and if not what are the possible
''weak limits'' of Willmore surfaces ?
\item[iv)] How stable is the Willmore equation ? that means : following a sequence of ''almost Willmore'' surfaces ''solving more and more''
the Willmore equation - Willmore Palais Smale sequences for instance - does such a sequence converges to a Willmore surface ?
\item[v)] Can one apply fundamental variational principles such as {\it Ekeland's variational Principles} or {\it Mountain pass lemma} to the Willmore
functional ?
\item[vi)] Is there a weak notion of Willmore immersions and, if yes, is any such a weak solution smooth ?
\end{itemize}

We will devote the rest of the course to these questions which are very much related to another - as an experienced non-linear analysts
could anticipate ! -. To this aim we have to find a suitable framework for developing calculus of variation questions for Willmore functional.

The first step consists naturally in trying to confront  the  Euler Lagrange Willmore equation  (\ref{VI.44}) we obtained to the questions i)$\cdots$vi).

\medskip

In codimension 3 the Schadow's-Thomsen equation of Willmore surfaces is particularly attractive because of it's apparent simplicity :
\begin{itemize}
\item[i)] The term $\Delta_gH$ is the application of a somehow {\bf classical linear elliptic operator} - the Laplace Beltrami Operator - on the mean-curvature $H$.
\item[ii)] The nonlinear terms $2H\ (H^2-K)$ is an {\bf algebraic function of the principal curvatures}.
\end{itemize}
Despite it's elegance, Schadow-Thomsen's equation is however very complex for an analysis approach and for the previous questions i)$\cdots$vi) we posed. Indeed 
\begin{itemize}
\item[i)] The term $\Delta_gH$ is in fact the application of the Laplace Beltrami Operator to the mean curvature but this operator {\bf it depends on the metric $g$}
which itself is varying depending of the immersion $\vec{\Phi}$. Hence $\Delta_g$ could, in the minimization procedures we aim to follow, strongly degenerate
as the immersion $\vec{\Phi}$ degenerates.
\item[ii)] Already in codimension 1 where the problem admits a simpler formulation, the non-linear term $2H\ (H^2-K)$ is somehow {\bf supercritical} with respect to the
Lagrangian :  indeed the Willmore Lagrangian ''controls'' the $L^2$ norm of the mean curvature and the $L^2$ norm of the second fundamental form for a given closed surface for instance.
However the non-linearity in the Eular Lagrange equation in the form (\ref{VI.44}) is {\bf cubic} in the second fundamental form.  Having some weak notion of immersions with only $L^2$-bounded second fundamental form would then be insufficient to write the equation though the Lagrangian from which this equation is deduced would make sense for such a weak immersion  !   This provides some apparent {\bf functional analysis paradox}.
\end{itemize}

\subsubsection{The conservative form of Willmore surfaces equation.}

In \cite{Riv2} an alternative form to the Euler Lagrange equation of Willmore functional was proposed. We the following result plays a central
role in the rest of the course.
\begin{Th}
\label{th-VI.8}
Let $\vec{\Phi}$ be a smooth immersion of a two dimensional manifold $\Sigma^2$ into ${\R}^m$ then the following identity holds
\be
\label{VI.64}
\begin{array}{c}
\ds \Delta_\perp\vec{H}+\ti{A}(\vec{H})-2|\vec{H}|^2\,\vec{H}\\[5mm]
\ds=\frac{1}{2}d^{\ast_g}\lf[d\vec{H}-3\pi_{\vec{n}}(d\vec{H})+\star(\ast_g d\vec{n}\wedge\vec{H})\rg]
\end{array}
\ee
where $\vec{H}$ is the mean curvature vector of the immersion $\vec{\Phi}$, $\Delta_\perp$ is the negative covariant laplacian on the normal bundle to the immersion, $\ti{A}$ is the linear
map given by (\ref{VI.44a}), $\ast_g$ is the Hodge operator associated to the pull-back metric $\vec{\Phi}^\ast g_{{\R}^m}$ on $\Sigma^2$, $d^{\ast_g}=-\ast_g\,d\,\ast_g$ is the adjoint operator 
with respect to the metric $g$ to the exterior differential $d$, $\vec{n}$ is the Gauss map to the immersion, $\pi_{\vec{n}}$ is the orthogonal projection onto the normal space to the tangent
space $\vec{\Phi}_\ast T\Sigma^2$ and $\star$ is the Hodge operator from $\wedge^p{\R}^m$ into $\wedge^{m-p}{\R}^m$ for the canonical metric in ${\R}^m$.
\hfill$\Box$
\end{Th}
A straightforward but important consequence of theorem~\ref{th-VI.8} is the following conservative form of Willmore surfaces equations.
\begin{Co}
\label{co-VI.8}  
An immersion $\vec{\Phi}$ of a 2-dimensional manifold $\Sigma^2$ is Willmore if and only if the 1-form given by
\[
\ast_g\lf[d\vec{H}-3\pi_{\vec{n}}(d\vec{H})+\star(\ast_g d\vec{n}\wedge\vec{H})\rg]
\]
is closed.\hfill $\Box$
\end{Co}
The analysis questions for Willmore immersions we raised in the previous subsection  can be studied basically from two point of views
\begin{itemize}
\item[i)] By working with the maps $\vec{\Phi}$ themseves.
\item[ii)] By working with the immersed surface : the image $\vec{\Phi}(\Sigma^2)\subset{\R}^m$.
\end{itemize}
The drawback of the first approach is the huge invariance group of the problem containing the positive diffeomorphism group of $\Sigma^2$.
The drawback of the second approach comes from the fact that the Euler Lagrange equation is defined on an unknown object $\vec{\Phi}(\Sigma^2)$.
In this course we shall take the first approach but by trying to ''break'' as much as we can the symmetry invariance given by the action of positive diffeomorphisms of $\Sigma^2$.
From Gauge theory and in particular from Yang-Mills theory a natural way to ''break'' a symmetry group is to look for {\it Coulomb gauges}.
As we will explain the concept of  {\it Coulomb gauge} transposed to the present setting of immersions of 2-dimensional manifolds is given by
the Isothermal coordinates (or conformal parametrization). We shall make an intensive use of these conformal parametrization and therefore
we shall make an intensive use of the conservative form of Willmore surfaces equation written in isothermal coordinates.
A further corollary of theorem~\ref{th-VI.8} giving the Willmore surfaces equations in isothermal coordinates is the following.
\begin{Co}
\label{co-VI.9}
A conformal immersion $\vec{\Phi}$ of the flat disc $D^2$ is Willmore if and only if
\be
\label{VI.65aa}
div\lf[\nabla\vec{H}-3\pi_{\vec{n}}(\nabla\vec{H})+\star(\nabla^\perp\vec{n}\wedge\vec{H})\rg]=0
\ee
where the operators $div$, $\nabla$ and $\nabla^\perp$ are taken with respect to the flat metric\footnote{$div\,X=\p_{x_1}X_1+\p_{x_2}X_2$,
$\nabla f=(\p_{x_1}f,\p_{x_2}f)$ and $\nabla^\perp f=(-\p_{x_2}f,\p_{x_1}f)$.} in $D^2$.
\hfill$\Box$
\end{Co}
In order to exploit analytically equation (\ref{VI.65aa}) we will need a more explicit expression of $\pi_{\vec{n}}(\nabla\vec{H})$.
Let  $\res$ be the interior multiplication between $p-$ and $q-$vectors  $p\ge q$ producing $p-q-$vectors in ${\R}^m$ 
such that (see
\cite{Fe} 1.5.1 combined with 1.7.5) : for every choice of $p-$, $q-$ and $p-q-$vectors,
respectively $\al$, $\beta$ and $\gamma$ the following holds
\[
\lf<\al\res\beta,\gamma\rg>=\lf<\al,\beta\wedge\gamma\rg>\quad.
\]
Let $(\vec{e}_1,\vec{e}_2)$ be an orthonormal basis of the orthogonal 2-plane to the $m-2$ plane given by $\vec{n}$
and positively oriented in such a way that
\[
\star(\vec{e}_1\wedge\vec{e}_2)=\vec{n}
\]
and let $(\vec{n}_1\cdots\vec{n}_{m-2})$ be a positively oriented orthonormal basis of the $m-2$-plane given by $\vec{n}$
satisfying $\vec{n}=\wedge_{\al}\vec{n}_\al$. One verifies easily that
\[
\lf\{
\begin{array}{l}
\vec{n}\res\vec{e}_i=0\\[5mm]
\vec{n}\res\vec{n}_\al= (-1)^{\al-1}\wedge_{\beta\ne\al}\vec{n}_\beta\\[5mm]
\vec{n}\res(\wedge_{\beta\ne\al}\vec{n}_\beta)=(-1)^{m+\al-2}\,\vec{n}_\al
\end{array}
\rg.
\]
We then deduce the following identity :
\be
\label{VI.65ab}
\forall \vec{w}\in {\R}^m\quad\quad \pi_{\vec{n}}(\vec{w})=(-1)^{m-1}\ \vec{n}\res(\vec{n}\res\vec{w})
\ee
From (\ref{VI.65ab}) we deduce in particular
\be
\label{VI.65ac}
\begin{array}{l}
\pi_{\vec{n}}(\nabla\vec{H})=\nabla\vec{H}-(-1)^{m-1}\ \nabla(\vec{n})\res(\vec{n}\res\vec{w})\\[5mm]
\quad\quad\quad\quad\quad\quad\quad-(-1)^{m-1}\ \vec{n}\res(\nabla(\vec{n})\res\vec{H})
\end{array}
\ee

\medskip

It is interesting to look at the equation (\ref{VI.65aa}) in the codimension 1 case ($m=3$). In this particular case we have proved the following equivalence :
Let $\vec{\Phi}$ be a conformal immersion from the 2 disc $D^2$ into ${\R}^3$ then
\be
\label{VI.65}
\begin{array}{c}
\ds\Delta_gH+2H\ (H^2-K)=0\\[5mm]
\ds\Updownarrow\\[5mm]
\ds div\lf[2\nabla\vec{H}-3H\,\nabla\vec{n}-\nabla^\perp\vec{n}\times\vec{H}\rg]=0
\end{array}
\ee
It is now clear, at least in codimension 1, that the conservative form of Willmore surface equation is now compatible with the Willmore lagrangian
in the sense that this 2nd equation in (\ref{VI.65}) has a distributional sense assuming only that 
\[
\begin{array}{c}
\ds\int_{D^2}|d\vec{n}|_g^2\ dvol_g=\int_{D^2}|\nabla\vec{n}|^2\ dx_1\, dx_2<+\infty\\[5mm]
\ds\Downarrow\\[5mm]
\ds\nabla\vec{H}\in H^{-1}(D^2)\quad,\\[5mm]
\ds 3H\,\nabla\vec{n}\in L^1(D^2)\quad\mbox{and}\quad\nabla^\perp\vec{n}\times\vec{H}\in L^1(D^2)
\end{array}
\]
Thus under the minimal assumption saying that the second fundamental form is in $L^2$ (that comes naturally from our variational problem $W(\vec{\Phi})<+\infty$), in
conformal coordinates, the quantity
\[
 div\lf[2\nabla\vec{H}-3H\,\nabla\vec{n}-\nabla^\perp\vec{n}\times\vec{H}\rg]
 \]
 is an ''honest'' distribution in ${\mathcal D}'(D^2)$ whereas, under such minimal assumption
 \[
 \Delta_gH+2H\ (H^2-K)
 \]
 has {\bf no distributional meaning} at all. This is why the conservative form of the Willmore surface equation is more suitable to solve the analysis questions
 $i)\cdots vi)$ we are asking. The same happens in higher codimension as well. Using (\ref{VI.65ac}),
 one sees that, under the assumption that $\nabla\vec{n}\in L^2$ one has
 \[
 \nabla\vec{H}-3\pi_{\vec{n}}(\nabla\vec{H})+\star(\nabla^\perp\vec{n}\wedge\vec{H})\in H^{-1}+L^1
 \]
 which is again an honest distribution.
 
 \medskip
 
 The second equation in (\ref{VI.65}) is in {\bf conservative-elliptic form} which is {\bf critical in 2 dimension }
  under the assumption that $\vec{n}\in W^{1,2}$. For a sake of clarity we present it in codimension 1 though this holds identically in arbitrary codimension.
  
   In codimension 1 we write the Willmore surface equation as follows
 \be
 \label{VI.65a}
 \Delta\vec{H}= div\lf[\frac{3}{2}H\,\nabla\vec{n}+\frac{1}{2}\nabla^\perp\vec{n}\times\vec{H}\rg]
 \ee
the right-hand-side  of (\ref{VI.65a}) is the flat Laplacian of $\vec{H}$ and the left-hand-side is the divergence
of a ${\R}^m$ vector-field wich is a bilinear map of the second fundamental form. Assuming hence $\vec{n}\in W^{1,2}$
we deduce 
$3/2\ H\,\nabla\vec{n}+1/2\ \nabla^\perp\vec{n}\times\vec{H}\in L^1(D^2)$. Adams result on Riesz potentials
\cite{Ad} implies that
\[
\Delta^{-1}_0div\lf[\frac{3}{2}H\,\nabla\vec{n}+\frac{1}{2}\nabla^\perp\vec{n}\times\vec{H}\rg]\in L^{2,\infty}
\]
where $\Delta^{-1}_0$ is the {\it Poisson Kernel} on the disc $D^2$. Inserting this information back in (\ref{VI.65a}) we obtain $\vec{H}\in L^{2,\infty}_{loc}(D^2)$ which is almost the information we started from\footnote{We will see in the next subsection that $\vec{n}\in W^{1,2}$ implies that the conformal factor is bounded
in $L^\infty$ and then, since again the parametrization is conformal, we have $$H\in L^{2,\infty}_{loc}(D^2)\quad\quad\Longrightarrow \quad\quad \nabla\vec{n}\in L^{2,\infty}_{loc}(D^2)$$.}. This phenomenon characterizes {\bf critical
elliptic systems} as we saw it already in the first section of this course while presenting the elliptic systems of quadratic growth in two dimension for $W^{1,2}$ norm.

 \medskip
 
 It remains now in this subsection to prove theorem~\ref{th-VI.8}.
 
\medskip

In order to make the proof of theorem~\ref{th-VI.8} more accessible we first give a proof of it in the codimension 1 setting
which is simpler  and then we will explain how to generalize it to higher codimension.

The result is a local one on $\Sigma^2$ therefore we can work locally in a disc-neighborhood of a point and use isothermal coordinates
on this disc. This means that we can assume $\vec{\Phi}$ to be a conformal immersion from the unit disc $D^2\subset {\R}^2$
into ${\R}^3$.

we will need the following general lemma for conformal immersions of the 2-disc in ${\R}^3$
\begin{Lm}
\label{lm-VI.10}
Let $\vec{\Phi}$ be a conformal immersion from $D^2$ into ${\R}^3$. Denote by $\vec{n}$ the Gauss map of the conformal immersion $\vec{\Phi}$ and denote by $H$ the mean curvature. Then the following identity holds
\be
\label{VI.66}
-2H\ \nabla\vec{\Phi}=\nabla\vec{n}+\vec{n}\times\nabla^\perp\vec{n}\quad
\ee
where $\nabla\cdot:=(\p_{x_1}\cdot,\p_{x_2}\cdot)$ and $\nabla^\perp\cdot:=(-\p_{x_2}\cdot,\p_{x_1}\cdot)$.\hfill$\Box$
\end{Lm}
{\bf Proof of lemma~\ref{lm-VI.10}.}
Denote $(\vec{e}_1,\vec{e}_2)$ the orthonormal basis of $\vec{\Phi}_\ast(T\Sigma^2)$ given by
$$
\vec{e}_i:=e^{-\la}\ \frac{\p \vec{\Phi}}{\p x_i}\quad,
$$
where $e^\la=|\p_{x_1}\vec{\Phi}|=|\p_{x_2}\vec{\Phi}|$. The oriented Gauss map $\vec{n}$ is then given by
\[
\vec{n}=e^{-2\la}\ \frac{\p\vec{\Phi}}{\p x_1}\times \frac{\p\vec{\Phi}}{\p x_2}\quad.
\]
We have
\[
\lf\{
\begin{array}{l}
<\vec{e}_1,\vec{n}\times\nabla^{\perp}\vec{n}>=-<\nabla^\perp\vec{n},\vec{e}_2>\\[5mm]
<\vec{e}_2,\vec{n}\times\nabla^{\perp}\vec{n}>=<\nabla^\perp\vec{n},\vec{e}_1>\quad.
\end{array}
\rg.
\]
From which we deduce
\[
\lf\{
\begin{array}{l}
\ds-\vec{n}\times\p_{x_2}\vec{n}=<\p_{x_2}\vec{n},\vec{e}_2>\ \vec{e}_1-<\p_{x_2}\vec{n},\vec{e}_1>\ \vec{e}_2\\[5mm]
\ds\vec{n}\times\p_{x_1}\vec{n}=-<\p_{x_1}\vec{n},\vec{e}_2>\ \vec{e}_1+<\p_{x_1}\vec{n},\vec{e}_1>\ \vec{e}_2
\end{array}
\rg.
\]
Thus
\[
\lf\{
\begin{array}{l}
\ds\p_{x_1}\vec{n}-\vec{n}\times\p_{x_2}\vec{n}=[<\p_{x_2}\vec{n},\vec{e}_2>+<\p_{x_1}\vec{n},\vec{e}_1>]\ 
\vec{e}_1\\[5mm]
\ds\p_{x_2}\vec{n}+\vec{n}\times\p_{x_1}\vec{n}=[<\p_{x_2}\vec{n},\vec{e}_2>+<\p_{x_1}\vec{n},\vec{e}_1>]\ \vec{e}_2
\end{array}
\rg.
\]
Since $H= -e^{-\la}\ 2^{-1}[<\p_{x_2}\vec{n},\vec{e}_2>+<\p_{x_1}\vec{n},\vec{e}_1>]$ we deduce (\ref{VI.66}) and 
Lemma~\ref{lm-VI.10} is proved.\hfill$\Box$

\medskip

\noindent {\bf Proof of theorem~\ref{th-VI.8} in codimension 1.}

We can again assume that $\vec{\Phi}$ is conformal. First take the divergence of (\ref{VI.66}) and multiply by $H$.
This gives
\be
\label{VI.67}
-2H^2\ \Delta\vec{\Phi}-2H\nabla H\cdot\nabla\vec{\Phi}=H\,div\lf[\nabla\vec{n}+\vec{n}\times\nabla^\perp\vec{n}\rg]\quad.
\ee
We replace $-2H\,\nabla\vec{\Phi}$ in (\ref{VI.67}) by the expression given by (\ref{VI.66}), moreover we also use
the expression of the mean curvature vector in terms of $\vec{\Phi}$ :
\be
\label{VI.67a}
\Delta\vec{\Phi}=2e^{2\la}\ \vec{H}\quad.
\ee
So (\ref{VI.67}) becomes
\be
\label{VI.68}
\begin{array}{l}
\ds-4H^2\ \vec{H}\ e^{2\la}+\nabla H\cdot\lf[\nabla\vec{n}+\vec{n}\times\nabla^\perp\vec{n}\rg]\\[5mm]
\ds\quad\quad= H\,div\lf[\nabla\vec{n}+\vec{n}\times\nabla^\perp\vec{n}\rg]\quad.
\end{array}
\ee
The definition of the Gauss curvature gives
\be
\label{VI.69}
K\ \vec{n}=-\frac{e^{-2\la}}{2}\ \nabla\vec{n}\times\nabla^\perp\vec{n}=-\frac{e^{-2\la}}{2}\ div\lf[\vec{n}\times\nabla^{\perp}\vec{n}\rg]
\ee
Inserting (\ref{VI.69}) in (\ref{VI.68}) gives
\be
\label{VI.70}
\begin{array}{l}
-4 e^{2\la}\ \vec{H}\ (H^2-K)+\nabla H\cdot\lf[\nabla\vec{n}+\vec{n}\times\nabla^\perp\vec{n}\rg]\\[5mm]
\ds\quad\quad= H\,div\lf[\nabla\vec{n}-\vec{n}\times\nabla^\perp\vec{n}\rg]\quad.
\end{array}
\ee
This becomes
\be
\label{VI.71}
\begin{array}{l}
-4 e^{2\la}\ \vec{H}\ (H^2-K)-2\Delta H\ \vec{n}\\[5mm]
\ds\quad\quad= div\lf[-2\nabla H\ \vec{n}+H\,\nabla\vec{n}-\vec{H}\times\nabla^\perp\vec{n}\rg]\quad.
\end{array}
\ee
Using now ''intrinsic notations'' (independent of the parametrization) on $(\Sigma^2,g)$, (\ref{VI.71}) says
\be
\label{VI.72}
\begin{array}{l}
4\ \vec{H}\ (H^2-K)+2\Delta_g H\ \vec{n}\\[5mm]
\ds\quad\quad= d^{\ast_g}\lf[-2d H\ \vec{n}+H\,d\vec{n}-\ast\vec{H}\times d\vec{n}\rg]\quad.
\end{array}
\ee
which is the desired identity (\ref{VI.64}) and theorem~\ref{th-VI.8} is proved in codimension 1.\hfill $\Box$

\medskip

Before to proceed to the proof of theorem~\ref{th-VI.8} in arbitrary codimension we first introduce some complex notations
that will be useful in the sequel. Assume $\vec{\Phi}$ is a conformal immersion into ${\R}^m$, one denotes
$z=x_1+ix_2$, $\p_z=2^{-1}(\p_{x_1}-i\p_{x_2})$, $\p_{\ov{z}}=2^{-1}(\p_{x_1}+i\p_{x_2})$.

\noindent Moreover we denote\footnote{Observe that the notation has been chosen in such a way that $\ov{\vec{e}_z}=\vec{e}_{\ov{z}}$.}
\[
\lf\{
\begin{array}{l}
\ds\vec{e}_{z}:=e^{-\la}\p_z\vec{\Phi}=2^{-1}(\vec{e}_1-i\vec{e}_2)\\[5mm]
\ds\vec{e}_{\ov{z}}:=e^{-\la}\p_{\ov{z}}\vec{\Phi}=2^{-1}(\vec{e}_1+i\vec{e}_2)
\end{array}
\rg.
\]
Observe that
\be
\label{VI.199a}
\lf\{
\begin{array}{l}
\ds\lf<\vec{e}_z,\vec{e}_z\rg>=0\\[5mm]
\ds\lf<\vec{e}_z,\vec{e}_{\ov{z}}\rg>=\frac{1}{2}\\[5mm]
\ds\vec{e}_z\wedge\vec{e}_{\ov{z}}=\frac{i}{2}\,\vec{e}_1\wedge\vec{e}_2
\end{array}
\rg.
\ee
Introduce moreover the {\it Weingarten Operator} expressed in our conformal coordinates $(x_1,x_2)$ :
\[
\vec{H}_0:=\frac{1}{2}\lf[\vec{\mathbb I}(e_1,e_1)-\vec{\mathbb I}(e_2,e_2)-2\,i\, \vec{\mathbb I}(e_1,e_2)\rg]\quad.
\]
With these notations the following lemma holds
\begin{Lm}
\label{lm-VI.11}
Let $\vec{\Phi}$ be a conformal immersion of $D^2$ into ${\R}^m$
\be
\label{z-VI.1}
\pi_T(\p_z\vec{H})-i\,\star(\p_z\vec{n}\wedge\vec{H})=-2\,\lf<\vec{H},\vec{H}_0\rg>\ \p_{\ov{z}}\vec{\Phi}
\ee
and hence
\be
\label{z-VI.2}
\begin{array}{l}
\ds\p_z\vec{H}-3\pi_{\vec{n}}(\p_z\vec{H})-i\,\star(\p_z\vec{n}\wedge\vec{H})\\[5mm]
\ds=-2\,\lf<\vec{H},\vec{H}_0\rg>\ \p_{\ov{z}}\vec{\Phi}-2\,\pi_{\vec{n}}(\p_z\vec{H})
\end{array}
\ee
\hfill $\Box$
\end{Lm}
\begin{Rm}
\label{rm-VI.11}
Observe that with these complex notations the identity (\ref{VI.64}),
which reads in conformal coordinates
\be
\label{VI.98}
\begin{array}{l}
\ds -\frac{e^{-2\la}}{2}\,div\lf[\nabla\vec{H}-3\pi_{\vec{n}}(\nabla \vec{H})+\star(\nabla^\perp\vec{n}\wedge\vec{H})\rg]\\[5mm]
\ds=\Delta_\perp\bH+\ti{A}(\vec{H})-2 |\bH|^2\ \bH\quad ,
\end{array}
\ee
becomes
\be
\label{z-VI.98}
\begin{array}{l}
\ds 4\,e^{-2\la}\ \Re\lf(\p_{\ov{z}}\lf[\pi_{\vec{n}}(\p_z\vec{H})+\lf<\vec{H},\vec{H}_0\rg>\ \p_{\ov{z}}\vec{\Phi}\rg]\rg)\\[5mm]
\ds =\Delta_\perp\bH+\ti{A}(\vec{H})-2\ |\bH|^2\ \bH\quad ,
\end{array}
\ee
\hfill $\Box$
\end{Rm}
{\bf Proof of lemma~\ref{lm-VI.11}.}
We denote by $(\bbe_1,\bbe_2)$ the orthonormal basis of $\vec{\Phi}_\ast(TD^2)$ given by
\[
\bbe_i=e^{-\la}\ \frac{\p \bP}{\p x_i}\quad .
\]
With these notations the second fundamental form $\bh$ which is a symmetric 2-form on $TD^2$  into  $(\vec{\Phi}_\ast TD^2)^\perp$
is given by
\be
\label{VI.74}
\begin{array}{l}
\bh=\sum_{\al,i,j}h^\al_{ij}\ \bn_\al\otimes(\bbe_i)^\ast\otimes(\bbe_j)^\ast\\[5mm]
\mbox{ with }\quad h^\al_{ij}=-e^{-\la}\,\lf(\frac{\p \bn_\al}{\p x_i},\bbe_j\rg)
\end{array}
\ee
We shall also denote
\[
\vec{h}_{ij}:=\vec{\mathbb I}(\bbe_i,\bbe_j)=\sum_{\al=1}^{m-2}h^\al_{ij}\ \vec{n}_\al
\]
In particular the mean curvature vector $\bH$ is given by
\be
\label{VI.75}
\bH=\sum_{\al=1}^{m-2} H^\al\,\bn_\al=\frac{1}{2}\sum_{\al=1}^{m-2}(h^\al_{11}+h^\al_{22})\, \bn_\al=\frac{1}{2}(\vec{h}_{11}+\vec{h}_{22})
\ee
Let $\bn$ be the $m-2$ vector of ${\R}^m$ given by $\bn=\bn_1\wedge\cdots\wedge\bn_2$. We identify vectors and $m-1$-vectors in ${\R}^m$ using the Hodge operator $\star$ of ${\R}^m$ for the canonical flat metric. Hence we have for instance
\be
\label{VI.76}
\star(\bn\wedge \bbe_1)=\bbe_2\quad\mbox{ and }\quad\star(\bn\wedge \bbe_2)=- \bbe_1
\ee
Since $\bbe_1,\bbe_2,\bn_1\cdots\bn_{m-2}$ is a basis of $T_{\bP(x_1,x_2)}{\R}^m$, we can write
for every $\al=1\cdots m-2$
\[
\nabla \bn_\al=\sum_{\beta=1}^{m-2}<\nabla \bn_\al,\bn_\beta>\, \bn_\beta+\sum_{i=1}^2<\nabla\bn_\al,\bbe_i>\,\bbe_i
\]
and consequently
\be
\label{VI.77}
\star(\bn\wedge\nabla^\perp\bn_\al)=<\nabla^\perp \bn_\al,\bbe_1>\ \bbe_2
-<\nabla^\perp \bn_\al,\bbe_2>\ \bbe_1
\ee
Hence
\[
\begin{array}{l}
\ds\star(\nabla^\perp\bn\wedge\vec{H})=-<\nabla^\perp \vec{H},\bbe_1>\ \bbe_2
+<\nabla^\perp \vec{H},\bbe_2>\ \bbe_1\\[5mm]
\quad=<\vec{H},\pi_{\vec{n}}(\nabla^\perp\bbe_1)>\ \bbe_2
-<\vec{H},\pi_{\vec{n}}(\nabla^\perp\bbe_2)>\ \bbe_1
\end{array}
\]
Using (\ref{VI.74}), we then have proved
\be
\label{z-VI.3}
\begin{array}{l}
\ds\star(\nabla^\perp\bn\wedge\vec{H})=\\[5mm]
\lf(
\begin{array}{c}
\ds-<\vec{H},\vec{h}_{12}>\ \p_{x_2}\vec{\Phi}\, +\,<\vec{H},\vec{h}_{22}>\ \p_{x_1}\vec{\Phi}\\[5mm]
\ds<\vec{H},\vec{h}_{11}>\ \p_{x_2}\vec{\Phi} \,-\,<\vec{H},\vec{h}_{12}>\ \p_{x_1}\vec{\Phi}
\end{array}
\rg)
\end{array}
\ee
The tangential projection of $\nabla\vec{H}$ is given by
\[
\begin{array}{l}
\pi_T(\nabla\vec{H})=<\nabla\vec{H},\vec{e}_1>\ \bbe_1+<\nabla\vec{H},\vec{e}_2>\ \bbe_2\\[5mm]
\quad=-<\vec{H},\pi_{\vec{n}}(\nabla \vec{e}_1)>\ \bbe_1-<\vec{H},\pi_{\vec{n}}(\nabla \vec{e}_2)>\ \bbe_2\quad.
\end{array}
\]
Hence
\be
\label{z-VI.4}
\begin{array}{l}
\pi_T(\nabla\vec{H})=\\[5mm]
\lf(
\begin{array}{c}
-<\vec{H},\vec{h}_{11}>\ \p_{x_1}\vec{\Phi}-<\vec{H},\vec{h}_{12}>\ \p_{x_2}\vec{\Phi}\\[5mm]
-<\vec{H},\vec{h}_{12}>\ \p_{x_1}\vec{\Phi}-<\vec{H},\vec{h}_{22}>\ \p_{x_2}\vec{\Phi}
\end{array}
\rg)
\end{array}
\ee
Combining (\ref{z-VI.3}) and (\ref{z-VI.4}) gives
\be
\label{z-VI.5}
\begin{array}{l}
\ds-\pi_T(\nabla\vec{H})-\star(\nabla^\perp\bn\wedge\vec{H})=\\[5mm]
\lf(
\begin{array}{c}
<\vec{H},\vec{h}_{11}-\vec{h}_{22}>\ \p_{x_1}\vec{\Phi}+2<\vec{H},\vec{h}_{12}>\ \p_{x_2}\vec{\Phi}\\[5mm]
2<\vec{H},\vec{h}_{12}>\ \p_{x_1}\vec{\Phi}+<\vec{H},\vec{h}_{22}-\vec{h}_{11}>\ \p_{x_2}\vec{\Phi}
\end{array}
\rg)
\end{array}
\ee
This last identity written with the complex coordinate $z$ is exactly (\ref{z-VI.1}) and lemma~\ref{lm-VI.11} is proved. \hfill $\Box$

\medskip

Before to move to the proof of theorem~\ref{th-VI.8} we shall need two more lemma. First we have
\begin{Lm}
\label{z-lm-VI.11}
Let $\vec{\Phi}$ be a conformal immersion of the disc $D^2$ into ${\R}^m$, denote $z:=x_1+ix_2$, $e^\la:=|\p_{x_1}\vec{\Phi}|=|\p_{x_2}\vec{\Phi}|$
denote
\be
\label{z-VI.200}
\vec{e}_i:=e^{-\la}\,\p_{x_i}\vec{\Phi}\quad,
\ee
and let $\vec{H}_0$ be the Weingarten Operator of the immersion expressed in the conformal coordinates $(x_1,x_2)$ :
\[
\vec{H}_0:=\frac{1}{2}\lf[{\mathbb I}(\vec{e}_1,\vec{e}_1)-{\mathbb I}(\vec{e}_1,\vec{e}_1)-2\,i\, {\mathbb I}(\vec{e}_1,\vec{e}_2)\rg]
\]
Then the following identities hold
\be
\label{VI.200}
\p_{\ov{z}}\lf[e^\la\, \vec{e}_{z}\rg]=\frac{e^{2\la}}{2}\vec{H}\quad,
\ee
and
\be
\label{VI.204}
\p_{z}\lf[e^{-\la}\vec{e}_z\rg]=\frac{1}{2}\, \vec{H}_0\quad.
\ee
\hfill$\Box$
\end{Lm}
\noindent{\bf Proof of lemma~\ref{z-lm-VI.11}.}
The first identity (\ref{VI.200}) comes simply from the fact that $\p_{\ov{z}}\p_z\vec{\Phi}=\frac{1}{4}\Delta\vec{\Phi}$, from  (\ref{z-VI.200}) and
the expression of the mean curvature vector in conformal coordinates that we have seen several times and which is given by
\[
\vec{H}=\frac{e^{-2\la}}{2}\,\Delta\vec{\Phi}\quad.
\]
It remains to prove the identity (\ref{VI.204}).
One has moreover
\be
\label{VI.201}
\p_{z}\lf[e^\la\vec{e}_z\rg]=\p_{z}\p_{z}\vec{\Phi}=\frac{1}{4}\lf[\p^2_{x_1^2}\vec{\Phi}-\p^2_{x_2^2}\vec{\Phi}-2\, i\ \p^2_{x_1x_2}\vec{\Phi}\rg]\quad.
\ee
In one hand the projection into the normal direction gives
\be
\label{VI.202}
\pi_{\vec{n}}\lf[\p^2_{x_1^2}\vec{\Phi}-\p^2_{x_2^2}\vec{\Phi}-2\, i\ \p^2_{x_1x_2}\vec{\Phi}\rg]=2\,{e^{2\la}}\, \vec{H}_0\quad.
\ee
In the other hand the projection into the tangent plane gives
\[
\begin{array}{l}
\ds\pi_{T}\lf[\p^2_{x_1^2}\vec{\Phi}-\p^2_{x_2^2}\vec{\Phi}-2\, i\ \p^2_{x_1x_2}\vec{\Phi}\rg]\\[5mm]
\ds =e^{-\la}\ \lf<\p_{x_1}\vec{\Phi},\lf[\p^2_{x_1^2}\vec{\Phi}-\p^2_{x_2^2}\vec{\Phi}-2\, i\ \p^2_{x_1x_2}\vec{\Phi}\rg]\rg>\ \vec{e}_1\\[5mm]
\ds +\,e^{-\la}\ \lf<\p_{x_2}\vec{\Phi},\lf[\p^2_{x_1^2}\vec{\Phi}-\p^2_{x_2^2}\vec{\Phi}-2\, i\ \p^2_{x_1x_2}\vec{\Phi}\rg]\rg>\ \vec{e}_2\quad.
\end{array}
\]
This implies after some computation
\be
\label{VI.203}
\begin{array}{l}
\ds\pi_{T}\lf[\p^2_{x_1^2}\vec{\Phi}-\p^2_{x_2^2}\vec{\Phi}-2\, i\ \p^2_{x_1x_2}\vec{\Phi}\rg]\\[5mm]
\ds =2\,e^\la\ \lf[\p_{x_1}\la-i\p_{x_2}\la\rg]\ \vec{e}_1-2\,e^\la\ \lf[\p_{x_2}\la+i\p_{x_1}\la\rg]\ \vec{e}_2\\[5mm]
\ds=8\ \p_{z}e^\la\ \vec{e}_z\quad.
\end{array}
\ee
The combination of (\ref{VI.201}), (\ref{VI.202}) and (\ref{VI.203}) gives
\[
\p_{z}\lf[e^\la\vec{e}_z\rg]=\frac{e^{2\la}}{2}\, \vec{H}_0+2\,\p_{z}e^\la\ \vec{e}_z\quad,
\]
which implies (\ref{VI.204}).\hfill$\Box$

\medskip

The last lemma we shall need in order to prove theorem~\ref{th-VI.8} is the  Codazzi-Mainardi identity that we recall
and prove below.

\begin{Lm}
\label{z-lm-VI.12}{\bf[Codazzi-Mainardi Identity.]}
Let $\vec{\Phi}$ be a conformal immersion of the disc $D^2$ into ${\R}^m$, denote $z:=x_1+ix_2$, $e^\la:=|\p_{x_1}\vec{\Phi}|=|\p_{x_2}\vec{\Phi}|$
denote
\be
\label{z-VI.200}
\vec{e}_i:=e^{-\la}\,\p_{x_i}\vec{\Phi}\quad,
\ee
and let $\vec{H}_0$ be the Weingarten Operator of the immersion expressed in the conformal coordinates $(x_1,x_2)$ :
\[
\vec{H}_0:=\frac{1}{2}\lf[{\mathbb I}(\vec{e}_1,\vec{e}_1)-{\mathbb I}(\vec{e}_1,\vec{e}_1)-2\,i\, {\mathbb I}(\vec{e}_1,\vec{e}_2)\rg]
\]
Then the following identity holds
\be
\label{z-VI.203}
e^{-2\la}\,\p_{\ov{z}}\lf(e^{2\la}\,<\vec{H},\vec{H}_0>\rg)=<\vec{H},\p_{z}\vec{H}>+<\vec{H}_0,\p_{\ov{z}}\vec{H}>\quad.
\ee
\hfill$\Box$
\end{Lm}
{\bf Proof of lemma~\ref{z-lm-VI.12}.}
Using (\ref{VI.204}) we obtain
\[
\begin{array}{l}
\ds<\p_{\ov{z}}\vec{H}_0, \vec{H}>=2\, \lf<\p_{\ov{z}}\lf[\p_z\lf(e^{-2\la}\,\p_z\vec{\Phi}\rg)\rg],\vec{H}\rg>\\[5mm]
\ds\quad=2\, \lf<\p_{{z}}\lf[\p_{\ov{z}}\lf(e^{-2\la}\,\p_z\vec{\Phi}\rg)\rg],\vec{H}\rg>\quad.
\end{array}
\]
Thus
\[
\begin{array}{l}
\ds<\p_{\ov{z}}\vec{H}_0, \vec{H}>\\[5mm]
\ds=-4\,\lf<\p_{{z}}\lf[\p_{\ov{z}}\la\ e^{-2\la}\ \p_z\vec{\Phi}\rg],\vec{H}\rg>+\lf<\p_z\lf[\frac{e^{-2\la}}{2}\ \Delta\vec{\Phi}\rg],\vec{H}\rg>\\[5mm]
\ds=-2\p_{\ov{z}}\la\ \lf<\vec{H}_0,\vec{H}\rg>+\lf<\p_z\vec{H},\vec{H}\rg>\quad.
\end{array}
\]
This last identity implies the Codazzi-Mainardi identity (\ref{z-VI.203}) and lemma~\ref{z-lm-VI.12} is proved.\hfill $\Box$

\medskip

\noindent{\bf Proof of theorem~\ref{th-VI.8}.}
Du to lemma~\ref{lm-VI.11}, as explained in remark~\ref{rm-VI.11}, it suffices to prove in conformal parametrization the identity
(\ref{z-VI.98}).
First of all we observe that
\be
\label{z-VI.205a}
\begin{array}{l}
4\,e^{-2\la}\,\Re\lf(\pi_{\vec{n}}\lf(\p_{\ov{z}}\lf[\pi_{\vec{n}}(\p_z\vec{H})\rg]\rg)\rg)\\[5mm]
\ds=e^{-2\la}\ \pi_{\vec{n}}\lf(div\lf[\pi_{\vec{n}}(\nabla\vec{H})\rg]\rg)\\[5mm]
=\Delta_\perp\vec{H}
\end{array}
\ee
The tangential projection gives
\be
\label{z-VI.205}
\begin{array}{l}
\ds 4\,e^{-2\la}\,\pi_{T}\lf(\p_{\ov{z}}\lf[\pi_{\vec{n}}(\p_z\vec{H})\rg]\rg)\\[5mm]
\ds =8\, e^{-2\la}\ \lf<\p_{\ov{z}}(\pi_{\vec{n}}(\p_z\vec{H})),\vec{e}_z\rg>\ \vec{e}_{\ov{z}}\\[5mm]
\ds+\,8\, e^{-2\la}\ \lf<\p_{\ov{z}}(\pi_{\vec{n}}(\p_z\vec{H})),\vec{e}_{\ov{z}}\rg>\ \vec{e}_{{z}}
\end{array}
\ee
Using the fact that $\vec{e}_z$ and $\vec{e}_{\ov{z}}$ are orthogonal to the normal plane we have in one
hand using (\ref{VI.200})
\be
\label{z-VI.206}
\begin{array}{l}
\ds\lf<\p_{\ov{z}}(\pi_{\vec{n}}(\p_z\vec{H})),\vec{e}_z\rg>=-e^{-\la}\,\lf<\pi_{\vec{n}}(\p_z\vec{H}),\p_{\ov{z}}\lf[ e^\la\,\vec{e}_z\rg]\rg>\\[5mm]
\ds\quad=-\frac{e^{\la}}{2}\,\lf<\p_z\vec{H},\vec{H}\rg>
\end{array}
\ee
and in the other hand using (\ref{VI.204})
\be
\label{z-VI.207}
\begin{array}{l}
\ds\lf<\p_{\ov{z}}(\pi_{\vec{n}}(\p_z\vec{H})),\vec{e}_{\ov{z}}\rg>=-e^{\la}\,\lf<\pi_{\vec{n}}(\p_z\vec{H}),\p_{\ov{z}}\lf[ e^{-\la}\,\vec{e}_{\ov{z}}\rg]\rg>\\[5mm]
\ds\quad=-\frac{e^{\la}}{2}\,\lf<\p_z\vec{H},\ov{\vec{H}_0}\rg>
\end{array}
\ee
Combining (\ref{z-VI.205}), (\ref{z-VI.206}) and (\ref{z-VI.207}) we obtain
\be
\label{z-VI.208}
\begin{array}{l}
\ds 4\,e^{-2\la}\,\pi_{T}\lf(\p_{\ov{z}}\lf[\pi_{\vec{n}}(\p_z\vec{H})\rg]\rg)\\[5mm]
=-4\ e^{-2\la}\lf[\lf<\p_z\vec{H},\vec{H}\rg>\,\p_{\ov{z}}\vec{\Phi}+\lf<\p_z\vec{H},\ov{\vec{H}_0}\rg>\,\p_z\vec{\Phi}\rg]
\end{array}
\ee
Putting (\ref{z-VI.205a}) and (\ref{z-VI.208}) together we obtain
\be
\label{z-VI.209}
\begin{array}{l}
\ds 4\,e^{-2\la}\,\Re\lf(\p_{\ov{z}}\lf[\pi_{\vec{n}}(\p_z\vec{H})\rg]\rg)\\[5mm]
=\Delta_\perp\vec{H}-4\,e^{-2\la}\Re\lf[\lf[\lf<\p_{{z}}\vec{H},\vec{H}\rg>+\lf<\p_{\ov{z}}\vec{H},{\vec{H}_0}\rg>\rg]\,\p_{\ov{z}}\vec{\Phi}\rg]
\end{array}
\ee
Using Codazzi-Mainardi identity (\ref{z-VI.203})   and using also again identity (\ref{VI.204}), (\ref{z-VI.209}) becomes
\be
\label{z-VI.209}
\begin{array}{l}
\ds 4\,e^{-2\la}\,\Re\lf(\p_{\ov{z}}\lf[\pi_{\vec{n}}(\p_z\vec{H})+<\vec{H},\vec{H}_0>\ \p_{\ov{z}}\vec{\Phi}\rg]\rg)\\[5mm]
=\Delta_\perp\vec{H}+2\,\Re\lf(\lf<\vec{H},\vec{H}_0\rg>\ \ov{\vec{H}_0}\rg)\quad.
\end{array}
\ee
The definition (\ref{VI.44a}) of $\ti{A}$ gives
\[
\ti{A}(\vec{H})=\sum_{i,j=1}^2<\vec{H},\vec{h}_{ij}>\ \vec{h}_{ij}\quad
\]
hence a short elementary computation gives
\[
\begin{array}{l}
\ds\ti{A}(\vec{H})-2|\vec{H}|^2\ \vec{H}\\[5mm]
=2^{-1}\,\lf<\vec{H},{\vec{h}_{11}-\vec{h}_{22}}\rg>\ (\vec{h}_{11}-\vec{h}_{22})+2<\vec{H},\vec{h}_{12}>\ \vec{h}_{12}
\end{array}
\]
Using $\vec{H}_0$ this expression becomes
\be
\label{z-VI.210}
\ti{A}(\vec{H})-2|\vec{H}|^2\ \vec{H}=2\Re\lf(\lf<\vec{H},\vec{H}_0\rg>\ \ov{\vec{H}_0}\rg)
\ee
Combining (\ref{z-VI.209}) and (\ref{z-VI.210}) gives (\ref{z-VI.98})
which is the desired inequality and theorem~\ref{th-VI.8} is proved.\hfill$\Box$

\newpage

\subsection{Construction of Isothermal Coordinates.}

In the previous subsection we discussed the difficulty to work with the Willmore surfaces equation
due to the huge invariance group given by the space of positive diffeomorphisms of the surface.

A classical way to by-pass this difficulty consists in ''breaking'' the symmetry group (or ''gauge group'') of 
coordinates by restricting to a special subclass satisfying the {\it Coulomb condition}. In the present subsection we will
explain why this choice corresponds to the {\it conformal condition}. We have seen in the previous subsection that, in such coordinates,
the Willmore surface equations can be written in a {\it conservative-elliptic} form which is critical with respect to the $L^2-$norm of the second fundamental form. This triggers the hope to give answers to the analysis questions we raised  
for Willmore surfaces. 

Breaking the symmetry group of coordinates by taking isothermal ones is however not enough per se. The uniformization theorem tells us that, taking the conformal class defined by the pull-back metric $\vec{\Phi}^\ast g_{{\R}^m}$ on $\Sigma^2$, there is a system of coordinates on the fundamental domain of either ${\C}\cup{\infty}$, ${\C}$ or the Poincar\'e half-plane corresponding to this class in which our immersion is conformal. However, in a minimization procedure for instance, taking a minimizing sequence of immersions $\vec{\Phi}_k$, assuming the conformal class defined by $\vec{\Phi}_k$ on $\Sigma^2$ is controlled - is not converging to the boundary of the moduli space -  there is a-priori no control of the conformal factor corresponding to the pull-back metric $\vec{\Phi}_k^\ast g_{{\R}^m}$ and the fact that we are
in conformal coordinates is not helping much. We need then to have 

\medskip

\centerline{ Conformal coordinates + estimates of the conformal factor.}

\medskip

In Gauge theory such as in Yang-Mills problem in 4 dimension for instance, the Coulomb choice of gauge is the one that provides estimates
of the connection that will be controlled by the gauge invariant $L^2-$norm of the curvature, provided this energy is below
some universal threshold. Similarly in the present 
situation the Coulomb Gauge - or conformal choice of coordinates - will permit to control the $L^\infty$ norm
of the pull-back metric (the conformal factor) with the help of the ''gauge invariant'' $L^2$-norm of the second
fundamental form provided it stays below the universal threshold : $\sqrt{8\pi/3}$. This will be the up-shot of the
present subsection. This $L^\infty$ control of the conformal factor will be an application of integrability of compensation and
Wente estimates more specifically.

\subsubsection{The Chern Moving Frame Method.}

We present a method originally due to S.S.Chern in order to construct local isothermal coordinates. It is based on the following observation.

Let $\vec{\Phi}$ be a conformal immersion of the disc $D^2$ introduce the following tangent frame :
\[
(\vec{e}_1,\vec{e}_2)=e^{-\la}\ (\p_{x_1}\vec{\Phi},\p_{x_2}\vec{\Phi})\quad,
\]
where $e^\la=|\p_{x_1}\vec{\Phi}|=|\p_{x_2}\vec{\Phi}|$. 

A simple computation shows
\be
\label{VI.100a}
\lf<\vec{e}_1,\nabla\vec{e}_2\rg>=-\nabla^\perp\la\quad,
\ee
and in particular it follows
\be
\label{VI.100}
div\lf<\vec{e}_1,\nabla\vec{e}_2\rg>=0\quad.
\ee
This identity can be writen independently of the parametrization as follows
\be
\label{VI.101}
d^{\ast_g}\lf<\vec{e}_1,d\vec{e}_2\rg>=0\quad.
\ee
It appears clearly as the {\it Coulomb condition} : cancelation of the codifferential of the connection on the $S^1-$tangent orthonormal frame bundle given by the 1-form $i\lf<\vec{e}_1,d\vec{e}_2\rg>$ taking value into the Lie algebra $i{\R}$. 
Sections of this bundles are given by  maps $(\vec{f}_1,\vec{f}_2)$ from $D^2$ into $\vec{\Phi}_\ast(TD^2\times TD^2)$
such that $(\vec{f}_1,\vec{f}_2)(x_1,x_2)$ realizes a positive orthonormal basis of $\vec{\Phi}_\ast (T_{(x_1,x_2)}D^2)$. 

The passage from one section $(\vec{e}_1,\vec{e}_2)$ to another section $(\vec{f}_1,\vec{f}_2)$
is realized through a change of gauge which corresponds to the action of an $SO(2)$ rotation $e^{i\theta}$ on the tangent space $\vec{\Phi}_\ast (T_{(x_1,x_2)}D^2)$ :
\[
\vec{f}_1+i\vec{f}_2=e^{i\theta}\ (\vec{e}_1+i\vec{e}_2)\quad.
\]
The expression of the same connection but in the new trivialization given by the section $(\vec{f}_1,\vec{f}_2)$
satisfies the classical gauge change formula for an $S^1$-bundle :
\be
\label{VI.101a}
i\lf<\vec{f}_1,d\vec{f}_2\rg>=i\lf<\vec{e}_1,d\vec{e}_2\rg>+id\theta\quad.
\ee
The curvature of this connection is given by
\be
\label{VI.102}
\begin{array}{l}
\ds i\,d\lf<\vec{e}_1,d\vec{e}_2\rg>\\[5mm]
\ds=i\lf[<D_{e_1}\vec{e}_1,D_{e_2}\vec{e}_2>-<D_{e_2}\vec{e}_1,D_{e_1}\vec{e}_2>\rg]\ e_1^\ast\wedge e_2^\ast\\[5mm]
\ds=i\lf[<\vec{\mathbb I}(\vec{e}_1,\vec{e}_1),\vec{\mathbb I}(\vec{e}_2,\vec{e}_2)>-|\vec{\mathbb I}(\vec{e}_1,\vec{e}_2)|^2\rg]\ \ e_1^\ast\wedge e_2^\ast\\[5mm]
\ds=i\ K\ dvol_g
\end{array}
\ee
where we recall the notations we already introduced : $e_i$ is the vector field on $D^2$ given by $d\vec{\Phi}\cdot e_i=\vec{e}_i$, 
$D_{e_i}\vec{e}_j:=\pi_{\vec{n}}(d\vec{e}_j\cdot e_i)$ and $K$ is the Gauss curvature of $(D^2,\vec{\Phi}^{\ast}g)$.
In the last identity we have made use of Gauss theorem (theorem 2.5 chap. 6 in \cite{doC2}).

\medskip

Combining (\ref{VI.100a}) and (\ref{VI.102}) gives the well known expression of the Gauss curvature 
in isothermal coordinates in terms of the conformal factor $\la$ :
\be
\label{VI.103}
-\Delta\la=<\nabla^\perp\vec{e}_1,\nabla\vec{e}_2>=e^{2\la}\ K\quad.
\ee

\medskip

We have seen how any conformal parametrization generate a Coulomb frame in the tangent bundle. 
S.S.Chern observed that this is in fact an exact matching, the reciproque is also true : starting from a Coulomb frame one can generate
isothermal coordinates.

Let $\vec{\Phi}$ be an immersion of the disc $D^2$ into ${\R}^m$ and let $(\vec{e}_1,\vec{e}_2)$ be a {\it Coulomb tangent
orthonormal moving frame} : a map from $D^2$ into $\vec{\Phi}_\ast(TD^2\times TD^2)$
such that $(\vec{e}_1,\vec{e}_2)(x_1,x_2)$ realizes a positive orthonormal basis of $\vec{\Phi}_\ast (T_{(x_1,x_2)}D^2)$
and such that condition (\ref{VI.101}) is satisfied.
 
Let $\la$ be the solution of
\be
\label{VI.104}
\lf\{
\begin{array}{l}
\ds d\la=\ast_g<\vec{e}_1,d\vec{e}_2>\\[5mm]
\ds \int_{\p D^2}\la=0
\end{array}
\rg.
\ee
Denote moreover $e_i:=d\vec{\Phi}^{-1}\cdot\vec{e}_i$ and $(e_1^\ast,e_2^\ast)$ to be the dual basis to $(e_1,e_2)$. The Cartan formula\footnote{The Cartan formula
for the exterior differential of a 1 form $\al$ on a differentiable manifolfd $M^m$ says that for any pair of vector fields $X,Y$ on this manifold the following identity holds
\be
\label{VI.105}
d\al(X,Y)=d(\al(Y))\cdot X-d(\al(X))\cdot Y-\al([X,Y])
\ee
see corollary 1.122 chapter I of \cite{GHL}.} for the exterior differential of a 1-form implies
\be
\label{VI.106}
\begin{array}{l}
\ds de_i^\ast(e_1,e_2)=d(e_i^\ast(e_2))\cdot e_1-d(e_i^\ast(e_1))\cdot e_2- e_i^\ast([e_1,e_2])\\[5mm]
\ds \quad=- e_i^\ast([e_1,e_2])\\[5mm]
\ds\quad=-((\vec{\Phi}^{-1})^\ast e_i^\ast)([d\vec{\Phi}\cdot e_1,d\vec{\Phi}\cdot e_2])\\[5mm]
\ds\quad=-((\vec{\Phi}^{-1})^\ast e_i)([\vec{e}_1,\vec{e}_2])
\end{array}
\ee
The Levi-Civita connection $\nabla$ on $\vec{\Phi}_\ast T D^2$ issued from the restriction to the tangent space to $\vec{\Phi}(D^2)$ of the canonical metric in ${\R}^m$ is
given by $\nabla_X\sigma:=\pi_T(d\sigma\cdot X)$ where $\pi_T$ is the orthogonal projection onto the tangent plane. The Levi-Civita connection moreover is symmetric\footnote{We recall that a connection $\nabla$ on the tangent bundle of a manifold $M^m$ is symmetric if for any pair of tangent fields $X$ and $Y$ one has
\[
T(X,Y):=\nabla_XY-\nabla_YX-[X,Y]=0\quad.
\]} (see \cite{doC2} theorem 3.6 chap. 2)
hence we have in particular
\be
\label{VI.107}
[\vec{e}_1,\vec{e}_2]=\nabla_{e_1}\vec{e}_2-\nabla_{e_2}\vec{e}_1=\pi_T(d\vec{e}_2\cdot e_1-d\vec{e}_1\cdot e_2)
\ee
Since $\vec{e}_1$ and $\vec{e}_2$ have unit length, the  tangential projection of $d\vec{e}_1$ (resp. $d\vec{e}_2$ ) are oriented along $\vec{e}_2$ (resp. $\vec{e}_1$).
So we have
\be
\label{VI.108}
\lf\{
\begin{array}{l}
\ds\pi_T(d\vec{e}_2\cdot e_1)=<d\vec{e}_2,\vec{e}_1>\cdot e_1\ \vec{e}_1\\[5mm]
\ds\pi_T(d\vec{e}_1\cdot e_2)=<d\vec{e}_1,\vec{e}_2>\cdot e_2\ \vec{e}_2
\end{array}
\rg.
\ee 
Combining (\ref{VI.106}), (\ref{VI.107}) and (\ref{VI.108}) gives then
\be
\label{VI.109}
\lf\{
\begin{array}{l}
\ds de_1^\ast(e_1,e_2)=-<d\vec{e}_2,\vec{e}_1>\cdot e_1\\[5mm]
 \ds de_2^\ast(e_1,e_2)=<d\vec{e}_1,\vec{e}_2>\cdot e_2
\end{array}
\rg.
\ee
Equation (\ref{VI.104}) gives
\be
\label{VI.110}
\lf\{
\begin{array}{l}
\ds -<d\vec{e}_2,\vec{e}_1>\cdot e_1=(\ast_gd\la,e_1)=-d\la\cdot e_2\\[5mm]
 \ds <d\vec{e}_1,\vec{e}_2>\cdot e_2=(\ast_gd\la,e_2)=d\la\cdot e_1
\end{array}
\rg.
\ee
Thus combining (\ref{VI.109}) and (\ref{VI.110}) gives then
\be
\label{VI.111}
\lf\{
\begin{array}{l}
\ds de_1^\ast=-d\la\cdot e_2\ e_1^\ast\wedge e_2^\ast= d\la\wedge e_1^\ast\\[5mm]
\ds de_2^\ast=d\la\cdot e_1\ e_1^\ast\wedge e_2^\ast=d\la\wedge e_2^\ast\quad.
\end{array}
\rg.
\ee
We have thus proved at the end
\be
\label{VI.112}
\lf\{
\begin{array}{l}
d\lf(e^{-\la}e_1^\ast\rg)=0\\[5mm]
d\lf(e^{-\la}e_2^\ast\rg)=0
\end{array}
\rg.
\ee
Introduce $(\phi_1,\phi_2)$ the functions with average 0 on the disc $D^2$ such that
\[
d\phi_i:=e^{-\la} e_i^\ast\quad.
\]
since rank$(d\phi_1,d\phi_2)=2$, $\phi:=(\phi_1,\phi_2)$ realizes a diffeomorphism from $D^2$ into $\phi(D^2)$. 

From the previous identity
we have
\[
e^{-\la\circ\phi^{-1}}\ g(e_j,\p_{y_i}\phi^{-1})=e^{-\la\circ\phi^{-1}}\ (e_j^\ast,\p_{y_i}\phi^{-1})=\delta_{ij}\quad.
\]
where $g:=\vec{\Phi}^\ast g_{{\R}^m}$. This implies
\be
\label{VI.113}
g(\p_{y_i}{\phi}^{-1},\p_{y_j}{\phi}^{-1})=e^{2\la\circ\phi^{-1}}\delta_{ij}\quad.
\ee
or in other words
\be
\label{VI.114}
<\p_{y_i}(\vec{\Phi}\circ\phi^{-1}),\p_{y_j}(\vec{\Phi}\circ\phi^{-1})>=e^{2\la\circ\phi^{-1}}\,\delta_{ij}\quad.
\ee
This says that $\vec{\Phi}\circ\phi^{-1}$ is a conformal immersion from $\phi(D^2)$ into ${\R}^m$. The Riemann Mapping theorem\footnote{See for instance \cite{Rud} chapter 14.}  gives the existence of a biholomorphic
diffeomorphism $h$ from $D^2$ into $\phi(D^2)$. Thus $\vec{\Phi}\circ\phi^{-1}\circ h$ realizes a conformal immersion from $D^2$ onto $\vec{\Phi}(D^2)$.

\subsubsection{The space of Lipschitz Immersions with $L^2-$bounded Second Fundamental Form.}

In the previous subsection we have seen the equivalence between {\it tangent Coulomb moving frames} and {\it isothermal coordinates}. It remains now
to construct {\it tangent Coulomb moving frames} in order to produce {\it isothermal coordinates} in which Willmore surfaces equation admits a nice {\it conservative
elliptic form}. For a purpose that will become clearer later in this book we are extending the framework of smooth immersions to a more general framework : 
the space of {\it Lipschitz immersions with $L^2-$bounded second fundamental form}.

\medskip

Let $\Sigma^2$ be a smooth compact oriented 2-dimensional manifold (with or without boundary). 
Let $g_0$ be a reference smooth metric on $\Sigma$. One defines the Sobolev spaces $W^{k,p}(\Sigma,{\R}^m)$ of measurable maps from $\Sigma$ into 
${\R}^m$ in the following way
\[
W^{k,p}(\Sigma^2,{\R}^m)=\lf\{f\ :\ {\Sigma^2}\rightarrow {\R}^m\ ;\ \sum_{l=0}^k\int_{\Sigma}|\nabla^l f|_{g_0}^p\ dvol_{g_0}<+\infty\rg\}
\]
Since $\Sigma^2$ is assumed to be compact it is not difficult to see that this space is independent of the choice we have made of $g_0$.

\medskip

A lipschitz immersion of $\Sigma^2$ into ${\R}^m$ is a  map $\vec{\Phi}$ in $W^{1,\infty}(\Sigma^2,{\R}^m)$ for which
\be
\label{VI.115}
\exists\ c_0>0\quad\mbox{ s.t. }\quad|d\vec{\Phi}\wedge d\vec{\Phi}|_{g_0}\ge c_0>0 \quad,
\ee
where $d\vec{\Phi}\wedge d\vec{\Phi}$ is a 2-form on $\Sigma^2$ taking values into 2-vectors from ${\R}^m$ and given in local coordinates
by $2\,\p_{x_1}\vec{\Phi}\wedge\p_{x_2}\vec{\Phi}\ dx_1\wedge dx_2$.         
The condition (\ref{I.1}) is again independent of the choice of the metric $g_0$ . This assumption implies that $g:=\vec{\Phi}^\ast g_{{\R}^m}$ defines
an $L^\infty$ metric comparable to the reference metric $g_0$ : there exists $C>0$ such that
\be
\label{VI.115a}
\begin{array}{l}
\ds\forall X\in T\Sigma^2\ \\[5mm]
 \quad\quad C^{-1}g_0(X,X)\le \vec{\Phi}^\ast g_{{\R}^m}(X,X)\le C\ g_0(X,X)\quad.
 \end{array}
\ee
For a Lipschitz immersion satisfying (\ref{VI.115}) we can define the Gauss map as being the following measurable map in $L^\infty(\Sigma)$
\[
\vec{n}_{\vec{\Phi}}:=\star\frac{\p_{x_1}\vec{\Phi}\wedge\p_{x_2}\vec{\Phi}}{|\p_{x_1}\vec{\Phi}\wedge\p_{x_2}\vec{\Phi}|}\quad.
\] 
for an arbitrary choice of local positive coordinates $(x_1,x_2)$. We then introduce the space ${\mathcal E}_{\Sigma}$ of Lipschitz immersions of $\Sigma$ 
with bounded second fundamental form as follows :
\[
\mathcal{E}_\Sigma:=\lf\{
\begin{array}{c}
\ds\vec{\Phi}\in W^{1,\infty}(\Sigma,{\R}^m)\quad\mbox{ s.t. } \vec{\Phi} \mbox{ satisfies }(\ref{VI.115})\\[5mm]
\ds\mbox{ and }\quad\int_{\Sigma}|d\vec{n}|_g^2\ dvol_g<+\infty
\end{array}
\rg\}\quad .
\]
When $\Sigma$ is not compact we extend the definition of ${\mathcal E}_\Sigma$ as follows : we require $\vec{\Phi}\in W^{1,\infty}_{loc}(\Sigma,{\R}^m)$, we require that (\ref{VI.115})
holds locally on any compact subset of $\Sigma^2$ and we still require the global $L^2$ control of the second fundamental form
\[
\int_{\Sigma}|d\vec{n}|_g^2\ dvol_g<+\infty\quad.
\]

\subsubsection{Energy controlled liftings of $W^{1,2}-$maps into the Grassman manifold $Gr_2({\R}^m)$.}

In the next subsection we will apply the Chern moving frame method in the context of lipschitz Immersions with $L^2-$bounded Second Fundamental Form.
To that aim we need first to
construct local Coulomb tangent moving frames with controlled $W^{1,2}$ energy. We shall do it in two steps. First we will explore the possibility to ''lift'' the Gauss map and to construct tangent moving frame with bounded $W^{1,2}-$energy.
This is the purpose of the present subsection. The following result
{\bf lifting theorem} proved by F.H\'elein in \cite{He}.

\begin{Th}
\label{th-VI.12}
Let $\vec{n}$ be a $W^{1,2}$ map from the disc $D^2$ into the Grassman manifold\footnote{The Grassman manifold $Gr_{m-2}({\R}^m)$ can be seen 
as being the sub-manifold of the euclidian space $\wedge^{m-2}{\R}^m$ of $m-2$-vectors in ${\R}^m$ made of unit simple $m-2$-vectors
and then one defines
\[
W^{1,2}(D^2,Gr_{m-2}({\R}^m)):=\lf\{\vec{n}\in  W^{1,2}(D^2,\wedge^{m-2}{\R}^m)\ ;\ \vec{n}\in Gr_{m-2}({\R}^m)\mbox{ a.e. }\rg\}\quad.
\] } of oriented $m-2$-planes
in ${\R}^m$ : $Gr_{m-2}({\R}^m)$.
There exists a constant $C>0$, such that, if one assumes that  
\be
\label{VI.116}
\int_{D^2}|\nabla\vec{n}|^2<\frac{8\pi}{3}\quad,
\ee
then there exists $\vec{e}_1$ and $\vec{e}_2$ in $W^{1,2}(D^2,S^{m-1})$ such that
\be
\label{VI.117}
\vec{n}=\star(\vec{e}_1\wedge\vec{e}_2)\quad,
\ee
and\footnote{The condition (\ref{VI.117}) together with the fact that the $\vec{e}_i$ are in $S^{m-1}$ valued imply, since $\vec{n}$ has norm one, that $\vec{e}_1$ and $\vec{e}_2$
are orthogonal to each other.}
\be
\label{VI.118}
\int_{D^2}\sum_{i=1}^2|\nabla\vec{e}_i|^2\ dx_1\,dx_2\le C\ \int_{D^2}|\nabla\vec{n}|^2\ dx_1,dx_2\quad.
\ee
\hfill$\Box$
\end{Th}
The requirement for the $L^2$ nom of $\nabla\vec{n}$ to be below a threshold, inequality (\ref{VI.116}), is necessary and it is conjectured in \cite{He} that 
$8\pi/3$ could be replaced by $8\pi$ and that this should be optimal.

\medskip

We can illustrate the need to have an energy restriction such as (\ref{VI.116}) with the following example :

\medskip

Let $\pi$ be the stereographic projection from $S^2$ into ${\C}\cup\{\infty\}$ which sends the north pole $N=(0,0,1)$ to $0$ and the south pole $S=(0,0,-1)$ to $\infty$.
For $\la$ sufficiently small we consider on $D^2$ the following map $\vec{n}_{\la}$ taking value into $S^2$ :
\[
\lf\{
\begin{array}{l}
\ds\vec{n}_\la(x):=\pi^{-1}(\la\, x) \quad\quad\mbox{ for }\quad |x|\le 1/2\\[5mm]
\ds\vec{n}_\la(x):=\frac{(1-r)\ \pi^{-1}(\la\, x)+(r-1/2)\ S}{|(1-r)\ \pi^{-1}(\la\, x)+(r-1/2)\ S|}\quad\mbox{ for } \frac{1}{2}<|x|\le1.
\end{array}
\rg.
\]
The map $\vec{n}_\la$ has been constructed in such a way that on $B_{1/2}(0)$ it covers the most part of $S^2$ conformally\footnote{The map $x\rightarrow\pi^{-1}(\la\, x)$ is conformal.} and the missing small
part which is a small geodesic ball centered at the south pole $S$ is covered in the annulus $B_1(0)\setminus B_{1/2}(0)$ using some simple interpolation
 between $\pi^{-1}(\la\, x)$ with the south pole composed with the reprojection onto $S^2$. $\vec{n}_\la$ is clearly surjective onto $S^2$, is sending $\p D^2$ to the south pole
 and since points of the north hemisphere admit exactly one preimage by $\vec{n}_\la$, $\vec{n}_\la$ is of degree one. Because of the conformality of $\vec{n}_\la$ on the major part of the image
 it is not difficult to verify that
 \be
 \label{VI.119}
 \int_{D^2}|\nabla \vec{n}_\la|^2=8\pi+o_\la(1)
\ee
where $o_\la(1)$ is a positive function which goes to zero as $\la$ goes to $+\infty$. 
Since $\vec{n}_\la$ is constant on $\p D^2$ equal to the south pole $S=(0,0,-1)$ we can extend $\vec{n}_\la$ by $S$ on the whole plane ${\R}^2$ and we still have
that the $L^2$ norm of $\nabla\vec{n}_\la$ on ${\R}^2$ is equal to $8\pi+o_\la(1)$.
\be
 \label{VI.119a}
 \int_{{\R}^2}|\nabla \vec{n}_\la|^2=8\pi+o_\la(1)
\ee
Let now $0<\rho<1$ and denote $\vec{n}_\la^\rho(x):\vec{n}_\la(x/\rho)$. Due to the conformal invariance of the Dirichlet energy we have
\be
 \label{VI.120}
 \int_{D^2}|\nabla \vec{n}_\la^\rho|^2=8\pi+o_\la(1)
\ee
Consider an orthonormal moving frame\footnote{Such an orthonormal moving frame exists because
the pull-back by $\vec{n}^\rho_\la$ of the frame bundle of $S^2$ over $D^2$ is a trivial bundle since $D^2$ is contractible.} $(\vec{e}_1,\vec{e}_2)$ such that 
$$
\vec{n}^\rho_\la=\star(\vec{e}_1\wedge\vec{e}_2)\quad.
$$
Identifying the horizontal plane with ${\R}^2$, $\vec{e}_1$ realizes a map from $\p D^2$
into the unit circle $S^1$ of ${\R}^2$. Assume the topological degree of $\vec{e}_i$  ($i=1,2$) would be zero, then $\vec{e}_i$ would be homotopic to a constant.  We would then realize this homotopy in the annulus $B_2(0)\setminus B_1(0)$ in such a way that both $\vec{n}^\rho_\la$ and $\vec{e}_i$ would be constant on $\p B_2(0)$. We could
then identify every points on $\p B_2$ in such a way that $\vec{n}_\la^\rho$ realizes a degree
one map from $S^2$ into $S^2$. Any degree one map from $S^2$ into $S^2$ is homotopic to the identity map and the pull-back bundle $((\vec{n}_\la^\rho)^{-1})^{-1}TS^2$ is then bundle equivalent to $TS^2$ (see theorem 4.7 in chapter 1 of \cite{Hus}).$\vec{e}_i$ would realize a global section of this bundle that would be trivial
which would contradicts Brouwer's theorem. Hence the restriction to $\p D^2$ of $\vec{e}_i$ has a non zero degree\footnote{This degree is in fact equal to 2 which is the Euler characteristic  of $S^2$. See theorem 11.16 of \cite{BoTu} and example 11.18.} and by homotopy this is also the case on any circle $\p B_r(0)$ for $1>r>\rho$. Since $\vec{e}_1$ has non zero degree
on each of these circles one has
\[
\forall\, \rho<r<1\quad\quad 2\pi\le\lf|\int_{\p B^2_r}(\vec{e}_1)^\ast d\theta\rg|\le\ (2\pi r)^{1/2}\ \lf[\int_{\p B^2_r}|\nabla\vec{e}_1|^2\rg]^{1/2}
\]
We deduce from this inequality for $i=1,2$
\be
\label{VI.121}
\int_{D^2}|\nabla\vec{e}_i|^2\ dx_1\,dx_2\ge 2\pi\log\frac{1}{\rho}\rightarrow +\infty\ \mbox{ as }\ \rho\rightarrow 0\quad.
\ee
By taking $\la\rightarrow +\infty$ and $\rho\rightarrow +\infty$, we can deduce the following lemma. 
\begin{Lm}
\label{lm-VI.13}
There exists a sequence $\vec{n}_k$ in $W^{1,2}(D^2,Gr_{m-2}({\R}^m))$ such that
\[
\int_{D^2}|\nabla\vec{n}_k|^2\ dx_1\,dx_2\longrightarrow 8\pi
\]
and
\[
\inf\lf\{
\begin{array}{l}
\ds \int_{D^2}\sum_{i=1}^2|\nabla\vec{e}_i|^2\ dx_1\,dx_2\ \mbox{ s.t. }\\[5mm]
 \vec{e}_i\in W^{1,2}(D^2,S^{m-1})\mbox{ and }\vec{e}_1\wedge\vec{e}_2=\star\vec{n}
\end{array}
\rg\}\rightarrow +\infty
\]
\hfill $\Box$
\end{Lm}
This lemma says that it is necessary to stay strictly below the threshold $8\pi$ for the Dirichlet energy of maps into $Gr_{m-2}({\R}^m)$ 
in order to hope to construct energy controlled liftings.
However, if one removes the requirement to control the energy, every map $\vec{n}\in W^{1,2}(D^2, Gr_{m-2}({\R}^m))$ admits a $W^{1,2}$
lifting. The following theorem is proved in \cite{He} chapter 5.2.
\begin{Th}
\label{th-VI.14}
Let  $\vec{n}$ in $W^{1,2}(D^2,Gr_{m-2}({\R}^m))$, then there exists $\vec{e}_1,\vec{e}_2\in W^{1,2}(D^2,S^{m-1})$ such that
\[
\vec{e}_1\wedge\vec{e}_2=\star\vec{n}\quad.
\]
\hfill $\Box$
\end{Th}
We can then make the following observation. Let $\vec{n}\in W^{1,2}(D^2,S^{m-1})$ and $\vec{e}_1,\vec{e}_2\in W^{1,2}(D^2,S^{m-1})$ given by the previous theorem.
Consider $\theta\in W^{1,2}$ to be the unique solution of
\[
\lf\{
\begin{array}{l}
\ds\Delta\theta=-div<\vec{e}_1,\nabla\vec{e}_2>\quad\quad\mbox{ in }\quad D^2\\[5mm]
\ds \frac{\p\theta}{\p\nu}=-\lf<\vec{e}_1,\frac{\p\vec{e}_2}{\p\nu}\rg>
\end{array}
\rg.
\]
then the lifting $\vec{f}=\vec{f}_1+i\vec{f}_2$ given by
\[
\vec{f}_1+i\vec{f}_2=e^{i\theta}\ (\vec{e}_1+i\vec{e}_2)\quad,
\]
is Coulomb. Let $\la$ such that $-\nabla^\perp\la=<\vec{f}_1,\nabla\vec{f}_2>$ and such that $\int_{\p D^2}\la=0$,
one easily sees that $\la$ satisfies
\be
\label{VI.122}
\lf\{
\begin{array}{l}
-\Delta\la=<\nabla^\perp\vec{f}_1,\nabla\vec{f}_2>\quad\quad\mbox{ in }\quad D^2\\[5mm]
\la=0\quad\quad\mbox{ on }\quad \p D^2
\end{array}
\rg.
\ee
Observe that since $\nabla\vec{e}_1$ is perpendicular to $\vec{e}_1$ and since $\nabla\vec{e}_2$ is perpendicular to $\vec{e}_2$ one has
\be
\label{VI.123}
<\nabla^\perp\vec{f}_1,\nabla\vec{f}_2>=<\pi_{\vec{n}}( \nabla^\perp\vec{f}_1), \pi_{\vec{n}}( \nabla\vec{f}_2)>
\ee
We make now use of (\ref{VI.65ab}) and we deduce
\be
\label{VI.124}
\begin{array}{rl}
\ds\pi_{\vec{n}}( \nabla^\perp\vec{f}_i)&=(-1)^{m-1}\ \vec{n}\res(\vec{n}\res\nabla^\perp\vec{f}_i)\\[5mm]
 &\ds=\nabla^\perp(\pi_{\vec{n}}(\vec{f}_i))+(-1)^{m-1}\ \nabla^\perp\vec{n}\res(\vec{n}\res\vec{f}_i)\\[5mm]
 &\ds\quad+(-1)^{m-1}\ \vec{n}\res(\nabla^\perp\vec{n}\res\vec{f}_i)
 \end{array}
\ee
Using the fact that $\pi_{\vec{n}}(\vec{f}_i)\equiv 0$ we obtain from (\ref{VI.124}) that
\be
\label{VI.125}
\int_{D^2}|\pi_{\vec{n}}( \nabla\vec{f}_i)|^2\ dx_1\,dx_2\le 2\int_{D^2}|\nabla\vec{n}|^2\ dx_1\,dx_2
\ee
Combining (\ref{VI.123}) and (\ref{VI.125}) we obtain
\be
\label{VI.126}
\int_{D^2}|<\nabla^\perp\vec{f}_1,\nabla\vec{f}_2>|\ dx_1\,dx_2\le 2\int_{D^2}|\nabla\vec{n}|^2\ dx_1\,dx_2
\ee
This estimate together with standard elliptic estimates (see for instance \cite{Ad}) give
\be
\label{VI.127}
\begin{array}{l}
\ds\|<\vec{f}_1,\nabla\vec{f}_2>\|_{L^{2,\infty}(D^2)}=\|\nabla\la\|_{L^{2,\infty}(D^2)}\\[5mm]
\ds\quad\quad\le C\ \int_{D^2}|\nabla\vec{n}|^2\ dx_1\,dx_2\quad.
\end{array}
\ee
From (\ref{VI.125}) and (\ref{VI.127}) we deduce the following theorem\footnote{Similarly it would be interesting to explore the possibility to construct {\it global gauges with estimates} 
in non abelian gauge theory. For instance one can ask the following question : for a given curvature  of a $W^{1,2}$ 
$SU(n)-$connection over the 4-dimensional ball, can one construct a gauge in which
the $L^{4,\infty}$ norm of the connection is controlled by the $L^2-$norm of the curvature ?}  

.
\begin{Th}
\label{th-VI.15}
Let  $\vec{n}$ in $W^{1,2}(D^2,Gr_{m-2}({\R}^m))$, then there exists $\vec{e}_1,\vec{e}_2\in W^{1,2}(D^2,S^{m-1})$ such that
\[
\vec{e}_1\wedge\vec{e}_2=\star\vec{n}\quad,
\]
\[
div<\vec{e}_1,\nabla\vec{e}_2>=0\quad,
\]
and satisfying
\[
\sum_{i=1}^2\|\nabla\vec{e}_i\|_{L^{2,\infty}(D^2)}\le C\ \|\nabla\vec{n}\|_{L^2(D^2)}\ \lf[1+\|\nabla\vec{n}\|_{L^2(D^2)}\rg]\quad.
\]
\hfill $\Box$
\end{Th}

\subsubsection{Conformal parametrizations for lipschitz Immersions with $L^2-$bounded Second Fundamental Form.}

Starting from a lipschitz immersion $\vec{\Phi}$ of a surface $\Sigma$ with $L^2-$bounded Second Fundamental Formwe can then cover $\Sigma$ by disks in such a way that on each of these
discs, in some coordinates system, the $L^2$ norm of $\nabla\vec{n}$ is below $\sqrt{8\pi/3}$. Due to theorem~\ref{th-VI.12}, there exists a
$W^{1,2}$ frame $(\vec{e}_1,\vec{e}_2)$ with controlled energy. In order to produce a Coulomb frame on each of these discs one minimizes
\[
\min\lf\{\int_{D^2}|<\vec{f}_1,d\vec{f}_2>|_g^2\ dvol_g\quad;\quad \vec{f}_1+i\vec{f}_2=e^{i\theta}\, (\vec{e}_1+i\vec{e}_2) \rg\}
\]
where $g=\vec{\Phi}^\ast g_{{\R}^m}$. Using (\ref{VI.101a}) the previous problem corresponds to minimize the following energy
\[
\int_{D^2} |d\theta+<\vec{e}_1,d\vec{e}_2>)|^2_g\ dvol_g
\]
among all $\theta\in W^{1,2}(D^2,{\R})$. This Lagrangian is convex on the Hilbert space $W^{1,2}(D^2,{\R})$ and 
goes to $+\infty$ as $\|\theta\|_{W^{1,2}}\rightarrow +\infty$. Then there exists a unique minimum satisfying
\[
\lf\{
\begin{array}{l}
\ds d^{\ast_g}\lf[d\theta+<\vec{e}_1,d\vec{e}_2>\rg]=0\quad\quad\mbox{ in }D^2\\[5mm]
\ds\iota_{\p D^2}^\ast (\ast_g\lf[d\theta+<\vec{e}_1,d\vec{e}_2>\rg])=0\quad\quad\mbox{ on }\p D^2
\end{array}
\rg.
\]
where $\iota_{\p D^2}$ is the canonical inclusion of $\p D^2$ into $\ov{D^2}$. Then $\vec{f}:=\vec{f}_1+i\vec{f}_2$
given by $\vec{f}=e^{i\theta} \vec{e}$ is Coulomb :
\be
\label{VI.128}
\lf\{
\begin{array}{l}
\ds d^{\ast_g}\lf[<\vec{f}_1,d\vec{f}_2>\rg]=0\quad\quad\mbox{ in }D^2\\[5mm]
\ds\iota_{\p D^2}^\ast (\ast_g\lf[<\vec{f}_1,d\vec{f}_2>\rg])=0\quad\quad\mbox{ on }\p D^2
\end{array}
\rg.
\ee
We are now in position to start the Chern moving frame method in order to produce a conformal parametrization of $\vec{\Phi}$ on this disc. This however has to be done with the additional difficulty of keeping track of the regularity of the different actors at each step of the construction.

\medskip

First we introduce the function $\la\in W^{1,2}$ satisfying (\ref{VI.104}). The second equation of (\ref{VI.128}) implies
that the restriction to the boundary of $D^2$ of the one form $d\la$ is equal to zero. Hence this last fact combined 
with the second equation of (\ref{VI.104}) implies that $\la$ is identically equal to zero on $\p D^2$.

We have then
\be
\label{VI.129}
\lf\{
\begin{array}{l}
d\ast_gd\la=-<d\vec{e}_1,d\vec{e}_2>\quad\quad\mbox{ on }D^2\\[5mm]
\la=0\quad\quad\quad\mbox{ on }\p D^2
\end{array}
\rg.
\ee
which reads in the canonical coordinates of $D^2$
\[
\label{VI.130}
\lf\{
\begin{array}{l}
\ds\frac{\p}{\p x_i}\lf[\frac{g^{ij}}{\sqrt{det\,g}}\frac{\p\la}{\p{x_j}}\rg]=<\frac{\p\vec{e}_1}{\p x_1},\frac{\p\vec{e}_2}{\p x_2}>-<\frac{\p\vec{e}_1}{\p x_2},\frac{\p\vec{e}_2}{\p x_1}>\mbox{ on }D^2\\[7mm]
\la=0\quad\quad\quad\mbox{ on }\p D^2
\end{array}
\rg.
\]
where we are using an implicit summation in $i$ and $j$ and where $g^{ij}$ are the coefficient to the inverse matrix
to $g_{ij}:=<\p_{x_i}\vec{\Phi},\p_{x_j}\vec{\Phi}>$.
We are now in position to make use of the following generalization of Wente's theorem due to
S.Chanillo and Y.Y. Li \cite{ChLi}.
\begin{Th}
\label{th-VI.16}
Let $a$ and $b$ be two functions in $W^{1,2}(D^2,{\R})$. Let $(a^{ij})_{1\le i,j\le 2}$ be a $2\times 2$ symmetric matrix
valued map in $L^\infty(D^2)$ such that there exists $C>0$ for which
\[
\forall\, \xi=(\xi_1,\xi_2)\in{\R}^2\ \forall x\in D^2\quad\quad C^{-1}\ |\xi|^2\le a^{ij}(x)\xi_i\xi_j\le C\ |\xi|^2
\]
Let $\varphi$ be the solution in $W^{1,p}(D^2,{\R})$ for any $1\le p<2$ of the following equation
\be
\label{VI.131}
\lf\{
\begin{array}{l}
\ds \frac{\p}{\p x_i}\lf[a^{ij}\,\frac{\p\la}{\p{x_j}}\rg]=\frac{\p a}{\p x_1}\,\frac{\p b}{\p x_2}-\frac{\p a}{\p x_2}\frac{\p b}{\p x_1}\quad\quad\mbox{ on }D^2\\[7mm]
\ds\varphi=0\quad\quad\quad\mbox{ on }\p D^2\quad.
\end{array}
\rg.
\ee
Then $\varphi\in L^\infty\cap W^{1,2}(D^2,{\R})$ and there exists $C>0$ independent of $a$ and $b$ such that
\be
\label{VI.132}
\|\varphi\|_{L^\infty}+\|\nabla\varphi\|_{L^2}\le C\ \|\nabla a\|_{L^2}\ \|\nabla b\|_{L^2}\quad.
\ee
\hfill $\Box$
\end{Th}
Using theorem~\ref{th-VI.16} we obtain that $\la\in L^\infty(D^2,{\R})$. Let $e_i$ be the frame on $D^2$ given by
\[
d\vec{\Phi}\cdot e_i=\vec{e}_i
\]
and denote $e_i=e_i^1\p_{x_1}+e_i^2\p_{x_2}$. We have
\[
\sum_{k=1}^2 e_i^k\ g_{kj}=<\vec{e}_i,\p_{x_j}\vec{\Phi}>\quad.
\]
from which we deduce
\[
e^k_i=\sum_{j=1}^2g^{kj}\,<\vec{e}_i,\p_{x_j}\vec{\Phi}>
\]
Because of (\ref{VI.115a}) the maps $g^{kj}$ are in $L^\infty$. Thus
\be
\label{VI.133}
\begin{array}{l}
\ds e^\ast_i=\sum_{j=1}^2\lf[g^{1j}\,<\vec{e}_i,\p_{x_j}\vec{\Phi}>\,dx_1\rg.\\[5mm]
\ds\quad\quad\quad\lf. +g^{2j}\,<\vec{e}_i,\p_{x_j}\vec{\Phi}>\,dx_2\rg]\, \in L^\infty(D^2)\ .
\end{array}
\ee
Denote by $(f_1,f_2)$ the frame on $D^2$ such that $d\vec{\Phi}\cdot f_i=\vec{f}_i$, we have $f^ast=f_1^\ast+if_2^\ast=e^{-i\theta}\, (e_1^\ast+ie_2^\ast)$ which is
in $L^\infty(D^2)$ due to (\ref{VI.133}). Hence the map $\phi=(\phi_1,\phi_2)$ which is given by
\[
d\phi_i:=e^{-\la} f_i^\ast\quad.
\]
is a bilipschitz diffeomorphism between $D^2$ and $\phi(D^2)$. Equation (\ref{VI.114}) says that $\vec{\Phi}\circ\phi^{-1}$ is a conformal lipshitz immersion from $\phi(D^2)$ into ${\R}^m$. The Riemann Mapping theorem  gives the existence of a biholomorphic
diffeomorphism $h$ from $D^2$ into $\phi(D^2)$. Thus $\vec{\Phi}\circ\phi^{-1}\circ h$ realizes a conformal immersion from $D^2$ onto $\vec{\Phi}(D^2)$ which is in $W^{1,\infty}_{loc}(D^2,{\R}^m)$. We have then established the following theorem.

\begin{Th}
\label{th-VI.17} {\bf [Existence of a smooth conformal structure]}
Let $\Sigma^2$ be a closed smooth 2-dimensional manifold. Let $\vec{\Phi}$ be an element of ${\mathcal E}_\Sigma$ : a Lipshitz immersion with $L^2-$bounded second fundamental form.
Then there exists a finite covering of $\Sigma^2$ by discs $(U_i)_{i\in I}$ and Lipschitz diffeomorphisms $\psi_i$ from $D^2$ into $U_i$ such that $\vec{\Phi}\circ\psi_i$ realizes a
lipschitz conformal immersion of $D^2$. Since $\psi_j^{-1}\circ\psi_i$ on $\psi_i^{-1}(U_i\cap U_j)$ is conformal and positive (i.e. holomorphic) the system of charts $(U_i,\psi_i)$ defines
a smooth conformal structure $c$ on $\Sigma^2$ and in particular there exists a constant scalar curvature metric $g_c$ on $\Sigma^2$ and a lipshitz diffeomorphism $\psi$ of $\Sigma^2$
such that $\vec{\Phi}\circ\psi$ realizes a conformal immersion of the riemann surface $(\Sigma^2,g_c)$.\hfill $\Box$
\end{Th}

\subsubsection{Weak Willmore immersions}

Let $\Sigma$ be a smooth compact oriented 2-dimensional manifold and let $\vec{\Phi}$ be a Lipschitz immersion
with $L^2-$bounded second fundamental form : $\vec{\Phi}$ is an element of ${\mathcal E}_\Sigma$. Because of the
previous subsection we know that for any smooth disc $U$ included in $\Sigma$ there exists a Lipschitz diffeomorphism
from $D^2$ into $U$ such that $\vec{\Phi}\circ\Psi$ is a conformal Lipschitz immersion of $D^2$. In this chart the $L^2-$ norm of the second fundamental form
which is assumed to be finite is given by
\[
\int_{U}|\vec{\mathbb I}|_g^2\ dvol_g=\int_{D^2}|\nabla\vec{n}|^2\ dx_1\,dx_2\quad.
\]
where $\vec{\mathbb I}$ and $\vec{n}$  denote respectively the second fundamental form and the Gauss map of the conformal immersion $\vec{\Phi}\circ\Psi$.
The mean curvature vector of $\vec{\Phi}\circ\Psi$, that we simply denote $\vec{H}$,
is given by
\[
\vec{H}:=2^{-1}\,e^{-2\la}\ \lf[\vec{\mathbb I}(\p_{x_1},\p_{x_1})+\vec{\mathbb I}(\p_{x_2},\p_{x_2})\rg]\in L^2(D^2)\quad.
\]
Hence we have that
\[
\nabla \vec{H}\in H^{-1}(D^2)\quad,
\]
moreover, we have also
\[
\star(\nabla^\perp\vec{n}\wedge\vec{H})\in L^1(D^2)\quad.
\]
Using the expression (\ref{VI.65ac}), we also deduce that
\[
\pi_{\vec{n}}(\nabla\vec{H})\in H^{-1}+L^1(D^2)\quad.
\]
Hence for any immersion $\vec{\Phi}$ in ${\mathcal E}_\Sigma$ the quantity
\[
\nabla\vec{H}-3\pi_{\vec{n}}(\nabla\vec{H})+\star(\nabla^\perp\vec{n}\wedge\vec{H})\in H^{-1}+L^1(D^2)
\]
defines a distribution in ${\mathcal D}'(D^2)$. 

\medskip

Using corollary~\ref{co-VI.9} one can then generalize the notion of Willmore immersion that we defined for smooth immersions to immersions in ${\mathcal E}_\Sigma$ in the following way.
\begin{Dfi}
\label{df-VI.17}
Let $\Sigma$ be a smooth compact oriented 2-dimensional manifold and let $\vec{\Phi}$ be a Lipschitz immersion
with $L^2-$bounded second fundamental form. 

$\vec{\Phi}$ is called a weak Willmore immersion  if, in any lipschitz conformal chart $\Psi$ from $D^2$ into $(\Sigma,\vec{\Phi}^\ast g_{{\R}^m})$, the following holds
\[
div\lf[\nabla\vec{H}-3\pi_{\vec{n}}(\nabla\vec{H})+\star(\nabla^\perp\vec{n}\wedge\vec{H})\rg]=0\ \mbox{ in }{\mathcal D}'(D^2)\ .
\]
where $\vec{H}$ and $\vec{n}$ denote respectively the mean curvature vector and the Gauss map of $\vec{\Phi}$ in the chart $\Psi$ and the operators $div$, $\nabla$ and $\nabla^\perp$ are taken with respect to the flat metric\footnote{$div\,X=\p_{x_1}X_1+\p_{x_2}X_2$,
$\nabla f=(\p_{x_1}f,\p_{x_2}f)$ and $\nabla^\perp f=(-\p_{x_2}f,\p_{x_1}f)$.} in $D^2$.
\hfill$\Box$
\end{Dfi}

Having defined weak Willmore immersion that will naturally come in our minimization procedure it is a fair question to ask whether solutions to this elliptic non-linear system for which we have seen that it is critical in two dimensions for the Willmore energy are smooth. Or in other words we are asking the following question :

\medskip

\noindent{\it Are weak Willmore immersions smooth Willmore immersions ?}

\medskip 

We will answer positively to that question in the next sections.

\subsubsection{Isothermal Coordinates with Estimates.}

In the previous subsections we explained how to produce locally and globally conformal parametrizations, Coulomb gauges, for Lipschitz immersions with $L^2-$bounded second fundamental forms. This has been obtained in a qualitative way, without any care of establishing estimates. The goal of the present subsection is to remedy to it in 
giving $L^\infty$ controls of the metric in conformal parametrization  by the mean of several quantities ($L^2-$norm of the second fundamental form, area of the image, distance of the images of 2 distinct points) and under the assumption that the $L^2-$norm of the second fundamental form is below some threshold. Precisely we shall prove the following result.

\begin{Th}
\label{th-VI.18} {\bf [Control of local isothermal coordinates]}
Let $\vec{\Phi}$ be a lipschitz conformal immersion\footnote{A Lipschitz conformal map from $D^2$ into ${\R}^m$ satisfying (\ref{VI.115})} from the disc $D^2$ into ${\R}^m$. Assume 
\be
\label{VI.134}
\begin{array}{l}
\ds \int_{D^2}|\nabla\vec{n}_{\vec{\Phi}}|^2<8\pi/3\quad.
\end{array}
\ee
Denote $e^\la:=|\p_{x_1}\vec{\Phi}|=|\p_{x_2}\vec{\Phi}|$. Then for any $0<\rho<1$ there exists a constant $C_\rho$ independent of $\vec{\Phi}$ such that
\be
\label{VI.135}
\sup_{p\in B^2_\rho(0)}e^\la(p)\le C_\rho\ \lf[ \mbox{Area}(\vec{\Phi}(D^2))\rg]^{1/2}\ \exp\lf(C\ \int_{D^2}|\nabla\vec{n}_{\vec{\Phi}}|^2\rg)\quad.
\ee
Moreover, for two given distinct points $p_1$ and $p_2$ in the interior of $D^2$ and again for $0<\rho<1$ there exists a constant $C>0$
independent of $\vec{\Phi}$ such that
\be
\label{VI.136}
\begin{array}{rl}
\ds\|\la\|_{L^\infty(B^2_\rho(0))}\le & \ds C_\rho\ \int_{D^2}|\nabla\vec{n}_{\vec{\Phi}}|^2+C_\rho\ \lf|\log\frac{|\vec{\Phi}(p_1)-\vec{\Phi}(p_2)|}{|p_2-p_1|}\rg|\\[5mm]
 &\ds+C_\rho\ \log^+\lf[C_\rho\ \mbox{Area}(\vec{\Phi}(D^2))\rg]\quad.
 \end{array}
\ee
where $\log^+:=\max\{\log,0\}$. \hfill$\Box$
\end{Th}
{\bf Proof of theorem~\ref{th-VI.18}.}
Denote
\[
(\vec{f}_1,\vec{f}_2)=e^{-\la}\ (\p_{x_1}\vec{\Phi},\p_{x_2}\vec{\Phi})\quad.
\]
we have seen that 
\[
-\nabla^\perp\la=<\vec{f}_1,\nabla\vec{f}_2>\quad,
\]
from which we deduce
\be
\label{VI.137}
-\Delta\la=<\p_{x_1}\vec{f}_1,\p_{x_2}\vec{f}_2>-<\p_{x_2}\vec{f}_1,\p_{x_1}\vec{f}_2>\quad.
\ee
Let $(\vec{e}_1,\vec{e}_2)$ be the frame given by theorem~\ref{th-VI.12}. There exists $\theta$ such that
such that $e^{i\theta} (\vec{e}_1+i\vec{e}_2)=\vec{f}_1+i\vec{f}_2$ and hence
\[
<\nabla^\perp\vec{f}_1,\nabla\vec{f}_2>=<\nabla^\perp\vec{e}_1,\nabla\vec{e}_2>\quad.
\]
and $\la$ satisfies
\be
\label{VI.138}
-\Delta\la=<\p_{x_1}\vec{e}_1,\p_{x_2}\vec{e}_2>-<\p_{x_2}\vec{e}_1,\p_{x_1}\vec{e}_2>\quad.
\ee
Let $\mu$ be the solution of 
\be
\label{VI.139}
\lf\{
\begin{array}{l}
\ds-\Delta\mu=<\p_{x_1}\vec{e}_1,\p_{x_2}\vec{e}_2>-<\p_{x_2}\vec{e}_1,\p_{x_1}\vec{e}_2>\quad\\[5mm]
\ds\mu=0\quad\quad\quad\mbox{ on }\p D^2
\end{array}
\rg.
\ee
We are now in position to apply Wente's theorem~\ref{th-III.1} and we obtain the $\mu\in L^\infty(D^2)$ together with the estimate
\be
\label{VI.140}
\begin{array}{l}
\ds\|\mu\|_{L^\infty(D^2)}\le (2\pi)^{-1}\ \lf[\int_{D^2}|\nabla \vec{e}_1|^2\rg]^\frac{1}{2}\ \lf[\int_{D^2}|\nabla \vec{e}_2|^2\rg]^\frac{1}{2}\\[5mm]
\ds\quad\quad\quad\quad\ \le C\ \int_{D^2}|\nabla\vec{n}_{\vec{\Phi}}|^2
\end{array}
\ee
$\mu$ has been chosen in such a way that $\nu:=\la-\mu$ is harmonic. We deduce that
\[
\Delta e^{\nu}=\ |\nabla\nu|^2\ e^{\nu}\ge 0\quad.
\]
Harnack inequality (see \cite{GT}  theorem 2.1 for instance) implies that for any $0<\rho<1$ there exists a constant $C_\rho>0$ such that
\be
\label{VI.141}
\begin{array}{rl}
\ds\sup_{p\in D^2_\rho}e^{\nu(p)}&\ds\le C_\rho\ \lf[\int_{D^2}e^{2\nu}\rg]^\frac{1}{2}\\[5mm]
&\ds\le C_\rho\ \lf[\int_{D^2}e^{2\la}\rg]^\frac{1}{2}\ e^{\|\mu\|_\infty}\quad.
\end{array}
\ee
Combining (\ref{VI.140}), (\ref{VI.141}) and the fact that $$\int_{D^2} e^{2\la}=\mbox{Area}(\vec{\Phi}(D^2))$$ we obtain (\ref{VI.135}).

\medskip

Let $p_1\ne p_2$ be two distinct points in $D^2$ such that $\vec{\Phi}(p_1)\ne\vec{\Phi}(p_2)$. We have
\[
|\vec{\Phi}(p_1)-\vec{\Phi}(p_2)|\le\int_{[p_1,p_2]}e^\la\ d\sigma\quad.
\]
where $[p_1,p_2]$  is the segment joining the two points $p_1$ and $p_2$ in $D^2$ and $d\sigma$ is the length element along this segment.

The mean-value theorem implies then that there exists a point $p_0\in[p_1,p_2]$ such that 
\[
\la(p_0)\ge \log\frac{|\vec{\Phi}(p_1)-\vec{\Phi}(p_2)|}{|p_2-p_1|}\quad.
\]
Thus in particular we have
\be
\label{VI.142}
\begin{array}{rl}
\ds\nu(p_0)=(\la-\mu)(p_0)\ge&\ds \log\frac{|\vec{\Phi}(p_1)-\vec{\Phi}(p_2)|}{|p_2-p_1|}\\[5mm]
 &\ds -\ C\ \int_{D^2}|\nabla\vec{n}_{\vec{\Phi}}|^2\quad.
\end{array}
\ee
Let $\rho<1$ such that $B_\rho(0)$ contains the segment $[p_1,p_2]$. Let $r:=(1+\rho)/2$. The explicit expression of
the Poisson Kernel\footnote{See for instance \cite{GT} theorem 2.6.} gives
\[
\nu(p_0)=\frac{r^2-|p_0|^2}{2\pi r}\int_{\p B^2_r(0)}\frac{\nu(z)}{|z-p_0|}\ d\sigma(z)\quad,
\]
where $d\sigma(z)$ is the volume form on $\p B^2_r(0)$. Let $\nu^+=\sup\{0,\nu\}$ and $\nu^-=\inf\{0,\nu\}$. We then have
\be
\label{VI.143}
\begin{array}{l}
\ds\int_{\p B^2_r(0)}\frac{\nu^-(z)}{|z-p_0|}\ d\sigma(z)\ge \frac{2\pi r}{r^2-|p_0|^2}\ \nu(p_0)\\[5mm]
\ds \quad\quad\quad-\int_{\p B^2_r(0)}\frac{\nu^+(z)}{|z-p_0|}\ d\sigma(z)
\end{array}
\ee
We apply (\ref{VI.135}) on $B^2_r(0)$ and we have
\be
\label{VI.145}
\nu^+\le \frac{1}{2}\log^+\lf[C_r\,\int_{D^2}e^{2\la}\rg]+\|\mu\|_\infty\quad.
\ee
Combining this inequality with (\ref{VI.140}) and (\ref{VI.143}) gives
\be
\label{VI.144}
\begin{array}{l}
\ds\int_{\p B^2_r(0)}\frac{\nu^-(z)}{|z-p_0|}\ d\sigma(z)\ge -C_r\,\lf|\log\frac{|\vec{\Phi}(p_1)-\vec{\Phi}(p_2)|}{|p_2-p_1|}\rg| \\[5mm]
\ds \quad\quad\quad-C_r \log^+\lf[C_r\,\int_{D^2}e^{2\la}\rg]-C_r\,\int_{D^2}|\nabla\vec{n}_{\vec{\Phi}}|^2\quad.
\end{array}
\ee
We deduce from this inequality and from (\ref{VI.145}) combined with (\ref{VI.140}) that
\be
\label{VI.146}
\begin{array}{l}
\ds\int_{\p B^2_r(0)}|\nu|\ d\sigma\le C_r\,\lf|\log\frac{|\vec{\Phi}(p_1)-\vec{\Phi}(p_2)|}{|p_2-p_1|}\rg| \\[5mm]
\ds \quad\quad\quad\quad+C_r \log^+\lf[C_r\,\int_{D^2}e^{2\la}\rg]+C_r\,\int_{D^2}|\nabla\vec{n}_{\vec{\Phi}}|^2\ .
\end{array}
\ee
Using again the explicit expression of the Poisson Kernel, we have that for any $p$ in $B^2_\rho(0)$
\[
\nu(p)=\frac{r^2-|p|^2}{2\pi r}\int_{\p B^2_r(0)}\frac{\nu(z)}{|z-p|}\ d\sigma(z)\quad,
\]
For any point $p$ in $B^2_\rho(0)$ and any point $z$ in $\p B^2_r(0)$, $$(r^2-|p|^2)/2\pi r\,|z-p|$$ is bounded from above and from below
by constants which only depend on $\rho$. Thus there exists $C_\rho>0$ such that
\be
\label{VI.147}
\begin{array}{l}
\ds\|\nu\|_{L^\infty(B^2_\rho(0))}\le C_\rho\,\lf|\log\frac{|\vec{\Phi}(p_1)-\vec{\Phi}(p_2)|}{|p_2-p_1|} \rg|\\[5mm]
\ds \quad\quad\quad\quad+C_\rho \log^+\lf[C_r\,\int_{D^2}e^{2\la}\rg]+C_\rho\,\int_{D^2}|\nabla\vec{n}_{\vec{\Phi}}|^2 
\end{array}
\ee
The combination of (\ref{VI.140}) and (\ref{VI.147}) gives the inequality (\ref{VI.136}) and theorem~\ref{th-VI.18} is proved. \hfill $\Box$

\subsection{Conformal Willmore Surfaces.}

\subsubsection{The problem of passing to the limit in the Willmore surface equations.}

One of our ambition is to pass to the limit while following a sequence of Willmore immersions into ${\R}^m$ of closed surfaces with uniformly bounded energy, area and topology. Considering the current of integration in ${\R}^m$ along
the corresponding immersed surfaces, Federer Fleming theorem\footnote{See a presentation of this funding result of the Geometric Measure Theory in \cite{Mor} for instance.} asserts that from such a sequence one can always
extract a subsequence that weakly converges\footnote{This weak convergence holds in fact for the flat distance - see \cite{Fe}.} to a limiting rectifiable cycle. What can be said about this cycle
is one of the main question we raise in this part of the course.

\medskip

Before to look at the problem globally, it is worth to first look at a sequence of Willmore immersions of the disc $D^2$ assuming that the area is uniformly bounded and that the $L^2$ norm of the second fundamental
form stays below some small value. We choose this value to be less or equal to $\sqrt{8\pi/3}$ in such a way that, due to theorem~\ref{th-VI.18}, we have a conformal parametrization
$\vec{\Phi}_k$ of our immersion in which the pull-back metric $\vec{\Phi}_k^\ast g_{{\R}^m}$ does not  degenerate on the interior of $D^2$, having assumed also the {\it non collapsing  + non expanding contitions} : 

\medskip

{\it there exists $p_1\ne p_2$
such that $\log |\vec{\Phi}_k(p_1)-\vec{\Phi}_k(p_2)|$ is uniformly bounded.}

\medskip

For this sequence of conformal parametrization $\vec{\Phi}_k$ from $D^2$ into ${\R}^m$ the following holds :
\begin{itemize}
\item[i)]
\[
\int_{D^2}|\nabla\vec{n}_{\vec{\Phi}_k}|^2\le 8\pi/3\quad,
\]
\item[ii)]
\[
\limsup_{k\rightarrow +\infty}\mbox{Area}(\vec{\Phi}_k(D^2))<+\infty
\]
\item[iii)]
\[
\forall\,\rho<1\ \quad\limsup_{k\rightarrow +\infty}\|\log|\nabla\vec{\Phi}_k|\|_{L^\infty(B^2_\rho(0))}<+\infty\quad, 
\]
\item[iv)]
\[
div\lf[\nabla \vec{H}_k-3\pi_{\vec{n}_k}(\nabla\vec{H}_k)+\star(\nabla^\perp\vec{n}_k\wedge\vec{H}_k)\rg]=0\quad.
\]
\end{itemize}
Moreover we can assume that the immersion does not shift to infinity by taking
\be
\label{VI.147a}
\vec{\Phi}_k(0)=0\quad.
\ee
Since the immersions $\vec{\Phi}_k$ are conformal we have that
\[
\Delta\vec{\Phi}_k=2\,e^{2\la_k}\ \vec{H}_k
\]
where $e^{\la_k}=|\p_{x_1}\vec{\Phi}_k|=|\p_{x_2}\vec{\Phi}_k|$.

\medskip

The conditions i) and iii) imply then the following
\be
\label{VI.148}
\forall\,\rho<1\ \quad\limsup_{k\rightarrow +\infty}\|\Delta\vec{\Phi}_k\|_{L^2(B^2_\rho(0))}<+\infty\quad.
\ee
Moreover condition ii) implies
\be
\label{VI.149}
\begin{array}{l}
\ds\limsup_{k\rightarrow +\infty}\int_{D^2}|\nabla\vec{\Phi}_k|^2=\limsup_{k\rightarrow +\infty}2\,\int_{D^2}e^{2\la_k}\ \\[5mm]
\ds\quad\quad=\limsup_{k\rightarrow +\infty}2\,\mbox{Area}(\vec{\Phi}_k(D^2))<+\infty\quad.
 \end{array}
\ee
Combining (\ref{VI.147a}), (\ref{VI.148}) and (\ref{VI.149}) we obtain
\be
\label{VI.150}
\forall\,\rho<1\ \quad\limsup_{k\rightarrow +\infty}\|\vec{\Phi}_k\|_{W^{2,2}(B^2_\rho(0))}<+\infty
\ee
We can then extract a subsequence that we keep denoting $\vec{\Phi}_k$ such that there exists a map $\vec{\Phi}_\infty\in W^{2,2}_{loc}(D^2,{\R}^m)$ for which
\[
\vec{\Phi}_k\rightharpoonup\vec{\Phi}_\infty\quad\quad\mbox{ weakly in }W^{2,2}_{loc}(D^2,{\R}^m)\quad.
\]
Using Rellich Kondrachov theorem\footnote{See Chapter 6 of \cite{AdFo}.} we deduce that 
\[
\vec{\Phi}_k\longrightarrow\vec{\Phi}_\infty\quad\quad\mbox{ strongly in }W^{1,p}_{loc}(D^2,{\R}^m)\quad\forall\, p<+\infty\quad.
\]
Because of the strong convergence of the gradient of $\vec{\Phi}_k$ in $L^p$ for any $p<+\infty$, 
$$
\nabla\vec{\Phi}_k\quad\mbox{ converges almost everywhere towards }\quad\nabla\vec{\Phi}_\infty
$$ and then we can pass to the limit in the conformality conditions 
\[
\lf\{
\begin{array}{l}
\ds |\p_{x_1}\vec{\Phi}_k|=|\p_{x_2}\vec{\Phi}_k|=e^{\la_k}\quad,\\[5mm]
\ds<\p_{x_1}\vec{\Phi}_k,\p_{x_2}\vec{\Phi}_k>=0
\end{array}
\rg.
\]  
in order to deduce
\[
\lf\{
\begin{array}{l}
\ds |\p_{x_1}\vec{\Phi}_\infty|=|\p_{x_2}\vec{\Phi}_\infty|=e^{\la_\infty}\quad,\\[5mm]
\ds<\p_{x_1}\vec{\Phi}_\infty,\p_{x_2}\vec{\Phi}_\infty>=0
\end{array}
\rg.
\]  
The passage to the limit in the condition iii) gives
\be
\label{VI.151}
\la_\infty=\log|\p_{x_1}\vec{\Phi}_\infty|=\log|\p_{x_1}\vec{\Phi}_\infty|\in\ L^\infty_{loc}(D^2)\quad.
\ee
Hence $\vec{\Phi}_\infty$ realizes a conformal lipschitz immersion of the disc $D^2$. 

\medskip

Because of the pointwise convergence of $\nabla \vec{\Phi}_k$ towards $\nabla\vec{\Phi}_\infty$ we have that
\[
\p_{x_1}\vec{\Phi}_k\wedge\p_{x_2}\vec{\Phi}_k\longrightarrow \p_{x_1}\vec{\Phi}_\infty\wedge\p_{x_2}\vec{\Phi}_\infty\quad\quad\mbox{ almost everywhere,}
\]
 and, because of iii) $|\p_{x_1}\vec{\Phi}_k\wedge\p_{x_2}\vec{\Phi}_k|=e^{2\la_k}$  is bounded from below by a positive constant on each compact set included in the open disc $D^2$.
 Therefore
 \[
\vec{n}_{\vec{\Phi}_k}= e^{-2\la_k}\p_{x_1}\vec{\Phi}_k\wedge\p_{x_2}\vec{\Phi}_k\longrightarrow e^{-2\la_\infty}\p_{x_1}\vec{\Phi}_\infty\wedge\p_{x_2}\vec{\Phi}_\infty=\vec{n}_{\vec{\Phi}_\infty}\quad\mbox{ a. e.}
 \]
 The assumption i) implies that, modulo extraction of a subsequence, $\vec{n}_{\vec{\Phi}_k}$ converges weakly in $W^{1,2}(D^2, \wedge^{m-2}{\R}^m)$ to a limit $\vec{n}_\infty\in W^{1,2}(D^2, \wedge^{m-2}{\R}^m)$. Using again Rellich Kondrachov compactness result we know that this convergence is strong in $L^p(D^2)$ for any $p<+\infty$.
 The almost everywhere convergence of $\vec{n}_{\vec{\Phi}_k}$ towards $\vec{n}_{\vec{\Phi}_\infty}$ implies that
 \[
 \vec{n}_\infty=\vec{n}_{\vec{\Phi}_\infty}
 \]
 hence the limit is unique and the whole sequence $\vec{n}_{\vec{\Phi}_k}$ converges weakly in $W^{1,2}(D^2, \wedge^{m-2}{\R}^m)$ to $\vec{n}_{\vec{\Phi}_\infty}$. 
 From the lower semicontinuity of the $W^{1,2}$ norm, we deduce in particular that
 \be
 \label{VI.152}
\int_{D^2}|\nabla\vec{n}_{\vec{\Phi}_\infty}|^2\le 8\pi/3\quad.
 \ee
 Hence $\vec{\Phi}_\infty$ is a conformal Lipschitz immersion with $L^2-$bounded second fundamental form. It is natural to ask whether the equation iv)
 passes to the limit or in other words we are asking the following question :
 
 \medskip
 
\noindent {\it Does the weak limit $\vec{\Phi}_\infty$ define a weak Willmore immersion in the sense of definition~\ref{df-VI.17} ?}
 
 \medskip
 
Since $\nabla\vec{\Phi}_k$ converges strongly in $L^p_{loc}(D^2)$ to $\nabla\vec{\Phi}_\infty$ for any $p<+\infty$ and since $\inf_{p\in B^2_\rho(0)}|\nabla\vec{\Phi}_k|(p)$ is bounded away from zero uniformly in $k$, we have that
 $$
 \begin{array}{c}
 \ds2\ e^{-2\la_k}=|\nabla\vec{\Phi}_k|^{-2}\longrightarrow |\nabla\vec{\Phi}_\infty|^{-2}=2\ e^{-2\la_\infty}\\[5mm]
 \mbox{ stronlgy in }L^p_{loc}(D^2)\quad\quad \forall\,p<+\infty
 \end{array}
 $$
Since $\Delta\vec{\Phi}_k\rightharpoonup\Delta\vec{\Phi}_\infty$ weakly in $L^2_{loc}(D^2)$ we deduce that
\[
\ds\vec{H}_k=\frac{e^{-2\la_k}}{2}\Delta\vec{\Phi}_k\longrightarrow\frac{e^{-2\la_\infty}}{2}\Delta\vec{\Phi}_\infty\quad\mbox{ in }{\mathcal D}'(D^2)
\]
Since $|\vec{H}_k|^2\ e^{2\la_k}\le 2^{-1}\ |\nabla\vec{n}_k|^2$, because of the assumption i) $\vec{H}_k$ is uniformly bounded w.r.t. $k$ in $L^2(B^2_\rho(0))$ for any $\rho<1$. We can then deduce from the previous facts that 
\be
\label{VI.153}
\vec{H}_k\rightharpoonup \vec{H}_{\vec{\Phi}_\infty}\quad\mbox{ weakly in }L^2_{loc}(D^2)\quad.
\ee 
At this preliminary stage of our analysis of the passage to the limit inside the Willmore equation in conservative form
(\ref{VI.65aa}) it is not possible to identify the limits of the bilinearities such as
\[
\nabla^\perp\vec{n}_k\wedge\vec{H}_k\quad\longrightarrow\quad ?
\]
Indeed both $\nabla^\perp\vec{n}_k$ and $\vec{H}_k$ converge weakly in $L^2_{loc}$ but, because of these
weak convergences, one cannot a-priori identify the limit of the product $\nabla^\perp\vec{n}_k\wedge\vec{H}_k$
as being the product of the limit $\nabla^\perp\vec{n}_{\vec{\Phi}_\infty}\wedge\vec{H}_{\vec{\Phi}_\infty}$.
 
 \medskip
 
Before any more advanced study of the passage to the limit the best one can deduce at this stage is the existence of a locally Radon measure vector fields taking values in ${\R}^m$ : $\mu=\mu_1\,\p_{x_1}+\mu_2\,\p_{x_2}$ on $D^2$
such that\footnote{Each $\mu_1$ and $\mu_2$ are ${\R}^m$-valued Radon measures on $D^2$.}
\be
\label{VI.154}
div\lf[\nabla \vec{H}_\infty-3\pi_{\vec{n}_\infty}(\nabla\vec{H}_\infty)+\star(\nabla^\perp\vec{n}_\infty\wedge\vec{H}_\infty)\rg]=div\,\mu
\ee
where $\vec{H}_\infty$ and $\vec{n}_\infty$ stand for $\vec{H}_{\vec{\Phi}_\infty}$ and $\vec{n}_{\vec{\Phi}_\infty}$. In order to understand the possible limits of Willmore discs with small $L^2-$norm of the second fundamental form, it remains to identify the nature of $\mu$. This is the purpose of the following 2 subsections.

\subsubsection{Conservation laws for Willmore surfaces.}

As we saw in the first part of the course critical non-linear elliptic systems and equations do not always pass to the limit.
For the systems issued from conformally invariant Lagrangians we discovered conservation laws which were the key
objects for passing in the limit in these systems. We shall here also, for Willmore surfaces, find divergence free quantities
that will help us to describe the passage to the limit in Willmore equation.
\begin{Th}
\label{th-VI.19}
Let $\vec{\Phi}$ be a Lipschitz conformal immersion of the disc $D^2$ with $L^2-$bounded second fundamental form.
Assume $\vec{\Phi}$ is a weak Willmore immersion then there exists $\vec{L}\in L^{2,\infty}_{loc}(D^2)$ such that
\be
\label{VI.154a}
\nabla^\perp\vec{L}=\nabla\vec{H}-3\pi_{\vec{n}}(\nabla\vec{H})+\star(\nabla^\perp\vec{n}\wedge\vec{H})\quad,
\ee
where $\vec{n}$ and $\vec{H}$ denote respectively the Gauss map and the mean curvature vector associed to the immersion $\vec{\Phi}$.

Moreover the following conservation laws are satisfied
\be
\label{VI.155}
div<\vec{L},\nabla^\perp\vec{\Phi}>=0\quad,
\ee
and
\be
\label{VI.156}
div\lf[\vec{L}\wedge\nabla^\perp\vec{\Phi}+2\ (\star(\vec{n}\res\vec{H}))\res\nabla^\perp\vec{\Phi}\rg]=0\quad.
\ee
where $\star$ is the ususal Hodge operator on multivectors for the canonical scalar product in ${\R}^m$ and $\res$ is the operations between $p-$ and $q-$ vectors ($p\ge q$)  satisfying for any $\al\in\wedge^p{\R}^m$, $\beta\in\wedge^q{\R}^m$ and $\gamma\in\wedge^{p-q}{\R}^m$
\[
<\al\res\beta,\gamma>=<\al,\beta\wedge\gamma>\quad.
\]
\hfill $\Box$
\end{Th} 
{\bf Proof of theorem~\ref{th-VI.19}.}
Let 
\[
\vec{F}:=div\lf(\frac{1}{2\pi}\log\, r\star \chi\ \lf[\nabla^\perp\vec{H}-3\pi_{\vec{n}}(\nabla^\perp\vec{H})-\star(\nabla\vec{n}\wedge\vec{H})\rg]\rg)
\]
where $\chi$ is the characteristic function of the disc $D^2$. Under the assumptions of the theorem we have seen that 
\[
\lf[\nabla^\perp\vec{H}-3\pi_{\vec{n}}(\nabla^\perp\vec{H})-\star(\nabla\vec{n}\wedge\vec{H})\rg]\in L^1_{loc}\cap H^{-1}_{loc}(D^2)\quad
\]
Hence a classical result on Riesz potentials applies (see \cite{Ad}) in order to deduce that
\[
\vec{F}\in L^{2,\infty}_{loc}(D^2)\quad.
\]
Since $\vec{F}$ has been chosen in order to have
\be
\label{VI.157}
\Delta\vec{F}=div\lf[\nabla^\perp\vec{H}-3\pi_{\vec{n}}(\nabla^\perp\vec{H})-\star(\nabla\vec{n}\wedge\vec{H})\rg]\quad,
\ee
We introduce the distribution
\[
\vec{X}:=\nabla^\perp F+\nabla\vec{H}-3\pi_{\vec{n}}(\nabla\vec{H})+\star(\nabla^\perp\vec{n}\wedge\vec{H})\quad.
\]
Combining (\ref{VI.57}) and the fact that $\vec{\Phi}$ is Willmore, which is equivalent to (\ref{VI.65aa}), we obtain that $\vec{X}$ satisfies
\[
\lf\{
\begin{array}{l}
div\, \vec{X}=0\\[5mm]
curl\, \vec{X}=0
\end{array}
\rg.
\]
Hence the components of $\vec{X}=(X_1,\cdots,X_m)$ realize harmonic vectorfields and there exists then an ${\R}^m-$valued harmonic map $\vec{G}=(G_1,\cdots, G_m)$ 
such that
\[
\vec{X}=\nabla^\perp \vec{G}\quad.
\]
Being harmonic, the map $\vec{G}$ is analytic in the interior of $D^2$, therefore\footnote{In fact (\ref{VI.154a}) is telling us
that $\nabla\vec{L}$ belongs to $H^{-1}+L^1(D^2)$ and using a more sophisticated result (see \cite{BoBr} theorem 4) one can infer that in fact
$\vec{L}\in L^2(D^2)$.} $\vec{L}:=\vec{G}-\vec{F}$ is in $L^{2,\infty}_{loc}(D^2)$ and satisfies (\ref{VI.154a}).

\medskip

We now establish the first conservation law (\ref{VI.155}). We have 
\[
<\nabla\vec{\Phi},\nabla^\perp\vec{L}>=<\nabla\vec{\Phi},\nabla\vec{H}>+<\nabla\vec{\Phi},\star(\nabla^\perp\vec{n}\wedge\vec{H})>\quad.
\]
Multiplying (\ref{VI.73}) by $H_\al$, summing over $\al=1\cdots m-2$ and projecting over $\vec{\Phi}_\ast TD^2$ using the tangential projection $\pi_T$ gives
\be
\label{VI.158}
\pi_T(\nabla\vec{H}-\star(\nabla^\perp\vec{n}\wedge\vec{H}))=-2\ |\vec{H}|^2\ \nabla\vec{\Phi}\quad.
\ee
Hence we have
\[
<\nabla\vec{\Phi},\star(\nabla^\perp\vec{n}\wedge\vec{H})>=<\nabla\vec{\Phi},\nabla\vec{H}>+2\ |\vec{H}|^2\ |\nabla\vec{\Phi}|^2\quad.
\]
Thus
\be
\label{VI.159}
<\nabla\vec{\Phi},\nabla^\perp\vec{L}>=2<\nabla\vec{\Phi},\nabla\vec{H}>+2\ |\vec{H}|^2\ |\nabla\vec{\Phi}|^2\quad.
\ee
We have in one hand, since $<\vec{H},\nabla\vec{\Phi}>=0$,
\[
2<\nabla\vec{\Phi},\nabla\vec{H}>=-2<\Delta\vec{\Phi},\vec{H}>=-4\, e^{2\la}\ |\vec{H}|^2\quad,
\]
and in the other hand
\[
2\ |\vec{H}|^2\ |\nabla\vec{\Phi}|^2=4\ |\vec{H}|^2\ e^{2\la}\quad.
\]
Inserting these two last identities in (\ref{VI.159}) gives (\ref{VI.155}).

\medskip

Finally we establish the conservation law (\ref{VI.156}). Since $\nabla\vec{\Phi}\wedge\nabla\vec{\Phi}=0$, multiplying (\ref{VI.73}) by $H_\al$, summing over $\al=1\cdots m-2$
and wedging with $\nabla\vec{\Phi}$ gives
\be
\label{VI.160}
\begin{array}{l}
\ds\nabla\vec{\Phi}\wedge\star(\nabla^\perp\vec{n}\wedge\vec{H})=\nabla\vec{\Phi}\wedge\sum_{\al=1}^{m-2}H_\al\,\nabla\vec{n}_\al\\[5mm]
 \ds\quad\quad-\nabla\vec{\Phi}\wedge\sum_{\al,\beta=1}^{m-2}\vec{n}_\beta\ <\nabla\vec{n}_\al,\vec{n}_\beta>\ H_\al
 \end{array}
 \ee
We observe that 
\[
\sum_{\al,\beta=1}^{m-2}\vec{n}_\beta\ <\nabla\vec{n}_\al,\vec{n}_\beta>\ H_\al=\pi_{\vec{n}}(\nabla\vec{H})-\sum_{\al=1}^{m-2}\nabla H_\al\ \vec{n}_\al\quad.
\]
Inserting this identity in (\ref{VI.160}) gives
\be
\label{VI.161}
\begin{array}{rl}
\ds\nabla\vec{\Phi}\wedge\star(\nabla^\perp\vec{n}\wedge\vec{H})&=\nabla\vec{\Phi}\wedge\nabla\vec{H}-\nabla\vec{\Phi}\wedge\pi_{\vec{n}}(\nabla\vec{H})\\[5mm]
 &\ds=\nabla\vec{\Phi}\wedge\pi_T(\nabla\vec{H})\quad.
 \end{array}
\ee
We have
\be
\label{VI.162}
\begin{array}{l}
\nabla\vec{\Phi}\wedge\pi_T(\nabla\vec{H})=e^{\la}\ \lf[<\vec{e_2},\p_{x_1}\vec{H}>-<\vec{e_1},\p_{x_2}\vec{H}>\rg]\ \vec{e}_1\wedge\vec{e}_2\\[5mm]
 \ds\quad\quad=e^\la\ \lf[<\pi_{\vec{n}}(\p_{x_2}\vec{e}_1-\p_{x_1}\vec{e}_2),\vec{H}>\rg]\ \vec{e}_1\wedge\vec{e}_2\\[5mm]
  \ds\quad\quad=e^{2\la}\ \lf[<\vec{\mathbb I}(\vec{e}_1,\vec{e}_2)-\vec{\mathbb I}(\vec{e}_2,\vec{e}_1),\vec{H}>\rg]\ \vec{e}_1\wedge\vec{e}_2\\[5mm]
 \quad\quad=0
\end{array}
\ee
Therefore combining this identity with (\ref{VI.161}) gives
\be
\label{VI.163}
\nabla\vec{\Phi}\wedge\star(\nabla^\perp\vec{n}\wedge\vec{H})=0
\ee
We deduce from this equality that
\[
\nabla\vec{\Phi}\wedge\nabla^\perp\vec{L}=\nabla\vec{\Phi}\wedge\nabla\vec{H}-3\nabla\vec{\Phi}\wedge\pi_{\vec{n}}(\nabla\vec{H})
\]
Combining this fact with (\ref{VI.162}) gives
\be
\label{VI.164}
\nabla\vec{\Phi}\wedge\nabla^\perp\vec{L}=-2\nabla\vec{\Phi}\wedge\nabla\vec{H}
\ee
It remains now to express $\nabla\vec{\Phi}\wedge\nabla\vec{H}$ in terms of a linear combination of jacobians in order to be able
to ''factorize'' the divergence operator in (\ref{VI.164}). 

\medskip

The definition of the contraction operation $\res$ gives
\[
\vec{n}\res\vec{H}=\sum_{\al=1}^{m-2}(-1)^{\al-1}\, H_\al\, \wedge_{\beta\ne\al}\vec{n}_\beta\quad.
\]
Hence we have
\be
\label{VI.165}
\star(\vec{n}\res\vec{H})=\vec{e}_1\wedge\vec{e}_2\wedge\vec{H}\quad.
\ee
We shall now compute $\nabla(\star(\vec{n}\res\vec{H}))\res\nabla^\perp\vec{\Phi}$. 

\medskip

To that purpose we first compute
\be
\label{VI.165a}
\begin{array}{l}
\ds \nabla(\vec{e}_1\wedge\vec{e}_2\wedge\vec{H})=\pi_{\vec{n}}(\nabla\vec{e}_1)\wedge\vec{e}_2\wedge\vec{H}+\vec{e}_1\wedge\pi_{\vec{n}}(\nabla\vec{e}_2)\wedge\vec{H}\\[5mm]
\ds \quad\quad+\vec{e}_1\wedge\vec{e}_2\wedge\nabla\vec{H}\quad.
\end{array}
\ee
Using elementary rule\footnote{The definition of $\res$ implies that for any choice of 4 vectors $\vec{a}$, $\vec{b}$, $\vec{c}$ and $\vec{d}$ one has
\[
\begin{array}{l}
\ds(\vec{a}\wedge\vec{b}\wedge\vec{c})\res\vec{d}=<\vec{a},\vec{d}>\ \vec{b}\wedge\vec{c}\ -\ <\vec{b},\vec{d}>\ \vec{a}\wedge\vec{c}\ +\ <\vec{c},\vec{d}>\ \vec{a}\wedge\vec{b}\quad.
\end{array}
\]
    } on the contraction operation $\res$ we compute the following : 
first we have
\be
\label{VI.166}
(\pi_{\vec{n}}(\nabla\vec{e}_1)\wedge\vec{e}_2\wedge\vec{H})\res\nabla^\perp\vec{\Phi}=e^\la\ \pi_{\vec{n}}(\p_{x_1}\vec{e}_1)\wedge\vec{H}\quad,
\ee
then we have
\be
\label{VI.167}
(\vec{e}_1\wedge\pi_{\vec{n}}(\nabla\vec{e}_2)\wedge\vec{H})\res\nabla^\perp\vec{\Phi}=e^\la\ \pi_{\vec{n}}(\p_{x_2}\vec{e}_2)\wedge\vec{H}\quad,
\ee
and finally we have
\be
\label{VI.168}
\begin{array}{l}
\ds(\vec{e}_1\wedge\vec{e}_2\wedge\nabla\vec{H})\res\nabla^\perp\vec{\Phi}=e^\la\ [\vec{e}_1\wedge\p_{x_1}\vec{H}+\vec{e}_2\wedge\p_{x_2}\vec{H}]\\[5mm]
\ds\quad\quad+<\nabla\vec{H},\nabla^\perp\vec{\Phi}>\ \vec{e}_1\wedge\vec{e}_2
\end{array}
\ee
Observe that,  since $<\vec{H},\nabla^\perp\vec{\Phi}>=0$,
$$
<\nabla\vec{H},\nabla^\perp\vec{\Phi}>=div<\vec{H},\nabla^\perp\vec{\Phi}>=0\quad.
$$ 
Hence (\ref{VI.168}) becomes
\be
\label{VI.169}
\begin{array}{l}
\ds(\vec{e}_1\wedge\vec{e}_2\wedge\nabla\vec{H})\res\nabla^\perp\vec{\Phi}=e^\la\ [\vec{e}_1\wedge\p_{x_1}\vec{H}+\vec{e}_2\wedge\p_{x_2}\vec{H}]\\[5mm]
\ds\quad\quad\quad=\nabla\vec{\Phi}\wedge\nabla\vec{H}
\end{array}
\ee
The combination of (\ref{VI.165}) with (\ref{VI.165a}), (\ref{VI.166}), (\ref{VI.167}) and (\ref{VI.169}) gives
\[
\begin{array}{l}
\ds\nabla(\star(\vec{n}\res\vec{H}))\res\nabla^\perp\vec{\Phi}=e^\la\ [\pi_{\vec{n}}(\p_{x_1}\vec{e}_1)+\pi_{\vec{n}}(\p_{x_2}\vec{e}_2)]\wedge\vec{H}\\[5mm]
\ds\quad\quad\quad+\nabla\vec{\Phi}\wedge\nabla\vec{H}\quad.
\end{array}
\]
Observe that $\pi_{\vec{n}}(\p_{x_1}\vec{e}_1)+\pi_{\vec{n}}(\p_{x_2}\vec{e}_2)=2\ e^\la\ \vec{H}$ hence finally we obtain
\be
\label{VI.170}
\nabla(\star(\vec{n}\res\vec{H}))\res\nabla^\perp\vec{\Phi}=\nabla\vec{\Phi}\wedge\nabla\vec{H}\quad.
\ee
Combining (\ref{VI.164}) and (\ref{VI.170}) gives
\be
\label{VI.171}
\nabla\vec{\Phi}\wedge\nabla^\perp{\vec{L}}=-2\,\nabla(\star(\vec{n}\res\vec{H}))\res\nabla^\perp\vec{\Phi}
\ee
This is exactly the conservation law (\ref{VI.156}) and theorem~\ref{th-VI.19} is proved. \hfill$\Box$

\medskip

Having now found 2 new conserved quantities, as we did for the first one $\nabla\vec{H}-3\pi_{\vec{n}}(\nabla\vec{H})+\star(\nabla^\perp\vec{n}\wedge\vec{H})$, we can apply Poincar\'e 
lemma in order to obtain ''primitives'' of these quantities. These ''primitive'' quantities will satisfy a very particular elliptic system.
Precisely we have the following theorem.

\begin{Th}
\label{th-VI.20}
Let $\vec{\Phi}$ be a conformal lipschtiz immersion of the disc $D^2$ with $L^2-$bounded second fundamental form. Assume there exists $\vec{L}\in L^{2,\infty}(D^2,{\R}^m)$
satisfying
\be
\label{VI.172}
\lf\{
\begin{array}{l}
\ds div<\vec{L},\nabla^\perp\vec{\Phi}>=0\\[5mm]
\ds div\lf[\vec{L}\wedge\nabla^\perp\vec{\Phi}+2\ (\star(\vec{n}\res\vec{H}))\res\nabla^\perp\vec{\Phi}\rg]=0\quad.
\end{array}
\rg.
\ee
where $\vec{H}$ and $\vec{n}$ denote respectively the mean-curvature vector and the Gauss map of the immersion $\vec{\Phi}$.
There exists\footnote{We denote by $W^{1,(2,\infty)}$ the space of distribution in $L^2$ with gradient in $L^{2,\infty}$.} $S\in W^{1,(2,\infty)}_{loc}(D^2,{\R})$ and $\vec{R}\in W^{1,(2,\infty)}_{loc}(D^2,\wedge^2{\R}^m)$  such that 
\be
\label{VI.172a}
\lf\{
\begin{array}{l}
\ds \nabla S=<\vec{L},\nabla\vec{\Phi}>\\[5mm]
\ds\nabla\vec{R}=\vec{L}\wedge\nabla\vec{\Phi}+2\ (\star(\vec{n}\res\vec{H}))\res\nabla\vec{\Phi}\quad.
\end{array}
\rg.
\ee
and the following equation holds
\be
\label{VI.173}
\lf\{
\begin{array}{rcl}
\ds\nabla S&=&\ds-<\star\vec{n},\nabla^\perp\vec{R}>\\[5mm]
\ds\nabla\vec{R}&=&\ds(-1)^m\ \star(\vec{n}\bullet\nabla^\perp\vec{R})+ (-1)^{m-1}\nabla^\perp S\ \star\vec{n}\quad,
\end{array}
\rg.
\ee
where $\bullet$ is the following contraction operation which to a pair of respectively $p-$ and $q-$vectors of ${\R}^m$ assigns
a $p+q-2-$vector of ${\R}^m$ such that
\[
\forall\,\vec{a}\in\wedge^p{\R}^m\quad\ \forall\,\vec{b}\in\wedge^1{\R}^m\quad\vec{a}\bullet\vec{b}:=\vec{a}\res\vec{b}
\]
and
\[
\begin{array}{l}
\ds\forall\,\vec{a}\in\wedge^p{\R}^m\quad\ \forall\,\vec{b}\in\wedge^r{\R}^m\quad\ \forall\,\vec{c}\in\wedge^s{\R}^m\quad\\[5mm]
\ds\vec{a}\bullet(\vec{b}\wedge\vec{c}):=(\vec{a}\bullet\vec{b})\wedge\vec{c}+(-1)^{r\,s}(\vec{a}\bullet\vec{c})\wedge\vec{b}
\end{array}
\]
\hfill$\Box$
\end{Th}
\begin{Rm}
\label{rm-VI.20}
In the particular case $m=3$ both $\vec{n}$ and $\vec{R}$ can be identified with vectors by the mean of the Hodge operator $\star$. Once this identification is made the systems (\ref{VI.172a}) and (\ref{VI.173}) become respectively
\be
\label{VI.172x}
\lf\{
\begin{array}{l}
\ds \nabla S=<\vec{L},\nabla\vec{\Phi}>\\[5mm]
\ds\nabla\vec{R}=\vec{L}\times\nabla\vec{\Phi}+2\ H\ \nabla\vec{\Phi}\quad.
\end{array}
\rg.
\ee
and
\be
\label{VI.173x}
\lf\{
\begin{array}{rcl}
\ds\nabla S&=&\ds-<\vec{n},\nabla^\perp\vec{R}>\\[5mm]
\ds\nabla\vec{R}&=&\ds\vec{n}\times\nabla^\perp\vec{R}+\nabla^\perp S\ \vec{n}\quad,
\end{array}
\rg.
\ee
\end{Rm}
{\bf Proof of theorem~\ref{th-VI.20}.}
The existence of $S$ and $\vec{R}$ satisfying (\ref{VI.172a}) is obtained exactly like for $\vec{L}$ in the beginning of the
proof of theorem~\ref{th-VI.19} taking successively the convolution of the divergence free quantities with $-(2\pi)^{-1}\log\ r$, then taking the 
$curl$ operator and finally subtracting some harmonic ${\R}$ (for $S$) or $\wedge^2{\R}$ (for $\vec{R}$) valued map.

\medskip

It remains to prove (\ref{VI.173}). Let $\vec{N}$ be a normal vector, exactly like for the particular case of $\vec{H}$ in (\ref{VI.165})
\be
\label{VI.174}
(-1)^{m-1}\,\star(\vec{n}\res\vec{N})=\vec{e}_1\wedge\vec{e}_2\wedge\vec{N}\quad.
\ee
We deduce from this identity
\be
\label{VI.175}
(-1)^{m-1}\,(\star(\vec{n}\res\vec{N}))\res\nabla\vec{\Phi}=(\vec{e}_1\wedge\vec{e}_2\wedge\vec{N})\res\nabla\vec{\Phi}=-\nabla^\perp\vec{\Phi}\wedge\vec{N}\quad.
\ee
Applying this identity to $\vec{N}:=\vec{H}$ implies
\be
\label{VI.175a}
\nabla\vec{R}=\vec{L}\wedge\nabla\vec{\Phi}-2\,(-1)^{m-1}\,\nabla^\perp\vec{\Phi}\wedge\vec{H} 
\ee
We take now the $\bullet$ contraction between $\vec{n}$ and $\nabla\vec{R}$ and we obtain
\be
\label{VI.176}
\begin{array}{rl}
\ds\vec{n}\bullet\nabla\vec{R}&\ds=(\vec{n}\res\vec{L})\wedge\nabla\vec{\Phi}+2\,(-1)^{m-1}\,(\vec{n}\res\vec{H})\wedge\nabla^\perp\vec{\Phi}\\[5mm]
 &\ds=(\vec{n}\res\pi_{\vec{n}}(\vec{L}))\wedge\nabla\vec{\Phi}+2\,(-1)^{m-1}\,(\vec{n}\res\vec{H})\wedge\nabla^\perp\vec{\Phi}\quad.
 \end{array}
\ee
For a normal vector $\vec{N}$ again a short computation gives
\be
\label{VI.177}
\star[(\vec{n}\res\vec{N})\wedge\nabla\vec{\Phi}]=(-1)^m\ \nabla^\perp\vec{\Phi}\wedge\vec{N}\quad,
\ee
from which we also deduce
\be
\label{VI.178}
\star[(\vec{n}\res\vec{N})\wedge\nabla^\perp\vec{\Phi}]=(-1)^{m-1}\ \nabla\vec{\Phi}\wedge\vec{N}\quad.
\ee
Combining (\ref{VI.176}), (\ref{VI.177}) and (\ref{VI.178}) gives then
\be
\label{VI.178a}
\star(\vec{n}\bullet\nabla\vec{R})=-\nabla^\perp\vec{\Phi}\wedge\pi_{\vec{n}}(\vec{L})+2\nabla\vec{\Phi}\wedge\vec{H}\quad.
\ee
from which we deduce
\be
\label{VI.179}
\star(\vec{n}\bullet\nabla^\perp\vec{R})=-\pi_{\vec{n}}(\vec{L})\wedge\nabla\vec{\Phi}+2\nabla^\perp\vec{\Phi}\wedge\vec{H}\quad.
\ee
Combining (\ref{VI.175}) and (\ref{VI.179}) gives 
\be
\label{VI.180}
(-1)^m\,\star(\vec{n}\bullet\nabla^\perp\vec{R})=\nabla\vec{R}+(-1)^m\,\pi_T(\vec{L})\wedge\nabla\vec{\Phi}\quad.
\ee
One verifies easily that
\be
\label{VI.181}
\begin{array}{rl}
\ds\pi_{T}(\vec{L})\wedge\nabla\vec{\Phi}&\ds=\nabla^\perp S\ \vec{e}_1\wedge\vec{e}_2\\[5mm]
 &\ds=\nabla^\perp S\ \star\vec{n}
 \end{array}
\ee 
The combination of (\ref{VI.180}) and (\ref{VI.181}) gives the second equation of(\ref{VI.173}). The first equation is obtained
by taking the scalar product between the first equation and $\star\vec{n}$ once one has observed that
\[
<\star\vec{n},  \star(\vec{n}\bullet\nabla^\perp\vec{R})>=0
\]
This later fact comes from (\ref{VI.178a}) which implies that $\star(\vec{n}\bullet\nabla^\perp\vec{R})$ is a linear combination
of wedges of tangent and normal vectors to $\vec{\Phi}_\ast TD^2$. Hence theorem~\ref{th-VI.20} is proved.\hfill $\Box$

\medskip

An important corollary of the previous theorem is the following.

\begin{Co}
\label{co-VI.21}
Let $\vec{\Phi}$ be a conformal lipschitz immersion of the disc $D^2$ with $L^2-$bounded second fundamental form. Assume there exists 
$\vec{L}$ in $L^{2,\infty}(D^2,{\R}^m)$ satisfying
\be
\label{VI.172bis}
\lf\{
\begin{array}{l}
\ds div<\vec{L},\nabla^\perp\vec{\Phi}>=0\\[5mm]
\ds div\lf[\vec{L}\wedge\nabla^\perp\vec{\Phi}+2\ (\star(\vec{n}\res\vec{H}))\res\nabla^\perp\vec{\Phi}\rg]=0\quad.
\end{array}
\rg.
\ee
where $\vec{H}$ and $\vec{n}$ denote respectively the mean-curvature vector and the Gauss map of the immersion $\vec{\Phi}$.
Let $S\in W^{1,(2,\infty)}_{loc}(D^2,{\R})$ and $\vec{R}\in W^{1,(2,\infty)}_{loc}(D^2,\wedge^2{\R}^m)$  such that 
\be
\label{VI.172abis}
\lf\{
\begin{array}{l}
\ds \nabla S=<\vec{L},\nabla\vec{\Phi}>\\[5mm]
\ds\nabla\vec{R}=\vec{L}\wedge\nabla\vec{\Phi}+2\ (\star(\vec{n}\res\vec{H}))\res\nabla\vec{\Phi}\quad.
\end{array}
\rg.
\ee
Then $(\vec{\Phi},S,\vec{R})$ satisfy the following system\footnote{In codimension 1 the systems reads
\be
\label{VI.182aa}
\lf\{
\begin{array}{rcl}
\ds\Delta S&=& -<\nabla\vec{n},\nabla^\perp\vec{R}>\\[5mm]
\ds\Delta\vec{R}&=&\nabla\vec{n}\times\nabla^\perp\vec{R}+\nabla^\perp S\,\nabla\vec{n}\\[5mm]
\ds\Delta\vec{\Phi}=\nabla^\perp S\,\nabla\vec{\Phi}+\nabla^\perp\vec{R}\times\nabla\vec{\Phi}
\end{array}
\rg.
\ee
Observe that the use of two different operations, $\bullet$
in arbitrary codimension where $\vec{R}$ is seen as a 2-vector and $\times$ in 3 dimension $\vec{R}$ is interpreted as a vector, generates formally different signs.
This might look first a bit confusing for the reader but we believed that the codimension 1 case which  more used in applications, deserved to be singled out.
}

\be
\label{VI.182}
\lf\{
\begin{array}{rcl}
\ds\Delta S&=&\ds-<\star\nabla\vec{n},\nabla^\perp\vec{R}>\\[5mm]
\ds\Delta\vec{R}&=&\ds(-1)^m\ \star(\nabla\vec{n}\bullet\nabla^\perp\vec{R})+\nabla^\perp S\ \star\nabla\vec{n}\\[5mm]
\ds\Delta\vec{\Phi}&=&2^{-1}\nabla^\perp S\cdot\nabla\vec{\Phi}-2^{-1}\nabla\vec{R}\res\nabla^\perp\vec{\Phi}
\end{array}
\rg.
\ee 
\end{Co}
\begin{Rm}
The spectacular fact in (\ref{VI.182}) is that we have deduced, from the Willmore equation, a system with quadratic non-linearities
which are made of linear combinations of jacobians. We shall exploit intensively this fact below in order to get the regularity of weak Willmore
immersions and pass to the limit in the system.
\end{Rm}
\noindent{\bf Proof of corollary~\ref{co-VI.21}.}
The two first identities of (\ref{VI.182}) are obtained by taking the divergence of (\ref{VI.173}).

\medskip

Using the definition of the operation $\bullet$ and the second line of (\ref{VI.172abis}), we have
\be
\label{VI.183}
\begin{array}{rl}
\ds\nabla^\perp\vec{\Phi}\bullet\nabla\vec{R}&\ds=<\nabla^\perp\vec{\Phi},\vec{L}>\cdot\nabla\vec{\Phi}-<\nabla^\perp\vec{\Phi},\nabla\vec{\Phi}>\ \vec{L}\\[5mm]
 &+2\nabla^\perp\vec{\Phi}\bullet\lf[(\star(\vec{n}\res\vec{H}))\res\nabla\vec{\Phi}\rg]
 \end{array}
\ee 
It is clear that
\be
\label{VI.183a}
<\nabla^\perp\vec{\Phi},\nabla\vec{\Phi}>=0
\ee
From (\ref{VI.165}) we compute
\be
\label{VI.184}
\star(\vec{n}\res\vec{H})\res\nabla\vec{\Phi}=(\vec{e}_1\wedge\vec{e}_2\wedge\vec{H})\res\nabla\vec{\Phi}=-\nabla^\perp\vec{\Phi}\wedge\vec{H}\quad.
\ee
Combining this identity with the definition of $\bullet$ we obtain
\be
\label{VI.185}
\nabla^\perp\vec{\Phi}\bullet\lf[(\star(\vec{n}\res\vec{H}))\res\nabla\vec{\Phi}\rg]=-2\ e^{2\la}\vec{H}=-\Delta\vec{\Phi}
\ee
Since $<\nabla^\perp\vec{\Phi},\vec{L}>=\nabla^\perp S$, combining (\ref{VI.183}), (\ref{VI.183a}) and (\ref{VI.185}) gives
\be
\label{VI.186}
\nabla\vec{R}\res\nabla^\perp\vec{\Phi}=\nabla^\perp\vec{\Phi}\bullet\nabla\vec{R}=\nabla^\perp S\cdot\nabla\vec{\Phi}-2\Delta\vec{\Phi}
\ee
which gives the last line of (\ref{VI.182}) and corollary~\ref{co-VI.21} is proved. \hfill$\Box$

\medskip

We shall now study a first consequence of the conservation laws (\ref{VI.172bis}) : the regularity of weak Willmore surfaces

\medskip

\subsubsection{The regularity of weak Willmore immersions.}

In the present subsection we prove that weak Willmore immersions are $C^\infty$ in conformal parametrization. This will be the consequence of the following theorem.

\begin{Th}
\label{th-VI.22}
Let $\vec{\Phi}$ be a conformal lipschitz immersion of the disc $D^2$ with $L^2-$bounded second fundamental form. Assume there exists 
$\vec{L}$ in $L^{2,\infty}(D^2,{\R}^m)$ satisfying
\be
\label{VI.187}
\lf\{
\begin{array}{l}
\ds div<\vec{L},\nabla^\perp\vec{\Phi}>=0\\[5mm]
\ds div\lf[\vec{L}\wedge\nabla^\perp\vec{\Phi}+2\ (\star(\vec{n}\res\vec{H}))\res\nabla^\perp\vec{\Phi}\rg]=0\quad.
\end{array}
\rg.
\ee
Then $\vec{\Phi}$ is $C^\infty$.
\hfill$\Box$.
\end{Th}
The combination of theorem~\ref{th-VI.19} and theorem~\ref{th-VI.22} gives immediately the following result
\begin{Co}
\label{co-VI.23}
Let $\vec{\Phi}$ be a lipschitz immersion of a smooth surface $\Sigma^2$ with $L^2-$bounded second fundamental form. Assume
moreover that $\vec{\Phi}$ is weak Willmore in the sense of definition~\ref{df-VI.17}.
Then $\vec{\Phi}$ is $C^\infty$ in conformal parametrization. \hfill $\Box$
\end{Co}
In order to prove theorem~\ref{th-VI.22} we will need the following consequence of Coifman Lions Meyer and Semmes result, theorem~\ref{th-III.3},
which is due to F.Bethuel (see \cite{Bet1}).
\begin{Th}
\label{th-VI.24}\cite{Bet1}
Let $a$ be a function such that $\nabla a\in L^{2,\infty}(D^2)$ and let $b$ be a function in $W^{1,2}(D^2)$. Let $\phi$ be the unique solution in $\cap_{p<2}W^{1,p}_0(D^2)$
of the equation\footnote{The Jacobian $\p_xa\,\p_yb-\p_ya\,\p_x b$ has to be understood in the weak sense
\[
 \p_xa\,\p_yb-\p_ya\,\p_x b:=div\lf[a\,\nabla^\perp b\rg]\quad.
 \]
 Since $\nabla a\in L^{2,\infty}(D^2)$ we have that $a\in L^q(D^2)$ for all $q<+\infty$ and hence $a\nabla^\perp b\in L^p(D^2)$ for all $p<2$.}
\be
\label{VI.188}
\lf\{
\begin{array}{l}
-\Delta\phi=\p_xa\,\p_yb-\p_ya\,\p_x b\quad\quad\quad\mbox{ in }D^2\\[5mm]
\phi=0\quad\quad\quad\mbox{ on }\p D^2\quad.
\end{array}
\rg.
\ee
Then $\phi$ lies in $W^{1,2}(D^2)$ and 
\be
\label{VI.189}
\|\nabla\phi\|_{L^2(D^2)}\le C\ \|\nabla a\|_{L^{2,\infty}(D^2)}\ \|\nabla b\|_{L^2(D^2)}\quad,
\ee
where $C>0$ is a constant independent of $a$ and $b$.\hfill $\Box$
\end{Th}
{\bf Proof of theorem~\ref{th-VI.24}.}
We assume first that $b$ is smooth and, once the estimate (\ref{VI.189}) will be proved 
we can conclude by a density argument. For a smooth $b$, classical elliptic theory tells us that $\nabla\phi$ is in $L^q(D^2)$ for any $q<+\infty$
and one has in particular
\[
\|\nabla\phi\|_{L^2(D^2)}=\sup_{\|X\|_{L^2}\le 1}\int_{D^2}X\cdot\nabla\phi\quad.
\]
For any such a vector field $X$ satisfying $\|X\|_{L^2}\le 1$ there exists a unique Hodge decomposition\footnote{$c$ is the minimizer of the following convex
problem
\[
\min_{c\in W^{1,2}_0(D^2)}\int_{D^2}|X-\nabla c|^2\quad.
\]}
 \[
X:=\nabla c+\nabla^\perp d\quad,
\]
where $c\in W^{1,2}_0$. Moreover, one easily verifies that
\be
\label{VI.190}
1=\|X\|^2_{L^2(D^2)}=\|\nabla c\|_{L^2(D^2)}^2+\|\nabla d\|_{L^2(D^2)}^2\quad.
\ee
Indeed one has
\[
\int_{D^2}\nabla^\perp d\cdot\nabla c=\int_{\p D^2}\frac{\p d}{\p\tau}\, c-\int_{D^2}div\lf[\nabla^\perp d\rg]\ c=0
\]
Similarly, replacing $c$ by $\phi$ the same argument gives
\be
\label{VI.191}
\int_{D^2}X\cdot\nabla\phi=\int_{D^2}\nabla c\cdot\nabla\phi\quad.
\ee
Using again the fact that $c=0$ on $\p D^2$, we have
\be
\label{VI.192}
\begin{array}{l}
\ds\int_{D^2}\nabla c\cdot\nabla\phi=-\int_{D^2}c\, \Delta\phi=-\int_{D^2}c\ div[a\ \nabla^\perp b]\\[5mm]
\ds \quad\quad=\int_{D^2}a\ \nabla^\perp b\cdot\nabla c\quad.
\end{array}
\ee
Let $\psi$ be the solution of 
\be
\label{VI.193}
\lf\{
\begin{array}{l}
\ds-\Delta\psi=\p_xb\,\p_yc-\p_yb\,\p_x c\quad\quad\quad\mbox{ in }D^2\\[5mm]
\ds\frac{\p\psi}{\p\nu}=0\quad\quad\quad\mbox{ on }\p D^2\quad.
\end{array}
\rg.
\ee
Using the Neuman version of theorem~\ref{th-III.3} together with the embedding of $W^{1,1}(D^2)$ into $L^{2,1}(D^2)$ we obtain the existence of a constant $C$ independent
of $b$ and $c$ such that
\be
\label{VI.194}
\|\nabla\psi\|_{L^{2,1}(D^2)}\le C\ \|\nabla b\|_{L^2(D^2)}\ \|\nabla c\|_{L^2(D^2)}\le C\ \|\nabla b\|_{L^2(D^2)}\quad.
\ee
The identity (\ref{VI.192}) becomes
\be
\label{VI.195}
\ds\int_{D^2}\nabla c\cdot\nabla\phi=\int_{D^2}a\ \Delta\psi=-\int_{D^2}\nabla a\cdot\nabla\psi\quad.
\ee 
Combining (\ref{VI.191}) with (\ref{VI.194}) and (\ref{VI.195}) gives
\be
\label{VI.196}
\|\nabla\phi\|_{L^2(D^2)}=\sup_{\|X\|_{L^2}\le 1}\int_{D^2}X\cdot\nabla\phi\le C\, \|\nabla a\|_{L^{2,\infty}}\ \|\nabla b\|_{L^2}\quad,
\ee
which is the identity (\ref{VI.189}) and the proof of theorem~\ref{th-VI.24} is complete.\hfill $\Box$

\medskip

It remains to prove theorem~\ref{th-VI.22} and the present subsection will be complete.

\medskip

\noindent{\bf Proof of theorem~\ref{th-VI.22}.}
Combining corollary~\ref{co-VI.21} and theorem~\ref{th-VI.22} gives in a straightforward way that $\nabla S$ and $\nabla \vec{R}$ given respectively
by the first and the second line of (\ref{VI.172x}) and satisfying the elliptic system (\ref{VI.182}) are in $L^2_{loc}(D^2)$.

Argueing exactly like in the proof of the regularity of solutions to CMC surfaces we obtain the existence of $\al>0$ such that
\[
\sup_{x_0\in B^2_{1/2}(0)\ ;\ r<1/4}r^{-\al}\ \int_{B^2_r(x_0)}|\Delta S|+|\Delta \vec{R}|<+\infty
\]
from which we deduce the existence of $p>2$ such that $\nabla S$ and $\nabla\vec{R}$ are in $L^p_{loc}(D^2)$. Bootstrapping this information 
again in the two first lines of (\ref{VI.182}), using lemma~\ref{lm-IV.ax}, gives that $\nabla S$ and $\nabla\vec{R}$ are in $L^p_{loc}(D^2)$ for any $p<+\infty$. Bootstarpping
the latest 
in the third equation of (\ref{VI.182}) gives that $\vec{\Phi}\in W^{2,p}_{loc}(D^2,{\R}^m)$ for any $p<+\infty$, from which we can deduce that $\vec{n}\in W^{1,p}_{loc}(D^2,Gr_{m-2}({\R}^m))$ for any $p<+\infty$, information that we inject back in the two first lines of the system (\ref{VI.182})...etc and one obtains after iterating again and again that $\vec{\Phi}$
is in $W^{k,p}_{loc}(D^2,{\R}^m)$ for any $k\in {\N}$ and $1\le p\le +\infty$. This gives that $\vec{\Phi}$ is in $C^\infty_{loc}(D^2)$ and theorem~\ref{th-VI.22}
is proved.\hfill $\Box$

\subsubsection{The conformal Willmore surface equation.}

As we have seen conformal immersions which are solutions to the conservation laws (\ref{VI.187})  for some $\vec{L}$ in $L^{2,\infty}(D^2)$ satisfy an elliptic system, the system (\ref{VI.182}) which
formally resembles to the CMC equation from which we deduced the smoothness of the immersion. Then comes naturally the question 

\centerline{\it Are solutions to the conservation laws (\ref{VI.187}) Willmore ?}

The answer to that question is {''\it almost''}. We shall see below that solutions to (\ref{VI.187})  for some $\vec{L}$ in $L^{2,\infty}(D^2)$ satisfy the
Willmore equation up to the addition of an holomorphic function times the Weingarten operator. Assuming that for some immersion $\vec{\Phi}$ of a given abstract surface $\Sigma^2$ (\ref{VI.187}) holds in any conformal chart,  then the union of these holomorphic functions can be ''glued together'' in order to produce an holomorphic quadratic differential of the riemann surface $\Sigma^2$ equipped with the conformal class given by $\vec{\Phi}^\ast g_{{\R}^m}$.
As we will explain below,the intrinsic equation obtained corresponds to the equation satisfied by critical points to the Willmore functional but with the constraint that the immersion realizes a  fixed conformal class. As we will see this holomorphic quadratic differential plays the role of a Lagrange multiplier.

\medskip

First we prove the following result.

\begin{Th}
\label{th-VI.25}
Let $\vec{\Phi}$ be a conformal lipschitz immersion of the disc $D^2$ with $L^2-$bounded second fundamental form. There exists 
$\vec{L}$ in $L^{2,\infty}(D^2,{\R}^m)$ satisfying
\be
\label{VI.197}
\lf\{
\begin{array}{l}
\ds div<\vec{L},\nabla^\perp\vec{\Phi}>=0\\[5mm]
\ds div\lf[\vec{L}\wedge\nabla^\perp\vec{\Phi}+2\ (\star(\vec{n}\res\vec{H}))\res\nabla^\perp\vec{\Phi}\rg]=0\quad.
\end{array}
\rg.
\ee
if and only if there exists an holomorphic function $f(z)$ such that\footnote{The operation $\Im$ assigns to a complex number
its imaginary part} 
\be
\label{VI.198}
div\lf[\nabla\vec{H}-3\pi_{\vec{n}}(\nabla\vec{H})+\star(\nabla^\perp\vec{n}\wedge\vec{H})\rg]=\Im\lf[f(z)\ \ov{\vec{H}_0}\rg]\quad,
\ee 
where $\ov{\vec{H}_0}$ is the expression in this chart of the conjugate to the Weingarten operator for the immersion $\vec{\Phi}$ :
\be
\label{VI.199}
\ov{\vec{H}_0}:=\frac{1}{2}\lf[\vec{\mathbb I}(e_1,e_1)-\vec{\mathbb I}(e_2,e_2)+2\,i\, \vec{\mathbb I}(e_1,e_2)\rg]\quad,
\ee
where $(e_1,e_2)$ is an arbitrary orthonormal frame of $(D^2,\vec{\Phi}^\ast g_{{\R}^m})$. \hfill $\Box$
\end{Th}
{\bf Proof of theorem~\ref{th-VI.25}.}
First we assume that $\vec{\Phi}$ satisfies the conservation laws (\ref{VI.197}). Let 
\[
\vec{e}_i:=e^{-\la}\,\p_{x_i}\vec{\Phi}\quad.
\]
where $e^\la=|\p_{x_1}\vec{\Phi}|=|\p_{x_2}\vec{\Phi}|$. And again we use the complex notations :$z=x_1+ix_2$, $\p_z=2^{-1}(\p_{x_1}-i\p_{x_2})$, $\p_{\ov{z}}=2^{-1}(\p_{x_1}+i\p_{x_2})$, and
\[
\lf\{
\begin{array}{l}
\ds\vec{e}_{z}:=e^{-\la}\p_z\vec{\Phi}=2^{-1}(\vec{e}_1-i\vec{e}_2)\\[5mm]
\ds\vec{e}_{\ov{z}}:=e^{-\la}\p_{\ov{z}}\vec{\Phi}=2^{-1}(\vec{e}_1+i\vec{e}_2)
\end{array}
\rg.
\]
Recall that from (\ref{VI.164}) the conservation laws (\ref{VI.197}) are equivalent to
\be
\label{VI.205}
\lf\{
\begin{array}{l}
\ds\lf<\nabla\vec{\Phi},\nabla^\perp\vec{L}\rg>=0\\[5mm]
\ds\nabla\vec{\Phi}\wedge\lf[\nabla^\perp\vec{L}+2\nabla\vec{H}\rg]=0
\end{array}
\rg.
\ee
Using the complex notations it becomes
\be
\label{VI.206}
\lf\{
\begin{array}{l}
\ds\Im\lf<\,\vec{e}_{\ov{z}}\,,\,\p_{z}\vec{L}\,\rg>=0\\[5mm]
\ds\Im\lf(\vec{e}_{\ov{z}}\wedge\lf[\p_z\vec{L}+2\,i\ \p_z\vec{H}\rg]\rg)=0
\end{array}
\rg.
\ee
Denote $$\p_z\vec{L}=A\, \vec{e}_{z}+B\,\vec{e}_{\ov{z}}+\vec{V}$$ where $A$ and $B$ are complex number and $\vec{V}:=\pi_{\vec{n}}(\p_z\vec{L})$
is a complex valued normal vector to the immersed surface. The first equation of (\ref{VI.206}), using (\ref{VI.199a}), is equivalent to
\be
\label{VI.207}
\Im\, A=0
\ee
Observe that if we write $$\p_z\vec{H}=C\ \vec{e}_z+D\ \vec{e}_{\ov{z}}+\vec{W}$$ where $\vec{W}=\pi_{\vec{n}}(\p_z\vec{H})$, one has, using (\ref{VI.201}) and the fact
that $\vec{H}$ is orthogonal to $\vec{e}_{\ov{z}}$
\be
\label{VI.208aa}
C=2\lf<\vec{e}_{\ov{z}},\p_z\vec{H}\rg>=-2\lf<\p_z(e^\la\ \vec{e}_{\ov{z}}),\vec{H}\rg>\ e^{-\la}=-e^\la\ |\vec{H}|^2\quad.
\ee
Hence we deduce in particular
\be
\label{VI.208}
\Im\, C=0\quad.
\ee
We have moreover using (\ref{VI.204})
\be
\label{VI.208ab}
D=2\lf<\vec{e}_{{z}},\p_z\vec{H}\rg>=-2\lf<\p_z(e^{-\la}\ \vec{e}_{{z}}),\vec{H}\rg>\ e^{\la}=-e^\la\ <\vec{H}_0,\vec{H}>\quad.
\ee
Thus combining (\ref{VI.208aa}) and (\ref{VI.208ab}) we obtain
\be
\label{VI.208ac}
\p_z\vec{H}=-|\vec{H}|^2\ \p_z\vec{\Phi}-<\vec{H}_0,\vec{H}>\ \p_{\ov{z}}\vec{\Phi}+\pi_{\vec{n}}(\p_z\vec{H})
\ee
The second line in the conservation law (\ref{VI.206}) is equivalent to
\be
\label{VI.209}
\lf\{
\begin{array}{l}
\ds\Im(i\,A-2C)=0\\[5mm]
\ds\Im\lf(\vec{e}_{\ov{z}}\wedge\lf[\vec{V}+2i\vec{W}\rg]\rg)=0
\end{array}
\rg.
\ee
We observe that $\vec{e}_1\wedge\lf[\vec{V}+2i\vec{W}\rg]$ and $\vec{e}_2\wedge\lf[\vec{V}+2i\vec{W}\rg]$ are linearly independent since $\lf[\vec{V}+2i\vec{W}\rg]$ is orthogonal to the tangent plane, moreover we combine (\ref{VI.208}) and (\ref{VI.209}) and we obtain that (\ref{VI.209}) is equivalent to
\be
\label{VI.210}
\lf\{
\begin{array}{l}
\ds\Im(i\,A)=0\\[5mm]
\ds\vec{e}_1\wedge\Im\lf(\vec{V}+2i\vec{W}\rg)=0\\[5mm]
\ds\vec{e}_2\wedge\Im\lf(i\ \lf[\vec{V}+2i\vec{W}\rg]\rg)=0
\end{array}
\rg.
\ee
Combining (\ref{VI.207}) and (\ref{VI.210}) we obtain that  the conservation laws (\ref{VI.197}) are equivalent to
\be
\label{VI.211}
\lf\{
\begin{array}{l}
A=0\\[5mm]
\vec{V}=-2i\,\vec{W}=-2i\,\pi_{\vec{n}}(\p_z\vec{H})
\end{array}
\rg.
\ee
Or in other words, for a conformal immersion $\vec{\Phi}$ of the disc into ${\R}^m$, there exists $\vec{L}$ from $D^2$ into ${\R}^m$ such that (\ref{VI.197}) holds if and only if there exists
a complex valued function $B$ and a map $\vec{L}$  from $D^2$ into ${\R}^m$ such that 
\be
\label{VI.212}
\p_z\vec{L}=B\,\vec{e}_{\ov{z}}-2i\,\pi_{\vec{n}}(\p_z\vec{H})\quad.
\ee
We shall now exploit the crucial fact that $\vec{L}$ is real valued by taking $\p_{\ov{z}}$ of (\ref{VI.212}). 
Let 
\be
\label{VI.212zz}
f:=\,e^\la B+2\,i\,e^{2\la}\,\lf<\vec{H},\vec{H}_0\rg>\quad.
\ee
With this notation (\ref{VI.212}) becomes\footnote{In real notations this reads also, after using the identity (\ref{z-VI.2}) in lemma~\ref{lm-VI.11}
\be
\label{VI.212bz}
\nabla^\perp\vec{L}=e^{-2\la}\lf(
\begin{array}{cc}
a & b\\[5mm]
-b & a
\end{array}\rg)\ \lf(
\begin{array}{c}
\p_{x_2}\vec{\Phi}\\[5mm]
\p_{x_1}\vec{\Phi}
\end{array}\rg)+\nabla\vec{H}-3\pi_{\vec{n}}(\nabla\vec{H})+\star\nabla^\perp\vec{n}\wedge\vec{H}
\ee
where $f=a+ib$.}
\be
\label{VI.212az}
\p_z\vec{L}=e^{-\la}\, f\,\vec{e}_{\ov{z}}-2i\lf<\vec{H},\vec{H}_0\rg>\ \p_{\ov{z}}\vec{\Phi}-2i\,\pi_{\vec{n}}(\p_z\vec{H})\quad.
\ee

We have using (\ref{VI.204})
\be
\label{VI.213}
\begin{array}{l}
\p_{\ov{z}}\p_z\vec{L}=\p_{\ov{z}}f\ e^{-\la}\,\vec{e}_{\ov{z}}+2^{-1}\,f\, \ov{\vec{H}_0}\\[5mm]
\ \quad\quad\quad-2i\,\p_{\ov{z}}\lf[\lf<\vec{H},\vec{H}_0\rg>\ \p_{\ov{z}}\vec{\Phi}+\pi_{\vec{n}}(\p_z\vec{H})\rg]\quad.
\end{array}
\ee
The fact that $\vec{L}$ is real valued implies that 
\[
\begin{array}{l}
0=\Im\lf(\p_{\ov{z}}f\ e^{-\la}\,\vec{e}_{\ov{z}}\rg)+2^{-1}\Im\lf(f\, \ov{\vec{H}_0}\rg)\\[5mm]
\quad-2\ \Re\lf(\p_{\ov{z}}\lf[\lf<\vec{H},\vec{H}_0\rg>\ \p_{\ov{z}}\vec{\Phi}+\pi_{\vec{n}}(\p_z\vec{H})\rg]\rg)
\end{array}
\]
Using identity (\ref{z-VI.98}), the previous equality is equivalent to
\be
\label{VI.214}
\begin{array}{l}
0=\Im\lf(2\,\p_{\ov{z}}f\ e^{-3\la}\,\vec{e}_{\ov{z}}\rg)+e^{-2\la}\, \Im\lf(f\, \ov{\vec{H}_0}\rg)\\[5mm]
\quad-\Delta_\perp\vec{H}-\ti{A}(\vec{H})+2|\vec{H}|^2\,\vec{H}
\end{array}
\ee
Decomposing this identity into the normal and tangential parts gives that (\ref{VI.214}) is equivalent to
\be
\label{VI.215}
\lf\{
\begin{array}{l}
\ds\Delta_\perp\vec{H}+\ti{A}(\vec{H})-2|\vec{H}|^2\,\vec{H}=e^{-2\la}\, \Im\lf(f\, \ov{\vec{H}_0}\rg)\\[5mm]
\ds\Im\lf(\p_{\ov{z}}f\ \,\vec{e}_{\ov{z}}\rg)=0\quad.
\end{array}
\rg.
\ee
The second line is equivalent to
\[
\p_{\ov{z}}f\ \,\vec{e}_{\ov{z}}-\p_z\ov{f}\,\vec{e}_z=0
\] 
Taking the scalar product respectively with $\vec{e}_z$ and $\vec{e}_{\ov{z}}$, using (\ref{VI.199a}), gives that (\ref{VI.215})
is equivalent to
\be
\label{VI.216}
\lf\{
\begin{array}{l}
\ds\Delta_\perp\vec{H}+\ti{A}(\vec{H})-2|\vec{H}|^2\,\vec{H}=e^{-2\la}\, \Im\lf(f\, \ov{\vec{H}_0}\rg)\\[5mm]
\ds\p_{\ov{z}}f=0\quad.
\end{array}
\rg.
\ee
We have then proved that (\ref{VI.197}) implies (\ref{VI.198}).

\medskip

Assuming now that (\ref{VI.198}) holds we can then go backwards in the equivalences in order to obtain 
\be
\label{VI.217}
\Im\lf[\p_{\ov{z}}\lf[B\,\vec{e}_{\ov{z}}-2i\,\pi_{\vec{n}}(\p_z\vec{H})\rg]\rg]=0
\ee
where $B$ is given by (\ref{VI.212zz}). Observe now that $$\Im[\p_{\ov{z}}\al]=0\quad\quad\Longleftrightarrow\quad\quad\p_{x_2}\al_1+\p_{x_1}\al_2=0$$
where $\al=\al_1+i\al_2$ and $\al_i\in{\R}$. Hence 
$$\Im[\p_{\ov{z}}\al]=0\quad\quad\Longleftrightarrow\quad\quad\exists\, a\in{\R}\ \mbox{ s.t. }\ \al=\p_za\quad.$$
Thus (\ref{VI.217}) is equivalent to the existence of a map $\vec{L}$ from $D^2$ into ${\R}^m$ such that  
\[
\p_z\vec{L}=B\,\vec{e}_{\ov{z}}-2i\,\pi_{\vec{n}}(\p_z\vec{H})\quad.
\]
This is exactly (\ref{VI.212}) for which we have proved that this is equivalent to (\ref{VI.197}). This finishes the proof of theorem~\ref{th-VI.25}.\hfill $\Box$

\medskip

Observe that we have just established the following lemma which gives some new useful formula.

\begin{Lm}
\label{lm-VI.26}
Let $\vec{\Phi}$ be a conformal immersion of the disc $D^2$. $\vec{\Phi}$ is {\bf conformal Willmore} on $D^2$ if and only  there exists 
a smooth map $\vec{L}$ from $D^2$ into ${\R}^m$ and an holomorphic function $f(z)$ such that
\be
\label{VI.218}
\p_z(\vec{L}-2i\vec{H})=2i\ |\vec{H}|^2\ \p_z\vec{\Phi}+[e^{-2\la}\ f(z)-4i\ <\vec{H},\vec{H}_0>]\ \p_{\ov{z}}\vec{\Phi}\quad.
\ee
In particular the following system holds
\be
\label{VI.219}
\lf\{
\begin{array}{l}
<\p_{\ov{z}}\vec{\Phi},\p_z(\vec{L}-2i\vec{H})>=i\ |\vec{H}|^2\ e^{2\la}\\[5mm]
<\p_{{z}}\vec{\Phi},\p_z(\vec{L}-2i\vec{H})>=2^{-1}\ f(z)-2i\ e^{2\la}\ <\vec{H},\vec{H}_0>
\end{array}
\rg.
\ee
\hfill $\Box$
\end{Lm}


\newpage

\begin{thebibliography}{99}
 \bibitem[Ad]{Ad} Adams, David R. ''A note on Riesz potentials.'' Duke Math. J. 42 (1975), no. 4, 765--778.
 \bibitem[AdFo]{AdFo} Adams, Robert A.; Fournier, John J. F. Sobolev spaces. Second edition. Pure and Applied Mathematics (Amsterdam), 140. Elsevier/Academic Press, Amsterdam, 2003.
 \bibitem[Bet1]{Bet1} Bethuel, Fabrice ''Un rŽsultat de rŽgularit\'e pour les solutions de l'Žquation de surfaces ˆ courbure moyenne prescrite.'' (French) [A regularity result for solutions to the equation of surfaces of prescribed mean curvature] C. R. Acad. Sci. Paris SŽr. I Math. 314 (1992), no. 13, 1003--1007.
 \bibitem[Bet2]{Bet2} Bethuel, Fabrice ''On the singular set of stationary harmonic maps.'' Manuscripta Math. 78 (1993), no. 4, 417--443. 
 \bibitem[Bla1]{Bla1} Blaschke, Wilhelm ''Vorlesungen \"uber differentialgeometrie und geometrische Grundlagen von Einsteins Relativit\"atstheorie''  I Die Grundlehren der mathematischen wissenschaften in einzeldarstellungen. Bd I, Elementare differentialgeometrie. 3. erweiterte aufl., bearbeitet und hrsg. von Gerhard Thomsen. 1930
 \bibitem[Bla3]{Bla3} Blaschke, Wilhelm ''Vorlesungen \"uber Differentialgeometrie und geometrische grundlagen von Einsteins relativit\"atstheorie'' III  Die Grundlehren der mathematischen wissenschaften in einzeldarstellungen. Bd. XXIX Differentialgeometrie der Kreise und Kugeln, bearbeitet von Gerhard Thomsen. 1929
 \bibitem[BoTu]{BoTu} Bott, Raoul; Tu, Loring W. Differential forms in algebraic topology. Graduate Texts in Mathematics, 82. Springer-Verlag, New York-Berlin, 1982.
 \bibitem[BoBr]{BoBr}  Bourgain, Jean; Brezis, Ha\"\i m ''New estimates for elliptic equations and Hodge type systems.'' J. Eur. Math. Soc. (JEMS) 9 (2007), no. 2, 277Ð315.
 \bibitem[doC1]{doC1} do Carmo, Manfredo Perdig\~ao ''Differential geometry of curves and surfaces.''  Prentice-Hall, Inc., Englewood Cliffs, N.J., 1976. viii+503 pp. 53-02.
 \bibitem[doC2]{doC2} do Carmo, Manfredo Perdig\~ao ''Riemannian geometry''.  Mathematics: Theory \& Applications. Birkh\"auser Boston, Inc., Boston, MA, 1992.
\bibitem[ChLi]{ChLi} Chanillo, Sagun; Li, Yan Yan Continuity of solutions of uniformly elliptic equations in R2. Manuscripta Math. 77 (1992), no. 4, 415Ð433.
 \bibitem[Che]{Che} Chen, Bang-Yen ``Some conformal invariants of submanifolds and their applications."  Boll. Un. Mat. Ital. (4) 10 (1974), 380--385.
 \bibitem[BYC2]{BYC2} Chen, Bang-yen ''On an inequality of T. J. Willmore.'' Proc. Amer. Math. Soc. 26 1970 473--479.
 \bibitem[Cho]{Cho} Chon\'e, Philippe ''A regularity result for critical points of conformally invariant functionals.'' Potential Anal. 4 (1995), no. 3, 269--296. 
 \bibitem[CLMS]{CLMS} Coifman, R.; Lions, P.-L.; Meyer, Y.; Semmes, S. "Compensated compactness and Hardy spaces".  J. Math. Pures Appl. (9) 72 (1993), no. 3, 247--286.
 \bibitem[DR1]{DR1} Da Lio Francesca and Rivi\`ere Tristan '' 3-commutator estimates and the regularity of
 1/2-harmonic maps into spheres.'' to appear in Analysis and PDE (2010).
 \bibitem[DR2]{DR2} Da Lio Francesca and Rivi\`ere Tristan ''Sub-criticality of non-local Schr\"odinger systems with antisymmetric potentials and applications'' preprint 2010.
  \bibitem[DHKW1]{DHKW1} Dierkes, Ulrich; Hildebrandt, Stefan; K\"uster, Albrecht; Wohlrab, Ortwin ''Minimal surfaces. I. Boundary value problems''. Grundlehren der Mathematischen Wissenschaften, 295. Springer-Verlag, Berlin, 1992.
 \bibitem[DHKW2]{DHKW2} Dierkes, Ulrich; Hildebrandt, Stefan; K\"uster, Albrecht; Wohlrab, Ortwin ''Minimal surfaces. II. Boundary regularity''. Grundlehren der Mathematischen Wissenschaften, 296. Springer-Verlag, Berlin, 1992.
 \bibitem[Ev]{Ev} Evans Craig "Partial regularity for stationary harmonic maps into spheres" Arch. Rat. Mech.
Anal. 116 (1991), 101-113.
 \bibitem[EvGa]{EvGa} Evans, Lawrence C.; Gariepy, Ronald F. Measure theory and fine properties of functions. Studies in Advanced Mathematics. CRC Press, Boca Raton, FL, 1992.
\bibitem[Fal]{Fal}  Falconer, K. J. The geometry of fractal sets. Cambridge Tracts in Mathematics, 85. Cambridge University Press, Cambridge, 1986.
\bibitem[Fe]{Fe}  Federer, Herbert ''Geometric measure theory.'' Die Grundlehren der mathematischen Wissenschaften, Band 153 Springer-Verlag New York Inc., New York 1969 
 \bibitem[Fre]{Fre} Frehse, Jens "A discontinuous solution of a midly nonlinear elliptic system".
Math. Z. 134 (1973), 229-230.
\bibitem[FMS]{FMS} Freire, Alexandre; M\"uller, Stefan; Struwe, Michael ''Weak convergence of wave maps from $(1+2)$-dimensional Minkowski space to Riemannian manifolds.'' Invent. Math. 130 (1997), no. 3, 589--617.
\bibitem[FJM]{FJM} Friesecke, Gero; James, Richard D.; M\"uller, Stefan A theorem on geometric rigidity and the derivation of nonlinear plate theory from three-dimensional elasticity. Comm. Pure Appl. Math. 55 (2002), no. 11, 1461--1506.
\bibitem[GHL]{GHL} Gallot, Sylvestre; Hulin, Dominique; Lafontaine, Jacques ''Riemannian geometry.'' Third edition. Universitext. Springer-Verlag, Berlin, 2004.
\bibitem[Ge]{Ge} Ge, Yuxin ``Estimations of the best constant involving the $L^2$ norm in Wente's inequality
and compact $H-$Surfaces in Euclidian space.'' C.O.C.V., 3, (1998), 263-300.
\bibitem[Gi]{Gi} Giaquinta, Mariano
''Multiple integrals in the calculus of variations and nonlinear elliptic systems.'' 
Annals of Mathematics Studies, 105. Princeton University Press, Princeton, NJ, 1983.
\bibitem[GT]{GT}  Gilbarg, David; Trudinger, Neil S. Elliptic partial differential equations of second order. Reprint of the 1998 edition. Classics in Mathematics. Springer-Verlag, Berlin, 2001.
\bibitem[Gr]{Gr} Gr\"uter, Michael ``Conformally invariant variational integrals and the removability of isolated singularities.''
  Manuscripta Math.  47  (1984),  no. 1-3, 85--104.
\bibitem[Gr2]{Gr2} Gr\"uter, Michael ``Regularity of weak $H$-surfaces''.  J. Reine Angew. Math.  329  (1981), 1--15.
\bibitem[GOR]{GOR}  Gulliver, R. D., II; Osserman, R.; Royden, H. L. ''A theory of branched immersions of surfaces.''  Amer. J. Math. 95 (1973), 750Ð812.
\bibitem[Hei1]{Hei1}ÊHeinz, Erhard "Ein RegularitŠtssatz fŸr schwache Lšsungen nichtlinearer elliptischer Systeme". (German) Nachr. Akad. Wiss. G\"ottingen Math.-Phys. Kl. II 1975, no. 1, 1--13. 
  \bibitem[Hei2]{Hei2} Heinz, Erhard "†ber die RegularitŠt schwacher Lšsungen nichtlinearer elliptischer Systeme". (German) [On the regularity of weak solutions of nonlinear elliptic systems] Nachr. Akad. Wiss. G\"ottingen Math.-Phys. Kl. II 1986, no. 1, 1--15. 
  \bibitem[He]{He} F.H\'elein ``Harmonic maps, conservation laws and moving frames'' Cambridge Tracts in Math. 150, Cambridge Univerity Press, 2002.
  \bibitem[Helf]{Helf} Helfrich, Wolfgang, Elastic properties of lipid bilayers - theory and possible experiments. Zeitschrift Fur Naturforschung C - A Journal Of Biosciences. 28. (1973) 693-703.
 \bibitem[Hil]{Hil} Hildebrandt, S. "Nonlinear elliptic systems and harmonic mappings."
 Proceedings of the 1980 Beijing Symposium on Differential Geometry and Differential
 Equations, vol 1,2,3 (Beijing, 1980), 481-615, Science Press, Beijing, 1982.
 \bibitem[Hil2]{Hil2} Hildebrandt, S. "Quasilinear elliptic systems in diagonal form". Systems of nonlinear partial differential equations (Oxford, 1982), 173-217, NATO Adv. Sci. Inst. Ser. C Math. Phys. Sci., 111,
 Reidel, Dordrecht, 1983.
 \bibitem[Hus]{Hus} Husemoller, Dale Fibre bundles. Third edition. Graduate Texts in Mathematics, 20. Springer-Verlag, New York, 1994.
 \bibitem[JaTa]{JaTa} Jaffe, Arthur; Taubes, Clifford Vortices and monopoles.  ''Structure of static gauge theories''. Progress in Physics, 2. Birkh\"auser, Boston, Mass., 1980.
 \bibitem[KaRu]{KaRu} Katzman, Dan; Rubinstein, Jacob ''Method for the design of multifocal optical elements'' United States Patent No.: US006302540B1, Oct. 16, 2001.
 \bibitem[Ki]{Ki} Kirchhoff, Gustav \"Uber das Gleichgewicht und die Bewegung einer elastischen Scheibe. J. Reine  Angew. Math. 40 (1850), 51Ð88.
 \bibitem[LaRi]{LaRi} Lamm, Tobias; Rivi\`ere Tristan ''Conservation laws for fourth order systems in four dimensions''. Comm.P.D.E., 33 (2008), no. 2, 245-262.
 \bibitem[LaLi]{LaLi} Landau, L. D.; Lifschitz, E. M. ''Theory of elasticity'' - Course of theoretical physics, volume 7 - third edition revised and enlarged by E.M.Lifshitz, A.M.Kosevich and L.P.Pitaevskii,
 Butterwoth-Heinemann 1986.
 \bibitem[LiYa]{LiYa} Li, P., Yau, S.-T., A New Conformal Invariant and its Applications to the
Willmore Conjecture and the First Eigenvalue on Compact Surfaces, Inventiones
Math. 69 (1982), 269-291.
\bibitem[MaNe]{MaNe} Marques, F.C.; Neves, A. ''Min-Max theory and the Willmore conjecture''.  arXiv:1202.6036 (2012).
\bibitem[Mor]{Mor}  Morgan, Frank ''Geometric measure theory. A beginner's guide''. Fourth edition. Elsevier/Academic Press, Amsterdam, 2009.
\bibitem[Mor1]{Mor1} Charles B. Morrey, Jr. ''The Problem of Plateau on a Riemannian Manifold''
Annals of Math., 49, No. 4 (1948), pp. 807-851. 
 \bibitem[Mul]{Mul} M\"uller, Stefan  "Higher integrability of determinants and weak convergence in $L\sp 1$". 
J. Reine Angew. Math. 412 (1990), 20--34.
\bibitem[Oss]{Oss}  Osserman, Robert ''A proof of the regularity everywhere of the classical solution to Plateau's problem''. Ann. of Math. (2) 91 1970 550Ð569.
\bibitem[OGR]{OGR} Gulliver, R. D., II; Osserman, R.; Royden, H. L. ''A theory of branched immersions of surfaces''.  Amer. J. Math. 95 (1973), 750Ð812. 
\bibitem[Poi]{Poi} Poisson, Sim\'eon Denis ''Extrait d'un m\'emoire sur les surfaces \'elastiques'' Correspondance de l'\'Ecole Polytechnique 3, 154 (1816).
\bibitem[Riv]{Riv}  Rivi\`ere, Tristan ''Everywhere discontinuous harmonic maps into spheres.'' Acta Math. 175 (1995), no. 2, 197Ð226.
 \bibitem[Riv1]{Riv1} Rivi\`ere, Tristan ''Conservation laws for conformally invariant variational problems.'' Invent. Math. 168 (2007), no. 1, 1--22.
\bibitem[Riv2]{Riv2} Rivi\`ere, Tristan ''Analysis aspects of Willmore surfaces,'' Inventiones Math., 174 (2008), no.1, 1-45.
\bibitem[Riv3]{Riv3} Rivi\`ere, Tristan  ''Sub-criticality of Schr\"odinger Systems with Antisymmetric Potentials'', preprint (2009).
\bibitem[Riv4]{Riv4} Rivi\`ere, Tristan ''Variational Principles for immersed Surfaces with L2-bounded Second Fundamental Form'', arXiv:1007.2997 (2010).
\bibitem[RiSt]{RiSt} Rivi\`ere, Tristan; Struwe, Michael ''Partial regularity for harmonic maps and related problems.'' 
Comm. Pure Appl. Math. 61 (2008), no. 4, 451--463. 
\bibitem[Rud]{Rud} Rudin, Walter ''Real and complex analysis''. Third edition. McGraw-Hill Book Co., New York, 1987.
\bibitem[Sha]{Sha} Shatah, Jalal
''Weak solutions and development of singularities of the ${\rm SU}(2)$ $\sigma$-model.'' 
Comm. Pure Appl. Math. 41 (1988), no. 4, 459--469. 
 \bibitem[ShS]{ShS} Shatah, Jalal; Struwe, Michael The Cauchy problem for wave maps. Int. Math. Res. Not. 2002, no. 11, 555--571.
\bibitem[Tao1]{Tao1} Tao, Terence ''Global regularity of wave maps. I. Small critical Sobolev norm in high dimension.'' Internat. Math. Res. Notices 2001, no. 6, 299--328.
\bibitem[Tao2]{Tao2} Tao, Terence ''Global regularity of wave maps. II. Small energy in two dimensions.'' Comm. Math. Phys. 224 (2001), no. 2, 443--544.
 \bibitem[Tar1]{Tar1} Tartar, Luc "Remarks on oscillations and Stokes' equation.  Macroscopic modelling of turbulent flows" (Nice, 1984), 24--31, Lecture Notes in Phys., 230, Springer, Berlin, 1985
\bibitem[Tar2]{Tar2} Tartar, Luc ''An introduction to Sobolev spaces and interpolation spaces.'' Lecture Notes of the Unione Matematica Italiana, 3. Springer, Berlin; UMI, Bologna, 2007.
\bibitem[Tho]{Tho} Thomsen, Gerhard ''\"Uber konforme Geometrie I ; Grundlagen der konformen Fl\"achentheorie." Ab. Math. Sem. Univ. Hamburg, 3, no 1 (1924) 31-56.
 \bibitem[To]{To} Topping, Peter ``The optimal constant in Wente's $L^\infty$ estimate'', Comm. Math. Helv. 72, (1997), 316-328.
 \bibitem[Uhl]{Uhl} Uhlenbeck, Karen K. ``Connections with $L\sp{p}$ bounds on curvature.''  Comm. Math. Phys.  83  (1982),
 no. 1, 31--42.
 \bibitem[Wei]{Wei} Weiner, Joel ``On a problem of Chen, Willmore, et al.'' Indiana U. Math. J., 27, no 1 (1978), 19-35.
 \bibitem[We]{We} Wente, Henry C.  "An existence theorem for surfaces of constant mean curvature". 
 \bibitem[Win]{Win} Wintgen, Peter On the total curvature of surfaces in E4. Colloq. Math. 39 (1978), no. 2, 289Ð296.
J. Math. Anal. Appl. 26 1969 318--344.


 

\end{thebibliography}
\end{document}